\documentclass[10pt]{article}

\usepackage[dvips]{epsfig} 
\usepackage{amssymb,amsmath,amsthm,mathrsfs} 
\usepackage{graphicx}
\usepackage{lineno}

\usepackage[authoryear,sort&compress]{natbib}
\usepackage{url}
\usepackage{enumerate}
\usepackage{colordvi,color,pspicture}
\usepackage{psfrag}

\DeclareMathOperator{\argsup}{argsup}

\DeclareMathOperator{\arginf}{arginf}

\DeclareMathOperator{\PAE}{PAE}
\DeclareMathOperator{\HLAE}{HLAE}

\DeclareMathOperator{\signum}{signum}

\newcommand{\E}[1]{\mathbf{E}#1}

\newcommand{\Var}[1]{\mathbf{Var}#1}
\newcommand{\Cov}[1]{\mathbf{Cov}#1}
\newcommand{\I}{\mathbf{I}}

\newcommand{\out}{\text{out}}

\newcommand{\C}{\mathcal{C}}

\newcommand{\U}{\mathcal{U}}

\newcommand{\G}{\Gamma}
\newcommand{\V}{\mathcal{V}}
\newcommand{\A}{\mathcal{A}}
\newcommand{\Y}{\mathcal{Y}}
\newcommand{\y}{\mathsf{y}}
\newcommand{\X}{\mathcal{X}}

\newcommand{\mS}{\mathcal{S}}

\newcommand{\mI}{\mathcal{I}}

\newcommand{\NCS}{N_{CS}}
\newcommand{\NPE}{N_{PE}}
\newcommand{\TY}{T(\Y_3)}
\newcommand{\N}{\mathcal{N}}
\newcommand{\R}{\mathbb{R}}

\newcommand{\RS}{\mathscr{R}_S}

\newcommand{\ve}{\varepsilon}

\theoremstyle{plain}
\newtheorem{theorem}{Theorem}[section]
\newtheorem{lemma}[theorem]{Lemma}

\newtheorem{corollary}[theorem]{Corollary}

\theoremstyle{definition}

\theoremstyle{remark}

\newtheorem{remark}[theorem]{Remark}

\begin{document}

\title{Technical Report \# KU-EC-10-3:\\
A Comparison of Two Proximity Catch Digraph Families in Testing Spatial Clustering}
\author{
Elvan Ceyhan\thanks{
Address:
Department of Mathematics, College of Sciences,
Ko\c{c} University, 34450 Sar{\i}yer, Istanbul, Turkey.
e-mail: elceyhan@ku.edu.tr, tel:+90 (212) 338-1845, fax: +90 (212) 338-1559.
}
}

\date{\today}
\maketitle

\begin{abstract}
\noindent
We consider two parametrized random digraph families,
namely, proportional-edge and central similarity \emph{proximity catch digraphs} (PCDs)
and
compare the performance of these two PCD families in testing spatial point patterns.
These PCD families are based on relative positions of data points from two classes
and the relative density of the PCDs is used as a statistic for testing segregation and association
against complete spatial randomness.
When scaled properly, the relative density of a PCD is a $U$-statistic.
We extend the distribution of the relative density of central similarity PCDs
for expansion parameter being larger than one.
We compare the asymptotic distribution of the statistic
for the two PCD families,
using the standard central limit theory of $U$-statistics.
We compare finite sample performance of the tests by Monte Carlo simulations
and prove the consistency of the tests under the alternatives.
The asymptotic performance of the tests under the alternatives
is assessed by Pitman's asymptotic efficiency.
We find the optimal expansion parameters of the PCDs
for testing each of the segregation and association alternatives
in finite samples and in the limit.
We demonstrate that in terms of empirical power (i.e., for finite samples)
relative density of central similarity PCD has better performance
(which occurs for expansion parameter values larger than one)
under segregation alternative,
while relative density of proportional-edge PCD has better performance
under association alternative.
The methods are illustrated in a real-life example from plant ecology.
\end{abstract}


\section{Introduction}
\label{sec:intro}
Spatial clustering has received considerable attention
in the statistical literature.
In recent years, a new clustering
approach has been developed which uses data-random proximity catch digraphs (PCDs)
and is based on the relative positions
of the data points from various classes.
A catch digraph is a directed graph whose vertices are pointed
sets (a pointed set is a pair $(S,p)$ where $S$ is a set and $p$ a distinguished point)
with an arc from vertex $(S_u,p_u)$ to vertex $(S_v,p_v)$ whenever $p_v \in S_u$.
Hence $S_u$ {\em catches} $p_v$.
\cite{priebe:2001} introduced the class cover catch digraphs (CCCDs) and gave the exact
and the asymptotic distribution of the domination number of the CCCD in $\R$.
For two classes, $\X$ and $\Y$, of points, let $\X$ be the class of interest
and $\Y$ be the reference class
and $\X_n$ and $\Y_m$ be samples of size $n$ and $m$ from classes $\X$ and $\Y$, respectively.
In the CCCD approach the points correspond to observations from class $\X$ and
the sets are defined to be (open) balls centered at the points with maximal radius
(relative to the other class $\Y$):
$S_x = B(x,r(x))$, where $r(x)=d(x,\Y_m)$ is the minimum
distance between the observation $x \in \X$ and
the observations of the other class, $\Y_m$.
The CCCD approach is extended to multiple dimensions by \cite{devinney:2002a},
\cite{marchette:2003}, \cite{priebe:2003b}, and \cite{priebe:2003a},
who demonstrated relatively good performance of it in classification by employing data reduction (condensing)
based on approximate minimum dominating sets as prototype sets
(since finding the exact minimum dominating set is an NP-hard problem ---in particular for CCCDs).

\cite{ceyhan:Phd-thesis} generalized CCCDs to PCDs.
In the PCD approach the points correspond to observations from class $\X$ and
the sets are defined to be (closed) regions (usually convex regions or simply triangles)
based on class $\X$ and $\Y$ points
and the regions increase as the distance of a class $\X$ point
from class $\Y$ points increases.
The (non-parametrized) central similarity proximity map and
parameterized proportional-edge proximity maps and the
associated random PCDs are introduced
in \cite{ceyhan:CS-JSM-2003} and \cite{ceyhan:TR-dom-num-NPE-spatial}, respectively.
In both cases, the space is partitioned by the Delaunay tessellation of class $\Y$ points
which is the Delaunay triangulation in $\R^2$.
In each triangle, a family of PCDs is constructed based on the relative
positions of the $\X$ points with respect to each other
and to $\Y$ points.
These proximity maps have the advantage that the calculations yielding the
asymptotic distribution of the relative density are analytically tractable.

Recently, the use of mathematical graphs has gained popularity in spatial analysis (\cite{roberts:2000})
providing a way to move beyond the usual Euclidean metrics for spatial analysis.
Graph theory is well suited to ecological applications
concerned with connectivity or movement,
although it is only recently introduced to landscape ecology (\cite{minor:2007}).
Conventional graphs reduce the utility of other geo-spatial information,
because they do not explicitly maintain geographic reference.
\cite{fall:2007} introduce spatial graphs that
preserve the relevant spatial information
by integrating a geometric reference system that ties patches and paths
to specific spatial locations and spatial dimensions.
However, usually the scale is lost
after a graph is constructed using spatial data
(see for instance, \cite{su:2007}).
Many concepts in spatial ecology depend on the idea of spatial adjacency which
requires information on the close vicinity of an object.
Graph theory conveniently can be adapted to express
and communicate adjacency information allowing one to compute
meaningful quantities related to a spatial point pattern.
Adding vertex and edge properties to graphs extends the problem domain to network modeling (\cite{keitt:2007}).
\cite{wu:2008} propose a new measure based on spatial interaction and graph theory,
which reflects intra-patch and inter-patch relationships by
quantifying contiguity within and among patches.
\cite{friedman:1983} also propose a graph-theoretic method
to measure multivariate association,
but their method is not designed to analyze spatial interaction
between two or more classes;
instead it is an extension of generalized correlation coefficient
(such as Spearman's $\rho$ or Kendall's $\tau$)
to measure multivariate (possibly nonlinear) correlation.

Intuitively, relative density should be useful
for testing association or segregation.
Under association, the observations from one class
tend to cluster around those of the other,
while under segregation they
tend to avoid observations from the other class.
For example, the pattern of spatial segregation has been investigated
for species (\cite{diggle:2003}), age classes of plants (\cite{hamill:1986})
and sexes of dioecious plants (\cite{nanami:1999}).
Under association, the defining proximity regions tend to be small,
and hence there should be fewer arcs;
while under segregation,
the proximity regions tend to be larger and cover many points,
resulting in many arcs.
Thus, the relative density
(number of arcs divided by the total number of possible arcs) is a reasonable
statistic to employ in this problem.
Unfortunately, in the case of the CCCD,
it is difficult to make precise calculations in multiple dimensions due to the geometry of the neighborhoods.
The domination number of the proportional-edge PCD with $r=3/2$
is used for testing segregation or association in \cite{ceyhan:dom-num-NPE-SPL}
and with general $r$ in \cite{ceyhan:dom-num-NPE-Spat2010}.

This is appropriate when both classes are comparably large.
\cite{ceyhan:arc-density-PE} used the relative density of the same proximity digraph for the same purpose
which is appropriate when only size of one of the classes is large.
The parameters of the PCDs expand the associated proximity region
as a function of the distance
from the point defining the proximity region
to the vertices or edges of the triangles
in which the point lies.

In this article,
we compare the two parameterized PCD families,
namely proportional-edge and central similarity PCDs in testing bivariate spatial patterns.
The graph invariant we use as a statistic is the relative density.
We also extend the (expansion) parameter of central similarity PCD
for values larger than one;
previously it was defined on for the range of $(0,1]$ (\cite{ceyhan:CS-JSM-2003,ceyhan:arc-density-CS}).
We compare the finite sample performance of the relative density of
these two PCD families
by empirical size and power analysis
based on extensive Monte Carlo simulations.
We also compare the asymptotic distributions
and asymptotic power performance of the tests under the alternatives.
We first consider the case of one triangle,
followed by the case of multiple triangles
(based on the Delaunay triangulation of four or more $\Y$ points).
We also propose a correction term for the proportion of
$\X$ points that lies outside the convex hull of $\Y$ points.

In Section \ref{sec:prox-map},
we provide a general definition of the proximity maps and the associated PCDs
and their relative density,
describe the two particular PCD families (namely, proportional-edge and central similarity PCDs).
We provide the asymptotic distribution of relative density of the PCDs
for uniform data in one and multiple triangles in Section \ref{sec:asy-dist-rho},
describe the alternative patterns of segregation and association,
provide the asymptotic normality under the alternatives,
present the standardized versions of the test statistics,
and prove their consistency in Section \ref{sec:alternatives}.
We present the empirical size performance of the PCDs in Section \ref{sec:emp-size-anal},
and empirical power analysis under the alternatives
in Section \ref{sec:emp-power-anal} by extensive Monte Carlo simulations.
The asymptotic performance of the tests is assessed by
comparison of Pitman asymptotic efficiency scores in Section \ref{sec:PAE}.
We propose a correction method for the $\X$ points outside the convex hull of $\Y_m$
in Section \ref{sec:conv-hull-correction},
illustrate the use of the tests in an ecological data set in Section \ref{sec:example}.
We present discussion and conclusions in Section \ref{sec:discussion}.
Derivations of some of the quantities and lengthy expressions are deferred to the
Appendix Sections.

\section{Proximity Maps and the Associated PCDs}
\label{sec:prox-map}
Our PCDs are based on the proximity maps which are defined in a
fairly general setting.
Let $(\Omega,\mathcal{M})$ be a measurable space and
consider a function $N:\Omega \times \wp(\Omega) \rightarrow \wp(\Omega)$,
where $\wp(\cdot)$ represents the power set function.
Then given $\Y_m \subseteq \Omega$,
the {\em proximity map}
$N(\cdot) = N(\cdot,\Y_m): \Omega \rightarrow \wp(\Omega)$
associates a {\em proximity region} $N(x) \subseteq \Omega$ with each point $x \in \Omega$.
The region $N(x)$ is defined in terms of the distance between $x$ and $\Y_m$.
If $\X_n:=\{X_1,X_2,\ldots,X_n\}$ is a set of $\Omega$-valued random variables,
then the $N(X_i),\; i=1,2,\ldots,n$, are random sets.
If the $X_i$ are independent and identically distributed (iid),
then so are the random sets $N(X_i)$.

Define the data-random PCD, $D$,
with vertex set $\V=\{X_1,X_2,\ldots,X_n\}$
and arc set $\A$ by
$(X_i,X_j) \in \A \iff X_j \in N(X_i)$.
The random digraph $D$ depends on
the (joint) distribution of the $X_i$ and on the map $N$.
The adjective {\em proximity}
--- for the catch digraph $D$ and for the map $N$ ---
comes from thinking of the region $N(x)$
as representing those points in $\Omega$ ``close'' to $x$.
An extensive treatment of the proximity graphs is presented in
\cite{toussaint:1980} and \cite{jaromczyk:1992}.

The \emph{relative density} of a digraph $D=(\V,\A)$
of order $|\V| = n$,
denoted $\rho(D)$,
is defined as
$$
\rho(D) = \frac{|\A|}{n(n-1)}
$$
where $|\cdot|$ stands for set cardinality (\cite{janson:2000}).
Thus $\rho(D)$ represents the ratio of the number of arcs
in the digraph $D$ to the number of arcs in the complete symmetric
digraph of order $n$, which is $n(n-1)$.

If $X_1,X_2,\ldots,X_n \stackrel{iid}{\sim} F$, then
the relative density
of the associated data-random PCD,
denoted $\rho(\X_n;h,N)$, is a U-statistic,
\begin{eqnarray}
\rho(\X_n;h,N) =
  \frac{1}{n(n-1)}
    \sum\hspace*{-0.1 in}\sum_{i < j \hspace*{0.25 in}}
      \hspace*{-0.1 in}h(X_i,X_j;N)
\end{eqnarray}
where
\begin{eqnarray}
h(X_i,X_j;N)&=& \I \{(X_i,X_j) \in \A\}+ \I \{(X_j,X_i) \in \A\} \nonumber \\
       &=& \I \{X_j \in N(X_i)\}+ \I \{X_i \in N(X_j)\}.
\end{eqnarray}
We denote $h(X_i,X_j;N)$ as $h_{ij}$ for brevity of notation.
Since the digraph is asymmetric, $h_{ij}$ is defined as
the number of arcs in $D$ between vertices $X_i$ and $X_j$, in order to produce a symmetric kernel with finite variance (\cite{lehmann:1988}).

The random variable $\rho_n := \rho(\X_n;h,N)$ depends on $n$ and $N$ explicitly
and on $F$ implicitly.
The expectation $\E[\rho_n]$, however, is independent of $n$
and depends on only $F$ and $N$:
\begin{eqnarray}
0 \leq \E[\rho_n] = \frac{1}{2}\E[h_{12}] \leq 1 \text{ for all $n\ge 2$}.
\end{eqnarray}
The variance $\Var[\rho_n]$ simplifies to
\begin{eqnarray}
\label{eqn:var-rho}
0 \leq
  \Var[\rho_n] =
     \frac{1}{2n(n-1)} \Var[h_{12}] +
     \frac{n-2}{n(n-1)} \Cov[h_{12},h_{13}]
  \leq 1/4.
\end{eqnarray}
A central limit theorem for $U$-statistics
(\cite{lehmann:1988})
yields
\begin{eqnarray}
\sqrt{n}(\rho_n-\E[\rho_n]) \stackrel{\mathcal{L}}{\longrightarrow} \N(0,\Cov[h_{12},h_{13}])
\end{eqnarray}
provided $\Cov[h_{12},h_{13}] > 0$
where $\N(\mu,\sigma^2)$ stands for the normal distribution with mean $\mu$
and variance $\sigma^2$.
The asymptotic variance of $\rho_n$, $\Cov[h_{12},h_{13}]$,
depends on only $F$ and $N$.
Thus, we need determine only
$\E[h_{12}]$
and
$\Cov[h_{12},h_{13}]$
in order to obtain the normal approximation
\begin{eqnarray}
\rho_n \stackrel{\text{approx}}{\sim}
\N\left(\E[\rho_n],\Var[\rho_n]\right) =
\N\left(\frac{\E[h_{12}]}{2},\frac{\Cov[h_{12},h_{13}]}{n}\right) \text{ for large $n$}.
\end{eqnarray}

\subsection{The Proximity Map Families}
\label{sec:two-proximity-map-families}
We now briefly define two proximity map families.
Let $\Omega = \R^d$
and let $\Y_m=\left \{\y_1,\y_2,\ldots,\y_m \right\}$ be $m$ points in
general position in $\R^d$ and $T_i$ be the $i^{th}$ Delaunay cell
for $i=1,2,\ldots,J_m$, where $J_m$ is the number of Delaunay cells.
Let $\X_n$ be a set of iid random variables from distribution $F$ in
$\R^d$ with support $\mS(F) \subseteq \C_H(\Y_m)$
where $\C_H(\Y_m)$ stands for the convex hull of $\Y_m$.
In particular, for illustrative purposes, we focus on $\R^2$ where
a Delaunay tessellation is a \emph{triangulation}, provided that no more
than three points in $\Y_m$ are cocircular (i.e., lie on the same circle).
Furthermore, for simplicity,
let $\Y_3=\{\y_1,\y_2,\y_3\}$ be three non-collinear points
in $\R^2$ and $\TY=T(\y_1,\y_2,\y_3)$ be the triangle
with vertices $\Y_3$.
Let $\X_n$ be a set of iid random variables from $F$ with
support $\mS(F) \subseteq \TY$.
Let $\U(T\left(\Y_3 \right))$ be the uniform distribution on $T\left(\Y_3 \right)$.
If $F=\U(\TY)$, a composition of translation,
rotation, reflections, and scaling
will take any given triangle $\TY$
to the basic triangle $T_b=T((0,0),(1,0),(c_1,c_2))$
with $0 < c_1 \le 1/2$, $c_2>0$,
and $(1-c_1)^2+c_2^2 \le 1$, preserving uniformity.
That is, if $X \sim \U(\TY)$ is transformed in the same manner to,
say $X'$, then we have $X' \sim \U(T_b)$.
In fact this will hold for data from any distribution $F$
up to scale.

\subsubsection{Proportional-Edge Proximity Maps and Associated Proximity Regions}
\label{sec:r-factor-PCD}
For the expansion parameter $r \in [1,\infty]$,
define $\NPE(x,r)$ to be the {\em proportional-edge} proximity map with expansion parameter $r$ as follows;
see also Figure \ref{fig:ProxMapDef} (left).
Using line segments from the center of mass of $\TY$ to the midpoints of its edges,
we partition $\TY$ into ``vertex regions" $R_V(\y_1)$, $R_V(\y_2)$, and $R_V(\y_3)$.
For $x \in \TY \setminus \Y_3$, let $v(x) \in \Y_3$ be the
vertex in whose region $x$ falls, so $x \in R_V(v(x))$.
If $x$ falls on the boundary of two vertex regions,
 we assign $v(x)$ arbitrarily to one of the adjacent regions.
Let $e(x)$ be the edge of $\TY$ opposite $v(x)$.
Let $\ell(x)$ be the line parallel to $e(x)$ through $x$.
Let $d(v(x),\ell(x))$ be the Euclidean distance from $v(x)$ to $\ell(x)$.
For $r \in [1,\infty)$, let $\ell_r(x)$ be the line parallel to $e(x)$
such that $d(v(x),\ell_r(x)) = rd(v(x),\ell(x))$ and $d(\ell(x),\ell_r(x)) < d(v(x),\ell_r(x))$.
Let $T_{PE}(x,r)$ be
the triangle similar to
and with the same orientation as $\TY$
having $v(x)$ as a vertex
and $\ell_r(x)$ as the opposite edge.
Then the {\em proportional-edge} proximity region
$\NPE(x,r)$ is defined to be $T_{PE}(x,r) \cap \TY$.
Notice that $r \ge 1$ implies $x \in \NPE(x,r)$.
Note also that
$\lim_{r \rightarrow \infty} \NPE(x,r) = \TY$
for all $x \in \TY \setminus \Y_3$,
so we define $\NPE(x,\infty) = \TY$ for all such $x$.
For $x \in \Y_3$, we define $\NPE(x,r) = \{x\}$ for all $r \in [1,\infty]$.
See \cite{ceyhan:TR-dom-num-NPE-spatial} for more detail.

\subsubsection{Central Similarity Proximity Maps and Associated Proximity Regions}
\label{sec:tau-factor-PCD}
\begin{figure} [h]
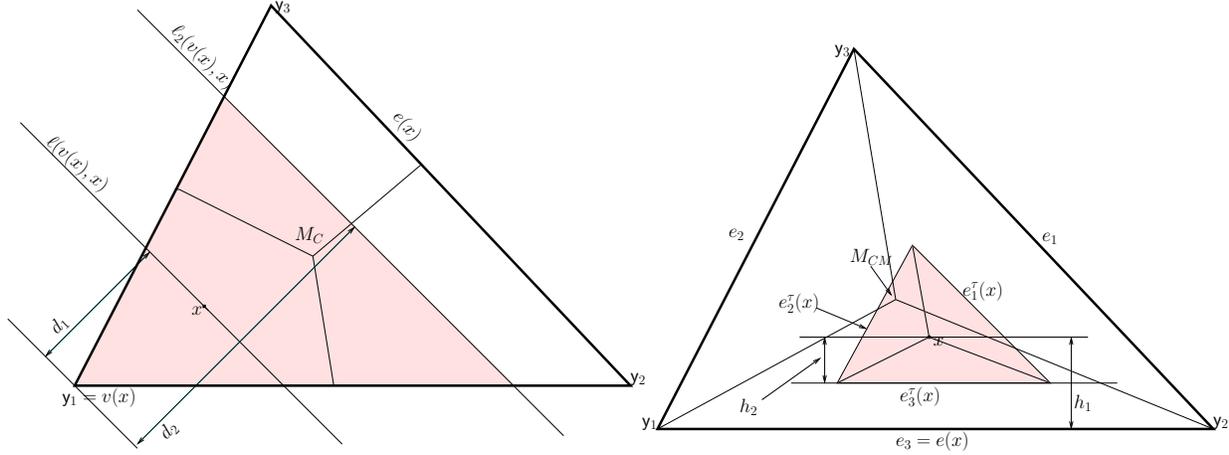

    \centering
    \scalebox{.35}{\input{Nofnu2.pstex_t}}
    \scalebox{.35}{\input{N_CSexample2.pstex_t}}
    \caption{
Plotted in the left is the illustration of the construction of proportional-edge proximity region, $\NPE(x,r=2)$ (shaded region)
for an $x \in R_V(\y_1)$
where $d_1=d(v(x),\ell(v(x),x))$
and
$d_2=d(v(x),\ell_2(v(x),x))=2\,d(v(x),\ell(v(x),x))$;
and in the right is the illustration of the
construction of central similarity proximity region, $\NCS(x,\tau=1/2)$ (shaded region)
for an $x \in R_E(e_3)$
where $h_2=d(x,e_3^\tau(x))=\frac{1}{2}\,d(x,e(x))$ and $h_1=d(x,e(x))$.}
\label{fig:ProxMapDef}
    \end{figure}

For the expansion parameter $\tau \in (0,\infty]$,
define $\NCS(x,\tau)$ to be the {\em central similarity proximity map} with expansion parameter $\tau$ as follows;
see also Figure \ref{fig:ProxMapDef} (right).
Let $e_j$ be the edge opposite vertex $\y_j$ for $j=1,2,3$,
and let ``edge regions'' $R_E(e_1)$, $R_E(e_2)$, $R_E(e_3)$
partition $\TY$ using line segments from the
center of mass of $\TY$ to the vertices.
For $x \in (\TY)^o$, let $e(x)$ be the
edge in whose region $x$ falls; $x \in R_E(e(x))$.
If $x$ falls on the boundary of two edge regions we assign $e(x)$ arbitrarily.
For $\tau > 0$, the central similarity proximity region
$\NCS(x,\tau)$ is defined to be the triangle $T_{CS}(x,\tau) \cap \TY$ with the following properties:
\begin{itemize}
\item[(i)]
For $\tau \in (0,1]$,
the triangle
$T_{CS}(x,\tau)$ has an edge $e_\tau(x)$ parallel to $e(x)$ such that
$d(x,e_\tau(x))=\tau\, d(x,e(x))$
and
$d(e_\tau(x),e(x)) \le d(x,e(x))$
and
for $\tau >1$,
$d(e_\tau(x),e(x)) < d(x,e_\tau(x))$
where $d(x,e(x))$ is the Euclidean distance from $x$ to $e(x)$,
\item[(ii)] the triangle $T_{CS}(x,\tau)$ has the same orientation as and is similar to $\TY$,
\item[(iii)] the point $x$ is at the center of mass of $T_{CS}(x,\tau)$.
\end{itemize}
Note that (i) implies the expansion parameter $\tau$,
(ii) implies ``similarity", and
(iii) implies ``central" in the name, (parametrized) {\em central similarity proximity map}.
Notice that $\tau>0$ implies that $x \in \NCS(x,\tau)$ and,
by construction, we have
$\NCS(x,\tau)\subseteq \TY$ for all $x \in \TY$.
For $x \in \partial(\TY)$ and $\tau \in (0,\infty]$, we define $\NCS(x,\tau)=\{x\}$.
For all $ x\in \TY^o$ the edges
$e_\tau(x)$ and $e(x)$ are coincident iff $\tau=1$.
Note also that
$\lim_{\tau \rightarrow \infty} \NCS(x,\tau) = \TY$
for all $x \in (\TY)^o$,
so we define $\NCS(x,\infty) = \TY$ for all such $x$.
Observe that the central similarity proximity maps in \cite{ceyhan:CS-JSM-2003} and  \cite{ceyhan:arc-density-CS}
are $\NCS(\cdot,\tau)$ with $\tau=1$ and $\tau \in (0,1]$, respectively.

\begin{remark}
Notice that $X_i \stackrel{iid}{\sim} F$,
with the additional assumption
that the non-degenerate two-dimensional
probability density function $f$ exists
with support$(f) \subseteq \TY$,
implies that the special case in the construction
of $\NPE(\cdot,r)$ ---
$X$ falls on the boundary of two vertex regions ---
occurs with probability zero;
similarly,
the special case in the construction
of $\NCS(\cdot,\tau)$ ---
$X$ falls on the boundary of two edge regions ---
occurs with probability zero.
$\square$
\end{remark}

\section{The Asymptotic Distribution of Relative Density for Uniform Data}
\label{sec:asy-dist-rho}

\subsection{The One Triangle Case}
\label{sec:asy-dist-rho-one-tri}
For simplicity, we consider $\X$ points iid uniform in one triangle only.
The null hypothesis we consider is a type of
{\em complete spatial randomness} (CSR); that is,
$$H_o: X_i \stackrel{iid}{\sim} \U(T\left(\Y_3 \right)) \text{ for } i=1,2,\ldots,n.$$
If it is desired to have the sample size be a random variable,
we may consider a spatial Poisson point process on $T\left(\Y_3 \right)$
as our null hypothesis.

We first present a ``geometry invariance" result
that will simplify our subsequent analysis by allowing us to
consider the special case of the equilateral triangle.

\begin{theorem}
\label{thm:geo-inv}
\textbf{(Geometry Invariance for Uniform Data)}
Let $\Y_3 = \{\y_1,\y_2,\y_3\} \subset \mathbb{R}^2$
be three non-collinear points.
For $i=1,2,\ldots,n$,
let $X_i \stackrel{iid}{\sim} F = \mathcal{U}(\TY)$.
Then
\begin{itemize}
\item[(i)]
for any $r \in [1,\infty]$
the distribution of relative density of proportional-edge PCDs,
$\rho_{_{PE}}(n,r)$,
is independent of $\Y_3$,
hence the geometry of $\TY$.
\item[(ii)]
for any $\tau \in (0,\infty]$
the distribution of relative density of central similarity PCDs,
$\rho_{_{CS}}(n,\tau)$,
is independent of $\Y_3$, hence the geometry of $\TY$.
\end{itemize}
\end{theorem}

The proof for (i) is are provided in \cite{ceyhan:arc-density-PE}
and the proof of (ii) for $\tau \in (0,1]$ is provided in \cite{ceyhan:arc-density-CS}
and the proof for $\tau > 1$ is similar.


In fact,
the geometry invariance of $\rho_{_{PE}}(n,\infty)$ for data from any continuous distribution on $\TY$ follows trivially,
since for $r=\infty$, $\rho_{_{PE}}(n,r)=1$ a.s. (i.e., it is degenerate).
Likewise,
the geometry invariance of $\rho_{_{CS}}(n,\infty)$ for data from any continuous distribution on $\TY$ follows trivially,
since for $\tau=\infty$, $\rho_{_{CS}}(n,\tau)=1$ a.s. (i.e., it is degenerate).

Based on Theorem \ref{thm:geo-inv} and our uniform null hypothesis,
we may assume that
$\TY$ is a standard equilateral triangle
with vertices $\Y_3 = \left\{(0,0),(1,0),\bigl( 1/2,\sqrt{3}/2 \bigr)\right\}$
henceforth.

\begin{remark}
Notice that, we proved the geometry invariance property for the relative density of PCDs
based on proportional-edge proximity regions
where vertex regions are defined with the lines joining $\Y_3$ to the center of mass $M_C$.
If we had used the orthogonal projections from $M_C$ to the edges,
the vertex regions (hence $\NPE(\cdot,r)$) would depend on
the geometry of the triangle.
That is, the orthogonal projections
from $M_C$ to the edges will not be mapped to the orthogonal
projections in the standard equilateral triangle.
Hence 
the exact and
asymptotic distribution of the relative density will depend on $c_1,c_2$ of $T_b$,
so one needs to do the calculations for each possible combination of $c_1,c_2$.
$\square$
\end{remark}

\subsection{Asymptotic Normality under the Null Hypothesis}
\label{sec:asy-norm-null}
By detailed geometric  probability calculations,
the means and the asymptotic variances of the relative
density of the proportional-edge and central similarity PCDs
can be calculated explicitly (\cite{ceyhan:arc-density-PE} and \cite{ceyhan:arc-density-CS}).

The central limit theorem for $U$-statistics
then establishes the asymptotic normality under the uniform null hypothesis.
For our proportional-edge proximity map and uniform null hypothesis, the
asymptotic null distribution of $\rho_{_{PE}}(n,r)$ can
be derived as a function of $r$.
Let $\mu_{_{PE}}(r):=\E[\rho_{_{PE}}(n,r)]$ and $\nu_{_{PE}}(r):=\Cov[h_{12},h_{13}]$.
Notice that $\mu_{_{PE}}(r)=\E[h_{12}]/2=P(X_2 \in \NPE(X_1,r))$ is the probability of an
arc occurring between any pair of vertices,
hence is called \emph{arc probability} also.
Similarly,
the asymptotic null distribution of $\rho_{_{CS}}(n,\tau)$
as a function of $\tau$ can be derived.
Let $\mu_{_{CS}}(\tau):=\E[\rho_{_{CS}}(n,\tau)]$,
then $\mu_{_{CS}}(\tau)=\E[h_{12}]/2=P\bigl(X_2 \in \NCS(X_1,\tau)\bigr)$
and
let $\nu_{_{CS}}(\tau):=\Cov[h_{12},h_{13}]$.
These results are summarized in the following theorems.

\begin{theorem}
\label{thm:asy-norm-PE}
For $r \in [1,\infty)$,
\begin{eqnarray}
 \frac{\sqrt{n}\,\bigl(\rho_{_{PE}}(n,r)-\mu_{_{PE}}(r)\bigr)}{\sqrt{\nu_{_{PE}}(r)}}
 \stackrel{\mathcal{L}}{\longrightarrow}
 \N(0,1)
\end{eqnarray}
where
\begin{eqnarray}
\label{eqn:PEAsymean}
\mu_{_{PE}}(r) =
 \begin{cases}
  \frac{37}{216}r^2                                 &\text{for} \quad r \in [1,3/2), \\
  -\frac{1}{8}r^2 + 4 - 8r^{-1} + \frac{9}{2}r^{-2}  &\text{for} \quad r \in [3/2,2), \\
  1 - \frac{3}{2}r^{-2}                             &\text{for} \quad r \in [2,\infty) ,\\
 \end{cases}
\end{eqnarray}
and
\begin{equation}
\label{eqn:PEAsyvar}
\nu_{_{PE}}(r) =\nu_1(r) \,\I(r \in [1,4/3)) + \nu_2(r) \,\I(r \in [4/3,3/2))+ \nu_3(r) \,\I(r \in [3/2,2)) + \nu_4(r) \,\I( r \in [2,\infty])
\end{equation}
with
{\small
\begin{align*}
  \nu_1(r) &=\frac{3007\,r^{10}-13824\,r^9+898\,r^8+77760\,r^7-117953\,r^6+48888\,r^5-24246\,r^4+60480\,r^3-38880\,r^2+3888}{58320\,r^4},\\
  \nu_2(r) &=\frac{5467\,r^{10}-37800\,r^9+61912\,r^8+46588\,r^6-191520\,r^5+13608\,r^4+241920\,r^3-155520\,r^2+15552}{233280\,r^4},  \\
  \nu_3(r) &=-[7\,r^{12}-72\,r^{11}+312\,r^{10}-5332\,r^8+15072\,r^7+13704\,r^6-139264\,r^5+273600\,r^4-242176\,r^3\\
& +103232\,r^2-27648\,r+8640]/[960\,r^6],\\
  \nu_4(r) &=\frac{15\,r^4-11\,r^2-48\,r+25}{15\,r^6}.
\end{align*}
}
For $r=\infty$, $\rho_{_{PE}}(n,r)$ is degenerate.
\end{theorem}

See \cite{ceyhan:arc-density-PE} for the proof.

\begin{theorem}
\label{thm:asy-norm-CS}
For $\tau \in (0,\infty)$,
\begin{eqnarray}
 \frac{\sqrt{n}(\rho_{_{CS}}(n,\tau)-\mu_{_{CS}}(\tau))}{\sqrt{\nu_{_{CS}}(\tau)} }
 \stackrel{\mathcal{L}}{\longrightarrow}
 \N(0,1)
7\end{eqnarray}
where
\begin{eqnarray}
\label{eqn:CSAsymean}
\mu_{_{CS}}(\tau) =
 \begin{cases}
  \tau^2/6                                 &\text{for} \quad \tau \in (0,1], \\
{\frac{\tau\, \left( 4\,\tau-1 \right) }{ 2 \left( 1+2\,\tau  \right)  \left( 2+\tau \right) }}             &\text{for} \quad \tau \in (1,\infty), \\
 \end{cases}
\end{eqnarray}
and
\begin{eqnarray}
\label{eqn:CSAsyvar}
\nu_{_{CS}}(\tau) =
 \begin{cases}
  \frac{\tau^4(6\,\tau^5-3\,\tau^4-25\,\tau^3+\tau^2+49\,\tau+14)}{45\,(\tau+1)(2\,\tau+1)(\tau+2)}         &\text{for} \quad \tau \in (0,1] ,\\
  \\
  {\frac{168\,\tau^7+886\,\tau^6+1122\,\tau^5+45\,\tau^4-470\,\tau^3-114\,\tau^2+48\,\tau+16}
  {5 \left( 2\,\tau+1 \right)^4 \left( \tau+2 \right)^4}}                                              &\text{for} \quad \tau \in (1,\infty). \\
 \end{cases}
\end{eqnarray}
For $\tau=0$, $\rho_{_{CS}}(n,\tau)$ is degenerate.
\end{theorem}

See \cite{ceyhan:arc-density-CS} for the derivation for $\tau \in (0,1]$
and Appendix 1 for $\tau > 1$.

\begin{figure}[h]
\centering
\psfrag{mu(r)}{{\Huge$\mu_{_{PE}}(r)$}}
\psfrag{r}{{\Huge$r$}}
\rotatebox{-90}{ \resizebox{2. in}{!}{\includegraphics{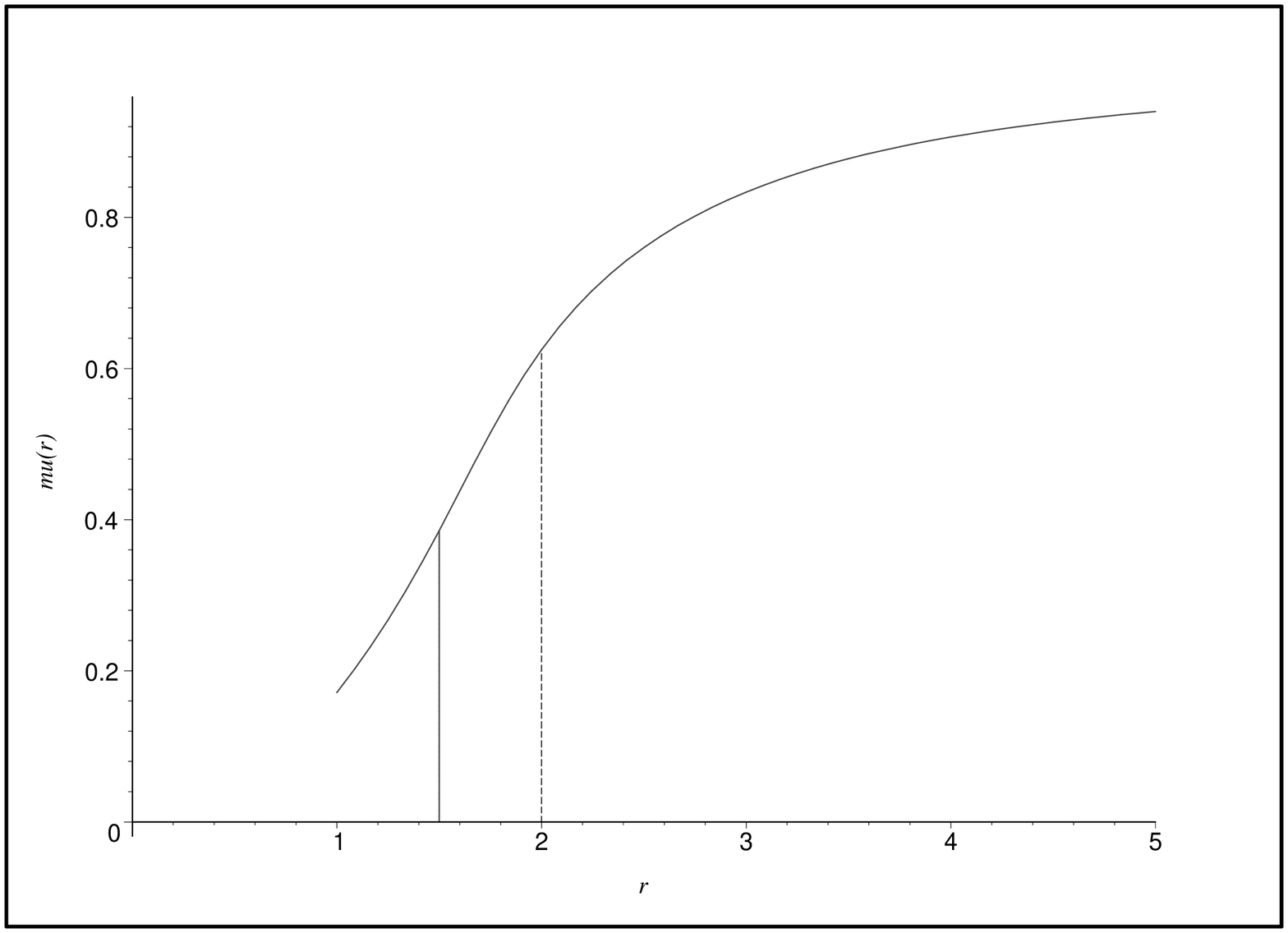}}}
\psfrag{m}{{\Huge$\mu_{_{CS}}(\tau)$}}
\psfrag{tau}{{\Huge$\tau$}}
\rotatebox{-90}{ \resizebox{2. in}{!}{\includegraphics{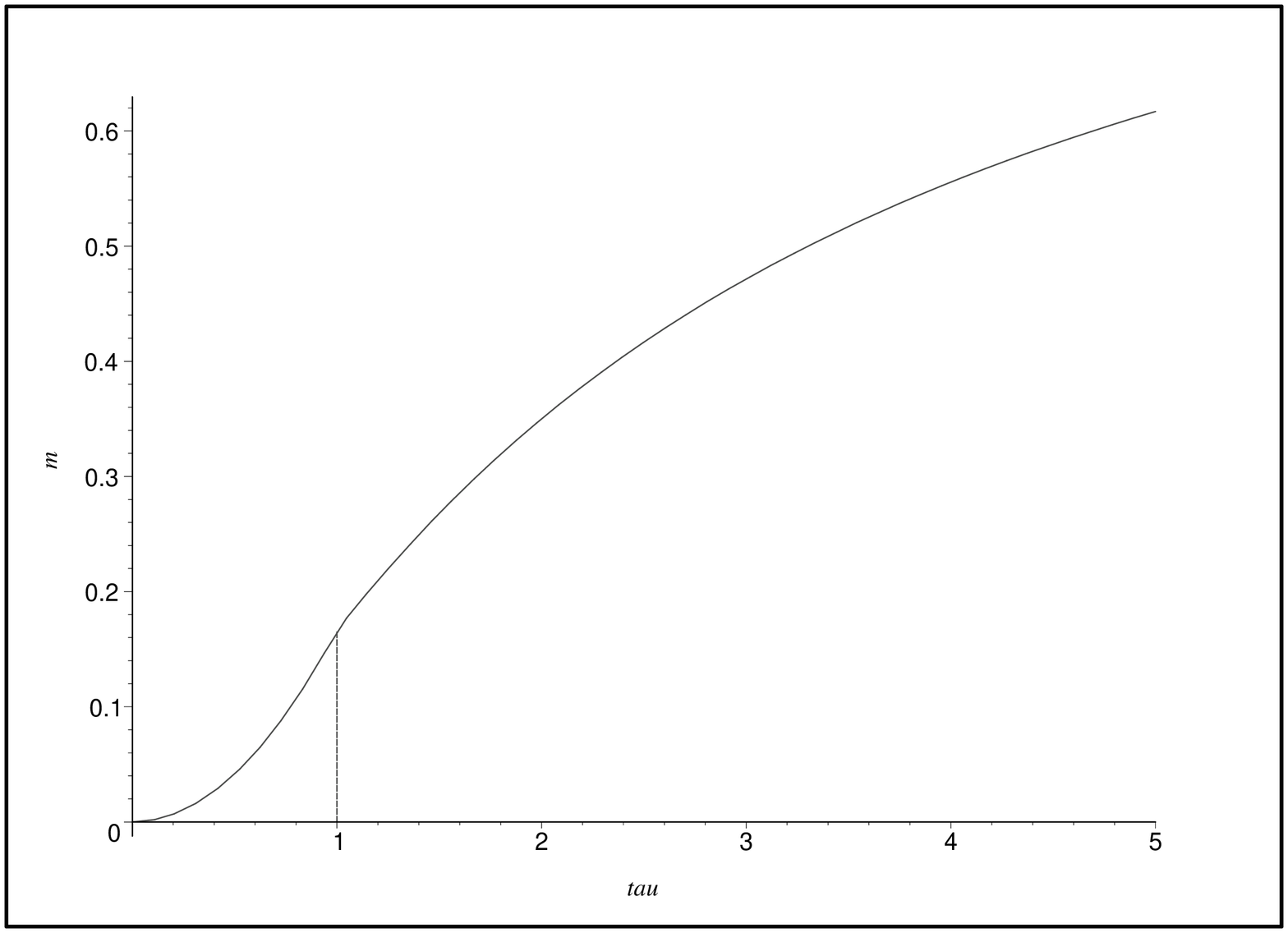}}}
\caption{
\label{fig:AsyMeans}
Asymptotic null means (i.e., arc probabilities)
$\mu_{_{PE}}(r)$ (left) and $\mu_{_{CS}}(\tau)$ (right)
as a function of the expansion parameters
from Theorems \ref{thm:asy-norm-PE}
and \ref{thm:asy-norm-CS}, respectively.
The vertical lines indicate the endpoints of the intervals in the piecewise definition of the functions.
Notice that the vertical axes are differently scaled.
}
\end{figure}

\begin{figure}[h]
\centering
\psfrag{mu}{{\Huge arc probability}}
\psfrag{t}{{\Huge expansion parameter}}
\rotatebox{-90}{ \resizebox{2. in}{!}{\includegraphics{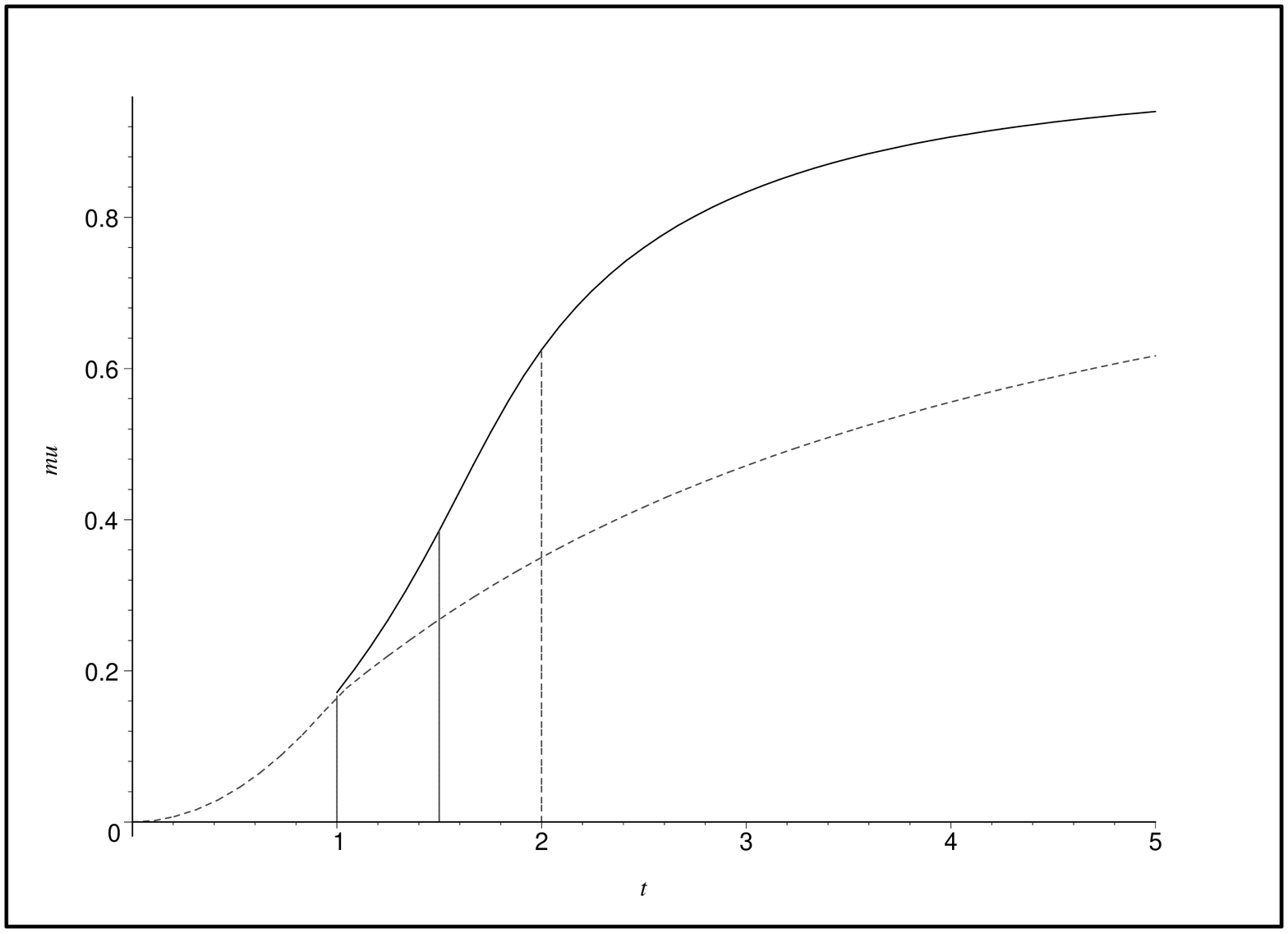}}}
\psfrag{cov}{{\Huge asymptotic variance}}
\rotatebox{-90}{ \resizebox{2. in}{!}{\includegraphics{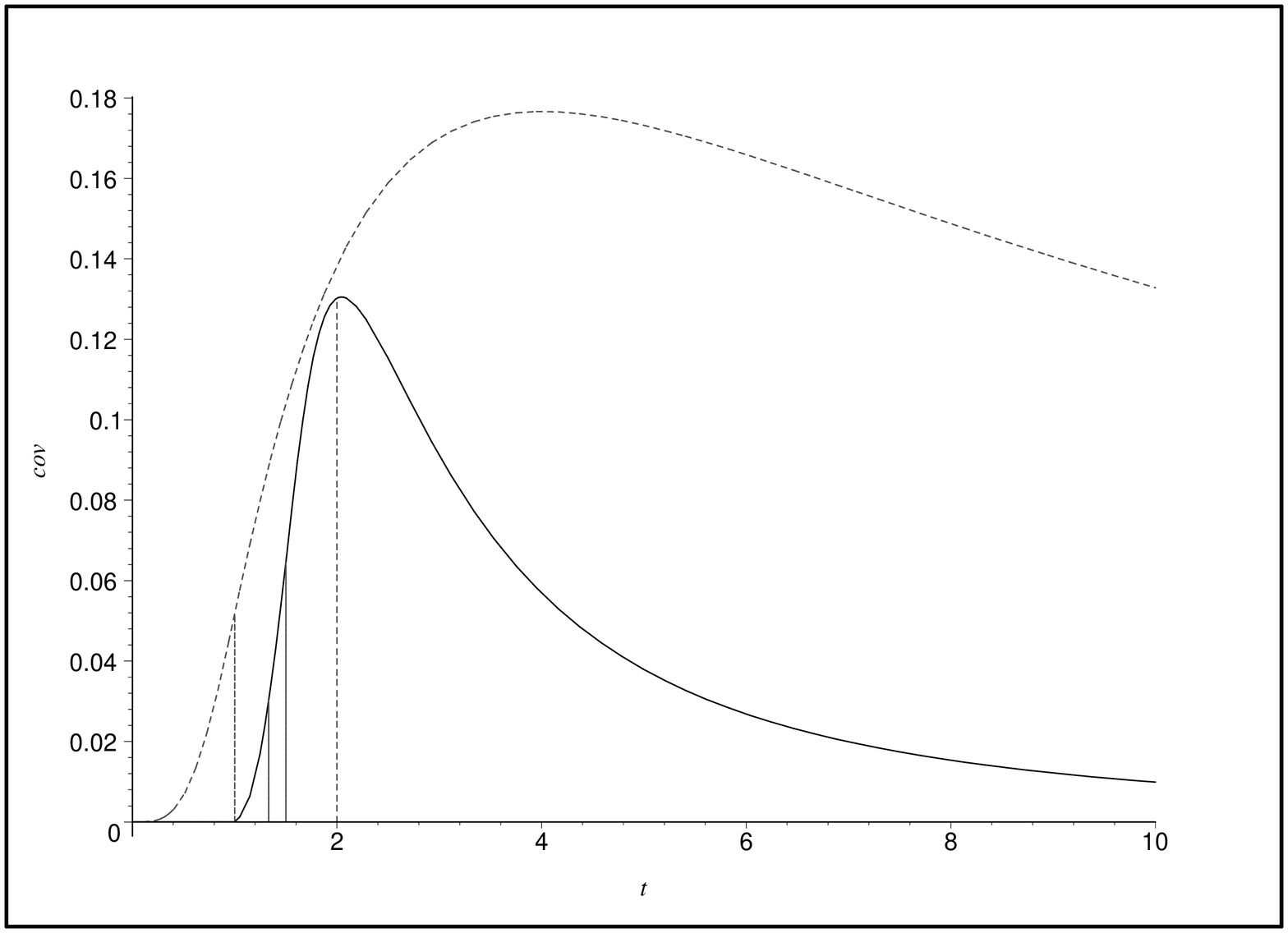}}}
\caption{
\label{fig:AsyMeansVars-Comp}
Asymptotic null means (i.e., arc probabilities) (left) and variances (right)
as a function of the expansion parameters
for relative density
of proportional-edge PCDs (solid line) and central similarity PCDs (dashed line).
The vertical lines indicate the endpoints of the intervals in the piecewise definition of the functions.
Notice that the vertical axes are differently scaled.
}
\end{figure}

Consider the forms of the mean functions,
which are depicted in Figure \ref{fig:AsyMeans}.
Note that $\mu_{_{PE}}(r)$ is monotonically increasing in $r$, since $\NPE(x,r)$
increases with $r$ for all $x \in R_V(\y_j) \setminus \RS(\NPE(\cdot,r),M_C)$,
where $\RS(\NPE(\cdot,r),M_C):=\{x \in \TY:\;\NPE(x,r)=\TY\}$.
In addition,
$\mu_{_{PE}}(r) \rightarrow 1$ as $r \rightarrow \infty$
(at rate $O\left( r^{-2} \right)$),
since the digraph becomes complete asymptotically,
which explains why $\rho_{_{PE}}(n,r)$ becomes degenerate, i.e., $\nu_{_{PE}}(r=\infty)=0$.
$\mu_{_{PE}}(r)$ is continuous, with the value at $r=1$,
$\mu_{_{PE}}(1) = 37/216 \approx .1713$.
Note also that $\mu_{_{CS}}(\tau)$ is monotonically increasing in $\tau$,
since $\NCS(x,\tau)$ increases with $\tau$ for all
$x \in R_E(e_j) \setminus \RS(\NCS(\cdot,\tau),M_C)$,
where $\RS(\NCS(\cdot,\tau),M_C):=\{x \in \TY:\;\NCS(x,\tau)=\TY\}$.
Note also that $\mu_{_{CS}}(\tau)$ is continuous in $\tau$
with $\mu_{_{CS}}(\tau=1)=1/6$ and
$\lim_{\tau \rightarrow 0}\mu_{_{CS}}(\tau)=0$.
In addition,
$\mu_{_{CS}}(\tau) \rightarrow 1$ as $\tau \rightarrow \infty$
(at rate $O\left(\tau^{-1} \right)$),
so $\rho_{_{CS}}(n,\tau)$ becomes degenerate as $\tau \rightarrow \infty$.
The asymptotic means
$\mu_{_{PE}}(r)$ and $\mu_{_{CS}}(\tau)$
are plotted together in Figure \ref{fig:AsyMeansVars-Comp} (left).
Observe that $\mu_{_{PE}}(r) > \mu_{_{CS}}(\tau)$
for all $r \in [1,\infty)$ and $\tau \in (0,\infty)$.

\begin{figure}[h]
\centering
\psfrag{nu(r)}{{\Huge$\nu_{_{PE}}(r)$}}
\psfrag{r}{{\Huge$r$}}
\psfrag{tau}{{\Huge$\tau$}}
\rotatebox{-90}{ \resizebox{2. in}{!}{\includegraphics{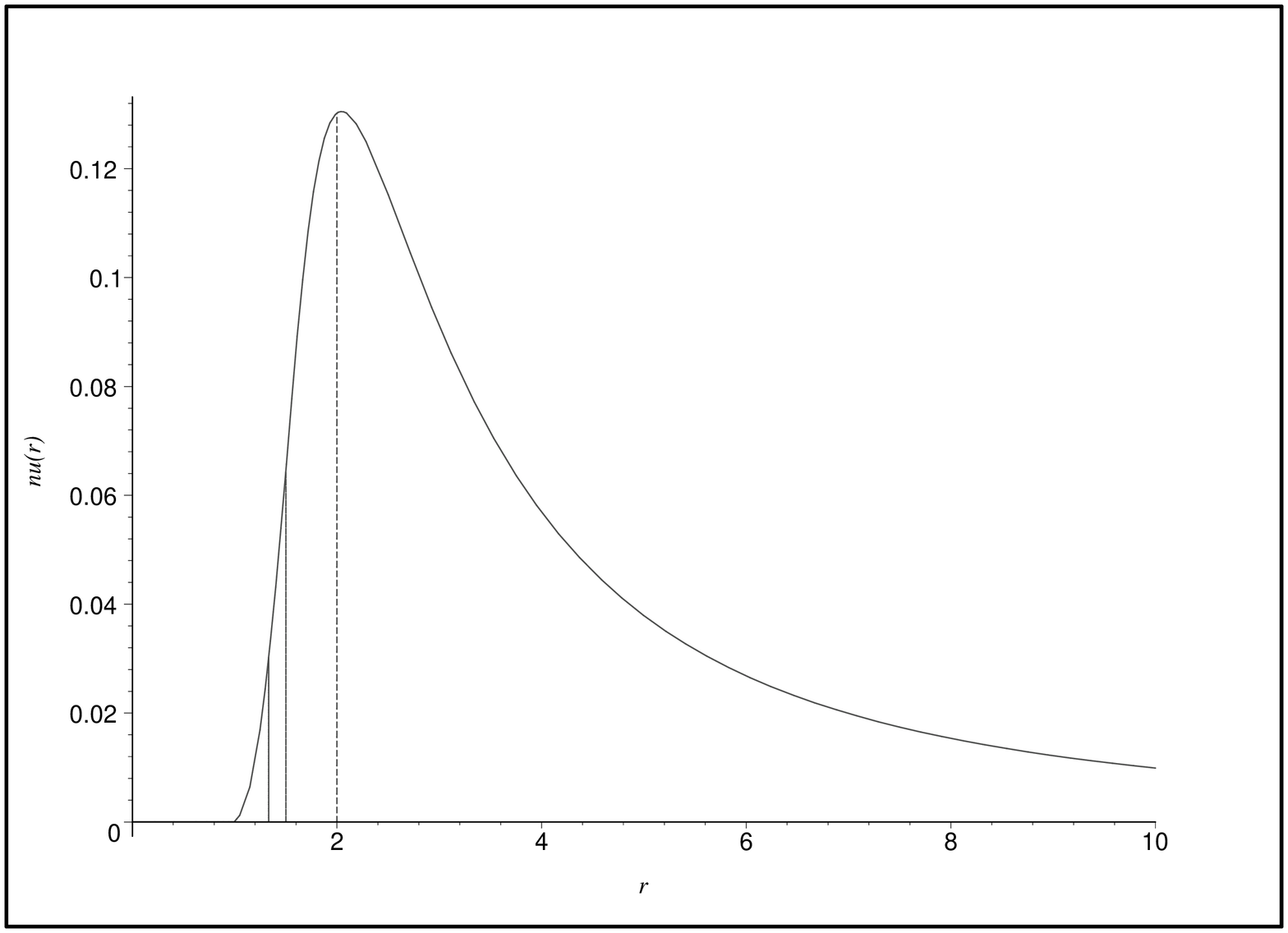}}}
\psfrag{cov}{{\Huge$\nu_{_{CS}}(\tau)$}}
\rotatebox{-90}{ \resizebox{2. in}{!}{\includegraphics{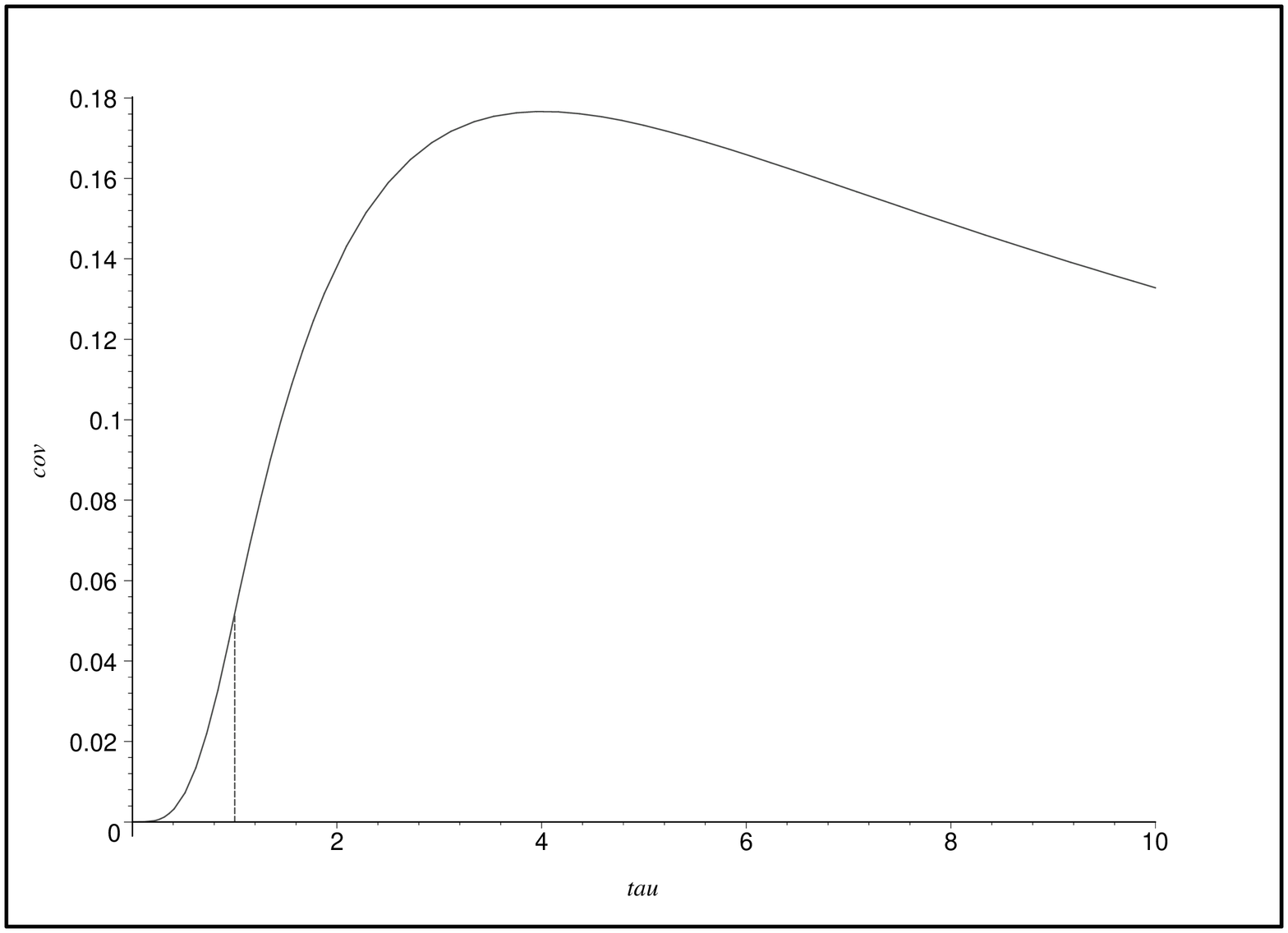}}}
\caption{
\label{fig:AsyVars}
Asymptotic null variances $\nu_{_{PE}}(r)$ (left) and $\nu_{_{CS}}(\tau)$ (right)
as a function of the expansion parameters
from Theorems \ref{thm:asy-norm-PE}
and \ref{thm:asy-norm-CS}, respectively.
The vertical lines indicate the endpoints of the intervals in the piecewise definition of the functions.
Notice that the vertical axes are differently scaled.
}
\end{figure}

The asymptotic variance functions
are depicted in Figure \ref{fig:AsyVars}.
Note that $\nu_{_{PE}}(r)$ is also
continuous in $r$ with $\lim_{r \rightarrow \infty} \nu_{_{PE}}(r) = 0$ and
$\nu_{_{PE}}(1) = 34/58320 \approx .000583$ and observe that
$\sup_{r \ge 1}\nu_{_{PE}}(r) \approx .1305$ which is attained at $r \approx 2.045$.
Note also that $\nu_{_{CS}}(\tau)$ is
continuous in $\tau$ with $\lim_{\tau \rightarrow \infty} \nu_{_{CS}}(\tau) = 0$
and $\nu(\tau=1)=7/135$ and $\lim_{\tau \rightarrow 0}\nu_{_{CS}}(\tau)=0$
---there are no arcs when $\tau=0$ a.s.--- which explains why
$\rho_n(\tau=0)$ is degenerate.
Moreover,
$\sup_{\tau > 0}\nu_{_{CS}}(\tau) \approx .1767$ which is attained at $\tau \approx 4.0051$.
The asymptotic variances $\nu_{_{PE}}(r)$ and $\nu_{_{CS}}(\tau)$
are plotted together in Figure \ref{fig:AsyMeansVars-Comp} (right).
Observe that $\nu_{_{CS}}(\tau) > \nu_{_{PE}}(r)$
for all $r \in [1,\infty)$ and $\tau \in (0,\infty)$.

To illustrate the limiting distribution, $r=2$ yields
$
\rho_{_{PE}}(n,2) \stackrel{\text{\scriptsize approx}}{\sim} \N\left(\frac{5}{8},\frac{25}{192n}\right)
$
or equivalently,
$$
\frac{\sqrt{n} \bigl( \rho_{_{PE}}(n,2) - \mu_{_{PE}}(2)\bigr)}{\sqrt{\nu_{_{PE}}(2)}}
=
\sqrt{\frac{192n}{25}} \left(\rho_{_{PE}}(n,2) - \frac{5}{8}\right)\stackrel{\mathcal{L}}{\longrightarrow} \N(0,1)$$
where $\stackrel{\mathcal L}{\longrightarrow}$
stands for convergence in law or distribution.

Similarly, $\tau=1$ yields
$
\rho_{_{CS}}(n,1) \stackrel{\text{approx}}{\sim} \N\left(\frac{1}{6},\frac{7}{135\,n}\right)
$
or equivalently,
$$
\frac{\sqrt{n} \bigl( \rho_{_{CS}}(n,1) - \mu_{_{CS}}(1) \bigr)}{\sqrt{\nu_{_{CS}}(1)}}
=
\sqrt{\frac{135\,n}{7}} \left(\rho_{_{CS}}(n,1) - \frac{1}{6}\right)\stackrel{\mathcal{L}}{\longrightarrow} \mathcal{N}(0,1).$$

The finite sample variance and skewness of $\rho_{_{PE}}(n,r)$ and $\rho_{_{CS}}(n,\tau)$
may be derived analytically in much the same way as was asymptotic variances.
In particular, the variance of $h_{12}$ for proportional-edge PCD is
\begin{multline*}
\omega_{_{PE}}(r)=\Var[h_{12}]= \omega^{1,1}_{_{PE}}(r)\,\I(r \in [1,4/3)) +\\
 \omega^{1,2}_{_{PE}}(r)\,\I(r \in [4/3,3/2))+\omega^{1,3}_{_{PE}}(r)\,\I(r \in [3/2,2))+\omega^{1,4}_{_{PE}}(r)\,\I(r \in [2,\infty))
\end{multline*}
where
{\small
\begin{align*}
\omega^{1,1}_{_{PE}}(r)&=\frac{-(1369\,r^8+4107\,r^7+902\,r^6-78084\,r^5+161784\,r^4-182736\,r^3-23328\,r^2+155520\,r-55296)}{11664\,(r+2)(r+1)r^2},\\
\omega^{1,2}_{_{PE}}(r)&=-\frac{1369\,r^7+4107\,r^6+9650\,r^5-98496\,r^4+132624\,r^3-79056\,r^2-57888\,r+72576}{11664\,(r+2)(r+1)r},\\
\omega^{1,3}_{_{PE}}(r)&=-\frac{r^{10}+3\,r^9-62\,r^8+968\,r^6-1704\,r^5-1824\,r^4+5424\,r^3-1168\,r^2-3856\,r+2208}{16\,(r+2)(r+1)r^4},\\
\omega^{1,4}_{_{PE}}(r)&=\frac{3\,r^3+3\,r^2+3\,r-13}{r^4(r+1)}.
\end{align*}
}

\begin{figure}[h]
\centering
\psfrag{var}{\huge{$\omega_{_{PE}}(r)$}}
\psfrag{r}{\huge{$r$}}
\rotatebox{-90}{ \resizebox{2. in}{!}{\includegraphics{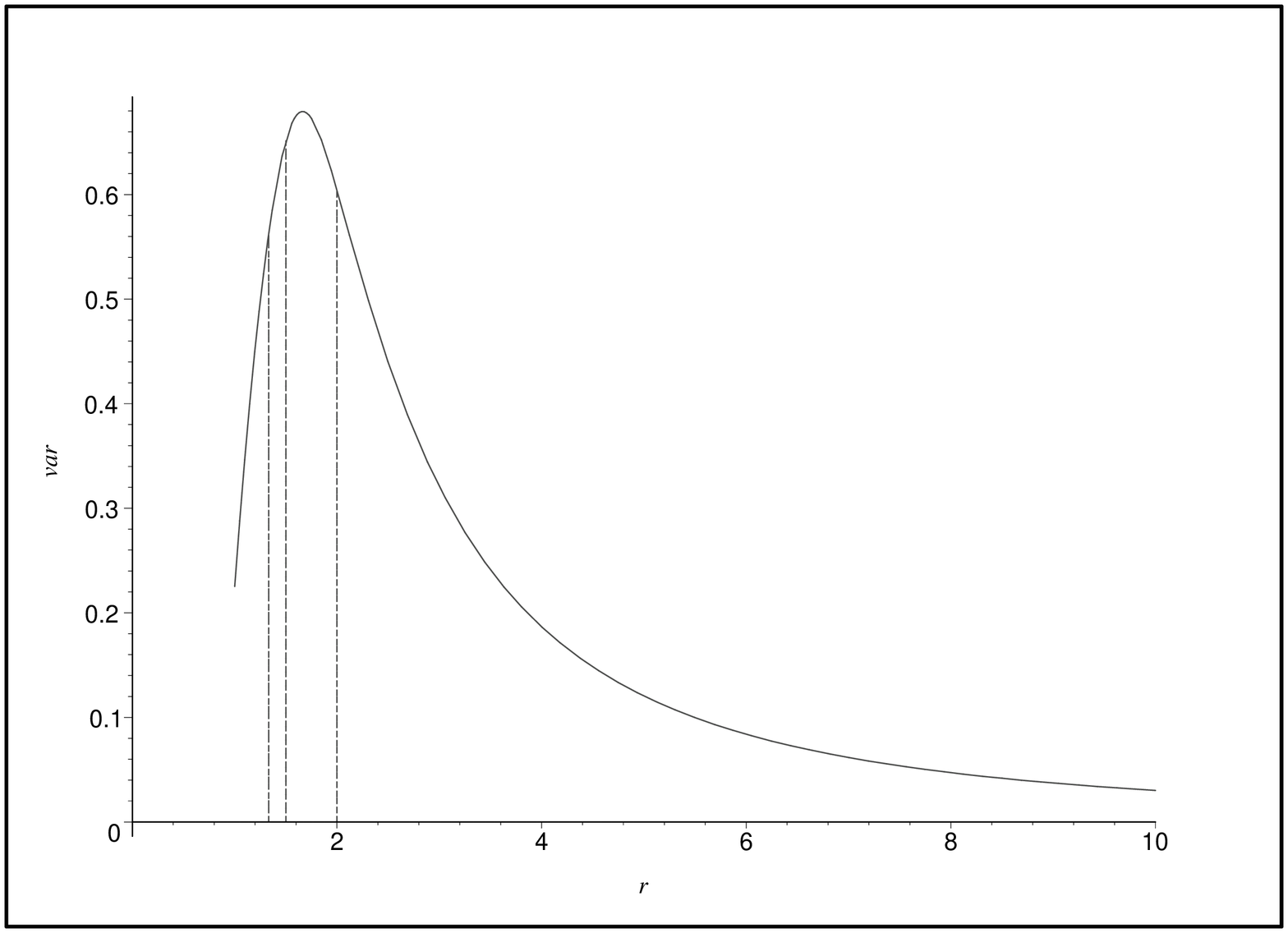}}}
\psfrag{var}{\huge{$\omega_{_{CS}}(\tau)$}}
\psfrag{tau}{\huge{$\tau$}}
\rotatebox{-90}{ \resizebox{2. in}{!}{\includegraphics{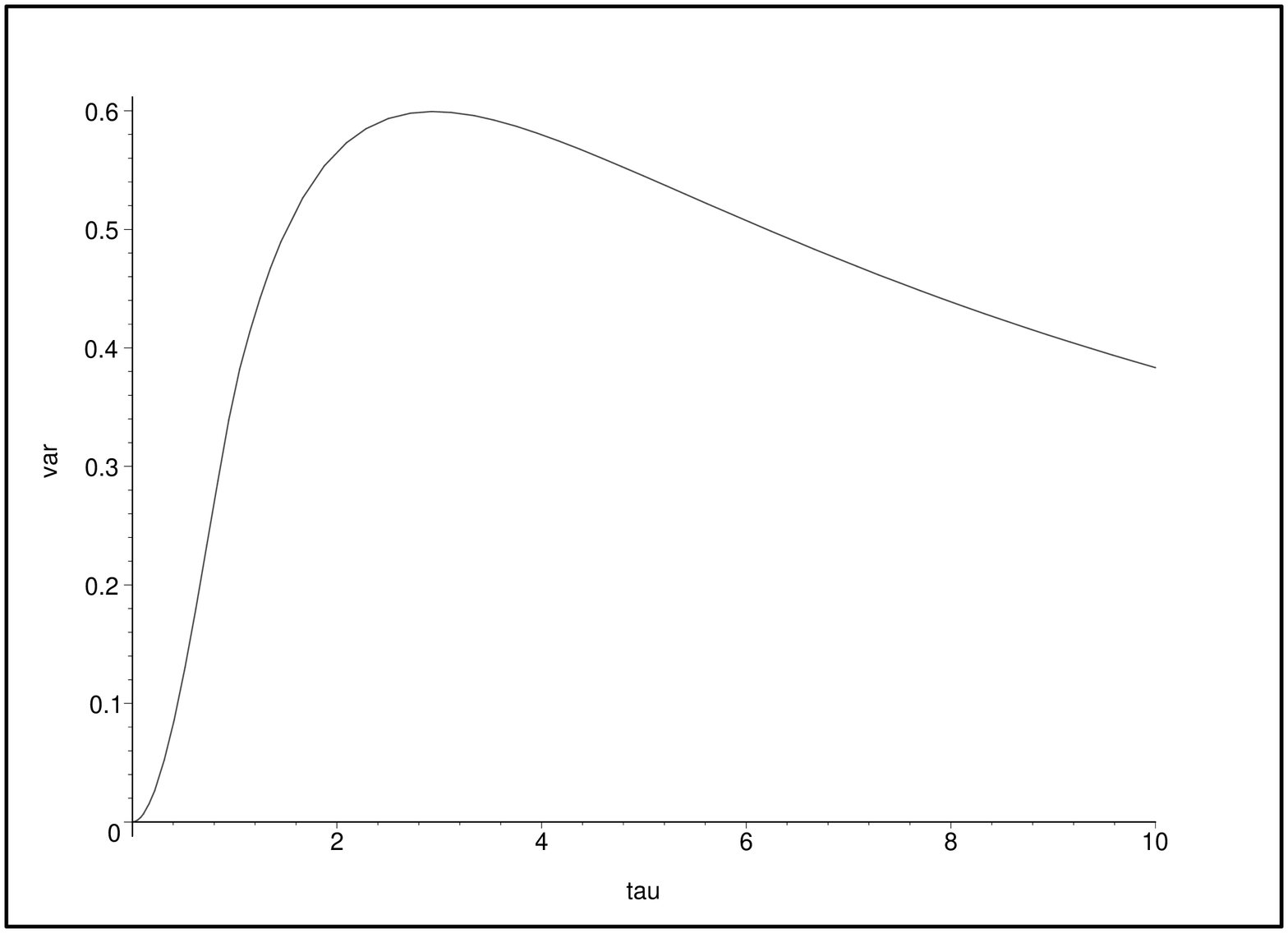}}}
\caption{
\label{fig:Var(h12)}
$\Var[h_{12}]=\omega_{_{PE}}(r)$ as a function of $r$ for $r \in [1,10]$ (left)
and
$\Var[h_{12}]=\omega_{_{CS}}(\tau)$ as a function of $\tau \in (0,10]$ (right).}
\end{figure}

In Figure \ref{fig:Var(h12)} (left) is the graph of $\omega_{_{PE}}(r)$ for $r \in [1,10]$.
Note that $\omega(r=1)=2627/11664\approx .2252$ and
$\lim_{r \rightarrow \infty}\omega_{_{PE}}(r)=0$ (at rate $O\left( r^{-2}\right)$),
$\sup_{r \in [1,\infty)} \omega_{_{PE}}(r) \approx .6796$
which is attained at $r \approx 1.66$.

\begin{figure}[h]
\centering
\psfrag{var}{{\Huge$\Var[h_{12}]$}}
\psfrag{t}{{\Huge expansion parameter}}
\rotatebox{-90}{ \resizebox{2. in}{!}{\includegraphics{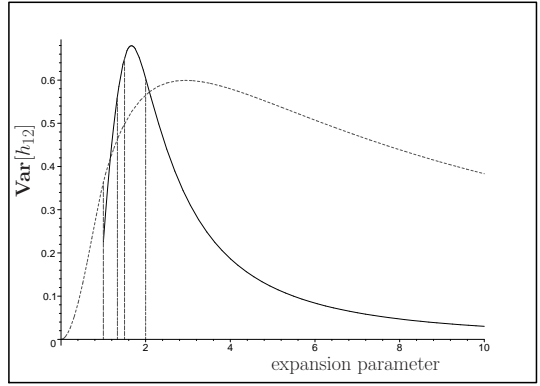}}}
\caption{
\label{fig:Varh12-Comp}
$\Var[h_{12}]$ terms for relative density
of proportional-edge PCDs (solid line) and central similarity PCDs (dashed line)
as a function of the expansion parameters.
The vertical lines indicate the endpoints of the intervals in the piecewise definition of the functions.
}
\end{figure}

Moreover, the variance of $h_{12}$ for central similarity PCDs is
\begin{equation*}
\omega_{_{CS}}(\tau) =
 \begin{cases}
  {\frac{ -\left( \tau^3+7\,\tau^2-5\,\tau-15 \right) \tau^2}{9\,\tau+3}}  &\text{for} \quad \tau \in [0,1/2), \\
  {\frac{ -\left( 2\,\tau^4+11\,\tau^3+9\,\tau^2-33\,\tau-81 \right) \tau^2}{ 9\,\left( \tau+3 \right)  \left( 2\,\tau+5 \right) }}  &\text{for} \quad \tau \in [1/2,1), \\
  {\frac{ 2\,\left( 22\,\tau^4+151\,\tau^3+244\,\tau^2+12\, \tau-15 \right) \tau}{ \left( \tau+2 \right)^2 \left( 2\,\tau+1
  \right)^2 \left( \tau+3 \right)  \left( 2\,\tau+5 \right) }}  &\text{for} \quad \tau \in [1,\infty).
 \end{cases}
\end{equation*}

In Figure \ref{fig:Var(h12)} (right) is the graph of $\omega_{_{CS}}(\tau)$ for $\tau \in [1,10]$.
Note that $\omega_{_{CS}}(\tau)$ is a continuous
function of $\tau$
with $\lim_{\tau \rightarrow 0}\omega(\tau)=0$
and $\omega(\tau=1)=23/63 \approx .3651$.
Furthermore,
$\lim_{\tau \rightarrow \infty}\omega_{_{CS}}(\tau)=0$ (at rate $O\left( \tau^{-2}\right)$),
$\sup_{\tau \in (0,\infty)} \omega_{_{CS}}(\tau) \approx .60$
which is attained at $\tau \approx 2.94$.
The variances $\Var[h_{12}]$, $\omega_{_{PE}}(r)$ and $\omega_{_{CS}}(\tau)$
are plotted together in Figure \ref{fig:Varh12-Comp}.
Observe that $\omega_{_{CS}}(t) > \omega_{_{PE}}(t)$
for $ 1 \le t \lesssim 1.165)$ and $t \gtrsim 2.09$;
and
$\omega_{_{PE}}(t) > \omega_{_{CS}}(t)$
for $1.165 \lesssim t \lesssim 2.09$

In fact,
the exact distribution of $\rho_{_{PE}}(n,r)$
is, in principle, available
by successively conditioning on the values of $X_i$.
Alas,
while the joint distribution of $h_{12},h_{13}$ is available,
the joint distribution of $\{h_{ij}\}_{1 \leq i < j \leq n}$,
and hence the calculation for the exact distribution of $\rho_{_{PE}}(n,r)$,
is extraordinarily tedious and lengthy for even small values of $n$.
The same holds for the the exact distribution of $\rho_{_{CS}}(n,\tau)$.

\begin{figure}[]
\centering
\psfrag{Density}{ \Huge{\bfseries{density}}}
\rotatebox{-90}{ \resizebox{1.84 in}{!}{ \includegraphics{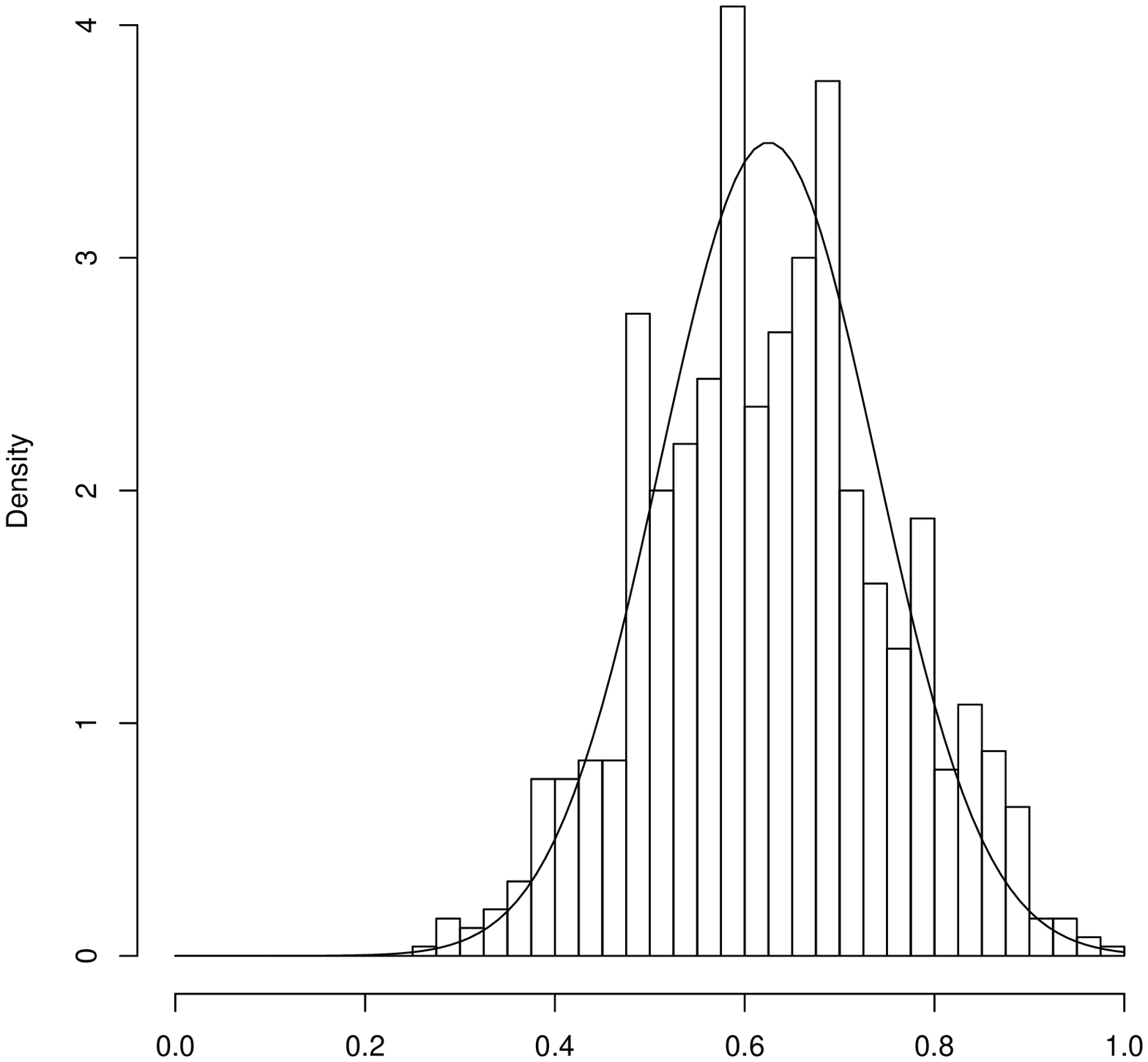}}}
\rotatebox{-90}{ \resizebox{1.84 in}{!}{ \includegraphics{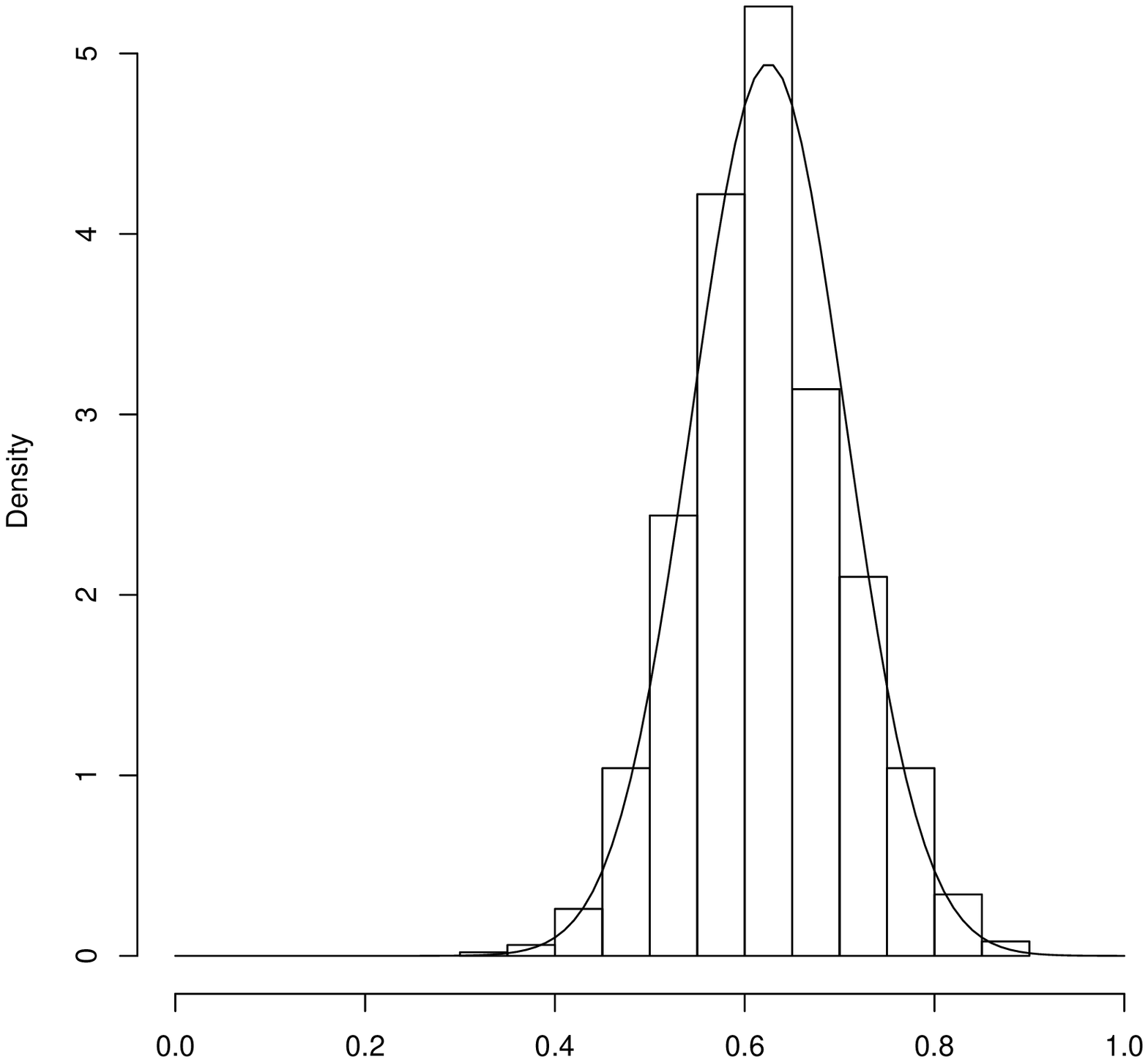}}}
\rotatebox{-90}{ \resizebox{1.84 in}{!}{ \includegraphics{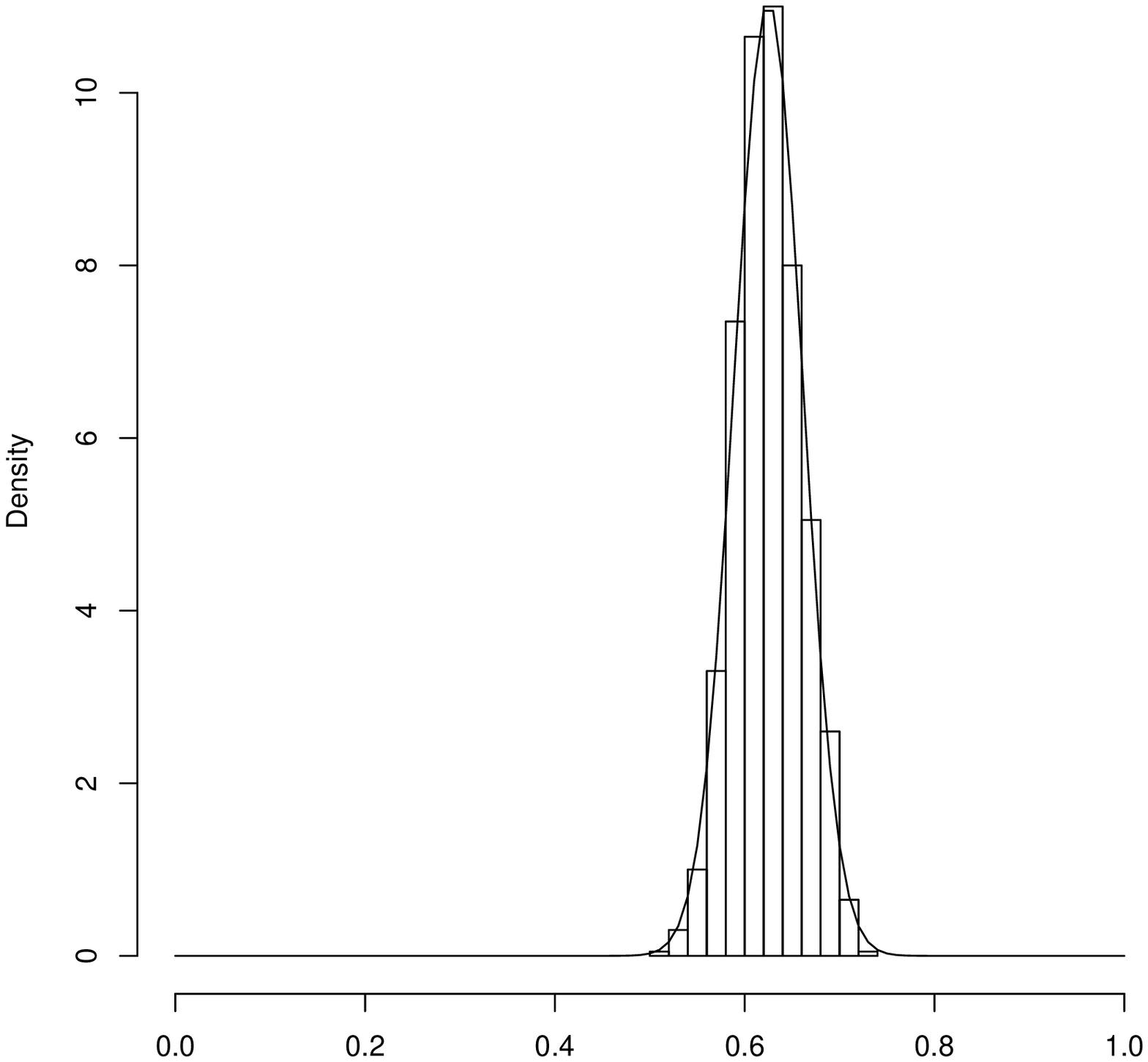}}}
\caption{
\label{fig:NormSkew}
Depicted are the distributions of
$\rho_{_{PE}}(n,2) \stackrel{\text{\scriptsize approx}}{\sim} \N\left(\frac{5}{8},\frac{25}{192n}\right)$
for $n=10,20,100$ (left to right).
Histograms are based on 1000 Monte Carlo replicates.
Solid curves represent the approximating normal densities given in Theorem \ref{thm:asy-norm-PE}.
Note that the vertical axes are differently scaled.
}
\end{figure}

\begin{figure}[]
\centering
\psfrag{Density}{ \Huge{\bfseries{density}}}
\rotatebox{-90}{ \resizebox{2.1 in}{!}{ \includegraphics{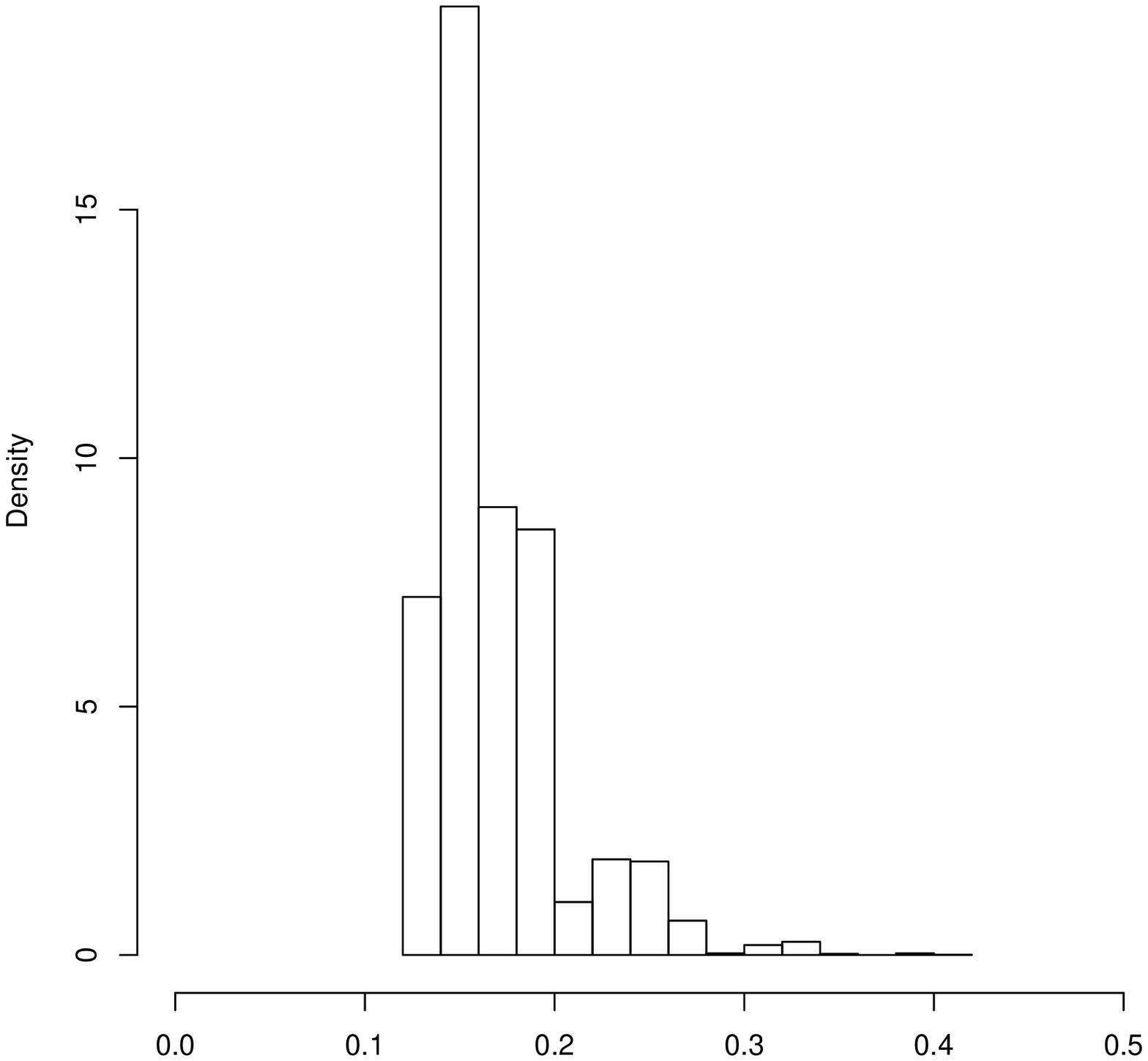}}}
\rotatebox{-90}{ \resizebox{2.1 in}{!}{ \includegraphics{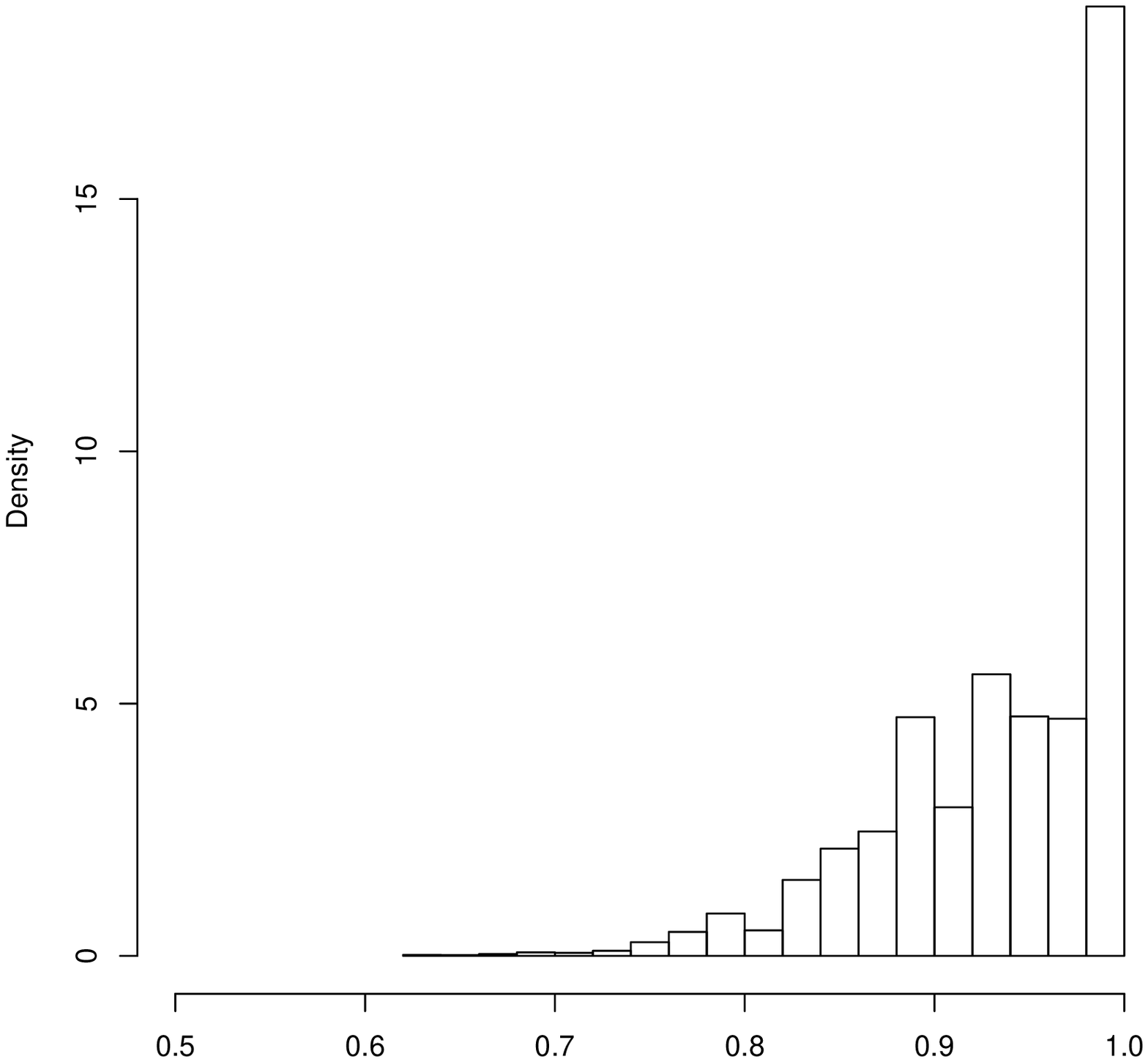}}}
\caption{
\label{fig:NormSkew1}
Depicted are the histograms of relative density for 10000 Monte Carlo replicates of
$\rho_{_{PE}}(10,1)$ (left) and $\rho_{_{PE}}(10,5)$ (right)
indicating severe small sample skewness for extreme values of $r$.}
\end{figure}

Figure \ref{fig:NormSkew}
indicates that, for $r=2$,
the normal approximation for the relative density of proportional-edge PCD
is accurate even for small $n$
(although kurtosis may be indicated for $n=10$).
Figure \ref{fig:NormSkew1} demonstrates,
however, that severe skewness obtains for small values of $n$ and extreme values of $r$.

\vspace*{0.5 in}

\begin{figure}[]
\centering
\psfrag{Density}{ \Huge{\bfseries{density}}}
\rotatebox{-90}{ \resizebox{1.8 in}{!}{ \includegraphics{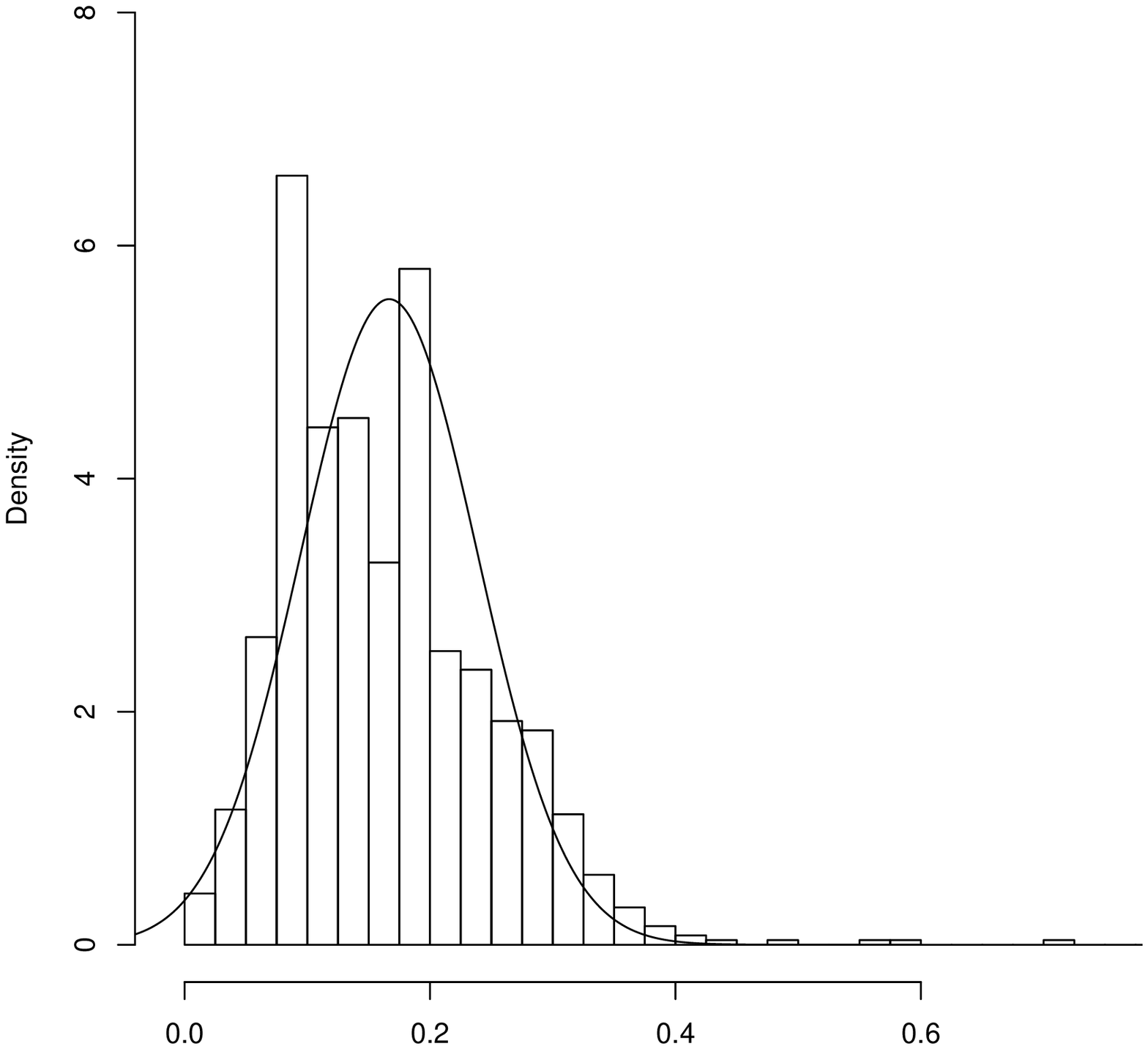} } }
\rotatebox{-90}{ \resizebox{1.8 in}{!}{ \includegraphics{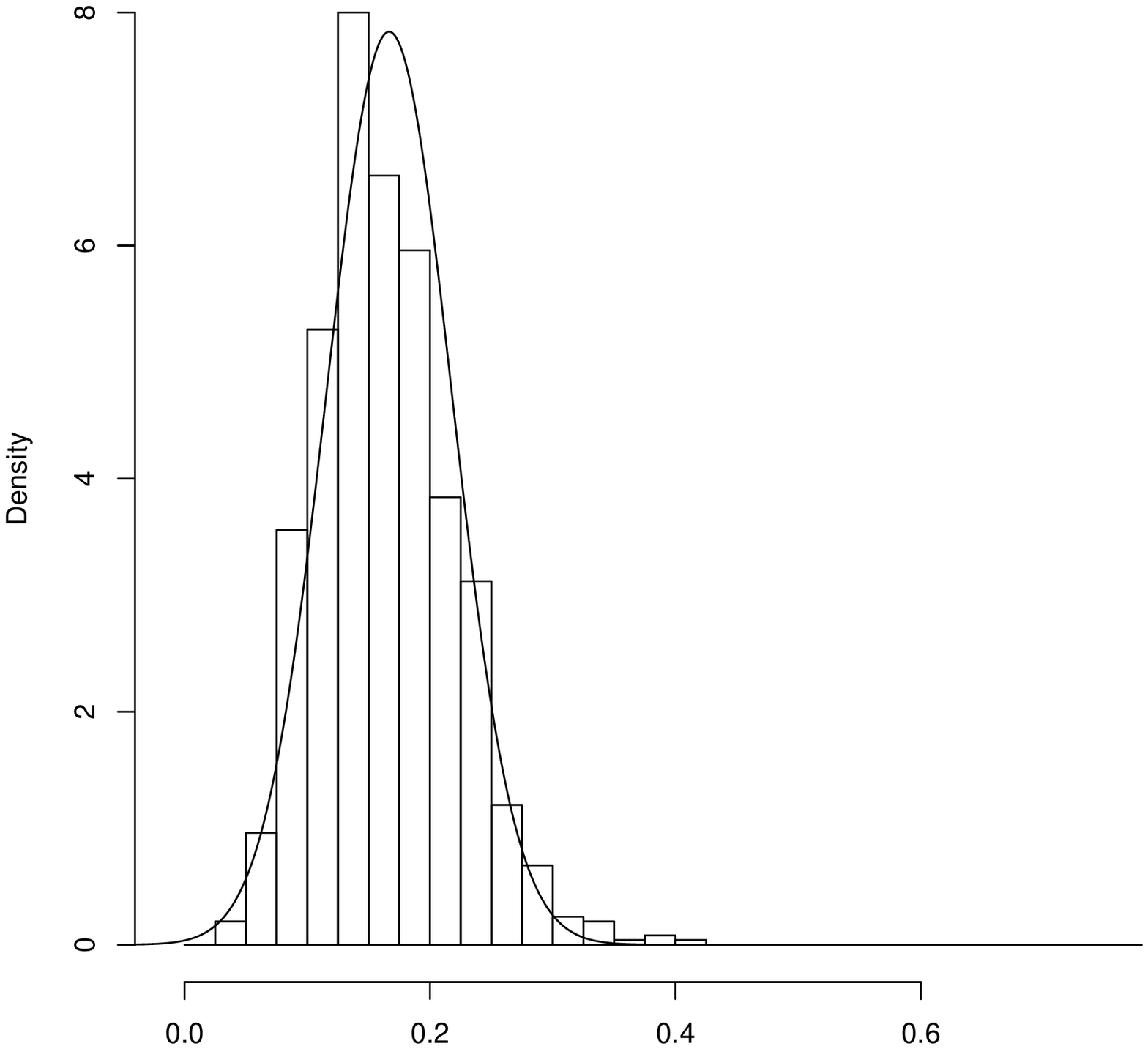} } }
\rotatebox{-90}{ \resizebox{1.8 in}{!}{ \includegraphics{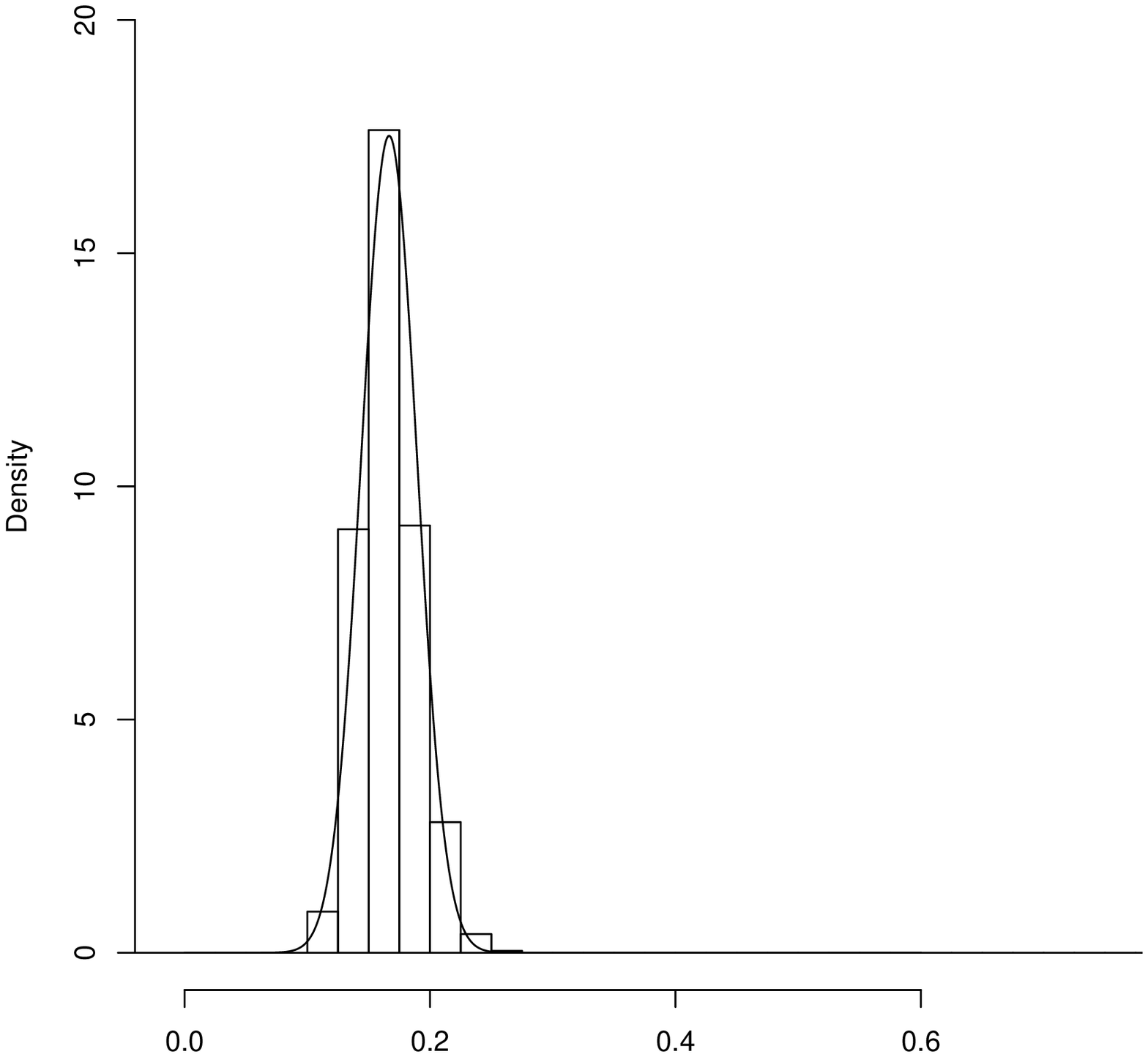} } }
\caption{
\label{fig:NormSkewCS}
Depicted are
$\rho_{_{CS}}(n,1) \stackrel{\text{\scriptsize approx}}{\sim} \mathcal{N}\left(\frac{1}{6},\frac{7}{135\,n}\right)$
for $n=10,\,20,\,100$ (left to right).
Histograms are based on 1000 Monte Carlo replicates.
Solid curves represent the approximating normal densities given in Theorem \ref{thm:asy-norm-CS}.
Note that the vertical axes are differently scaled.
}
\end{figure}

\begin{figure}[]
\centering
\psfrag{Density}{ \Huge{\bfseries{density}}}
\rotatebox{-90}{ \resizebox{1.8 in}{!}{ \includegraphics{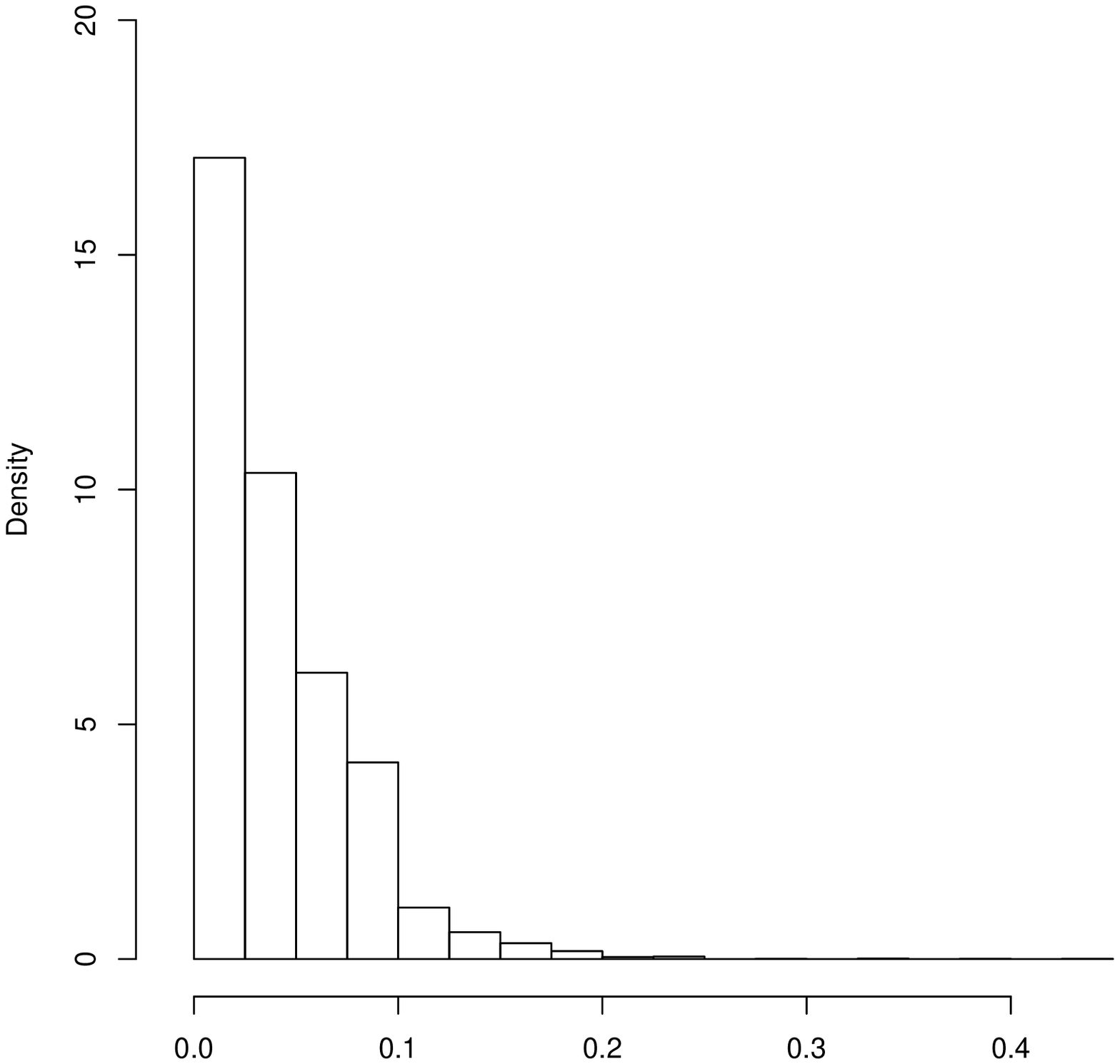} } }
\rotatebox{-90}{ \resizebox{1.8 in}{!}{ \includegraphics{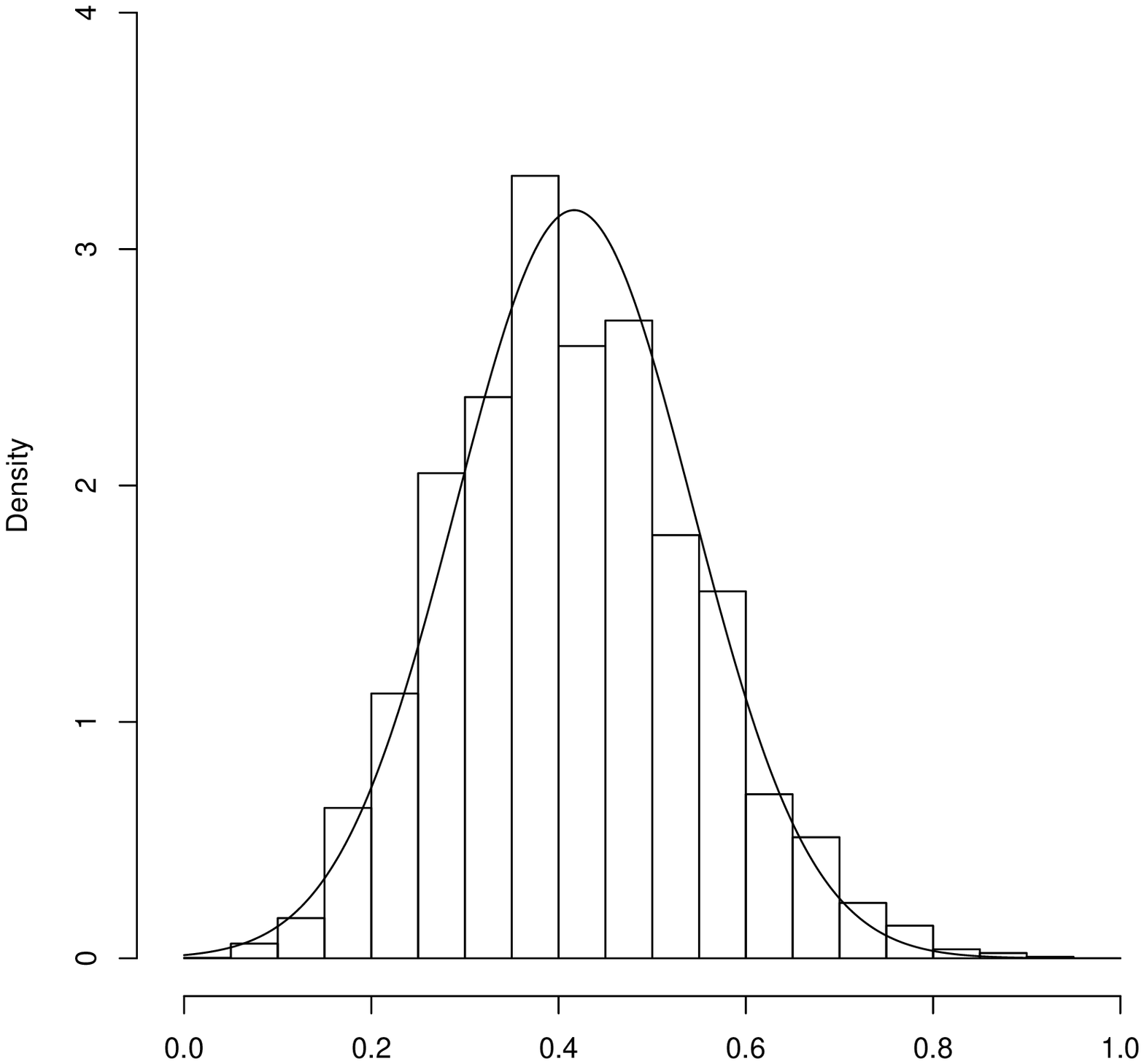} } }
\rotatebox{-90}{ \resizebox{1.8 in}{!}{ \includegraphics{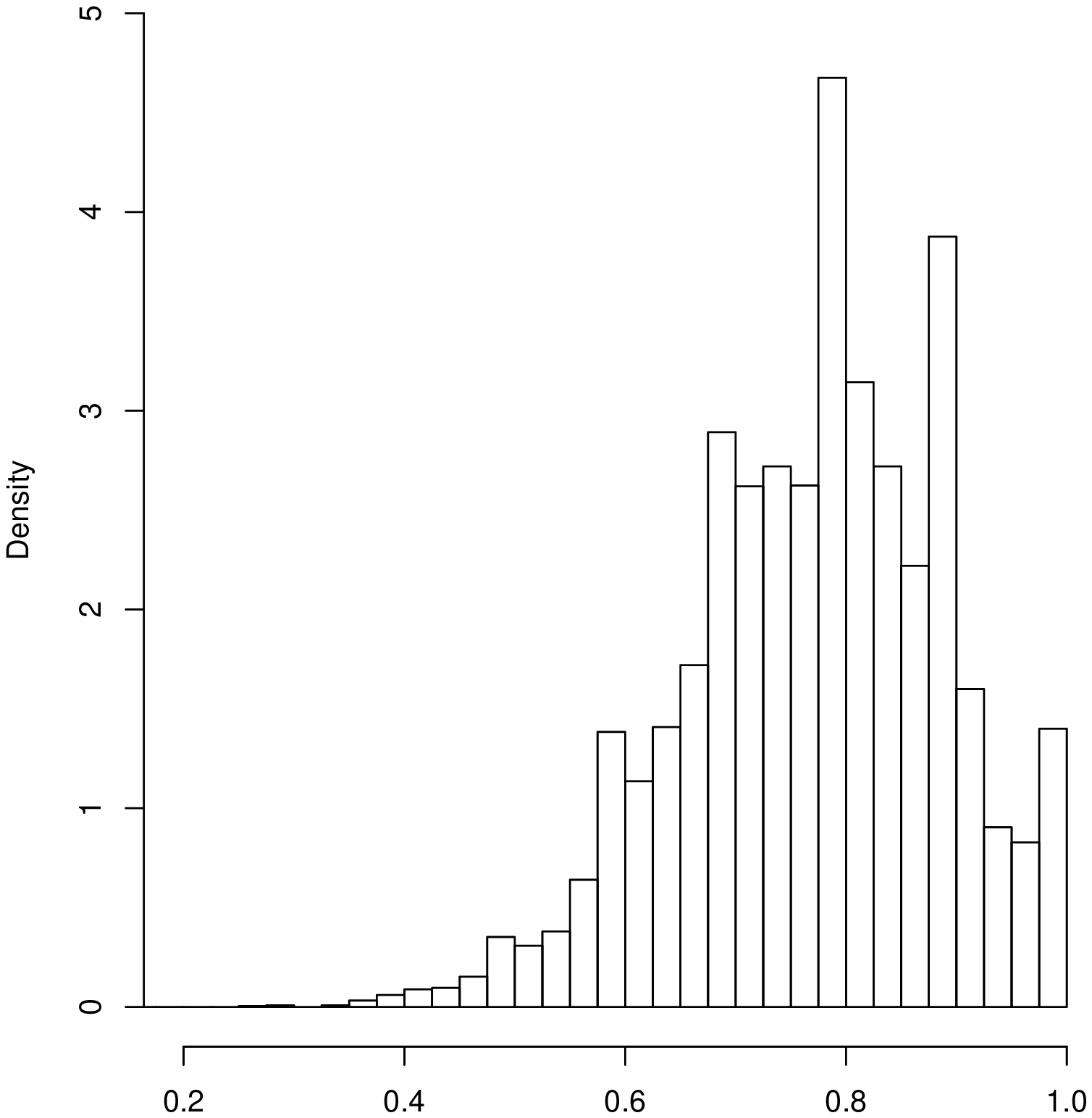} } }
\caption{
\label{fig:CSNormSkew1}
Depicted are the histograms for 10000 Monte Carlo replicates of
$\rho_{_{CS}}(10,1)$ (left), $\rho_{_{CS}}(10,2.5)$ (middle), and $\rho_{_{CS}}(10,10)$ (right)
indicating severe small sample skewness for extreme values of $\tau$
(i.e., $\tau=1$ or $\tau=10$).
}
\end{figure}

Figure \ref{fig:NormSkewCS}
indicates that, for $\tau=1$,
the normal approximation for the relative density of central similarity PCD is accurate even for small $n$
(although kurtosis and skewness may be indicated for $n=10,\,20$).
Figure \ref{fig:CSNormSkew1} demonstrates,
however, that the smaller the value of $\tau$,
the more severe the skewness of the probability density.

\subsection{The Multiple Triangle Case}
\label{sec:rho-mult-tri}
\begin{figure}[]
\centering
\epsfig{figure=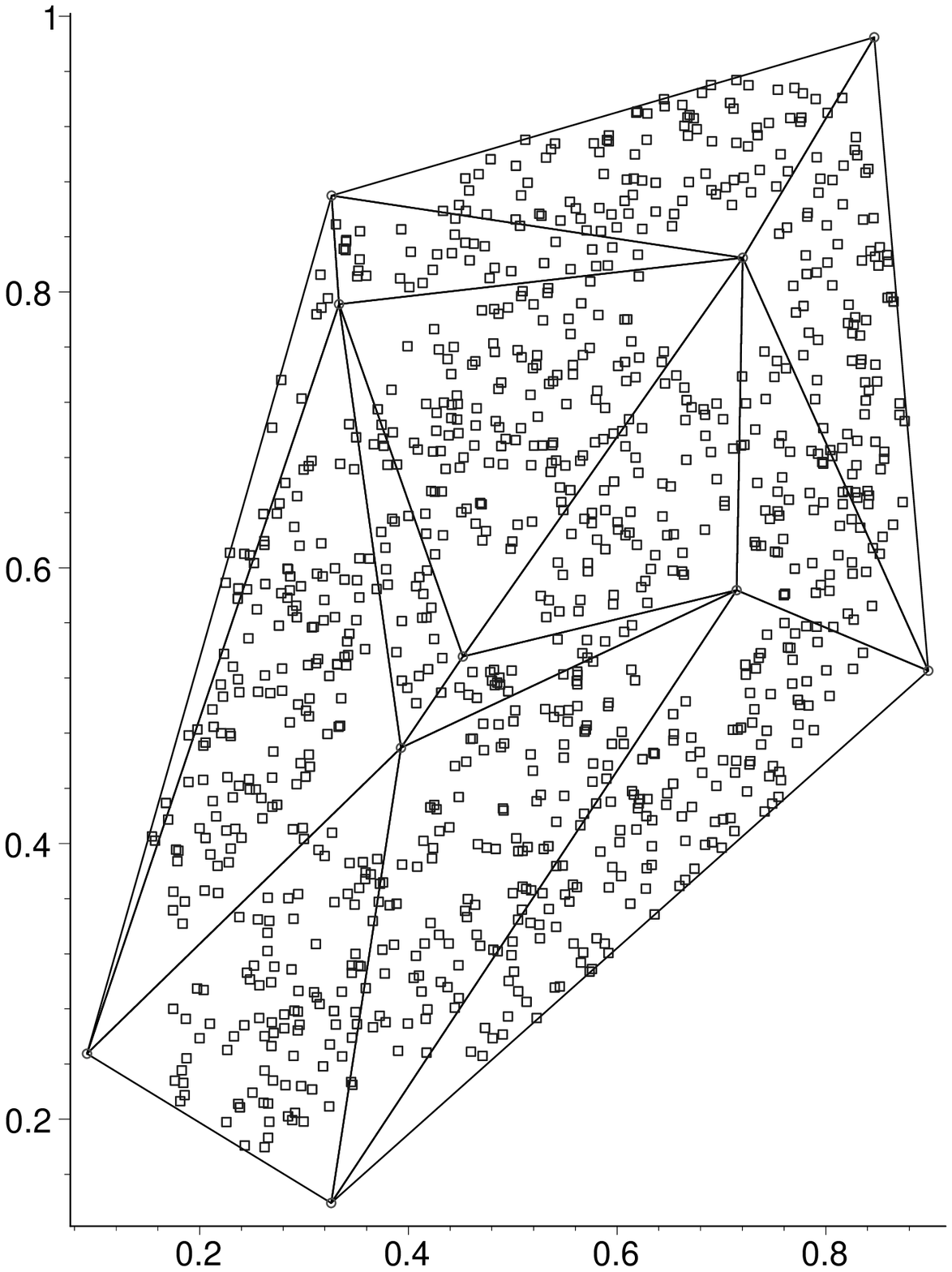, height=125pt, width=125pt}
\epsfig{figure=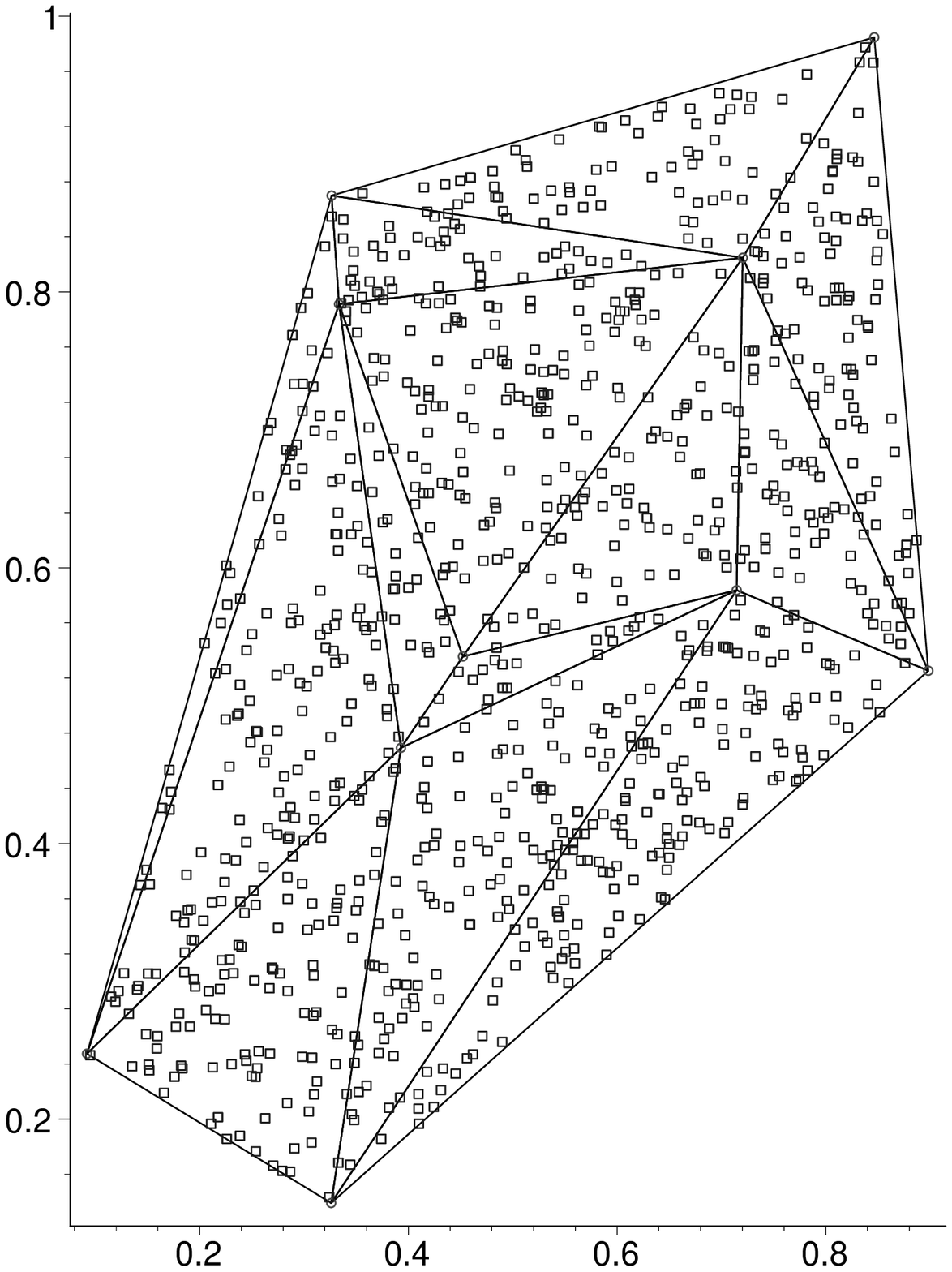, height=125pt, width=125pt}
\epsfig{figure=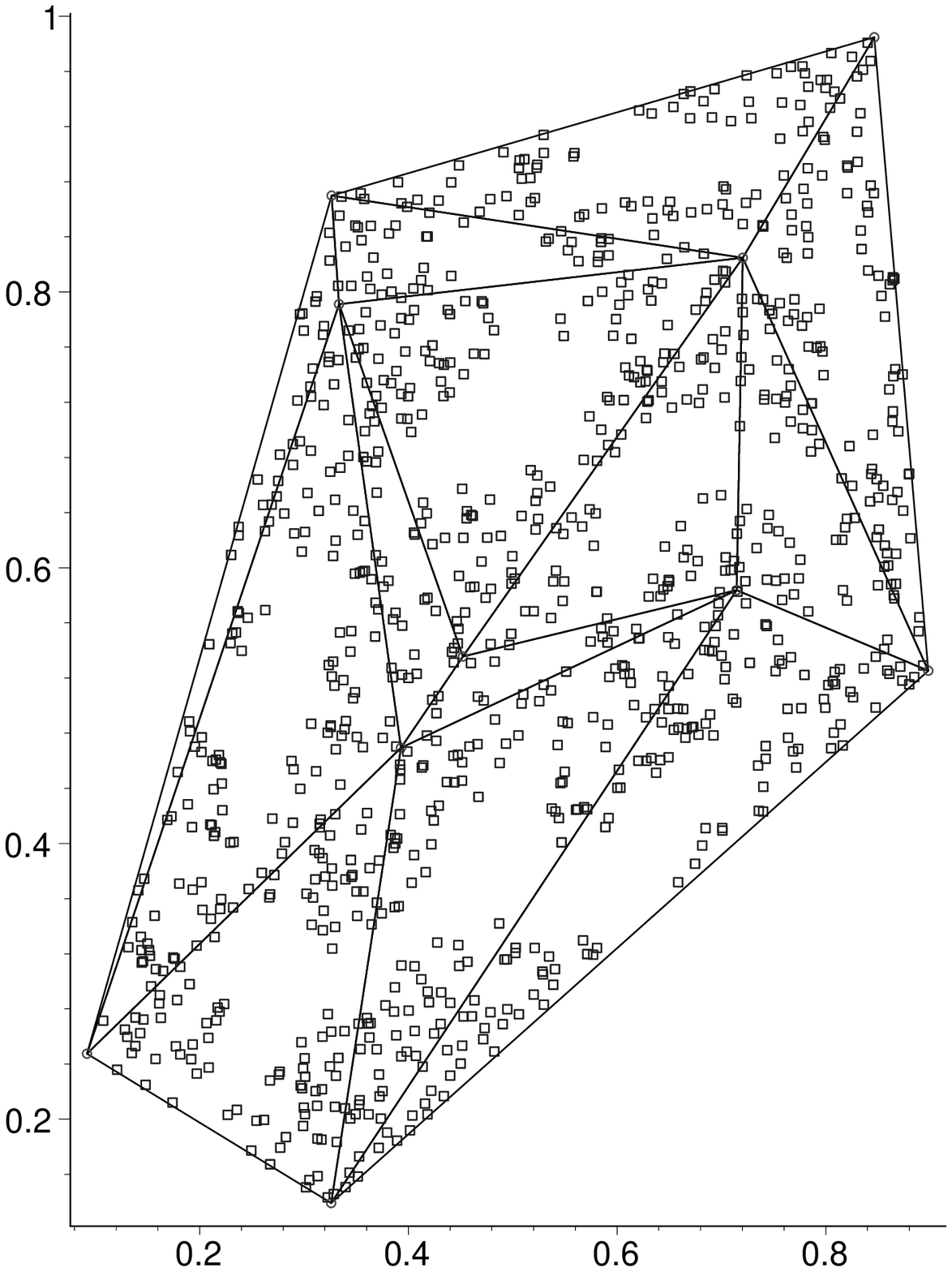, height=125pt, width=125pt}
\caption{\label{fig:deldata}
Realization of segregation (left), $H_o: CSR$ (middle), and association (right) for $|\Y_m|=10$.}
\end{figure}

In this section, we present the asymptotic
distribution of the relative density in multiple triangles.
Suppose $\Y_m=\{\y_1,\y_2,\ldots,\y_m\} \subset
\R^2$ be a set of $m$ points in general position with $m > 3$ and no
more than three points are cocircular.
As a result of the Delaunay triangulation of $\Y_m$ (\cite{okabe:2000}),
there are $J_m>1$ Delaunay triangles each of which is denoted as $T_j$.
The Delaunay triangles partition the convex hull of $\Y_m$.
We wish to investigate
\begin{equation}
\label{eqn:null-pattern-mult-tri}
H_o: X_i \stackrel{iid}{\sim} \U(C_H(\Y_m))  \text{ for } i=1,2,\ldots,n
\end{equation}
against segregation and association alternatives
(see Section \ref{sec:alternatives}).
Figure \ref{fig:deldata} (middle) presents a realization of 1000 observations
independent and identically distributed as $\U(C_H(\Y_m))$ for $m=10$ and $J_m=13$.


For $J_m>1$ (i.e., $m>3$),
as in Section \ref{sec:prox-map},
let $\widetilde \rho_{_{PE}}(n,m,r)=\left|\A\right|/(n\,(n-1))$
be the relative density for the proportional-edge PCD in the multiple triangle case.
Let $\widetilde \rho_{_{CS}}(n,m,\tau)$ and $\rho^{{}^{[i]}}_{_{CS}}(\tau)$ be
defined similarly for the central similarity PCD.
Let $n_i$ be the number of $\X$ points in $T_i$ for $i=1,2,\ldots,J_m$.
Letting $w_i = A(T_i) / A(C_H(\Y_m))$ with $A(\cdot)$ being the area function
and $\mathcal W = \{w_1,w_2,\ldots,w_{J_m}\}$,
we obtain the following as a corollary to Theorems \ref{thm:asy-norm-PE} and \ref{thm:asy-norm-CS}.

\begin{corollary}
\label{cor:MT-asy-norm}
For $r \in [1,\infty]$,
the asymptotic distribution for $\widetilde \rho_{_{PE}}(n,m,r)$ conditional on $\mathcal W$
is given by
\begin{equation}
\sqrt{n}\left(\widetilde \rho_{_{PE}}(n,m,r)-\widetilde \mu_{_{PE}}(m,r)\right)
\stackrel{\mathcal L}{\longrightarrow}\\
\mathcal{N}
 \left(
   0,
   4\,\widetilde \nu_{_{PE}}(m,r)
 \right),
\end{equation}
as $n \rightarrow \infty$,
where
$\widetilde \mu_{_{PE}}(m,r)=\mu_{_{PE}}(r) \left(\sum_{i=1}^{J_m}w_i^2\right)$
and
$$\widetilde \nu_{_{PE}}(m,r)=
\left[  \nu_{_{PE}}(r) \left(\sum_{i=1}^{J_m}w_i^3 \right)+
\left( \mu_{_{PE}}(r) \right)^2\left(\sum_{i=1}^{J_m}w_i^3-\left(\sum_{j=1}^{J_m}w_i^2 \right)^2\right) \right]$$
with $\mu_{_{PE}}(r)$ and $\nu_{_{PE}}(r)$ being as in Equations \eqref{eqn:PEAsymean} and \eqref{eqn:PEAsyvar},
respectively.
The asymptotic distribution of $\widetilde \rho_{CS}(n,m,\tau)$ with $\tau \in (0,\infty]$ is similar.
\end{corollary}

\noindent
\textbf{Proof:}
The expectation of $\widetilde \rho_{_{PE}}(n,m,r)$ is
$$ \E\left[\widetilde \rho_{_{PE}}(n,m,r)\right]=\frac{1}{n\,(n-1)}
\sum\hspace*{-0.1 in}\sum_{i < j \hspace*{0.25 in}}
\hspace*{-0.1 in}\,\E\left[h_{ij}(r)\right]=\E\left[ h_{12}(r) \right]/2=
P(X_2\in \NPE(X_1,r))=\widetilde \mu_{_{PE}}(r).$$
By definition of $\NPE(\cdot,r)$,
if $X_1$ and $X_2$ are in different triangles, then $P(X_2 \in \NPE(X_1,r))=0$.
So by the law of total probability
\begin{eqnarray*}
\widetilde \mu_{_{PE}}(r)&:=&P(X_2 \in \NPE(X_1,r))= \sum_{i=1}^{J_m}P(X_2 \in \NPE(X_1,r)\,|\,\{X_1,X_2\} \subset T_i)\,P(\{X_1,X_2\} \subset T_i)\\
&=& \sum_{i=1}^{J_m}\mu_{_{PE}}(r)\,P(\{X_1,X_2\} \subset T_i)
\text{ (since $P(X_2 \in \NPE(X_1,r)\,|\,\{X_1,X_2\} \subset T_i)=\mu_{_{PE}}(r)$)}\\
&=& \mu_{_{PE}}(r) \, \sum_{i=1}^{J_m}\left(\frac{A(T_i)}{\sum_{i=1}^{J_m}A(T_i)}\right)^2
\text{ (since $P(\{X_1,X_2\} \subset T_i)=\left(\frac{A(T_i)}{\sum_{i=1}^{J_m}A(T_i)}\right)^2$)}\\
&=& \mu_{_{PE}}(r) \, \left(\sum_{i=1}^{J_m}w_i^2\right).
\end{eqnarray*}
where $\mu_{_{PE}}(r)$ is given by Equation \eqref{eqn:PEAsymean}.

Likewise, we get $\widetilde \mu_{_{CS}}(\tau)=\mu_{_{CS}}(\tau)\,\left(\sum_{i=1}^{J_m}w_i^2\right)$
where $\mu_{_{CS}}(\tau)$ is given by Equation \eqref{eqn:CSAsymean}.

Furthermore, the asymptotic variance is
\begin{eqnarray*}
\widetilde \nu_{_{PE}}(m,r)&=&\E[h_{12}\,h_{13}]-\E[h_{12}]\,\E[h_{13}]\\
& = & P\bigl( \{X_2,X_3\} \subset \NPE(X_1,r) \bigr)+2\,P\bigl( X_2 \in \NPE(X_1,r), X_3 \in \G^{^{PE}}_1(X_1,r) \bigr)\\
& &+P\bigl(\{X_2,X_3\} \subset \G^{^{PE}}_1(X_1,r)\bigr)-4\,(\widetilde \mu_{_{PE}}(m,r))^2.
\end{eqnarray*}

Let $P_{_{PE}}^{2N}(r):=P\bigl( \{X_2,X_3\} \subset \NPE(X_1,r) \bigr)$,
$P_{_{PE}}^{2G}(r):=P\bigl(\{X_2,X_3\} \subset \G^{^{PE}}_1(X_1,r)\bigr)$, and
$P_{_{PE}}^{M}(r):=P\bigl( X_2 \in \NPE(X_1,r), X_3 \in \G^{^{PE}}_1(X_1,r)\bigr)$.
Then for $J_m>1$, we have
\begin{eqnarray*}
P\bigl( \{X_2,X_3\} \subset \NPE(X_1,r) \bigr)&=&\sum_{j=1}^{J_m}P\bigl( \{X_2,X_3\} \subset
\NPE(X_1,r)\,|\, \{X_1,X_2,X_3\} \subset T_j \bigr)\, P\bigl( \{X_1,X_2,X_3\} \subset T_j \bigr)\\
& = &\sum_{j=1}^{J_m}P_{_{PE}}^{2N}(r)\, \bigl( A(T_j) / A(C_H(\Y_m)) \bigr)^3 =
P_{_{PE}}^{2N}(r)\, \left(\sum_{j=1}^{J_m}w_j^3 \right).
\end{eqnarray*}
Similarly, $P\bigl( X_2 \in \NPE(X_1,r), X_3 \in \G^{^{PE}}_1(X_1,r) \bigr)=
P_{_{PE}}^{M}(r)\,\left(\sum_{j=1}^{J_m}w_j^3 \right) \text{  and }
P\bigl( \{X_2,X_3\} \subset \G^{^{PE}}_1(X_1,r) \bigr)=
P_{_{PE}}^{2G}(r)\,\left(\sum_{j=1}^{J_m}w_j^3 \right)$,
hence,
{\small
$$\widetilde \nu_{_{PE}}(m,r)=\bigl( P_{_{PE}}^{2N}(r)+2\,P^r_M+P_{_{PE}}^{2G}(r) \bigr)\,\left(\sum_{j=1}^{J_m}w_j^3 \right)-4\,\widetilde \mu_{_{PE}}(m,r)^2=
\nu_{_{PE}}(r)\,\left(\sum_{j=1}^{J_m}w_j^3 \right)+
4\,\mu_{_{PE}}(r)^2\,\left(\sum_{j=1}^{J_m}w_j^3-\left( \sum_{j=1}^{J_m}w_j^2
\right)^2\right),$$ }

Likewise, we get
$\widetilde \nu_{_{CS}}(\tau)=\nu_{_{CS}}(\tau)\,\left(\sum_{i=1}^{J_m}w_i^3\right)+
4\,\mu_{_{CS}}(\tau)^2\,\left(\sum_{i=1}^{J_m}w_i^3-\left(\sum_{i=1}^{J_m}w_i^2\right)^2\right).$

So, conditional on $\mathcal W$,
if $\widetilde \nu_{_{PE}}(r)>0$,
then
$\sqrt{n}\,\left(\widetilde \rho_{{}_{PE}}(n,m,r)-\widetilde \mu_{_{PE}}(r)\right)
\stackrel {\mathcal L}{\longrightarrow} \mathcal N\left(0,\widetilde \nu_{_{PE}}(r)\right)$.
A similar result holds for the relative density of the central similarity PCD.
$\blacksquare$

By an appropriate application of the Jensen's inequality,
we see that $\sum_{i=1}^{J_m}w_i^3 \ge \left(\sum_{i=1}^{J_m}w_i^2 \right)^2.$
So the covariance above is zero iff $\nu_{_{PE}}(r)=0$ and
$\sum_{i=1}^{J_m}w_i^3=\left(\sum_{i=1}^{J_m}w_i^2 \right)^2$,
so asymptotic normality may hold even though $\nu_{_{PE}}(r)=0$ in the multiple triangle case.
That is, $\widetilde \rho_{{}_{PE}}(n,m,r)$ has the asymptotic normality for $r = \infty$ also
provided that $\sum_{i=1}^{J_m}w_i^3 > \left(\sum_{i=1}^{J_m}w_i^2 \right)^2$.
The same holds for $\tau=\infty$ in the central similarity case.

\section{Alternative Patterns: Segregation and Association}
\label{sec:alternatives}
In a two class setting,
the phenomenon known as {\em segregation} occurs when members of
one class have a tendency to repel members of the other class.
For instance, it may be the case that one type of plant
does not grow well in the vicinity of another type of plant,
and vice versa.
This implies, in our notation,
that $X_i$ are unlikely to be located near elements of $\Y_m$.
Alternatively, association occurs when members of one class
have a tendency to attract members of the other class,
as in symbiotic species, so that $X_i$ will tend to
cluster around the elements of $\Y_m$, for example.
See, for instance, \cite{dixon:1994} and \cite{coomes:1999}.

These alternatives can be parametrized as follows.
In the one triangle case, without loss of generality
let $\Y_3=\left\{(0,0),(1,0),(c_1,c_2) \right\}$ and $T_b=T(\Y_3)$
with $\y_1=(0,0),\y_2=(1,0)$, and $\y_3=(c_1,c_2)$.
For the basic triangle $T_b$,
let $Q_{\theta}:=\{x \in T_b: d(x,\Y_3) \le \theta\}$ for $\theta \in (0,(c_1^2+c_2^2)/2]$
and $S(F)$ be the support of $F$.
Then consider
$$\mathscr H_S:=\{F: S(F) \subseteq T_b \text{ and } P_F(X \in Q_{\theta}) < P_U(X \in Q_{\theta}) \}$$
and
$$\mathscr H_A:=\{F: S(F) \subseteq T_b \text{ and } P_F(X \in Q_{\theta}) > P_U(X \in Q_{\theta}) \}$$
where $P_F$ and $P_U$ are probabilities with respect to distribution
function $F$ and the uniform distribution on $T_b$, respectively.
So if $X_i \stackrel{iid}{\sim} F \in \mathscr H_S$,
the pattern between class $\X$ and $\Y$ points is segregation,
but if $X_i \stackrel{iid}{\sim} F \in \mathscr H_A$,
the pattern between class $\X$ and $\Y$ points is association.
For example the distribution family
$$\mathscr F_S:=\{F: S(F) \subset T_b \text{ and the associated pdf
$f$ increases as $d(x,\Y_3)$ increases}\}$$
is a subset of $\mathscr H_S$ and
yields samples from the segregation alternatives.
Likewise, the distribution family
$$\mathscr F_A:=\{F: S(F) \subset T_b \text{ and the associated pdf
$f$ increases as $d(x,\Y_3)$ decreases}\}$$
is a subset of $\mathscr H_A$ and
yields samples from the association alternatives.

In the basic triangle, $T_b$,
we define the alternatives $H^S_{\ve}$ and $H^A_{\ve}$ with $\ve \in \left( 0,\sqrt{3}/3 \right)$,
for segregation and association alternatives, respectively.
Under $H^S_{\ve}$, $4 \ve^2/3\times 100$ \% of the area of $T_b$
is chopped off around each vertex so that the $\X$ points are restricted
to lie in the remaining region.
That is, for $\y_j \in \Y_3$,
let $e_j$ denote the edge of $T_b$ opposite vertex $\y_j$ for $j=1,2,3$,
and for $x \in T_b$,
let $\ell_j(x)$ denote the line parallel to $e_j$ through $x$.
Then define
$T_j(\ve) = \{x \in T_b: d(\y_j,\ell_j(x)) \le \ve_j\}$ where
$\displaystyle \ve_1=\frac{2\,c_2\,\ve}{3\sqrt{c_2^2+(1-c_1)^2}}$,
$\displaystyle \ve_2=\frac{2\,c_2\,\ve}{3\sqrt{c_1^2+c_2^2}}$, and
$\displaystyle \ve_3=\frac{2\,c_2\,\ve}{3}$.
Let $\mathcal T_\ve:=\bigcup_{j=1}^3 T_j(\ve)$.
Then under $H^S_{\ve}$,
we have
$X_i \stackrel{iid}{\sim} \U\left(T_b \setminus \mathcal T_\ve\right)$.
Similarly,
under $H^A_{\ve}$,
we have
$X_i \stackrel{iid}{\sim} \U\left(\mathcal T_{\sqrt{3}/3 - \ve} \right)$.
Thus the segregation model excludes the possibility of
any $X_i$ occurring around a $\y_j$,
and the association model requires
that all $X_i$ occur around $\y_j$'s.
The $\sqrt{3}/3 - \ve$ is used in the definition of the
association alternative so that $\ve=0$
yields $H_o$ under both classes of alternatives.
Thus,
we have the below parametrization of the
distribution families under the alternatives.

\begin{equation}
\label{eqn:eps-alt-1}
\mathscr U^S_{\ve}:=\{F: F = \U(T_b \setminus \mathcal T_\ve) \}
\text{ ~and~  }
\mathscr U^A_{\ve}:=\{F: F = \U(\mathcal T_{\sqrt{3}/3 - \ve}) \}.
\end{equation}
Clearly $\mathscr U^S_{\ve} \subsetneq \mathscr H_S$
and $\mathscr U_{\sqrt{3}/3 - \ve}^A \subsetneq \mathscr H_A$,
but
$\mathscr U^S_{\ve} \nsubseteq \mathscr F_S$
and $\mathscr U_{\sqrt{3}/3 - \ve}^A \nsubseteq \mathscr F_A$.

These alternatives $H^S_{\ve}$ and $H^A_{\ve}$ with $\ve \in \left( 0,\sqrt{3}/3 \right)$,
can be transformed into the equilateral triangle as in
\cite{ceyhan:arc-density-PE} and \cite{ceyhan:arc-density-CS}.

\begin{figure} []
\centering
\scalebox{.35}{\input{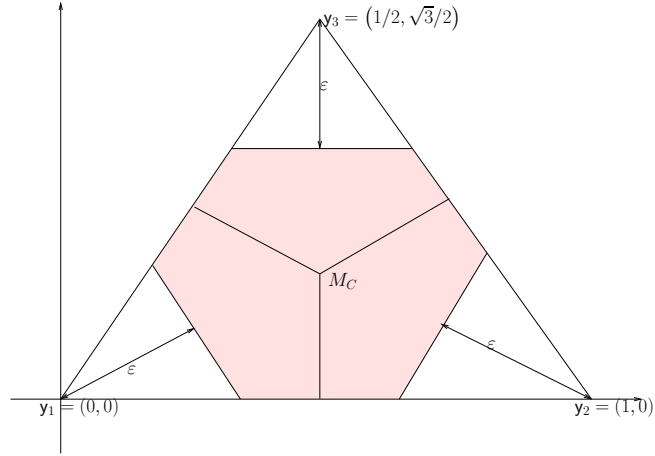}}
\caption{
\label{fig:seg-alt-1}
An example for the segregation alternative with a particular expansion parameter $\ve$ (shaded region),
and its complement is for the association alternative with expansion parameter $\sqrt{3}/3 - \ve$ (unshaded region)
on the standard equilateral triangle.}
\end{figure}

For the standard equilateral triangle,
in $T_j(\ve) = \{x \in T_e: d(\y,\ell_j(x)) \le \ve_j\}$,
we have $\ve_1=\ve_2=\ve_3=\ve$.
Thus $H^S_{\ve}$ implies
$X_i \stackrel{iid}{\sim} \U\left(T_e \setminus \mathcal T_\ve \right)$
and $H^A_{\ve}$ be the model under which
$X_i \stackrel{iid}{\sim} \U\left(\mathcal T_{\sqrt{3}/3 - \ve}\right)$.
See Figure \ref{fig:seg-alt-1} for a depiction of the above
segregation and the association alternatives in $T_e$.

\begin{remark}
The geometry invariance result of Theorem \ref{thm:geo-inv}
also holds under the alternatives $H^S_{\ve}$ and $H^A_{\ve}$
for both PCD families.
In particular, the segregation alternative with $\ve \in \left( 0,\sqrt{3}/4 \right)$ in the standard equilateral triangle
corresponds to the case that in an arbitrary triangle, $\kappa \times 100\%$ of the area is carved
away as forbidden from the vertices using line segments parallel to the opposite edge
where $\kappa = 4\ve^2$ (which implies $\kappa \in (0,3/4)$).
But the segregation alternative with $\ve \in \left( \sqrt{3}/4,\sqrt{3}/3 \right)$ in the standard equilateral triangle
corresponds to the case that in an arbitrary triangle, $\kappa \times 100\%$ of the area is carved
away as forbidden from each vertex using line segments parallel to the opposite edge
where $\kappa = 1-4 \left(1-\sqrt{3}\ve \right)^2$ (which implies $\kappa \in (3/4,1)$).
This argument is for the segregation alternative;
a similar construction is available for the association alternative. $\square$
\end{remark}

\begin{remark}
\label{rem:alt-multi-tri}
\textbf{The Alternatives in the Multiple Triangle Case:}
In the multiple triangle case,
the segregation and association alternatives,
$H^S_{\ve}$ and $H^A_{\ve}$ with $\ve \in \left( 0,\sqrt{3}/3 \right)$,
are defined as in the one-triangle case,
in the sense that,
when each triangle (together with the data in it)
is transformed to the standard equilateral triangle as in Theorem \ref{thm:geo-inv},
we obtain the same alternative pattern described above. 

%

Thus in the case of $J_m>1$, we have a (conditional) test of
$H_o: X_i \stackrel{iid}{\sim} \U(C_H(\Y_m))$
which once again
rejects against segregation for large values of $\rho_n(\tau,J)$ and
rejects against association for small values of $\rho_n(\tau,J)$.
The segregation (with $\kappa=1/16$, i.e., $\ve=\sqrt{3}/8$),
null, and association (with $\kappa=1/4$, i.e.,
$\ve=\sqrt{3}/12$) realizations (from left to right) are
depicted in Figure \ref{fig:deldata} with $n=1000$.
%
$\square$
\end{remark}

\subsection{Asymptotic Normality under the Alternatives}
\label{sec:asy-norm-alt}
Asymptotic normality of relative density of
the PCDs under both alternative hypotheses of
segregation and association can be established by the same method as
under the null hypothesis.
Let $\E^S_{\ve}[\cdot]$ ($\E^A_{\ve}[\cdot]$) be the expectation with respect to the
uniform distribution under the segregation ( association )
alternatives with $\ve \in \left( 0,\sqrt{3}/3 \right)$.

\begin{theorem}
\label{thm:asy-norm-alt-1}
\begin{itemize}
\item[]
\item[(i)]
Let $\mu^S_{_{PE}}(r,\ve)$ be the mean $\E^S_{\ve}[h_{12}]$ and
$\nu^S_{_{PE}}(r,\ve)$ be the covariance, $\Cov^S_{\ve}[h_{12},h_{13}]$ for
$r \in [1,\infty]$ and $\ve \in \bigl[ 0,\sqrt{3}/3 \bigr)$ under $H^S_{\ve}$.
Then as $n \rightarrow \infty$,
$\sqrt{n}\bigl( \rho_{_{PE}}(n,r)-\mu^S_{_{PE}}(r,\ve) \bigr)
\stackrel {\mathcal L}{\longrightarrow} \N(0,\nu^S_{_{PE}}(r,\ve))$ for
the values of $(r,\ve)$ for which $\nu^S_{_{PE}}(r,\ve)>0$.
A similar result holds under association.

\item[(ii)]
Let $\mu^S_{_{CS}}(\tau,\ve)$ be the mean $\E^S_{\ve}[h_{12}]$ and
$\nu^S_{_{CS}}(\tau,\ve)$ be the covariance, $\Cov^S_{\ve}[h_{12},h_{13}]$ for
$\tau \in (0,\infty]$ and $\ve \in \bigl[ 0,\sqrt{3}/3 \bigr)$ under $H^S_{\ve}$.
Then as $n \rightarrow \infty$,
$\sqrt{n}\bigl( \rho_{_{CS}}(n,\tau)-\mu^S_{_{PE}}(\tau,\ve) \bigr)
\stackrel {\mathcal L}{\longrightarrow} \N(0,\nu^S_{_{CS}}(\tau,\ve))$ for
the values of $(\tau,\ve)$ for which $\nu^S_{_{PE}}(\tau,\ve)>0$.
A similar result holds under association.
\end{itemize}
\end{theorem}

A sketch of the proof of part (i) is provided in (\cite{ceyhan:arc-density-PE}),
and of part (ii) for $\tau \in (0,1]$ is provided in (\cite{ceyhan:arc-density-CS}).
The proof of part (ii) for $\tau \in (1,\infty)$ is similar.

The explicit forms of $\mu^S_{_{PE}}(r,\ve)$ and $\mu^A_{_{PE}}(r,\ve)$
are given, defined piecewise, in (\cite{ceyhan:TR-rel-dens-NPE}).
Note that under $H^S_{\ve}$,
$$\nu^S_{_{PE}}(r,\ve)>0 \text{ for } (r,\ve) \in
\left[ 1,\sqrt{3}/(2 \ve) \right) \times \left(0,\sqrt{3}/4\right] \bigcup
\left[ 1,\sqrt{3}/\ve-2 \right) \times \left( \sqrt{3}/4,\sqrt{3}/3 \right),$$
and under $H^A_{\ve}$,
$$\nu^A_{_{PE}}(r,\ve) > 0\text{ for }(r,\ve)\in
(1,\infty) \times \left( 0,\sqrt{3}/3 \right) \bigcup \{1\} \times \left( 0,\sqrt{3}/12 \right).$$
The explicit forms of $\mu^S_{_{CS}}(\tau,\ve)$ and
$\mu^A_{_{CS}}(\tau,\ve)$ are given, defined piecewise, in (\cite{ceyhan:TR-CS-rel-dens}).
Note that under $H^S_{\ve}$,
$$\nu^S_{_{CS}}(\tau,\ve) > 0 \text{  for }(\tau,\ve) \in ( 0,\infty) \times \bigl( 0,3\,\sqrt{3}/10 \Bigr]
\bigcup \Biggl(\frac{2\,(\sqrt{3}-3\,\ve)}{4\,\ve-\sqrt{3}},\infty \Biggr]
\times \bigl( 3\,\sqrt{3}/10,\sqrt{3}/3 \bigr),$$
and under $H^A_{\ve}$,
$$\nu^A_{_{CS}}(\tau,\ve)>0 \text{ for }(\tau,\ve)\in (0,\infty] \times \bigl(0,\sqrt{3}/3 \bigr).$$

Notice that under association alternatives any $r \in [1,\infty)$ yields asymptotic normality
for relative density of proportional-edge PCD for all $\ve \in \left(0,\sqrt{3}/3 \right)$,
while under segregation alternatives,
only $r=1$ yields this universal asymptotic normality.
Furthermore, under association alternatives
any $\tau \in (0,\infty)$ yields asymptotic normality
for relative density of central similarity PCD
for all $\ve \in \bigl( 0,\sqrt{3}/3 \bigr)$.
The same holds under segregation alternatives.

The asymptotic normality also holds under the alternatives
in the multiple triangle case.
For example,
for the relative density of proportional-edge PCDs,
the asymptotic mean and variance are as in Corollary \ref{cor:MT-asy-norm}
with $\mu_{_{PE}}(r)$ ($\nu_{_{PE}}(r)$)  being replaced by $\mu^S_{_{PE}}(r,\ve)$ ($\nu^S_{_{PE}}(r,\ve)$) for segregation
and by $\mu^A_{_{PE}}(r,\ve)$ ($\nu^A_{_{PE}}(r,\ve)$) for association.

and

\subsection{The Test Statistics and Analysis}
\label{sec:test-stat-analysis}
The relative density of the PCD
is a test statistic for the segregation/association alternative;
rejecting for extreme values of $\rho_{_{PE}}(n,r)$ is appropriate
since under segregation we expect $\rho_{_{PE}}(n,r)$ to be large,
while under association we expect $\rho_{_{PE}}(n,r)$ to be small.

In the one triangle case,
using the standardized test statistic
\begin{equation}
\label{eqn:rho-PE-standardized}
R_{PE}(r) = \frac{\sqrt{n} \bigl( \rho_{_{PE}}(n,r) - \mu_{_{PE}}(r) \bigr)}{\sqrt{\nu_{_{PE}}(r)}},
\end{equation}
the asymptotic critical value
for the one-sided level $\alpha$ test against segregation
is given by
\begin{equation}
z_{\alpha} = \Phi^{-1}(1-\alpha)
\end{equation}
where $\Phi(\cdot)$ is the standard normal distribution function.
Against segregation, the test rejects for $R_{PE}(r) > z_{\alpha}$ and against association,
the test rejects for $R_{PE}(r) < z_{1-\alpha}$.
The same holds for the standardized test statistic
in the multiple triangle case,
$\widetilde R_{PE}(r) = \frac{\sqrt{n} \bigl( \widetilde \rho_{_{PE}}(n,r) - \widetilde \mu_{_{PE}}(r) \bigr)}{\sqrt{\widetilde \nu_{_{PE}}(r)}}$.

A similar construction is available for $\rho_{_{CS}}(n,\tau)$ with
\begin{equation}
\label{eqn:rho-CS-standardized}
R_{CS}(\tau) = \frac{\sqrt{n} \bigl(\rho_{_{CS}}(n,\tau) - \mu_{_{CS}}(\tau)\bigr)}{\sqrt{\nu_{_{CS}}(\tau)}}
\end{equation}
in the one triangle case,
and
with
$\widetilde R_{CS}(\tau) = \frac{\sqrt{n} \bigl(\widetilde \rho_{_{CS}}(n,\tau) - \widetilde \mu_{_{CS}}(\tau)\bigr)}{\sqrt{\widetilde \nu_{_{CS}}(\tau)}}$
in the multiple triangle case.

\subsection{Consistency}
\label{sec:consistency}

\begin{theorem}
\label{thm:consistency-1}
\begin{itemize}
\item[]
\item[(i)]
In the one triangle case,
the test against $H^S_{\ve}$ which rejects for $R_{PE}(r) > z_{1-\alpha}$
and
the test against $H^A_{\ve}$ which rejects for $R_{PE}(r) < z_{\alpha}$
are consistent for $r \in [1,\infty)$ and $\ve \in \left( 0,\sqrt{3}/3 \right)$.
The same holds in the multiple triangle case with $\widetilde R_{PE}(r)$.

\item[(ii)]
In the one triangle case,
the test against $H^S_{\ve}$ which rejects for $R_{CS}(\tau) > z_{1-\alpha}$
and
the test against $H^A_{\ve}$ which rejects for $R_{CS}(\tau) < z_{\alpha}$
are consistent for $\tau \in (0,\infty)$ and $\ve \in \left( 0,\sqrt{3}/3 \right)$.
The same holds in the multiple triangle case with $\widetilde R_{CS}(\tau)$.
\end{itemize}
\end{theorem}

For the one triangle case,
the proof of (i) is provided in (\cite{ceyhan:arc-density-PE})
and the proof of (ii) with $\tau \in (0,1]$ is provided in (\cite{ceyhan:arc-density-CS}).
The proofs for the multiple triangle cases
and for (ii) with $\tau>1$ are similar.

\section{Empirical Size Analysis}
\label{sec:emp-size-anal}

\subsection{Empirical Size Analysis for Proportional-Edge PCDs under CSR}
\label{sec:PE-emp-size-CSR}
\begin{figure}[ht]
\centering
\rotatebox{-90}{ \resizebox{2. in}{!}{\includegraphics{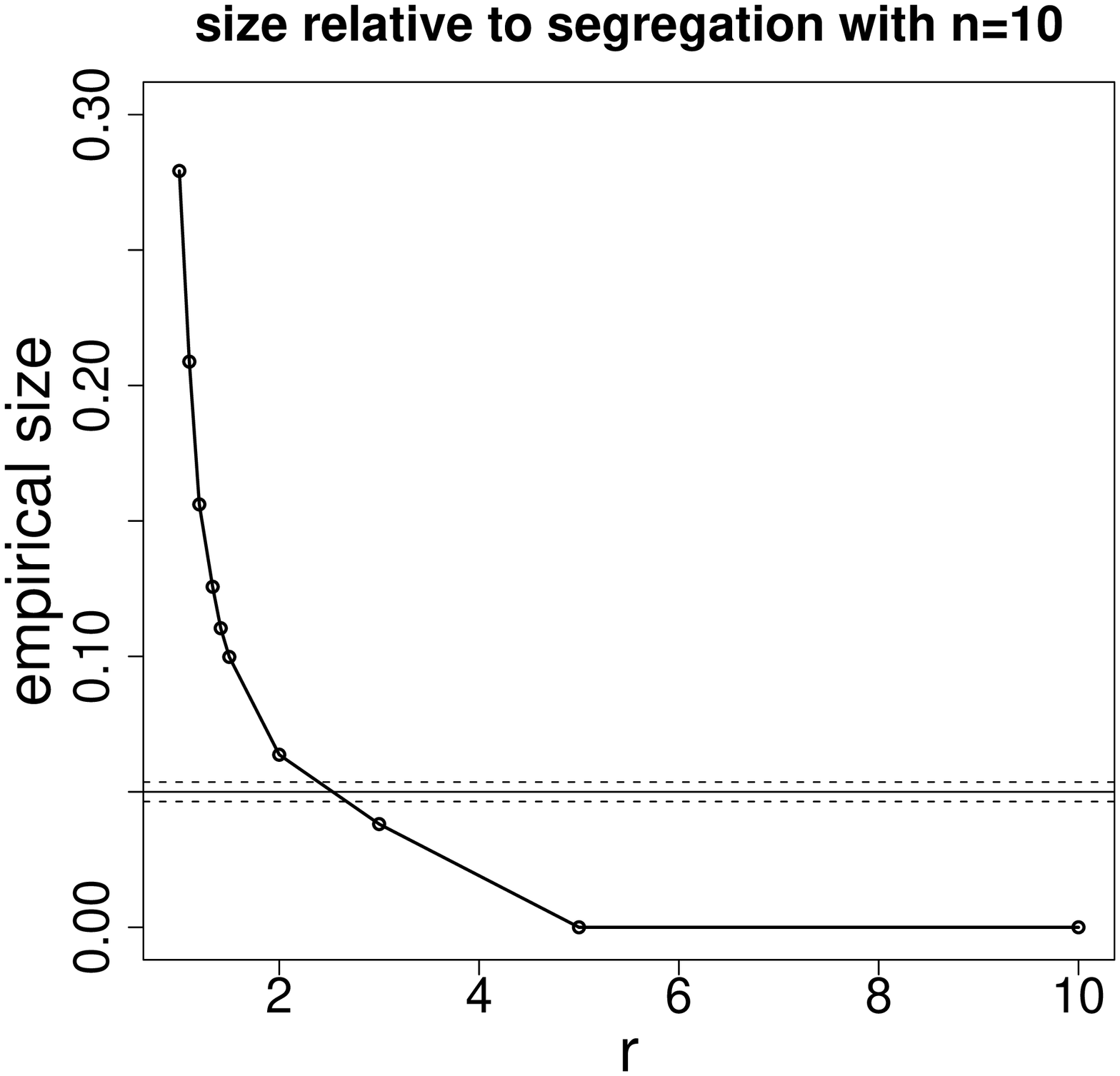} }}
\rotatebox{-90}{ \resizebox{2. in}{!}{\includegraphics{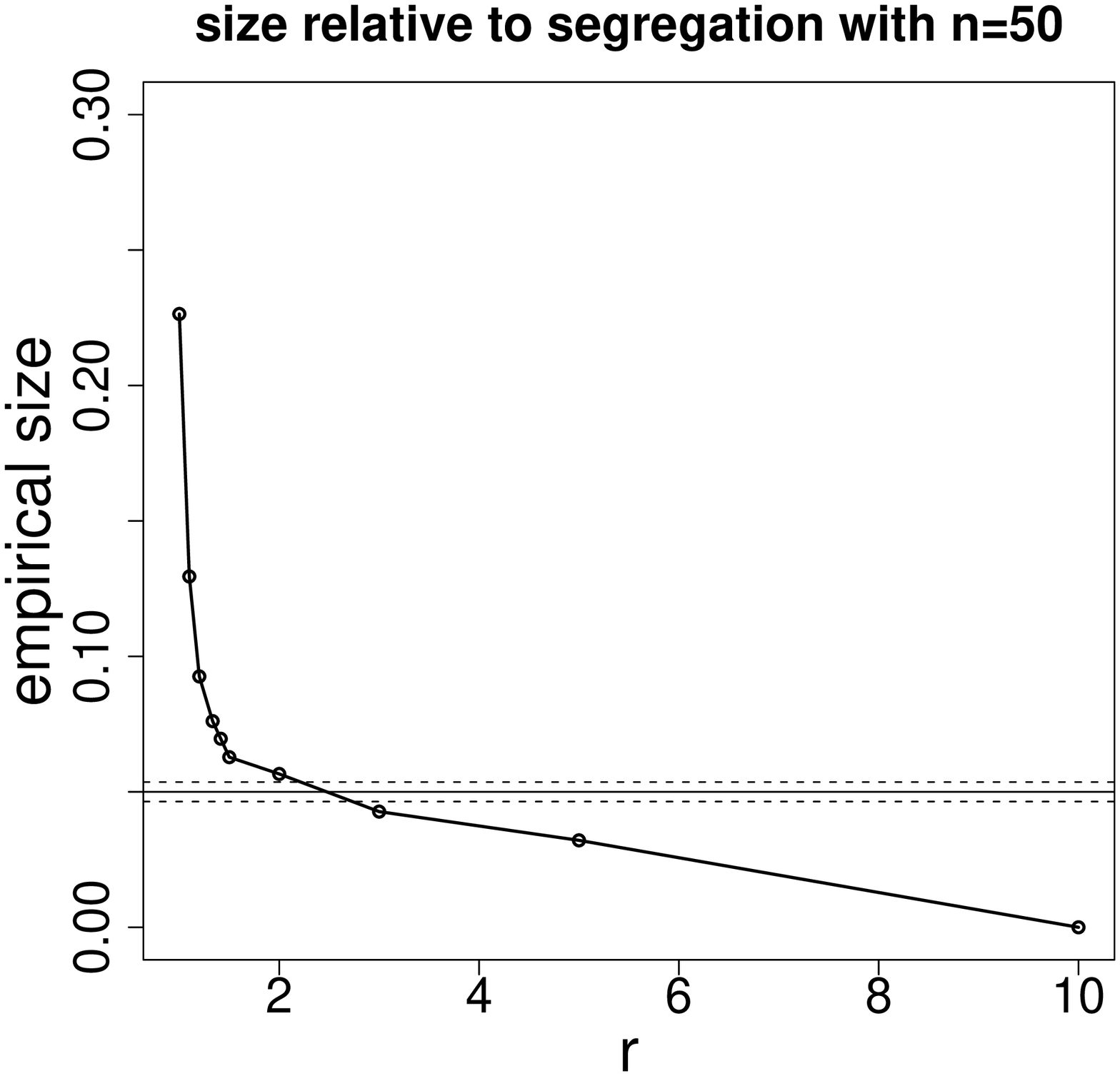} }}
\rotatebox{-90}{ \resizebox{2. in}{!}{\includegraphics{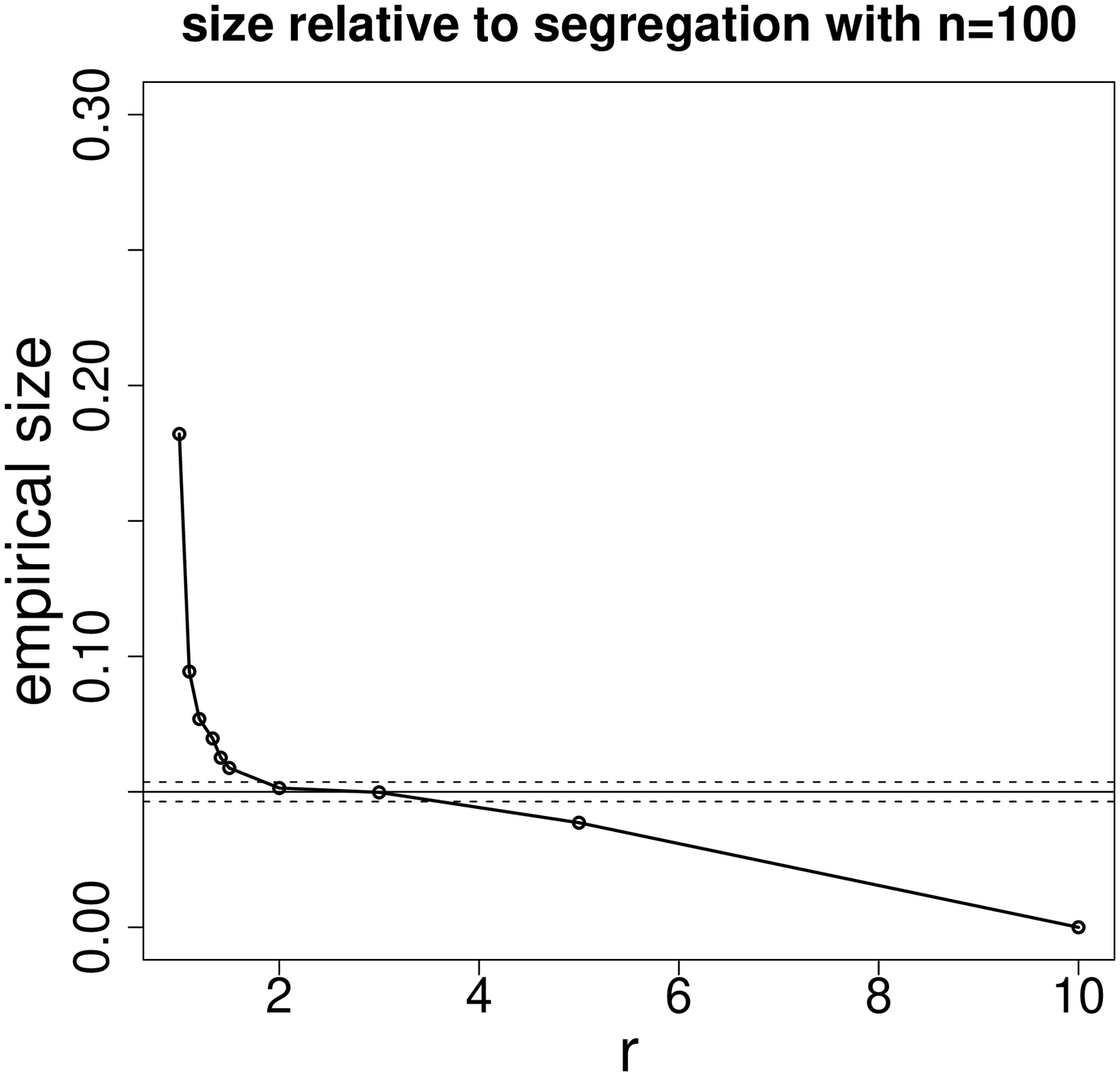} }}
\rotatebox{-90}{ \resizebox{2. in}{!}{\includegraphics{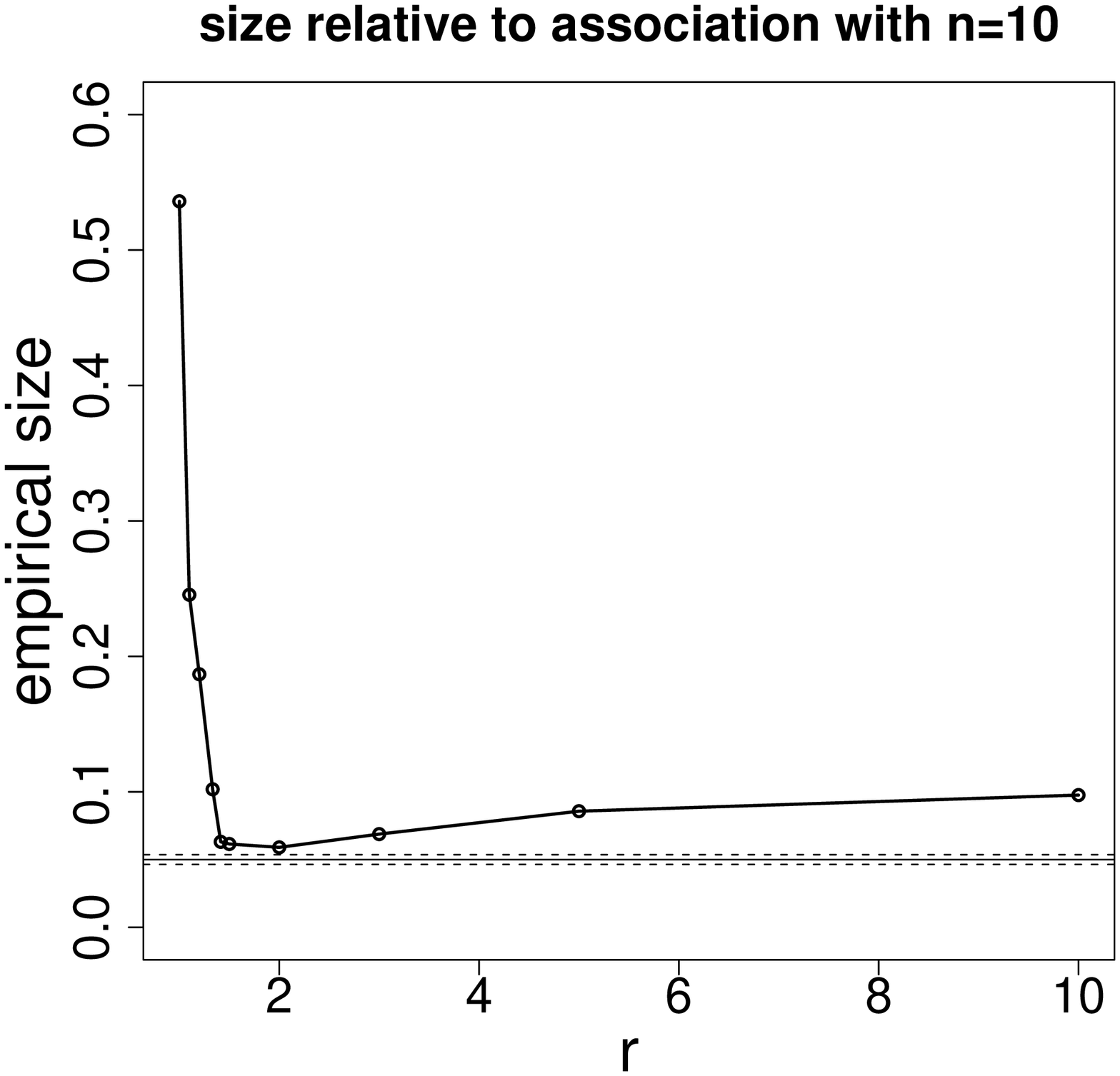} }}
\rotatebox{-90}{ \resizebox{2. in}{!}{\includegraphics{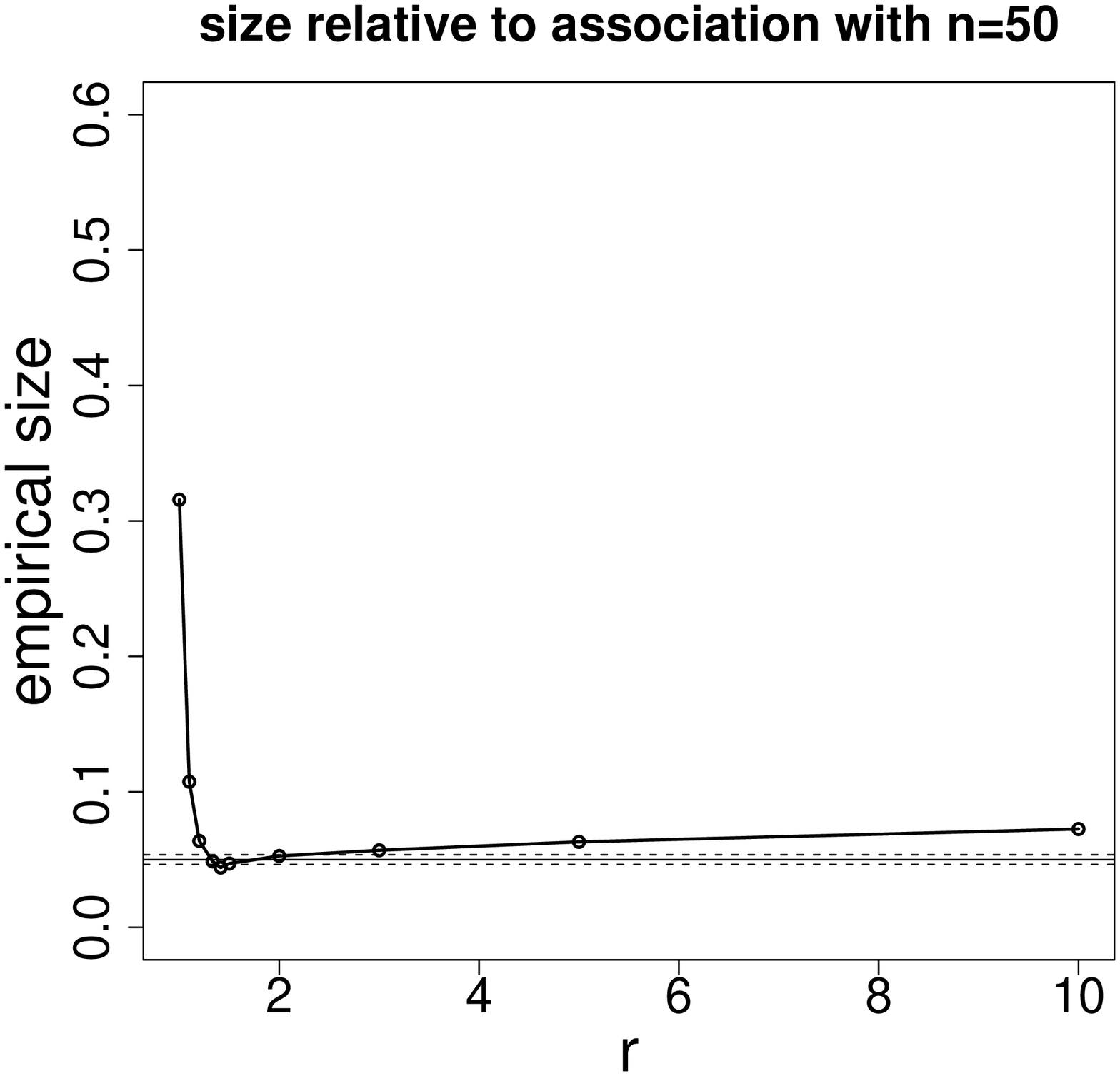} }}
\rotatebox{-90}{ \resizebox{2. in}{!}{\includegraphics{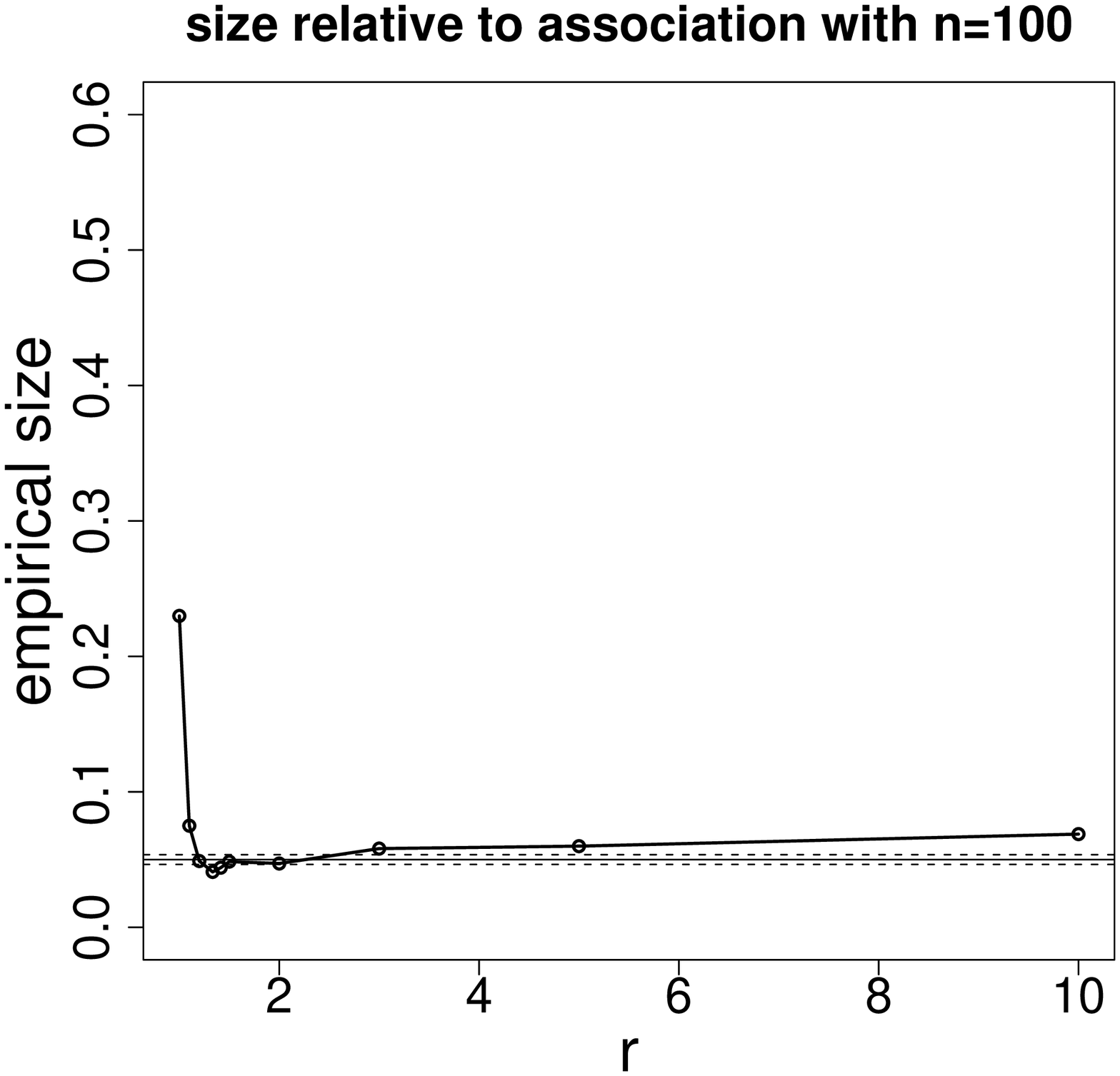} }}
\caption{
\label{fig:PE-emp-size-CSR}
\textbf{Empirical size for $R_{PE}(r)$ in the one triangle case:}
The empirical size estimates of the relative density of proportional-edge PCDs
in the one triangle case based on 10000 Monte Carlo replicates
for the left-sided alternative, i.e., relative to segregation (top)
and the right-sided alternative, i.e., relative to association (bottom)
with $n=10$ (left column), $n=50$ (middle column), and $n=100$ (right column) under the CSR pattern.
The horizontal lines are located at .0464 (upper threshold for conservativeness),
.050 (nominal level), and .0536 (lower threshold for liberalness).
}
\end{figure}

In one triangle case,
for the null pattern of CSR,
we generate $n$ $\X$ points iid $\U(T_e)$ where $T_e$ is the standard equilateral triangle.
We calculate the relative density of proportional-edge PCDs
for $r=1, 11/10, 6/5, 4/3, \sqrt{2},\\
3/2, 2, 3, 5, 10$ at each Monte Carlo replicate.
We repeat the Monte Carlo procedure $N_{mc}=10000$ times for each of $n=10, 50, 100$.
Using the critical values based on the normal approximation for the relative density,
we calculate the empirical size estimates for both right-sided (i.e., for segregation)
and left-sided (i.e., for association) tests as a function of the expansion parameter $r$.
Let $R_{PE}(r)(r,j):=\frac{\sqrt{n}\,\bigl( \rho_{_{PE}}(n,r,j)-\mu_{_{PE}}(r) \bigr)}{\sqrt{\nu_{_{PE}}(r)}}$
be the standardized relative density for Monte Carlo replicate $j$ with sample size $n$ for $j=1,2,\ldots,N_{mc}$.
For each $r$ value, the level $\alpha$ asymptotic critical value is
$\mu_{_{PE}}(r)+z_{(1-\alpha)} \cdot \sqrt{\nu_{_{PE}}(r)/n}$.
We estimate the empirical size against the segregation alternative as
$\frac{1}{N_{mc}}\sum_{j=1}^{N_{mc}}\I \left(R_{PE}(r)(r,j) > z_{1-\alpha} \right)$,
and
against the association alternative as
$\frac{1}{N_{mc}}\sum_{j=1}^{N_{mc}}\I\left( R_{PE}(r)(r,j) < z_{\alpha}\right)$.
The empirical sizes significantly smaller (larger) than .05 are deemed conservative (liberal).
The asymptotic normal approximation to proportions is used in determining the significance of
the deviations of the empirical sizes from .05.
For these proportion tests, we also use $\alpha=.05$ as the significance level.
With $N_{mc}=10000$, empirical sizes less than .0464 are deemed conservative,
greater than .0536 are deemed liberal at $\alpha=.05$ level.
The empirical sizes for the proportional-edge PCDs
together with upper and lower limits of liberalness and conservativeness
are plotted in Figure \ref{fig:PE-emp-size-CSR}.
Observe that as $n$ increases,
the empirical size gets closer to the nominal level of 0.05
(i.e., the normal approximation gets better).
For the right-sided tests (i.e., relative to segregation)
the size is close to the nominal level for $r \in (2,3)$,
for smaller $r$ values (i.e., $r<2$) the test seems to be liberal
with liberalness increasing as $r$ decreases;
and for larger $r$ values (i.e., $r>3$) the test seems to be conservative
with conservativeness increasing as $r$ increases.
For the left-sided tests (i.e., relative to association)
the size is close to the nominal level for $r \in (1.5,3)$,
for other $r$ values the test seems to be liberal
(more liberal for smaller $r$ values).
This is due to the fact that
very large and small values of $r$ require much larger sample
sizes for the normal approximation to hold.

\begin{figure}[ht]
\centering
\rotatebox{-90}{ \resizebox{2. in}{!}{\includegraphics{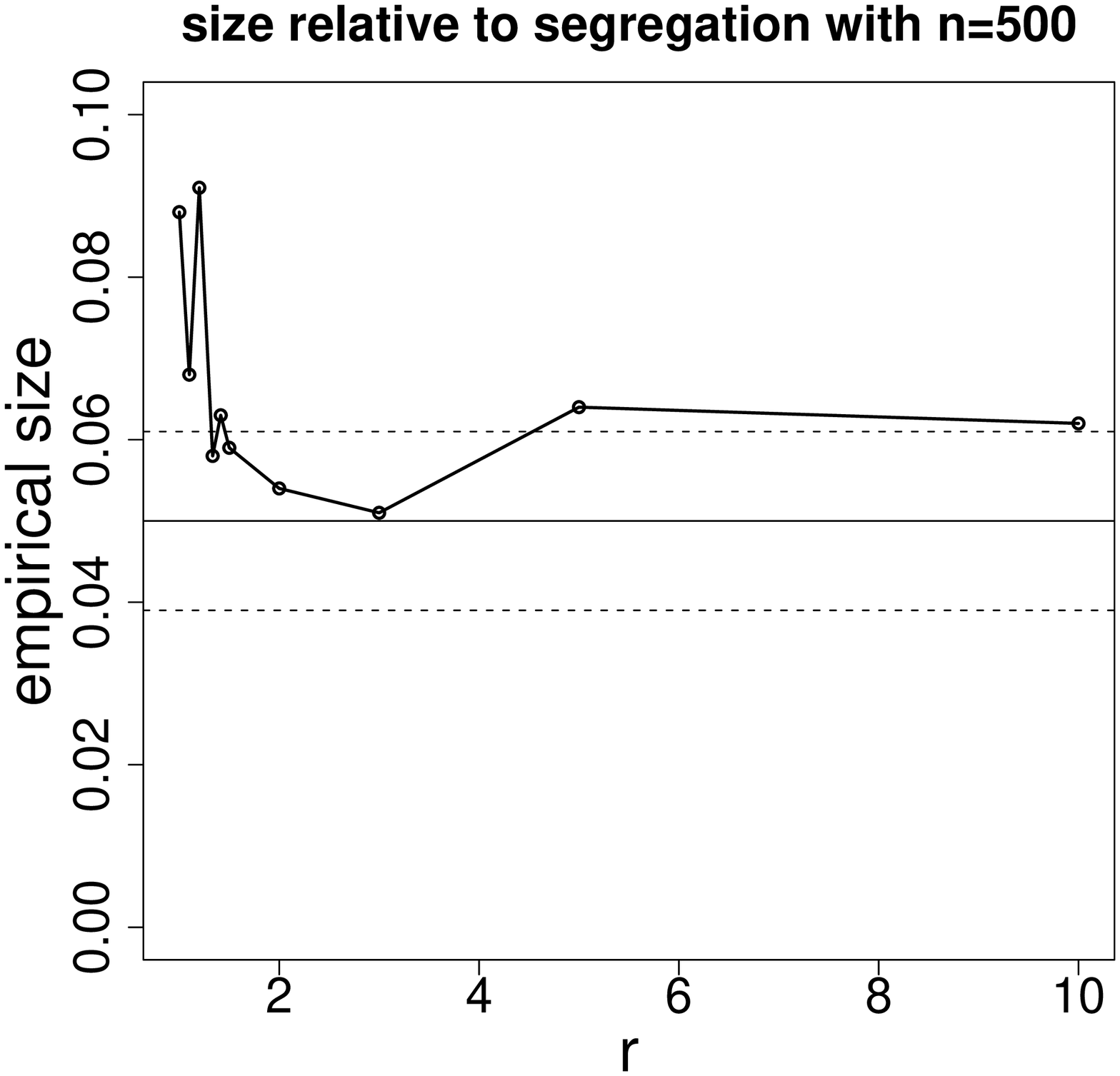} }}
\rotatebox{-90}{ \resizebox{2. in}{!}{\includegraphics{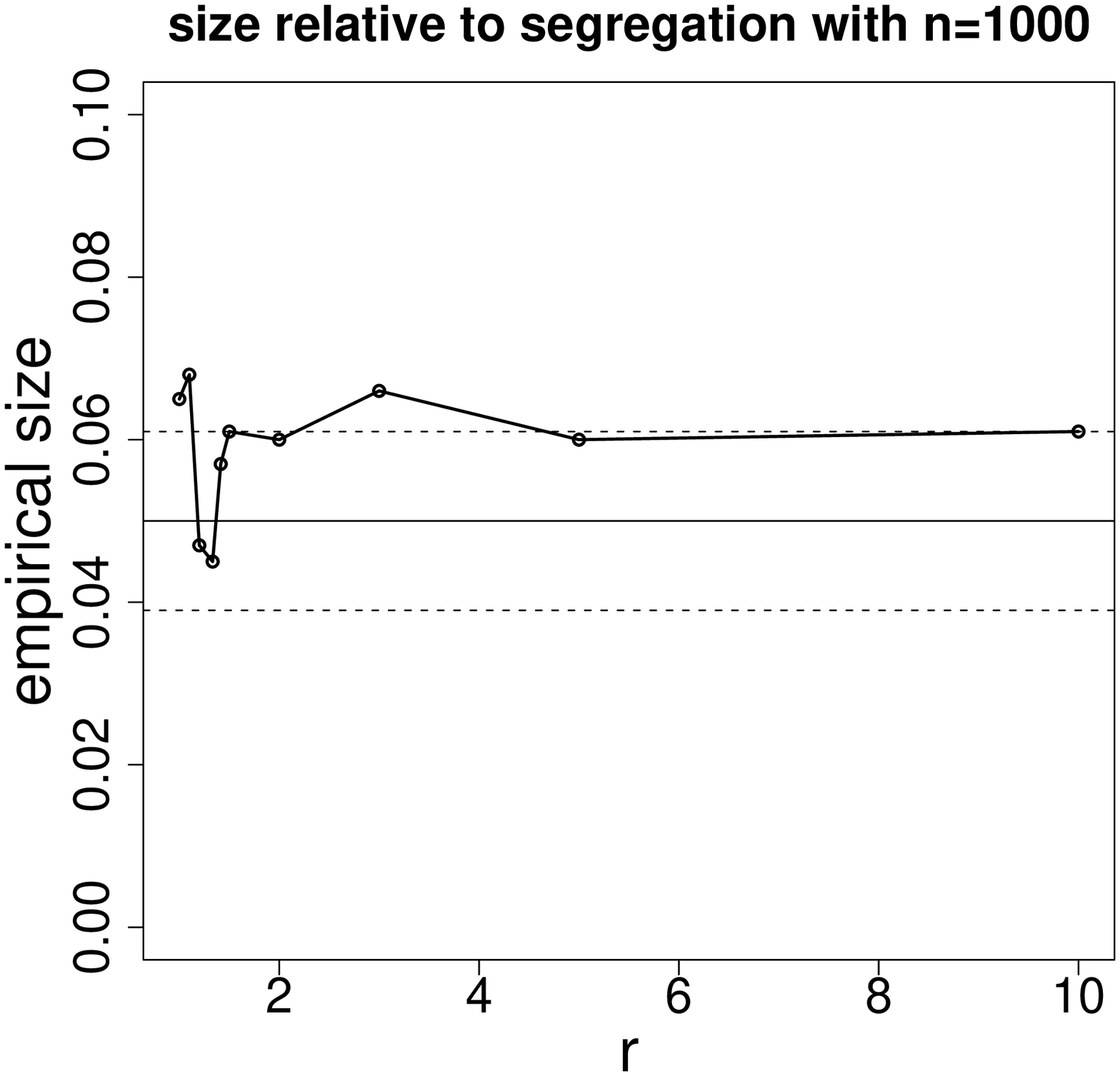} }}\\
\rotatebox{-90}{ \resizebox{2. in}{!}{\includegraphics{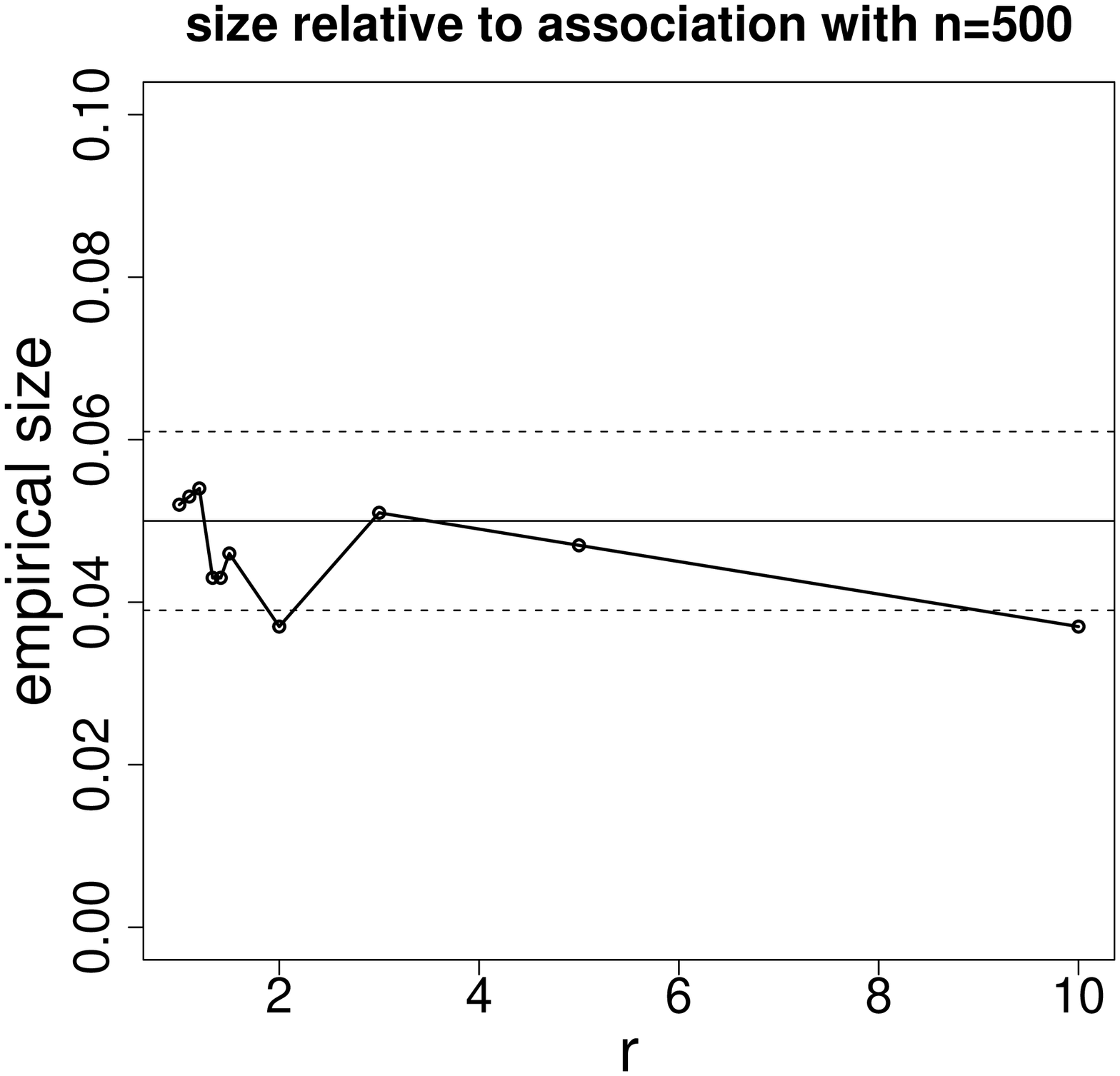} }}
\rotatebox{-90}{ \resizebox{2. in}{!}{\includegraphics{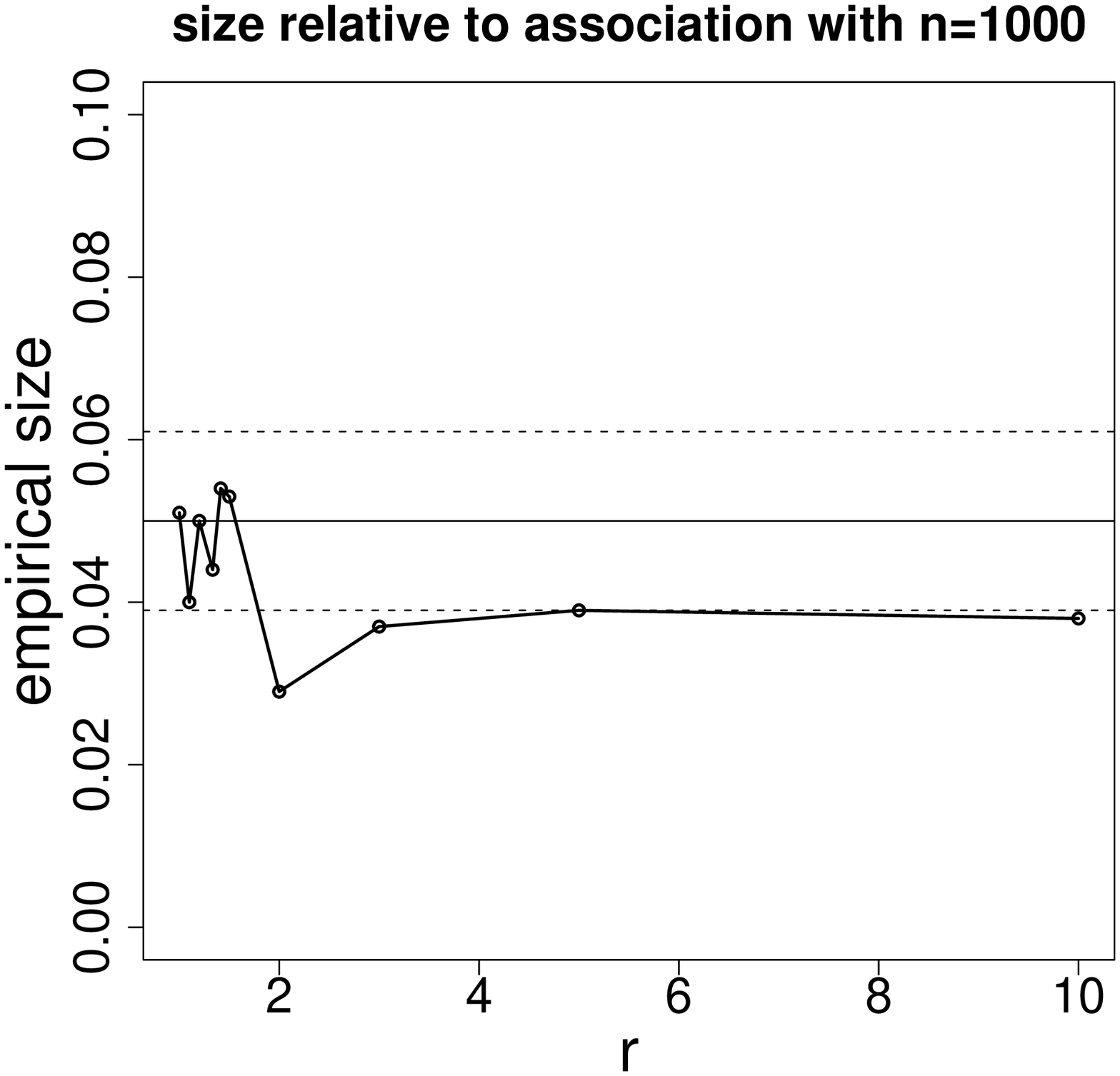} }}
\caption{
\label{fig:MT-PE-emp-size-CSR}
\textbf{Empirical size for $R_{PE}(r)$ in the multiple triangle case:}
The empirical size estimates of the relative density of proportional-edge PCDs
in the multiple triangle case based on 1000 Monte Carlo replicates
for the right-sided alternative (i.e., relative to segregation) (top)
and the left-sided alternative (i.e., relative to association) (bottom)
with $n=500$ (left column) and $n=1000$ (right column) under the CSR pattern.
The horizontal lines are located at .039 (upper threshold for conservativeness),
.050 (nominal level), and .061 (lower threshold for liberalness).
}
\end{figure}

In the multiple triangle case,
for the null pattern of CSR,
we generate $n$ $\X$ points iid $\U(C_H(\Y_{10}))$ where $\Y_{10}$ is the set of
the 10 class $\Y$ points given in Figure \ref{fig:deldata}.
With $N_{mc}=1000$, empirical sizes less than .039 are deemed conservative,
greater than .061 are deemed liberal at $\alpha=.05$ level.
The empirical sizes for the proportional-edge PCDs
together with upper and lower limits of liberalness and conservativeness
are plotted in Figure \ref{fig:MT-PE-emp-size-CSR}.
Observe that in the multiple triangle case
(which is more realistic than the one triangle case)
the empirical sizes are much closer to the nominal level
compared to the one triangle case.
For the right-sided alternative (i.e., against segregation),
the size is about the nominal level for $r \in (1.5,3)$,
and for the left-sided alternative (i.e., against association),
the size is about the nominal level for $r \in (1.1,2)$.
Furthermore,
although the empirical sizes for both right- and left-sided alternatives are about
the desired level for $r$ values between 1.5 and 2,
it seems that they are not very far from the nominal level for $r \in (1.5,10)$.
The test seems to be liberal for the right-sided alternative
and conservative for the left-sided alternative, if not at the desired level.

\subsection{Empirical Size Analysis for Central Similarity PCDs under CSR}
\label{sec:CS-emp-size-CSR}
\begin{figure}[ht]
\centering
\rotatebox{-90}{ \resizebox{2. in}{!}{\includegraphics{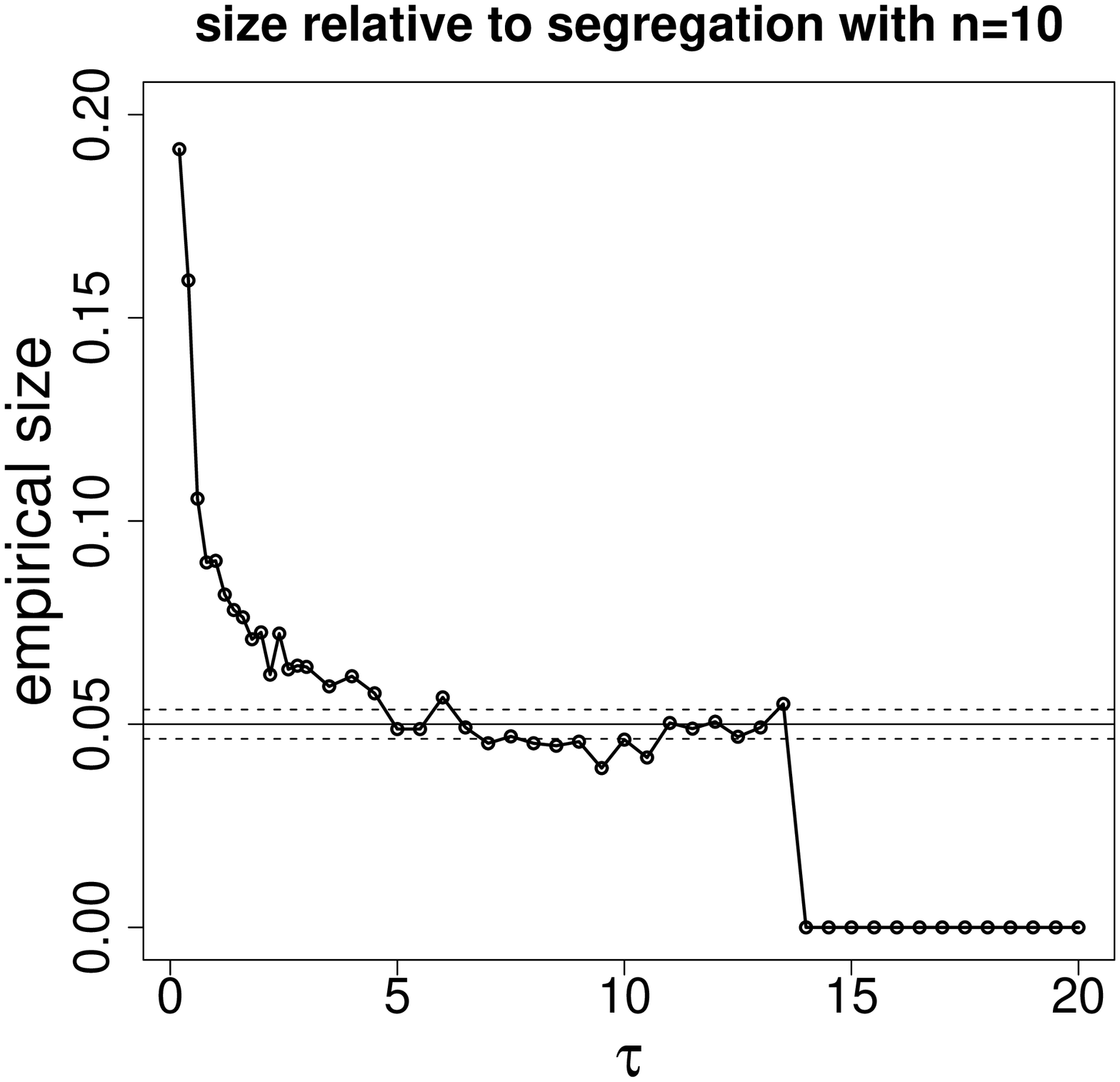} }}
\rotatebox{-90}{ \resizebox{2. in}{!}{\includegraphics{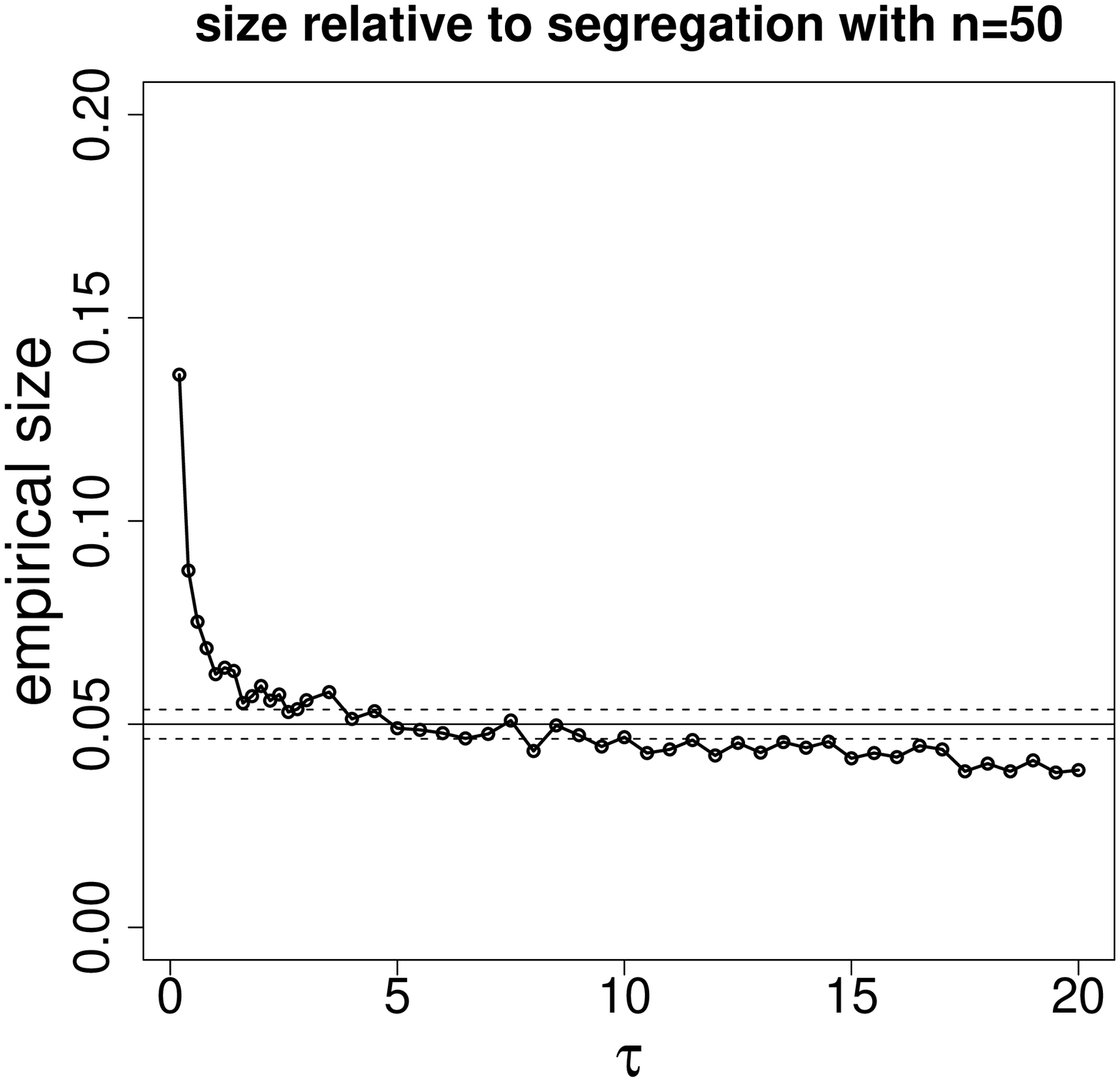} }}
\rotatebox{-90}{ \resizebox{2. in}{!}{\includegraphics{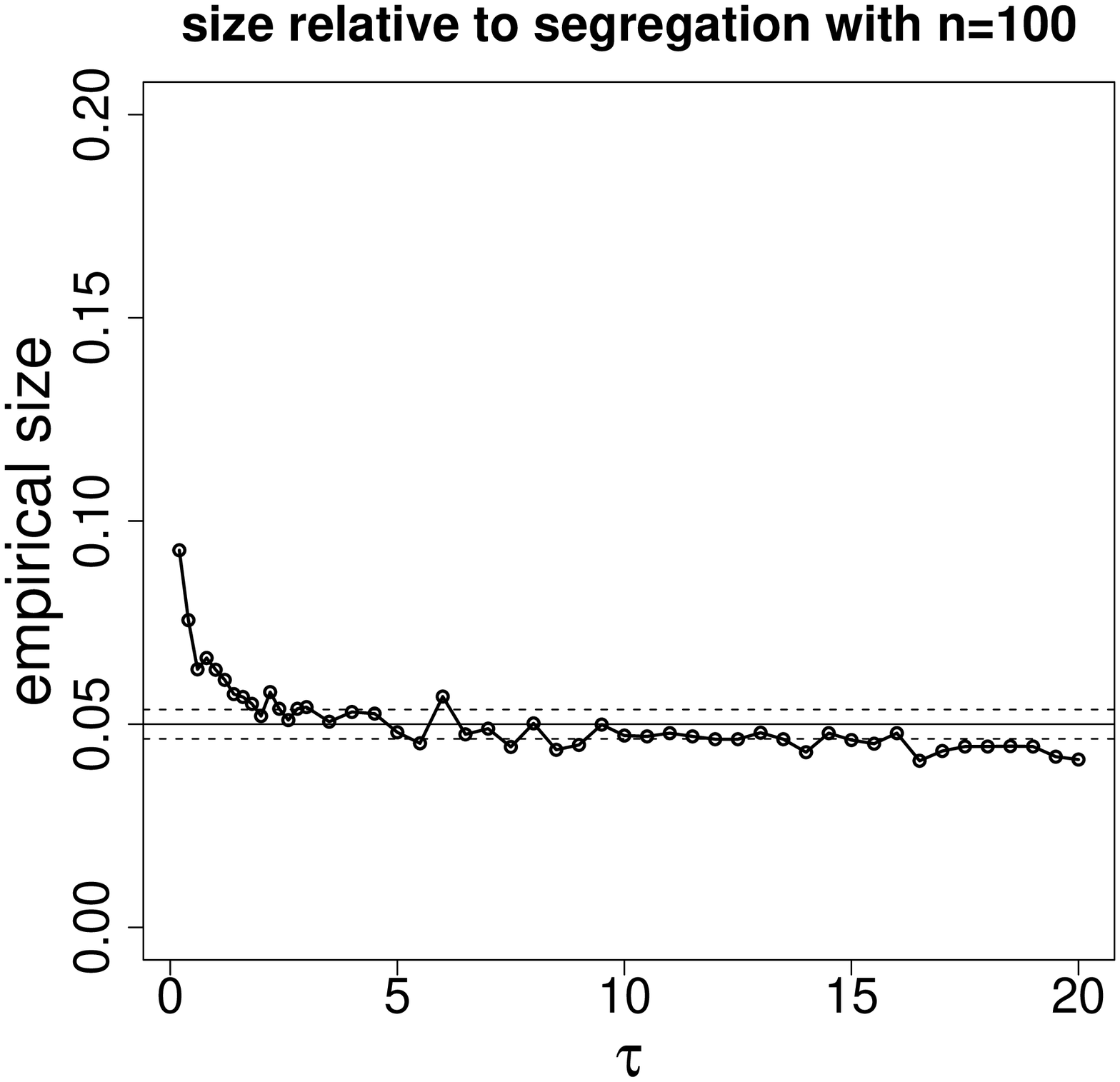} }}
\rotatebox{-90}{ \resizebox{2. in}{!}{\includegraphics{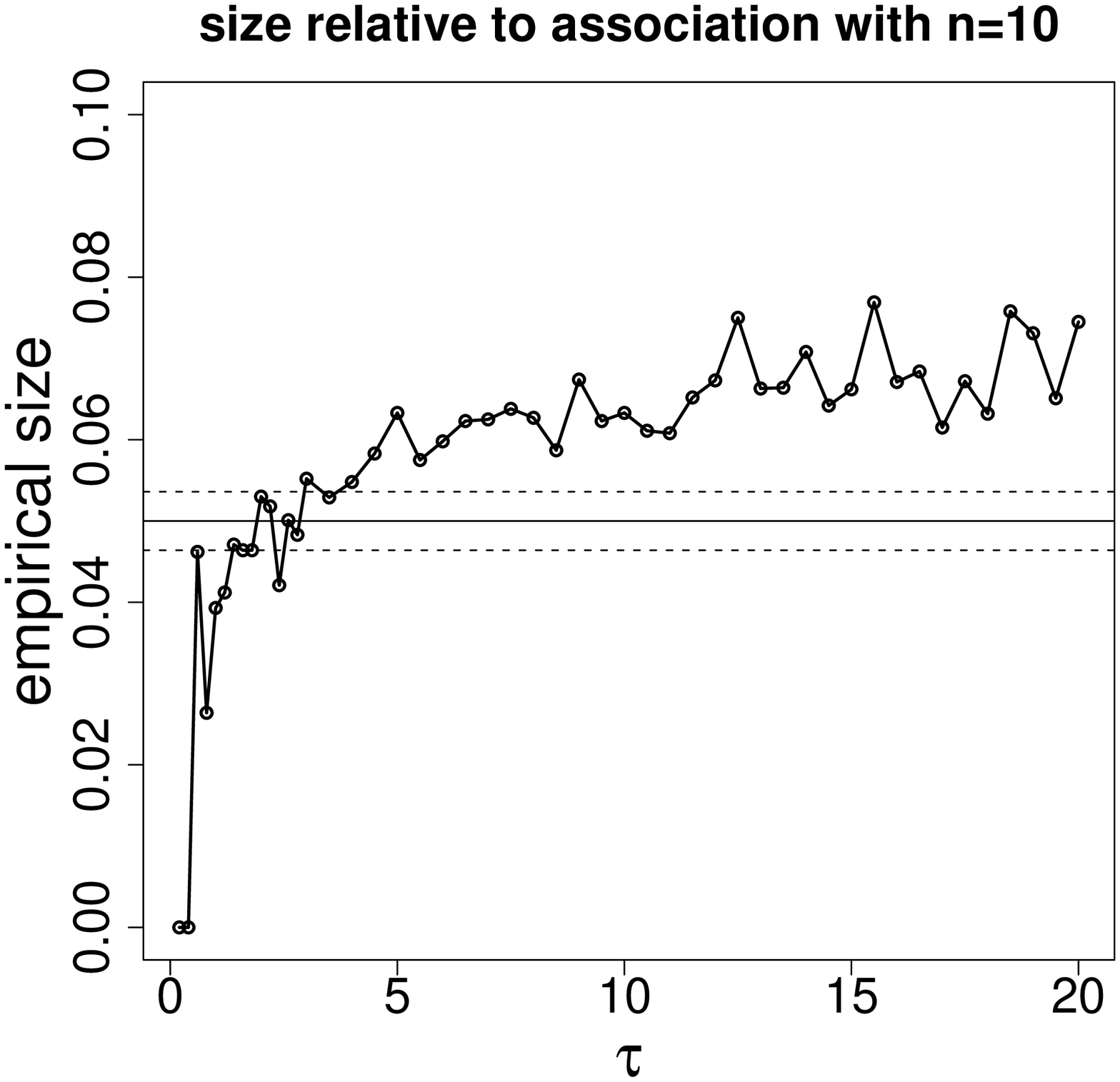} }}
\rotatebox{-90}{ \resizebox{2. in}{!}{\includegraphics{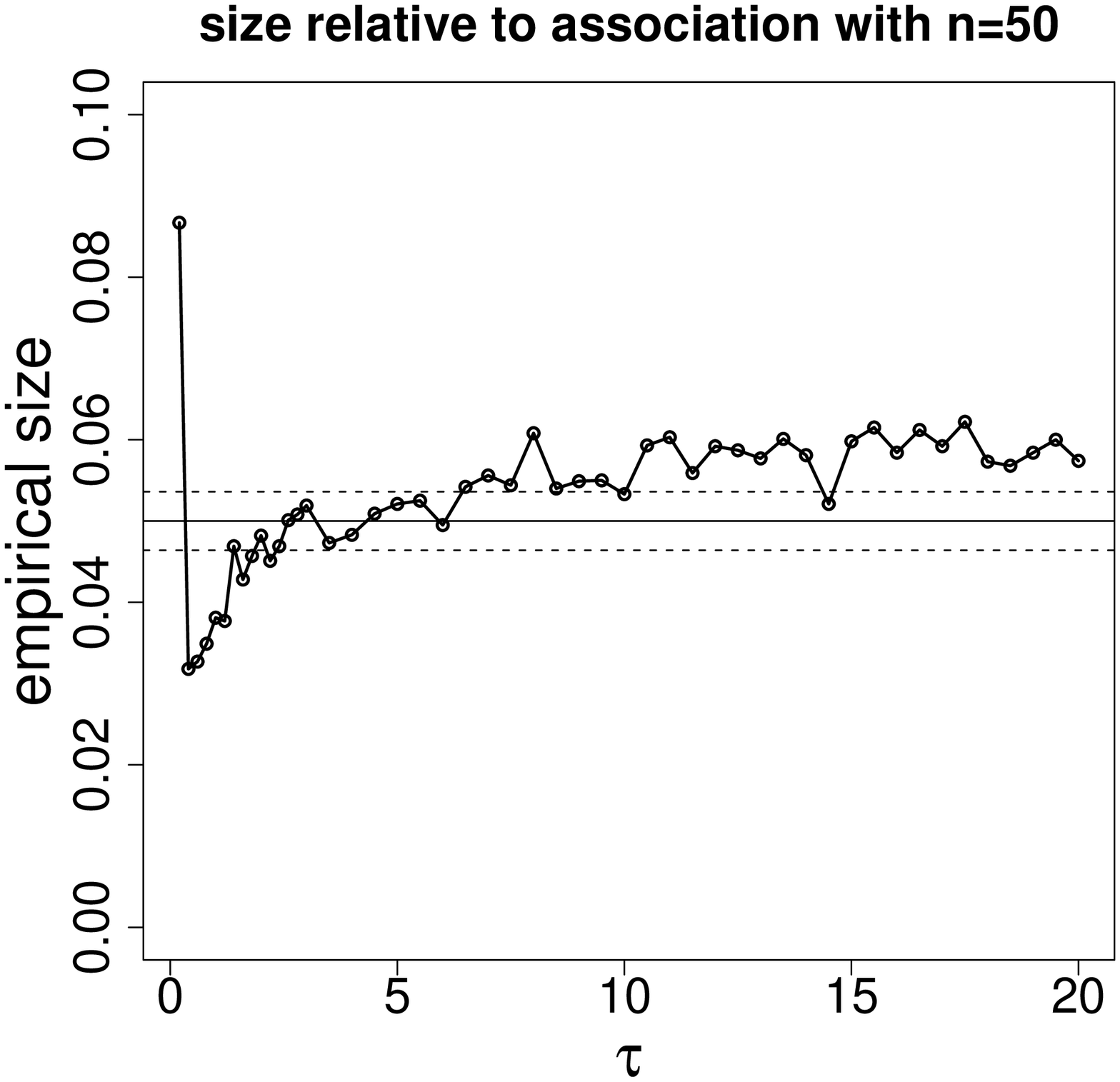} }}
\rotatebox{-90}{ \resizebox{2. in}{!}{\includegraphics{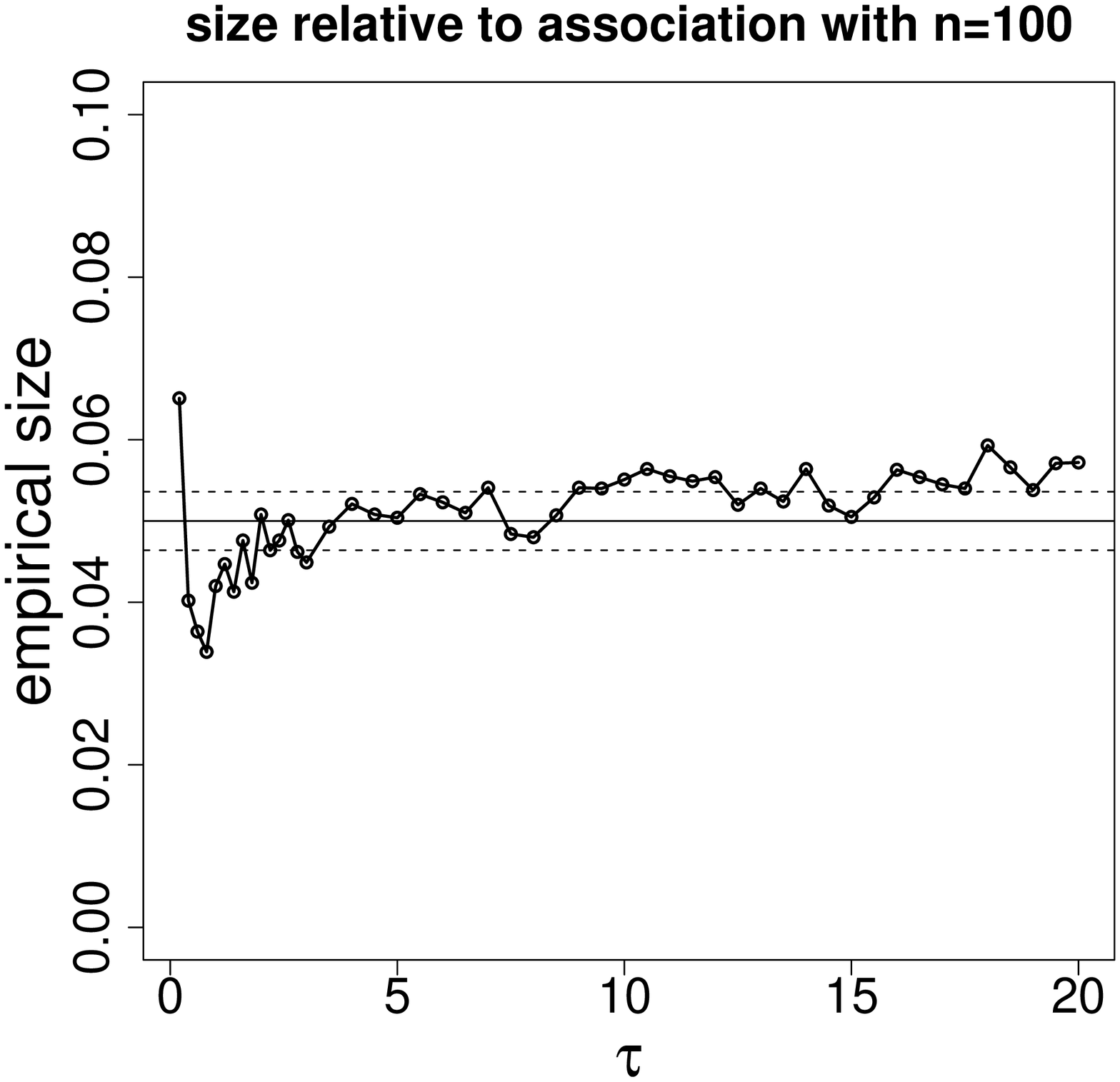} }}
\caption{
\label{fig:CS-emp-size-CSR}
\textbf{Empirical size for $R_{CS}(\tau)$ in the one triangle case:}
The empirical size estimates of the relative density of central similarity PCDs
in the one triangle case based on 10000 Monte Carlo replicates
for the left-sided alternative (i.e., relative to segregation) (top)
and the right-sided alternative (i.e., relative to association) (bottom)
with $n=10$ (left column), $n=50$ (middle column), and $n=100$ (right column) under the CSR pattern.
The horizontal lines are as in Figure \ref{fig:PE-emp-size-CSR}.
}
\end{figure}

In one and multiple triangle cases,
data generation is as in Section \ref{sec:PE-emp-size-CSR}
and we compute the relative density of central similarity PCDs
for $\tau= 0.2,0.4,0.6.\ldots,3.0,3.5,4.0,\ldots,20.0$ at each Monte Carlo replicate.
Let $R_{CS}(\tau)(\tau,j):=\frac{\sqrt{n}\,\bigl( \rho_{CS}(n,\tau,j)-\mu_{_{CS}}(\tau) \bigr)}{\sqrt{\nu_{_{CS}}(\tau)}}$
be the standardized relative density for Monte Carlo replicate $j$ with sample size $n$ for $j=1,2,\ldots,N_{mc}$.
For each $\tau$ value, the level $\alpha$ asymptotic critical value is
$\mu_{_{CS}}(\tau)+z_{(1-\alpha)} \cdot \sqrt{\nu_{_{CS}}(\tau)/n}$.
We estimate the empirical size against the segregation alternative as
$\frac{1}{N_{mc}}\sum_{j=1}^{N_{mc}}\I \left(R_{CS}(\tau)(\tau,j) > z_{1-\alpha} \right)$
and against the association alternative as
$\frac{1}{N_{mc}}\sum_{j=1}^{N_{mc}}\I \left(R_{CS}(\tau)(\tau,j) < z_{\alpha} \right)$.
In one triangle case,
the empirical sizes for the central similarity PCDs
together with upper and lower limits of liberalness and conservativeness
are plotted in Figure \ref{fig:CS-emp-size-CSR}.
Observe that as $n$ increases,
the empirical size gets closer to the nominal level of 0.05
(i.e., the normal approximation gets better).
For the right-sided tests,
the size is close to the nominal level for $\tau \in (5,14)$
and closest to 0.05 for $\tau \approx 5$ or $\tau \in (7,9)$ for all sample sizes;
for smaller $\tau$ values (i.e., $\tau \lesssim 4.5$)
the test seems to be liberal
with liberalness increasing as $\tau$ decreases;
and for $\tau \gtrsim 15$
the test is extremely conservative with size being virtually 0 for $n=10$
and the test is slightly conservative for $n=50$ and 100.
For larger $n$ (i.e., $n \ge 50$),
the test has the desired size for $\tau \ge 4$.
Considering all sample sizes,
we recommend $\tau \in (5,10)$ for testing against segregation.
For the left-sided tests
with $n=10$
the size is close to the nominal level for $\tau \in (2,4)$.
For $n=50$
the test has the desired size for $\tau \in (2,10)$
and for $n=100$
the test has the desired size for $\tau \in (2,15)$.
With all sample sizes,
the test seems to be conservative (slightly liberal)
for smaller (larger) $\tau$ values.
Considering all sample sizes,
we recommend $\tau \in (2.5,5)$ for testing against association.
The range of appropriate $\tau$ values gets wider with the increasing sample size
and
very large and small values of $\tau$ require much larger sample
sizes for the normal approximation to hold.

\begin{figure}[ht]
\centering
\rotatebox{-90}{ \resizebox{2. in}{!}{\includegraphics{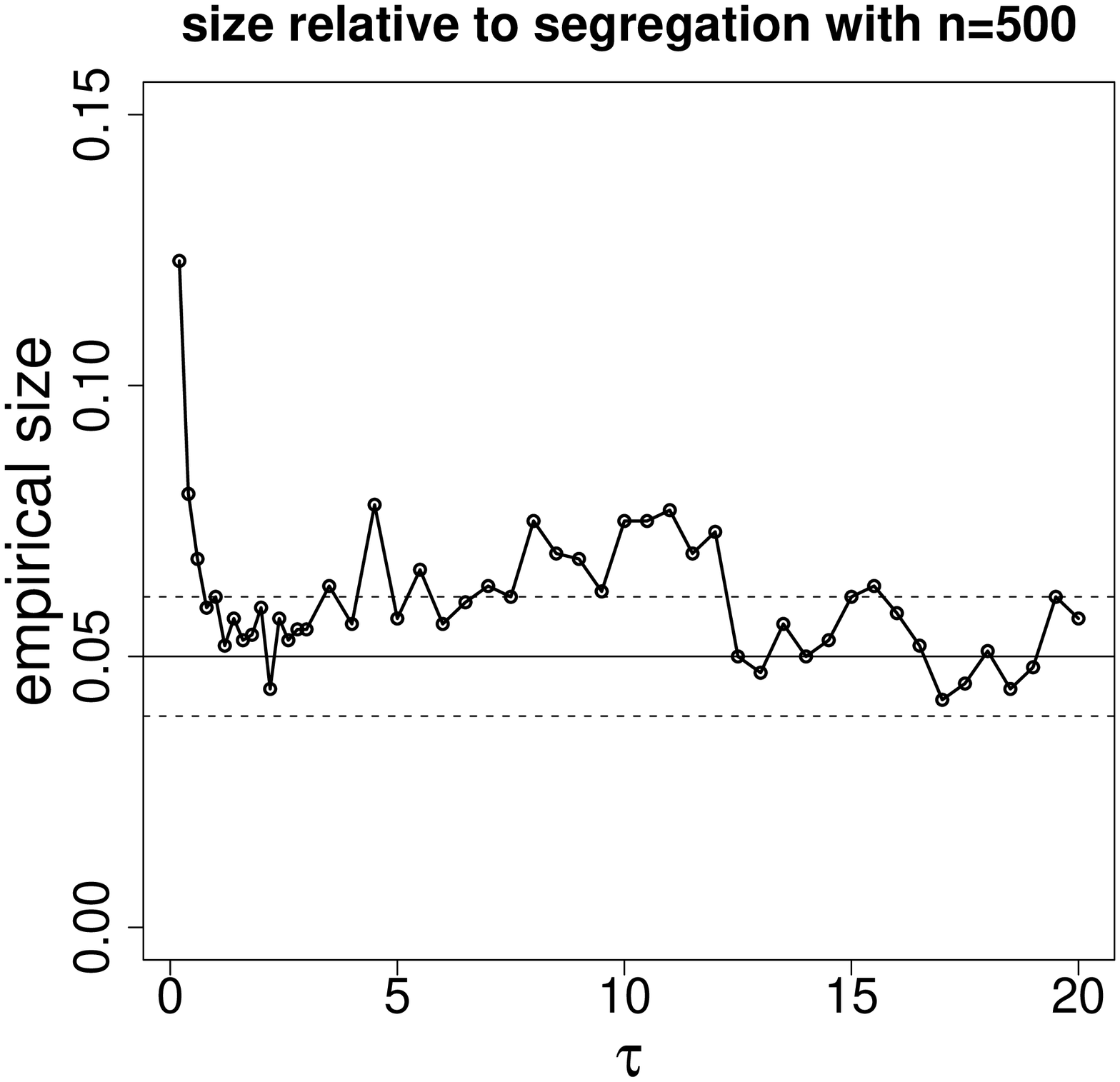} }}
\rotatebox{-90}{ \resizebox{2. in}{!}{\includegraphics{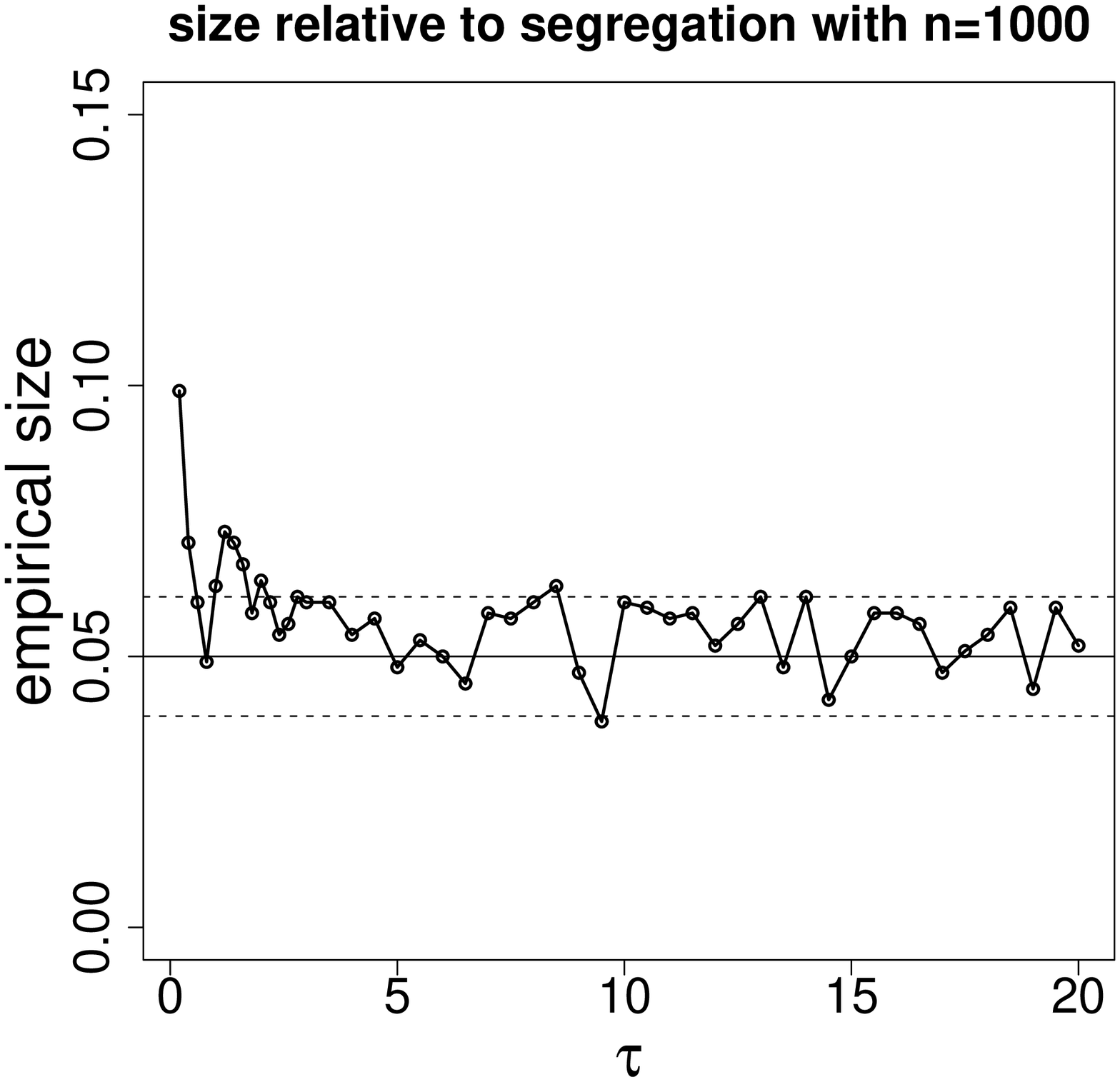} }}\\
\rotatebox{-90}{ \resizebox{2. in}{!}{\includegraphics{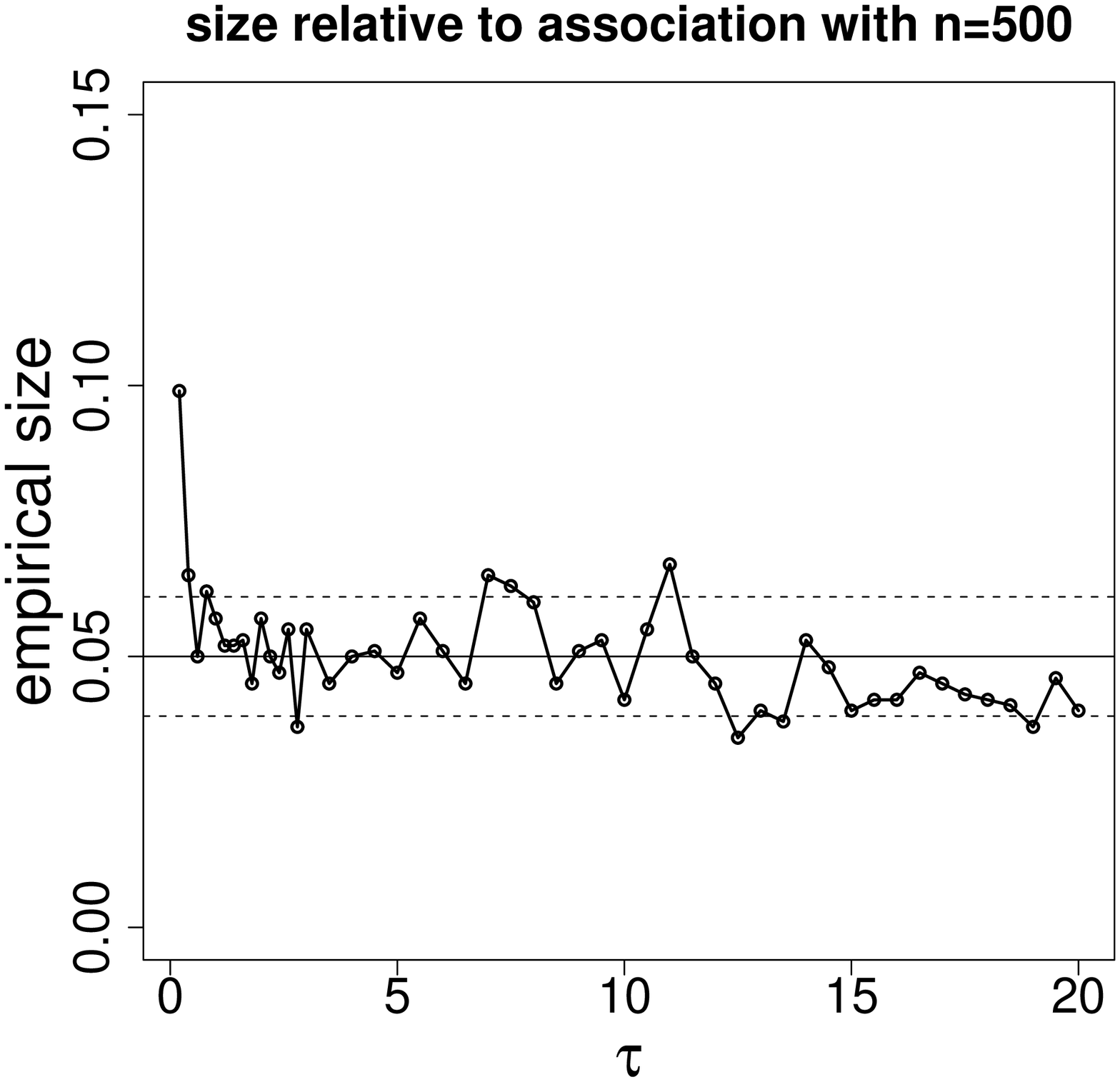} }}
\rotatebox{-90}{ \resizebox{2. in}{!}{\includegraphics{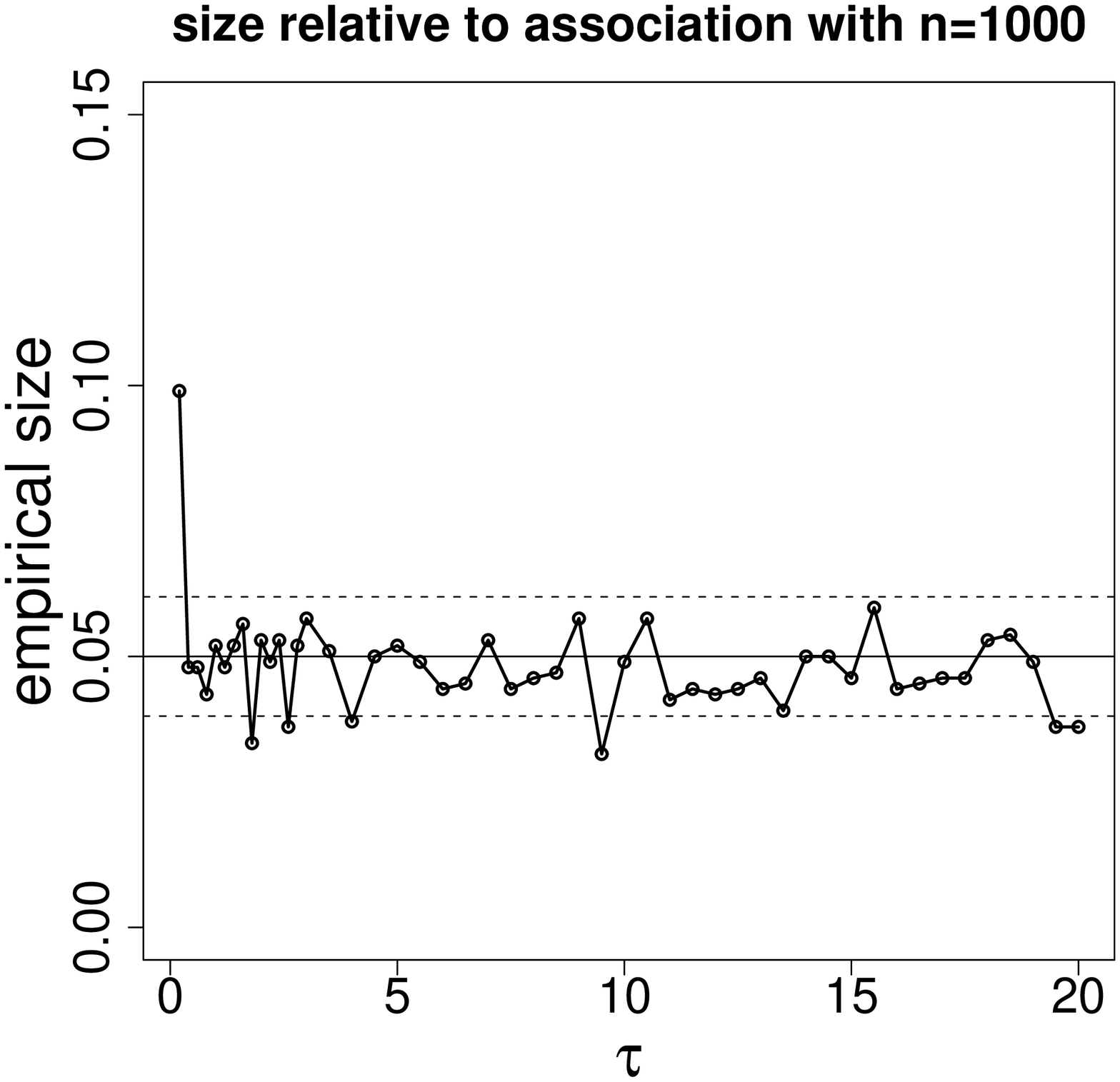} }}
\caption{
\label{fig:MT-CS-emp-size-CSR}
\textbf{Empirical size for $R_{CS}(\tau)$ in the multiple triangle case:}
The empirical size estimates of the relative density of central similarity PCDs
in the multiple triangle case based on 1000 Monte Carlo replicates
for the left-sided alternative (i.e., relative to segregation) (top)
and the right-sided alternative (i.e., relative to association) (bottom)
with $n=500$ (left column) and $n=1000$ (right column) under the CSR pattern.
The horizontal lines are as in Figure \ref{fig:MT-PE-emp-size-CSR}.
}
\end{figure}

In the multiple triangle case,
the empirical sizes for the central similarity PCDs
are plotted in Figure \ref{fig:MT-CS-emp-size-CSR}.
Observe that in the multiple triangle case
the empirical sizes are much closer to the nominal level
compared to the one triangle case.
Furthermore,
for the right-sided alternative with $n=500$,
the test has the desired level for $\tau \in (.8,4)$,
$\tau \approx 5$, and $\tau \in (12,20)$
and with $n=1000$ for $\tau \ge 2$.
Considering all sample sizes,
we recommend $\tau \in (2.5,8)$ for testing against segregation.
For the left-sided alternative,
with $n=500$,
$\tau \ge .5$ (except $\tau=7$ or 11)
seems to yield the appropriate level
and
with $n=1000$,
$\tau \ge .5$ seems to yield the appropriate level.
Considering all sample sizes,
we recommend $\tau \in (0.5,20)$ for testing against association.

\begin{remark}
\textbf{Empirical Size Comparison for the Two PCD Families:}
In the one triangle case,
the size estimates for the central similarity PCD
is close to the nominal level of 0.05 against the segregation alternative
for more of the expansion parameter values considered.
On the other hand,
the size estimates against association are close to the nominal level
for both PCD families,
but it seems that the size estimates for central similarity PCD
is closer to the nominal level.
In the multiple triangle case,
the size performance of the two PCD families is similar
and the size estimates are close to the nominal level
for both of the alternatives.
$\square$
\end{remark}

\section{Empirical Power Analysis under the Alternatives}
\label{sec:emp-power-anal}
To compare the power performance of the test statistics under the alternatives,
we generate $n$ $\X$ points uniformly in the corresponding support sets
as described in Section \ref{sec:alternatives}
and
provide the empirical power estimates of the tests under the segregation and
association alternatives.

\subsection{Empirical Power Analysis for Proportional-Edge PCDs under the Segregation Alternative}
\label{sec:PE-emp-power-seg}
In the one triangle case,
at each Monte Carlo replicate under segregation $H^S_{\ve}$,
we generate
$X_i \stackrel{iid}{\sim} \U\left(T_e \setminus \mathcal T_\ve\right)$,
for $i=1,2,\ldots,n$ for $n=10,50,100$.
At each Monte Carlo replicate,
we compute the relative density
of the proportional-edge PCDs.
We consider $r \in \{1,11/10,6/5,4/3,\sqrt{2},3/2,2,3,5,10\}$
for the proportional-edge PCD.
We repeat the above simulation procedure $N_{mc}=10000$ times.
We consider $\ve \in \{ \sqrt{3}/8, \sqrt{3}/4, 2 \, \sqrt{3}/7 \}$
(which correspond to 18.75 \%, 75 \%, and $4500/49\approx 91.84$ \%
of the triangle (around the vertices)
being unoccupied by the $\X$ points, respectively)
for the segregation alternatives.

Under segregation alternatives with $\ve>0$,
the distribution of $\rho_{_{PE}}(n,r)$ is degenerate
for large values of $r$.
For a given $\ve \in (0,\sqrt{3}/4)$,
the corresponding digraph is complete almost surely,
for $r \ge \frac{\sqrt{3}}{2\ve}$,
hence $\rho_{_{PE}}(n,r)=1$ a.s.
For $\ve \in (\sqrt{3}/4,\sqrt{3}/3)$,
the corresponding digraph is complete almost surely,
for $r \ge \frac{\sqrt{3}-2\ve}{\ve}$.
In particular,
for $\ve=\sqrt{3}/8$, $\rho_{_{PE}}(n,r)$ is degenerate for $r \ge 4$,
for $\ve=\sqrt{3}/4$, $\rho_{_{PE}}(n,r)$ is degenerate for $r \ge 2$,
and
for $\ve=2 \, \sqrt{3}/7$, $\rho_{_{PE}}(n,r)$ is degenerate for $r \ge 3/2$,

\begin{figure}[]
\centering
\psfrag{kernel density estimate}{ \Huge{\bfseries{kernel density estimate}}}
\psfrag{relative density}{ \Huge{\bfseries{relative density}}}
\rotatebox{-90}{ \resizebox{2. in}{!}{ \includegraphics{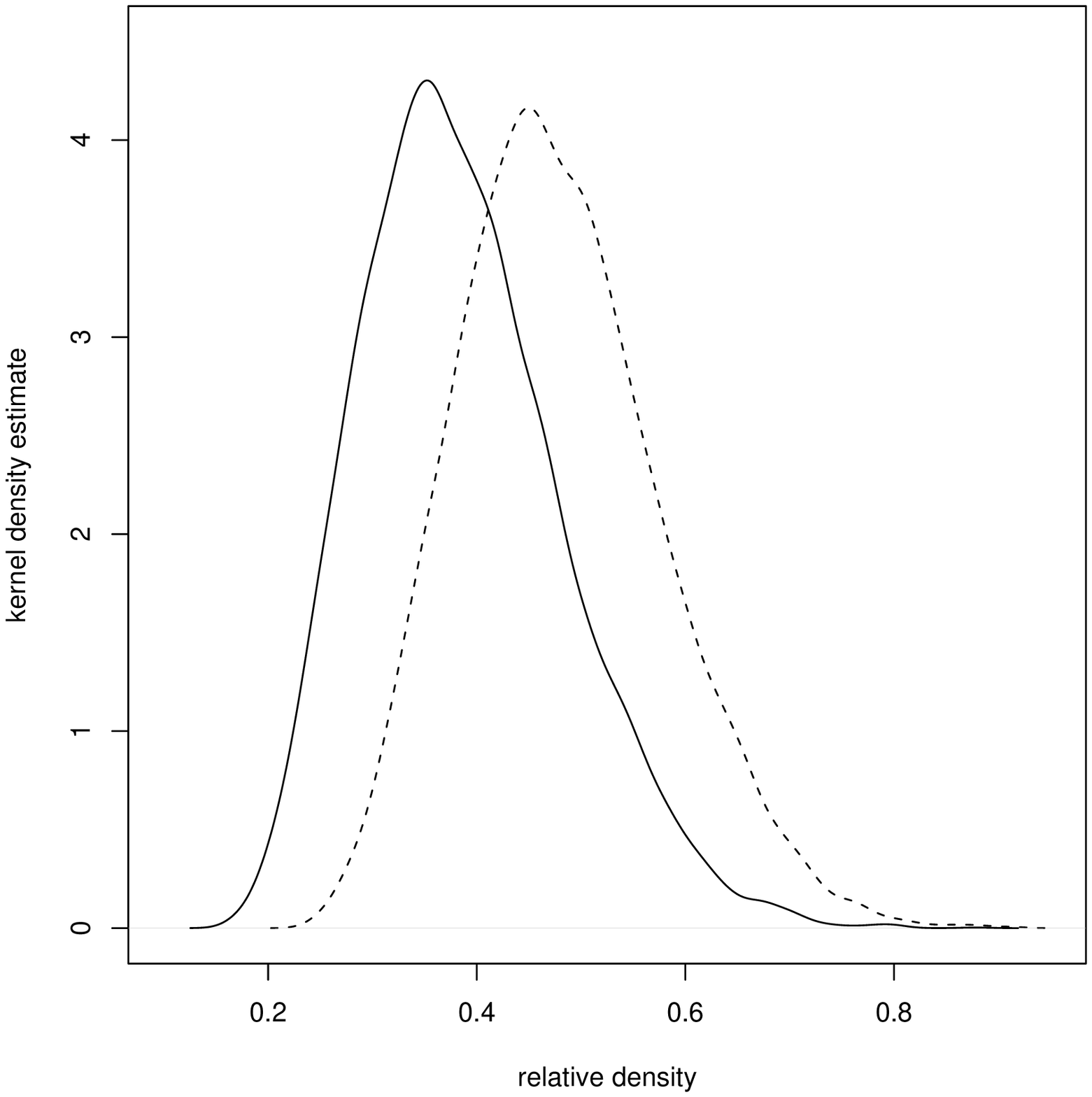}}}
\rotatebox{-90}{ \resizebox{2. in}{!}{ \includegraphics{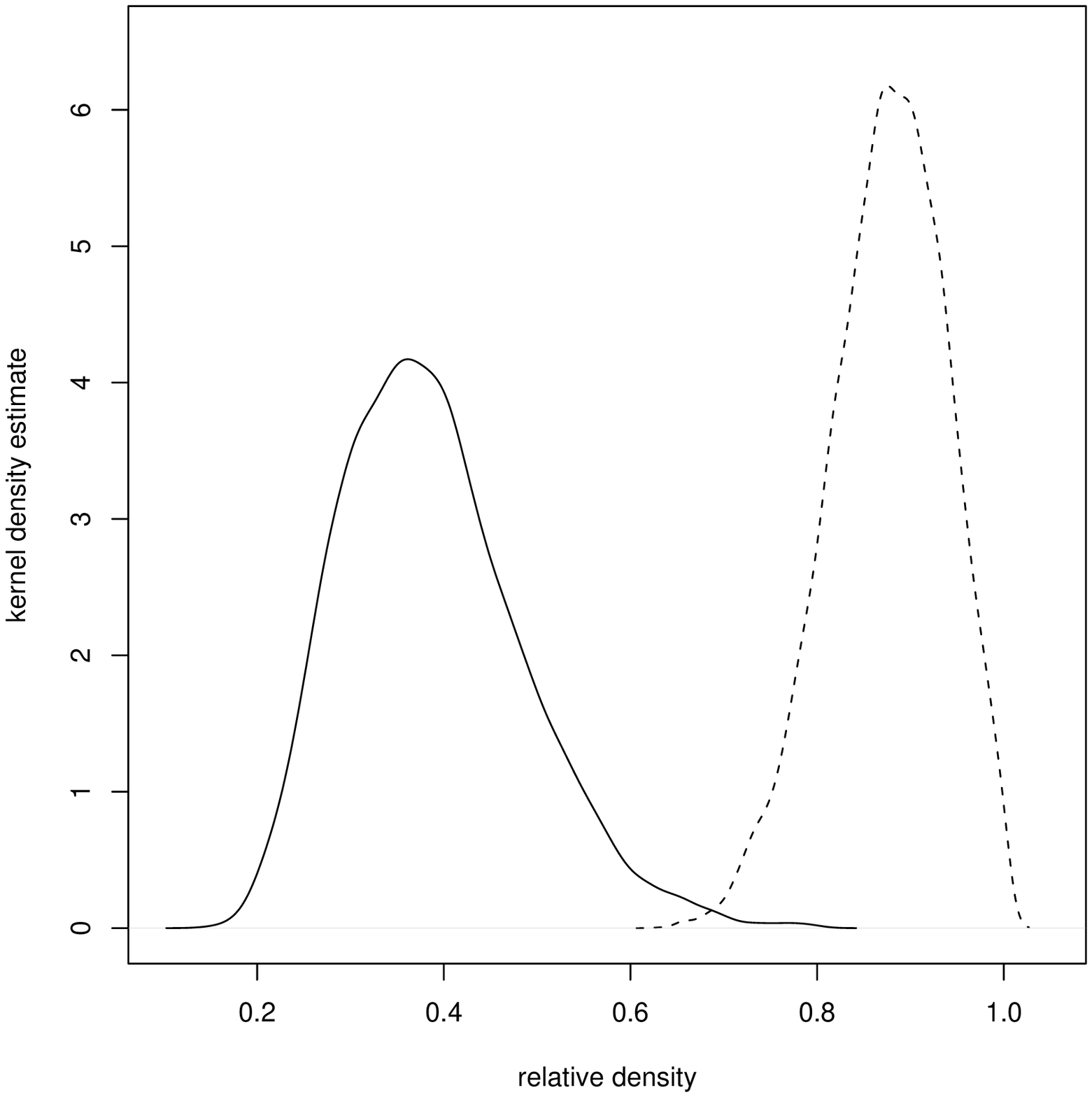}}}
\caption{\label{fig:segsim}
Kernel density estimates of the relative density of proportional-edge PCD,
$\rho_{_{PE}}(n,r)$,
under the null (solid line) and the segregation alternatives (dashed line) with
$H^S_{\sqrt{3}/8}$ (left) and $H^S_{\sqrt{3}/4}$ (right) for $r=3/2$
with $n=10$ based on $N_{mc}=10000$ replicates.
}
\end{figure}

\begin{figure}[ht]
\centering
\psfrag{kernel density estimate}{ \Huge{\bfseries{kernel density estimate}}}
\psfrag{relative density}{ \Huge{\bfseries{relative density}}}
\rotatebox{-90}{ \resizebox{2. in}{!}{ \includegraphics{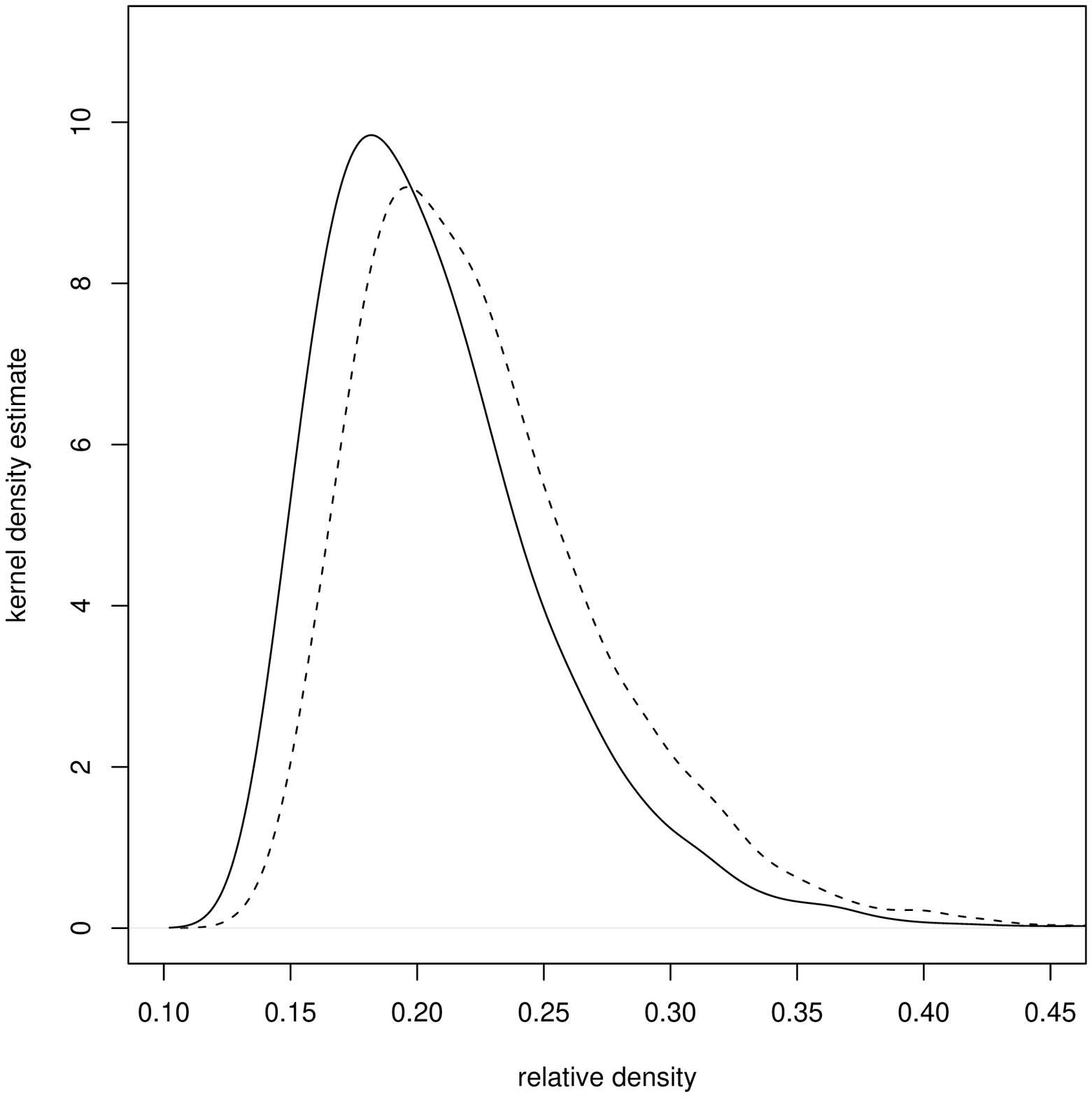}}}
\rotatebox{-90}{ \resizebox{2. in}{!}{ \includegraphics{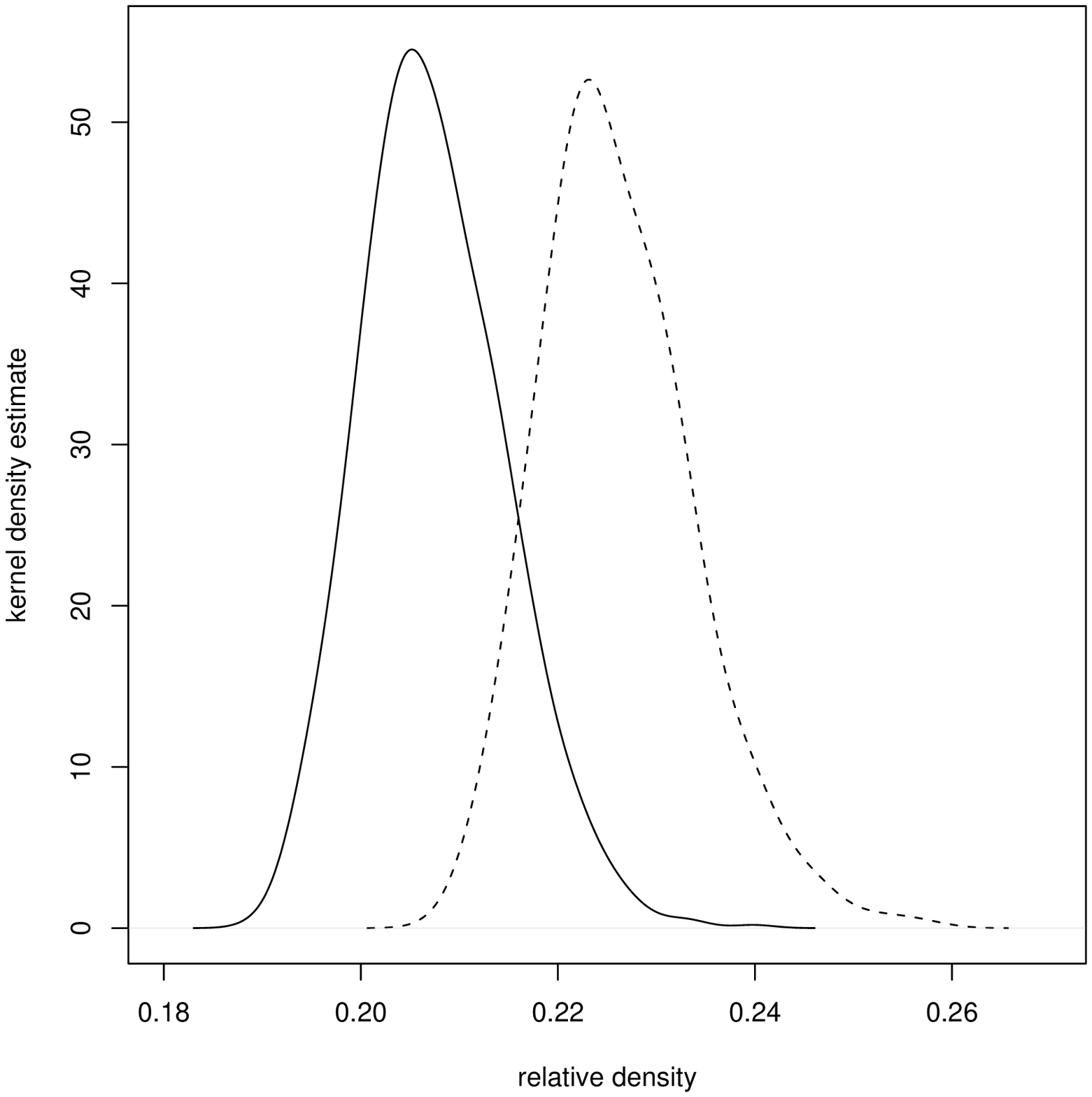}}}
\caption{ \label{fig:SegSimPowerPlots}
Depicted are kernel density estimates for $\rho_{_{PE}}(n,11/10)$ for
$n=10$ (left) and $n=100$ (right) under the null (solid line) and segregation alternative $H^S_{\sqrt{3}/8}$ (dashed line).
}
\end{figure}

In the one triangle case,
we plot
the kernel density estimates for the null
case and the segregation alternative with $\ve=\sqrt{3}/8$ and $\ve=\sqrt{3}/4$
with $n=10$ and $N_{mc}=10000$
in Figure \ref{fig:segsim}.
Observe that under
both $H_o$ and alternatives, kernel density estimates
are almost symmetric for $r=3/2$.
Moreover, there is much more separation between
the kernel density estimates of the null and alternatives
for $\ve=\sqrt{3}/4$ compared to $\ve=\sqrt{3}/8$,
implying more power for larger $\ve$ values.
In Figure \ref{fig:SegSimPowerPlots}, we present a Monte Carlo investigation
against the segregation alternative $H^S_{\sqrt{3}/8}$ for $r=11/10$, and $n=10$, $N_{mc}=10000$ (left), $n=100$, $N_{mc}=1000$ (right).
With $n=10$, the null and alternative kernel density functions
for $\rho_{10}(11/10)$ are very similar, implying small power.
With $n=100$,
there is more separation
between null and alternative kernel density functions
implying higher power.
Notice also that the probability density functions are more skewed for $n=10$,
while approximate normality holds for $n=100$.

\begin{figure}[ht]
\centering
\rotatebox{-90}{ \resizebox{1.7 in}{!}{ \includegraphics{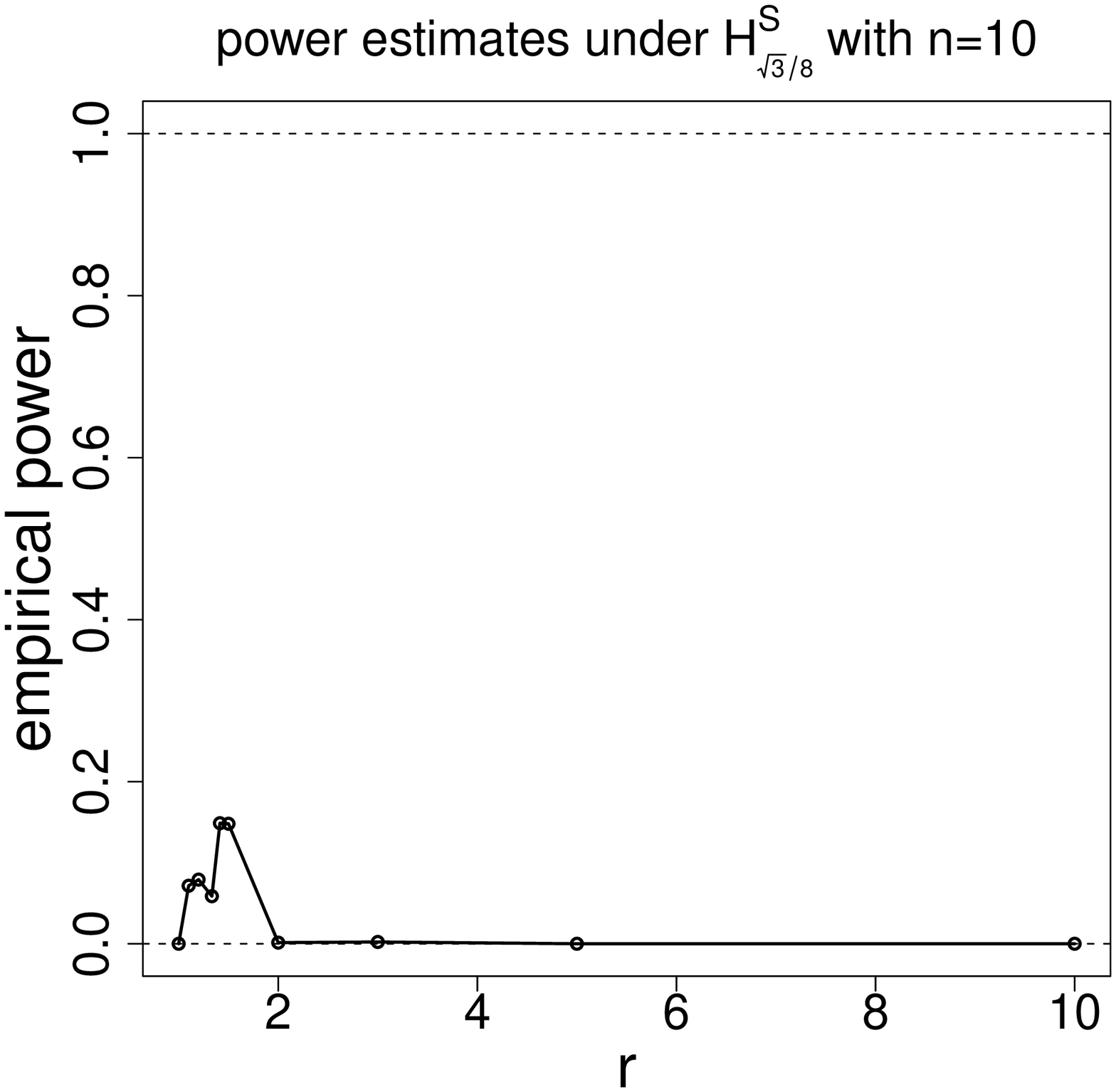}}}
\rotatebox{-90}{ \resizebox{1.7 in}{!}{ \includegraphics{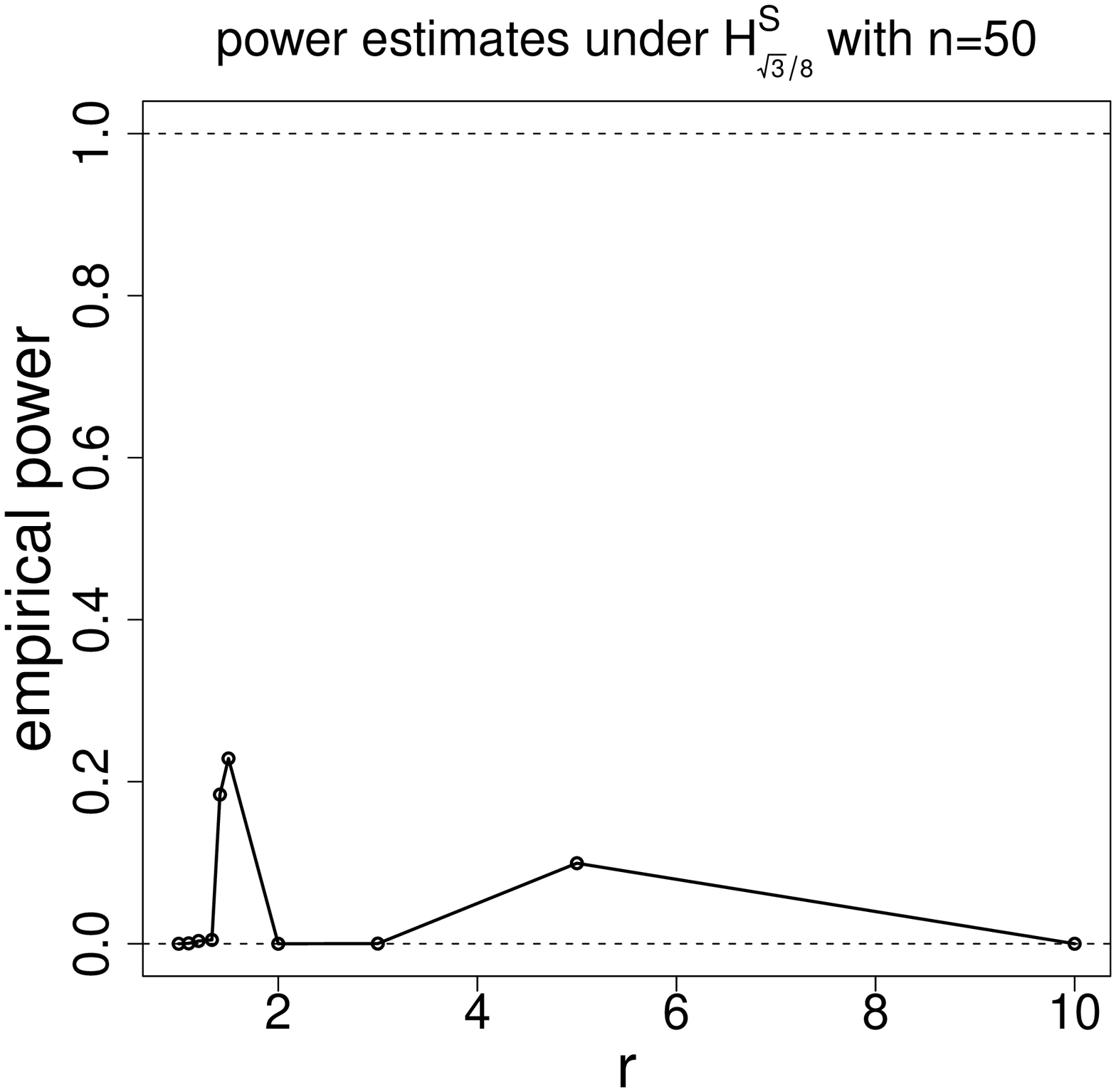}}}
\rotatebox{-90}{ \resizebox{1.7 in}{!}{ \includegraphics{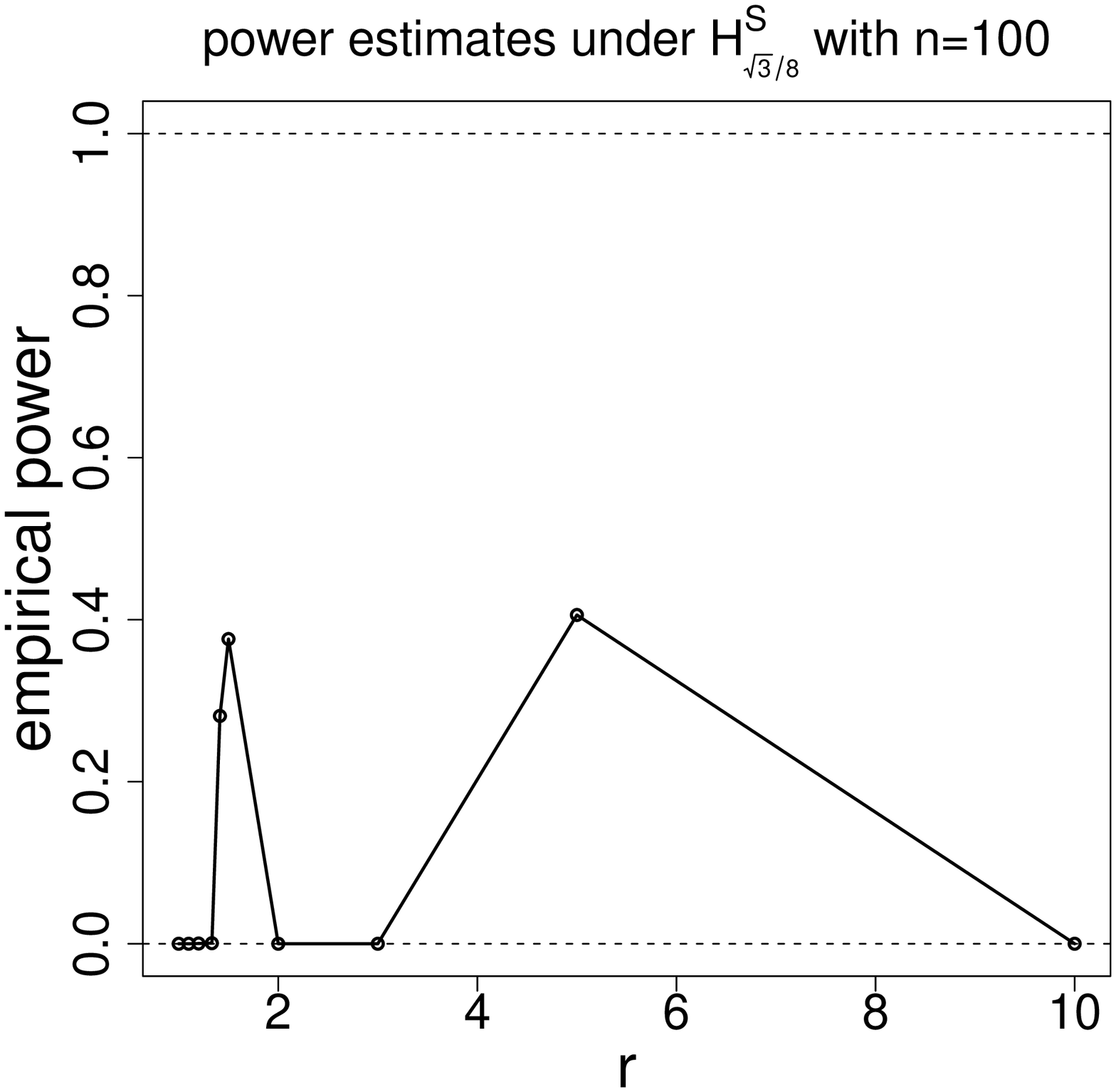}}}
\rotatebox{-90}{ \resizebox{1.7 in}{!}{ \includegraphics{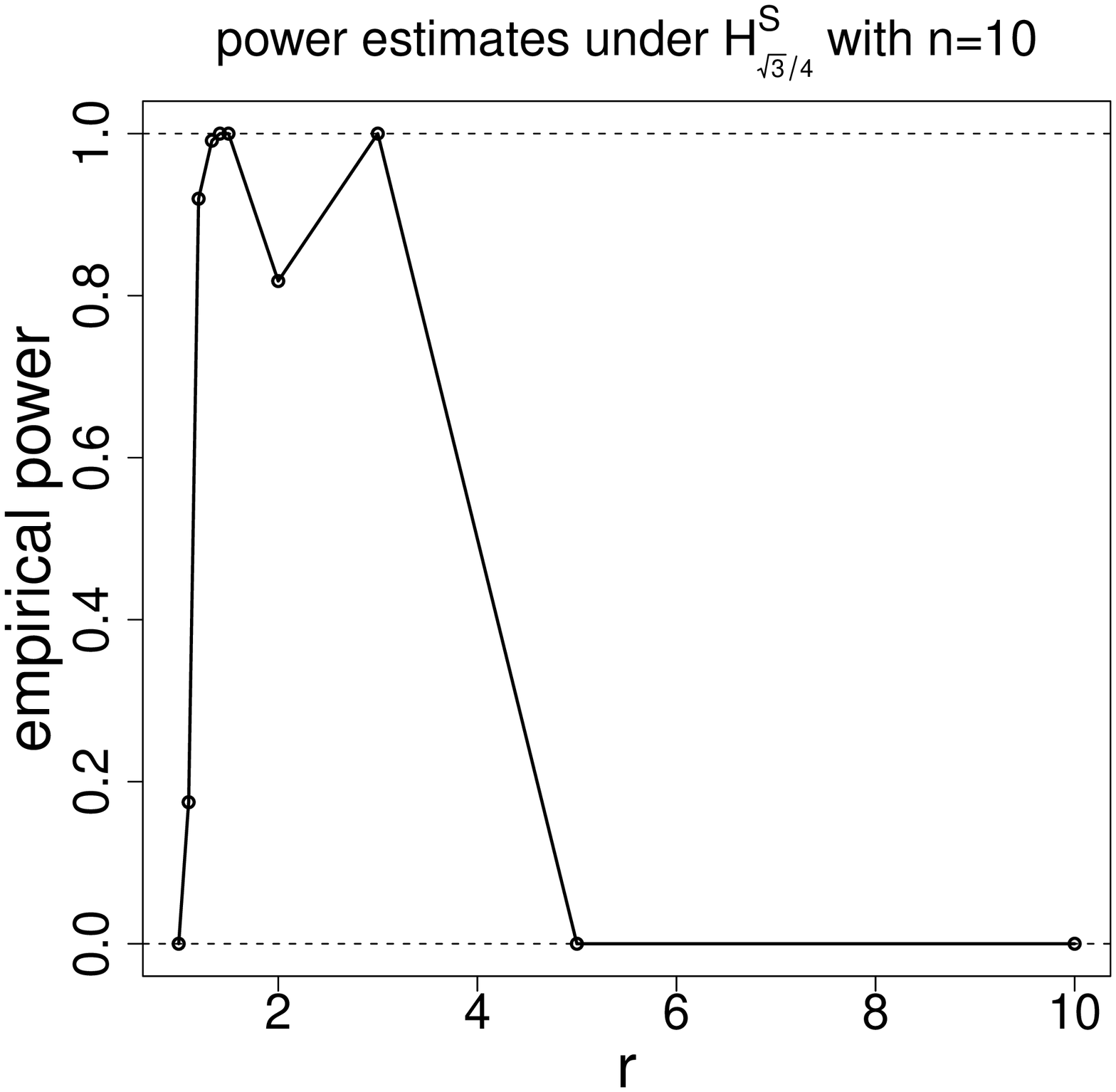}}}
\rotatebox{-90}{ \resizebox{1.7 in}{!}{ \includegraphics{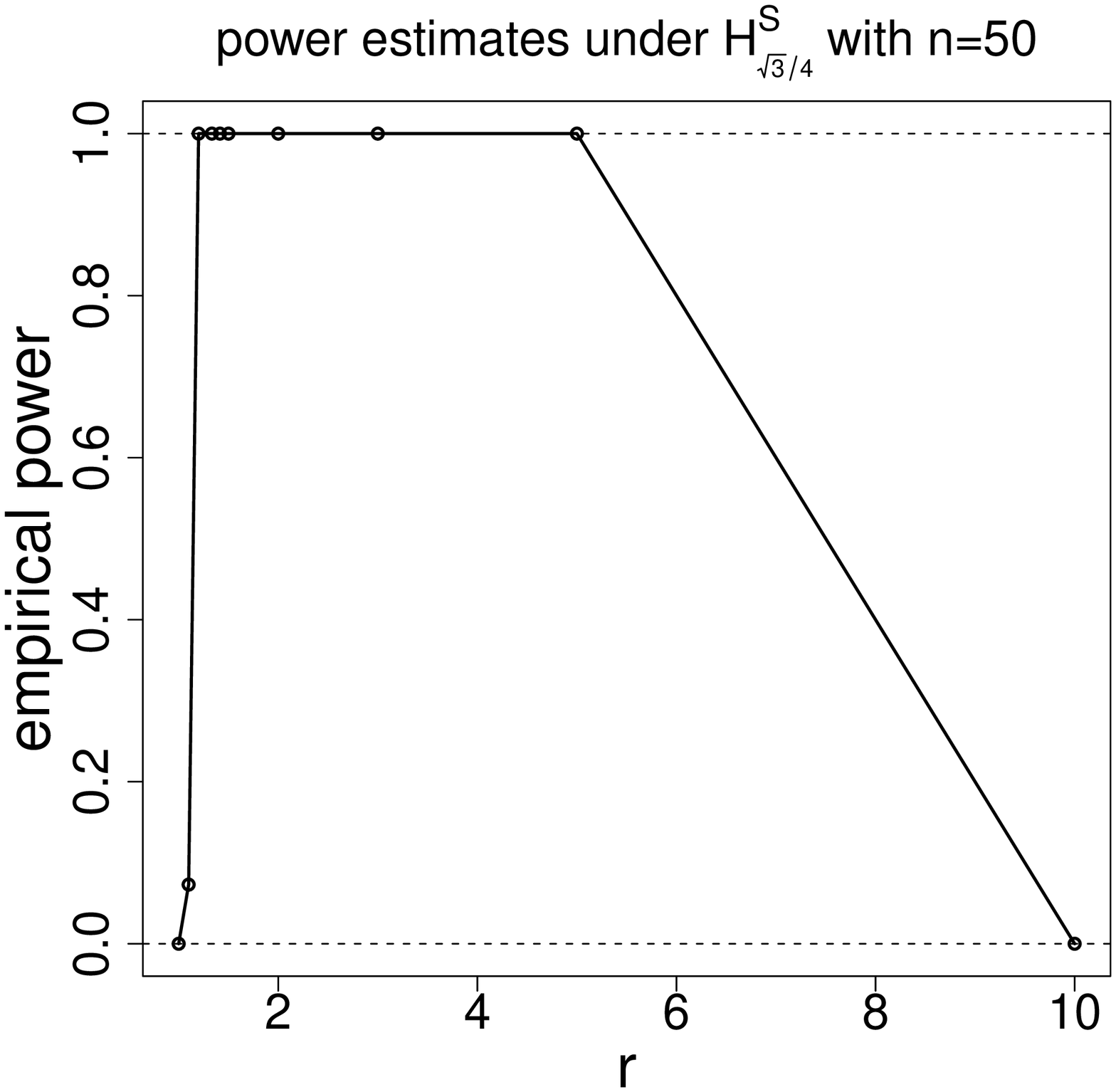}}}
\rotatebox{-90}{ \resizebox{1.7 in}{!}{ \includegraphics{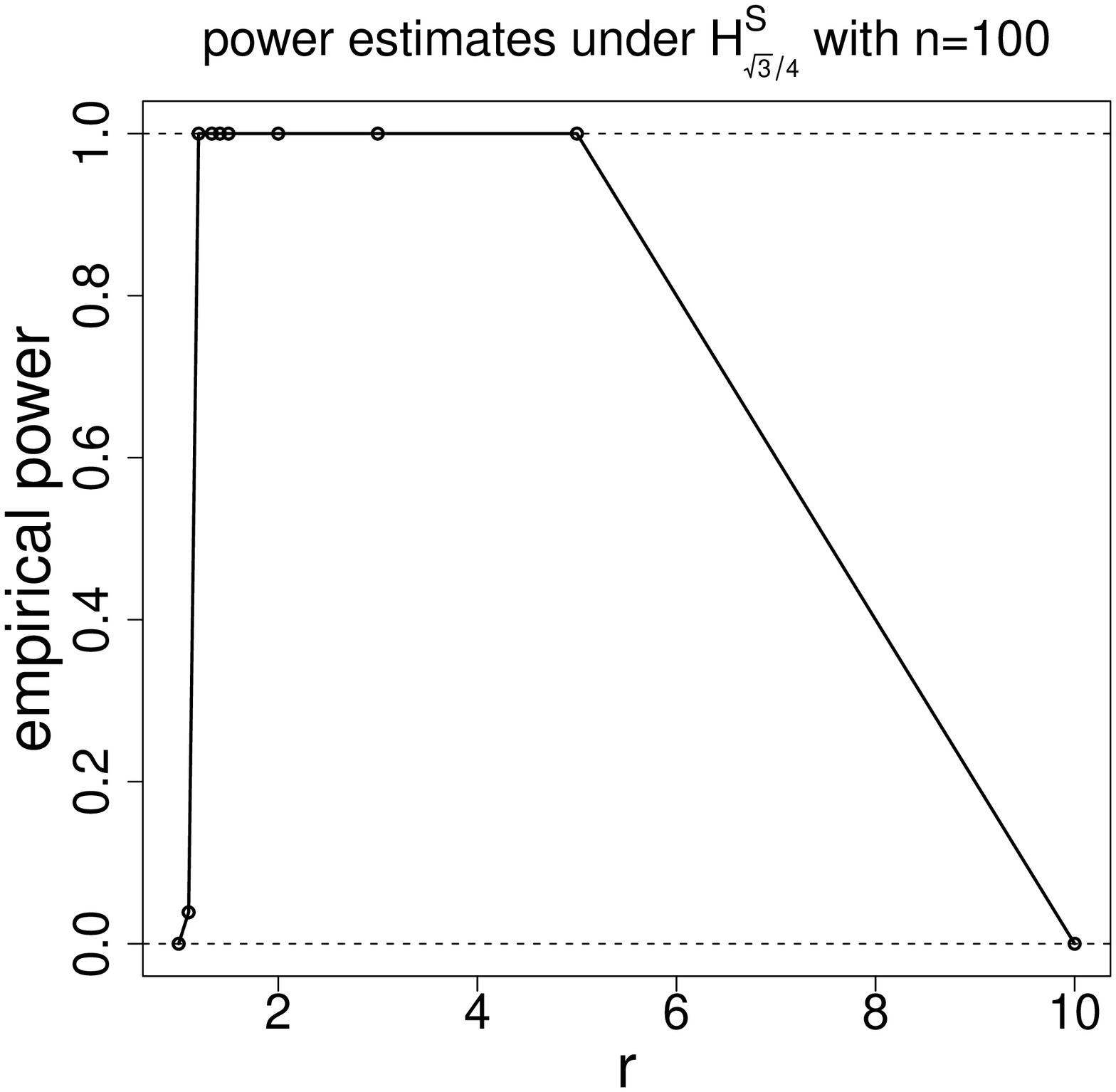}}}
\rotatebox{-90}{ \resizebox{1.7 in}{!}{ \includegraphics{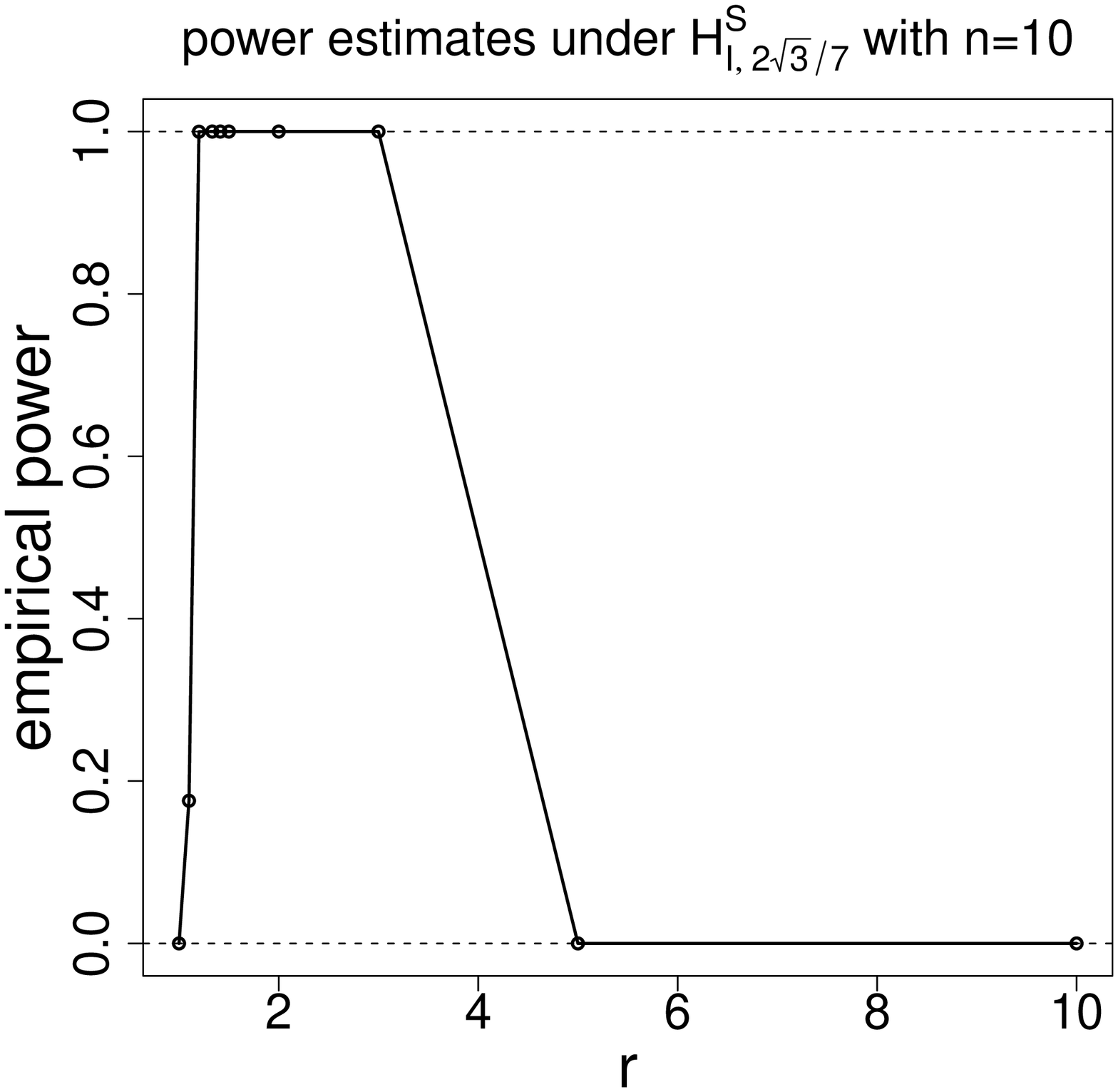}}}
\rotatebox{-90}{ \resizebox{1.7 in}{!}{ \includegraphics{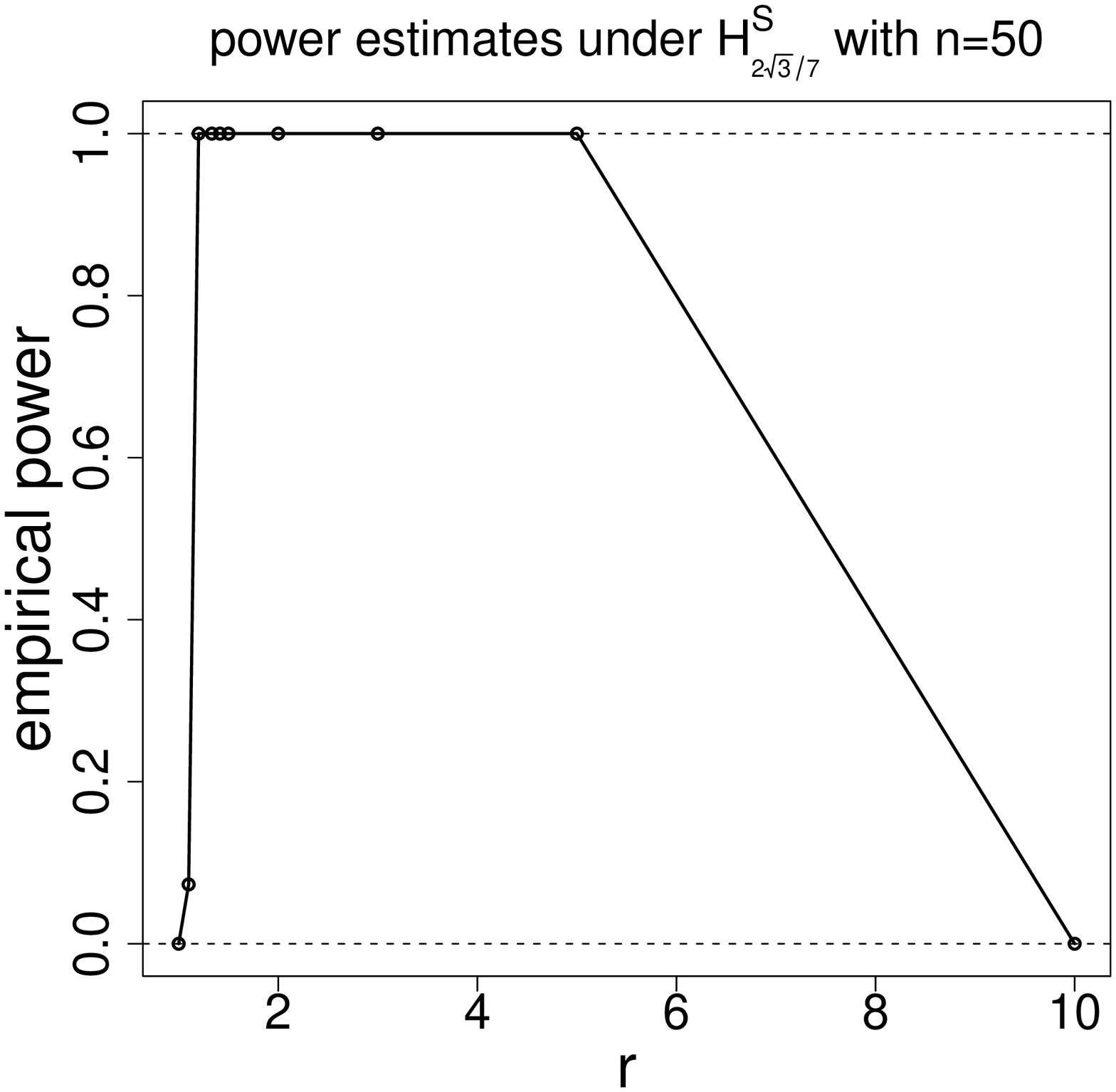}}}
\rotatebox{-90}{ \resizebox{1.7 in}{!}{ \includegraphics{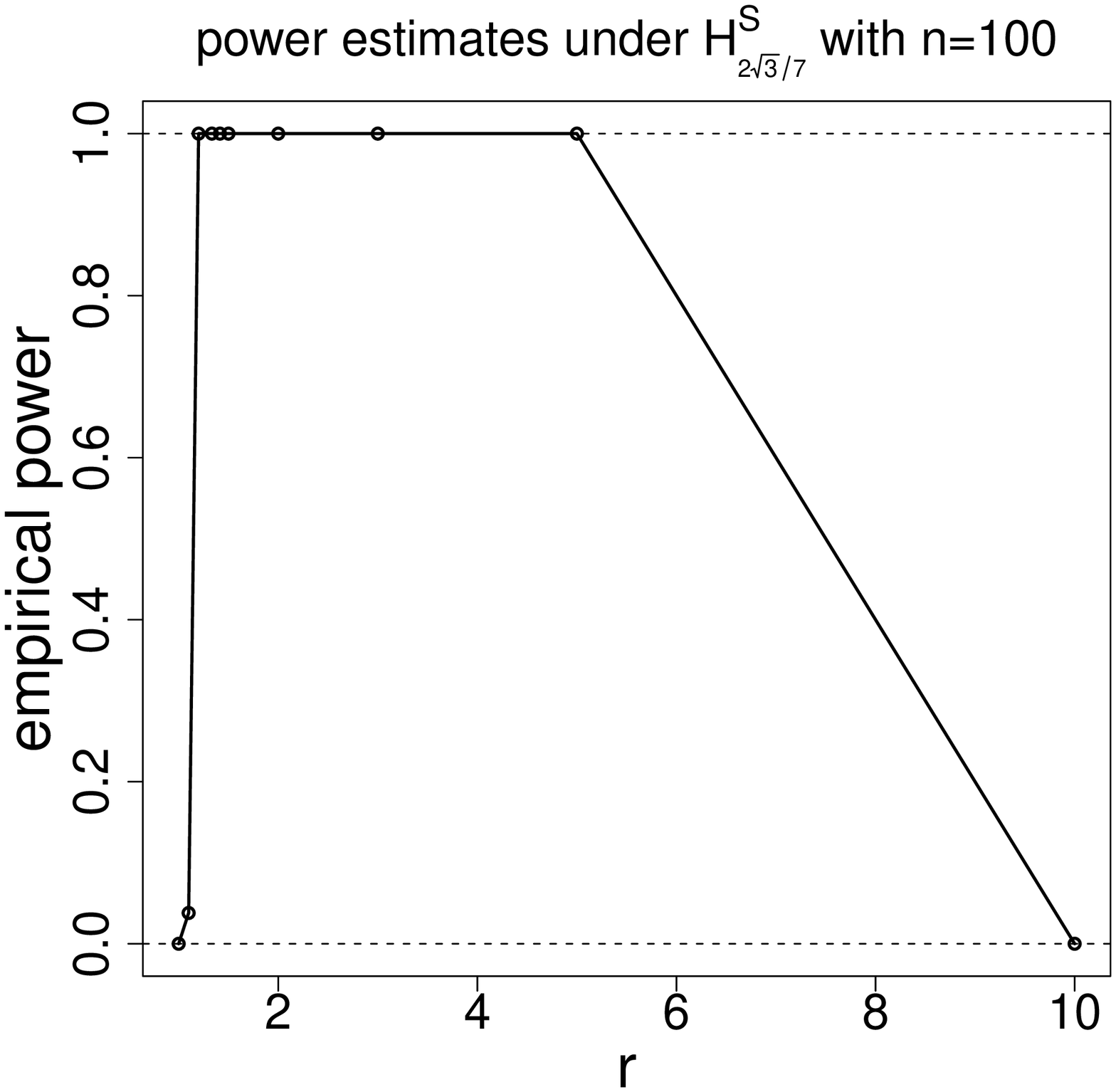}}}
\caption{\label{fig:PE-emp-power-seg}
\textbf{Empirical power for $R_{PE}(r)$ in the one triangle case:}
Monte Carlo power estimates for relative density of proportional-edge PCDs
in the one triangle case
using the asymptotic critical value against segregation alternatives
$H^S_{\sqrt{3}/8}$ (top row),
$H^S_{\sqrt{3}/4}$ (middle row),
and
$H^S_{2\,\sqrt{3}/7}$ (bottom row)
as a function of $r$, for $n=10$ (left column), $n=50$ (middle column), and $n=100$ (right column).
}
\end{figure}

For a given alternative and sample size, we analyze
the empirical power of the test based on $\rho_{_{PE}}(n,r)$ --- using the asymptotic critical value--- as
a function of the expansion parameter $r$.
We estimate the empirical power as
$\frac{1}{N_{mc}}\sum_{j=1}^{N_{mc}}\I \left(R_{PE}(r)(r,j) > z_{1-\alpha} \right)$.
In Figure \ref{fig:PE-emp-power-seg},
we present Monte Carlo power estimates for relative density of proportional-edge PCDs
in the one triangle case
against $H^S_{\sqrt{3}/8}$,
$H^S_{\sqrt{3}/4}$, and $H^S_{2\,\sqrt{3}/7}$ as a function of $r$ for $n=10,50,100$.
Notice that Monte Carlo power estimate
increases as $r$ gets larger and then decreases, due to the
magnitude of $r$ and $n$.
Because for small $n$ and large $r$, the
critical value is approximately 1 under $H_o$, as we get a complete
digraph with high probability.
Moreover, the more severe the segregation, the higher the power estimate at each $r$.
Under mild segregation (with $\ve=\sqrt{3}/8$),
$r$ around 1.5 or 5 yields the highest power
(for other $r$ values, the power performance is very poor).
Furthermore,
under moderate to severe segregation,
with $n=10$
the power estimate seems to be close to 1 for $r \in (1,4)$,
and
with $n=50$ or 100
the power estimate seems to be close to 1 for $r \in (1,5)$.
However, the power estimates are valid only for $r$ within $(2,3)$,
since the test has the desired size for this range of $r$ values
against the right-sided alternative.
So, for small sample sizes,
$r \approx 1.5$ is recommended,
and for larger sample sizes,
moderate values of $r$ (i.e., $r \in (2,3)$) are recommended for the segregation alternative
as they are more appropriate for normal
approximation and they yield the desired significance level.

\begin{figure}[]
\centering
\rotatebox{-90}{ \resizebox{1.7 in}{!}{ \includegraphics{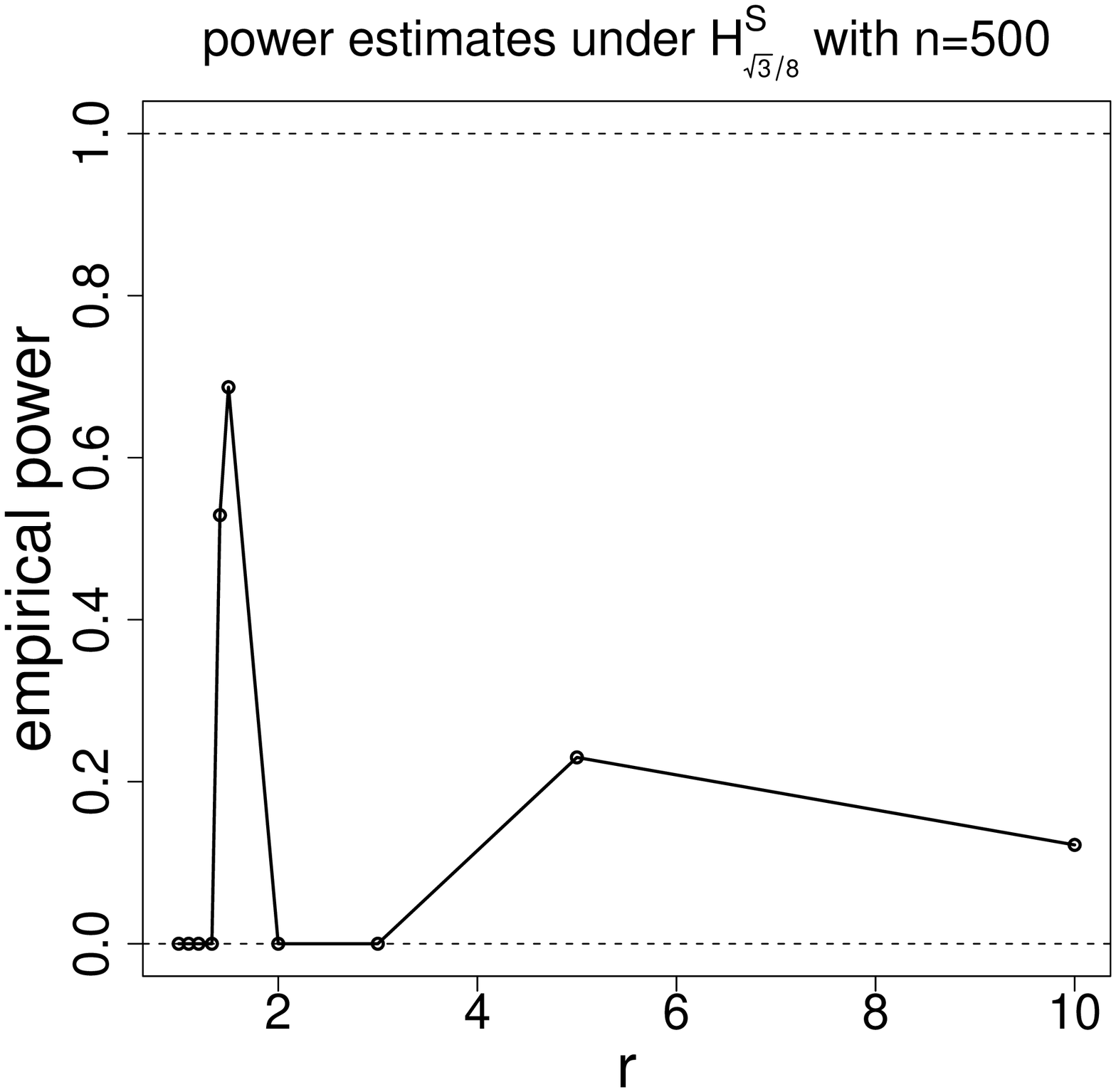}}}
\rotatebox{-90}{ \resizebox{1.7 in}{!}{ \includegraphics{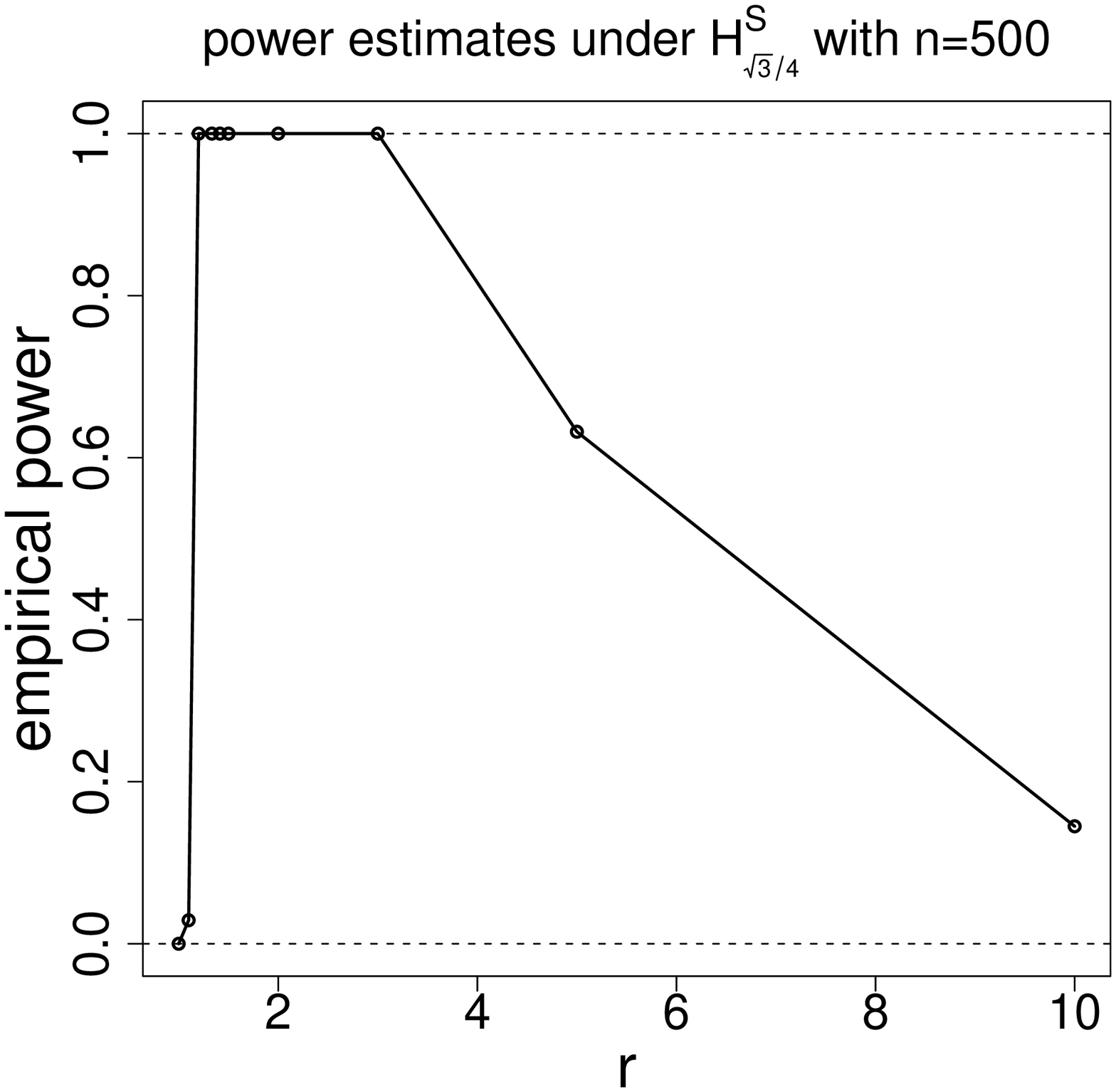}}}
\rotatebox{-90}{ \resizebox{1.7 in}{!}{ \includegraphics{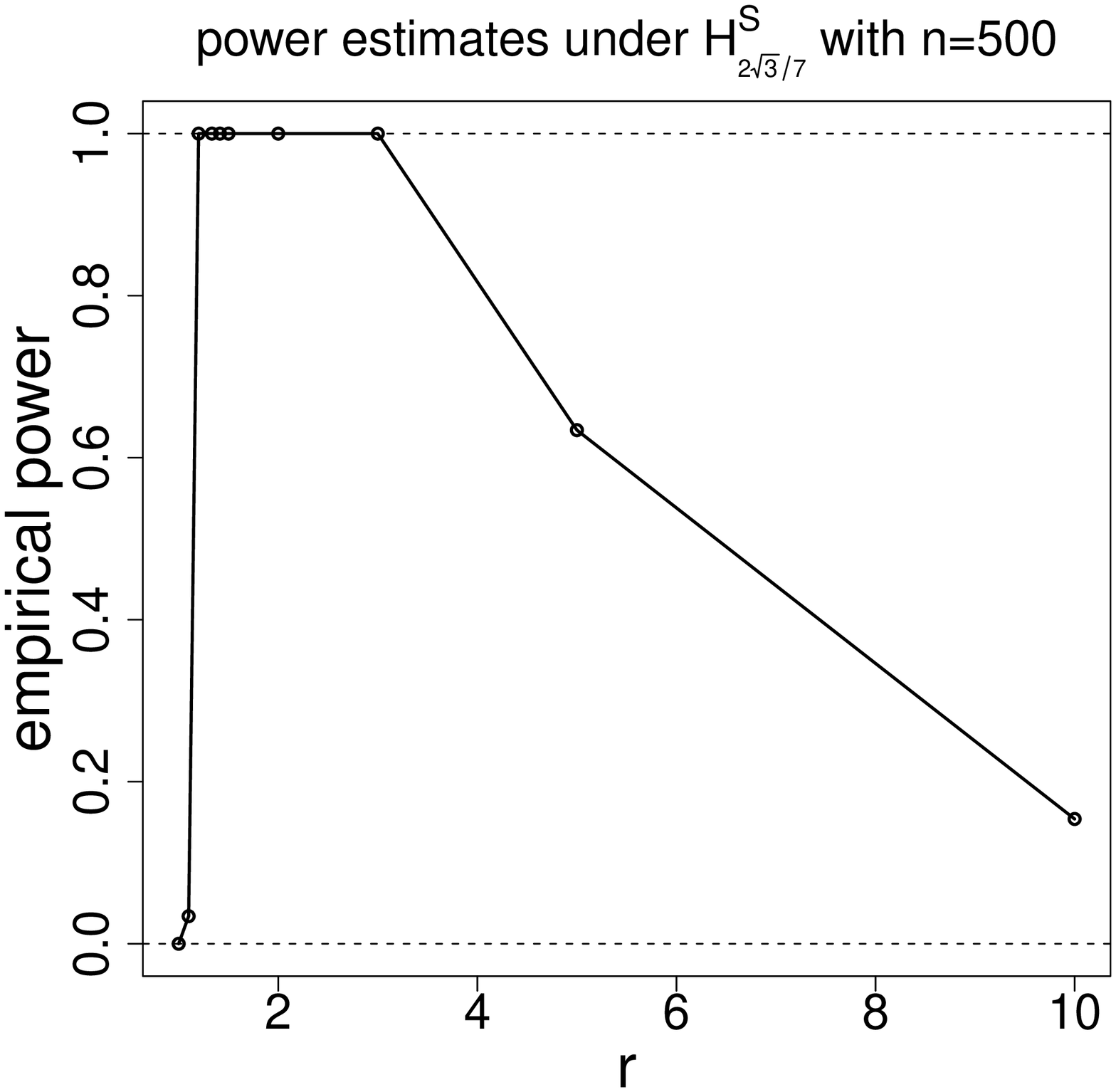}}}
\rotatebox{-90}{ \resizebox{1.7 in}{!}{ \includegraphics{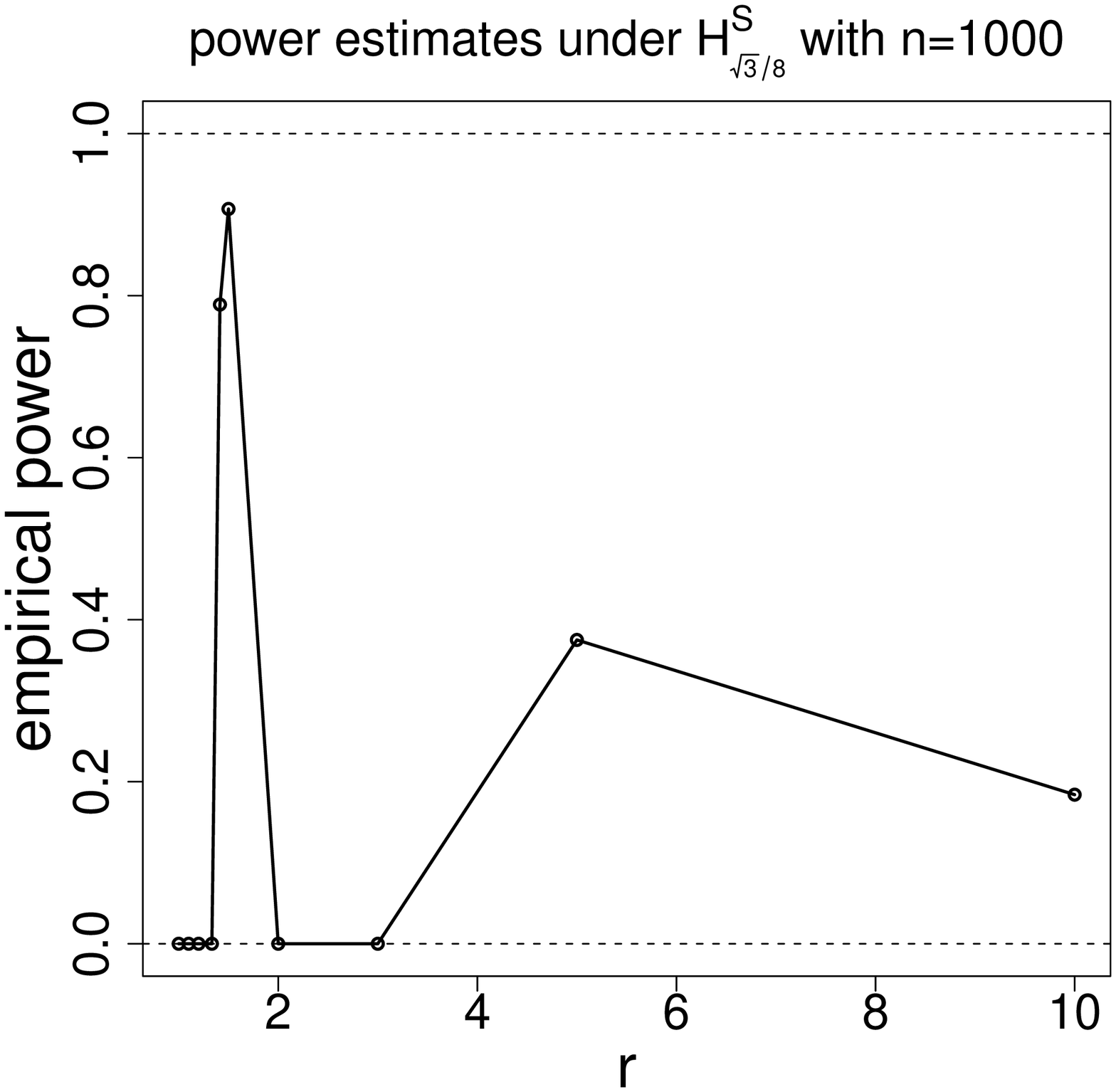}}}
\rotatebox{-90}{ \resizebox{1.7 in}{!}{ \includegraphics{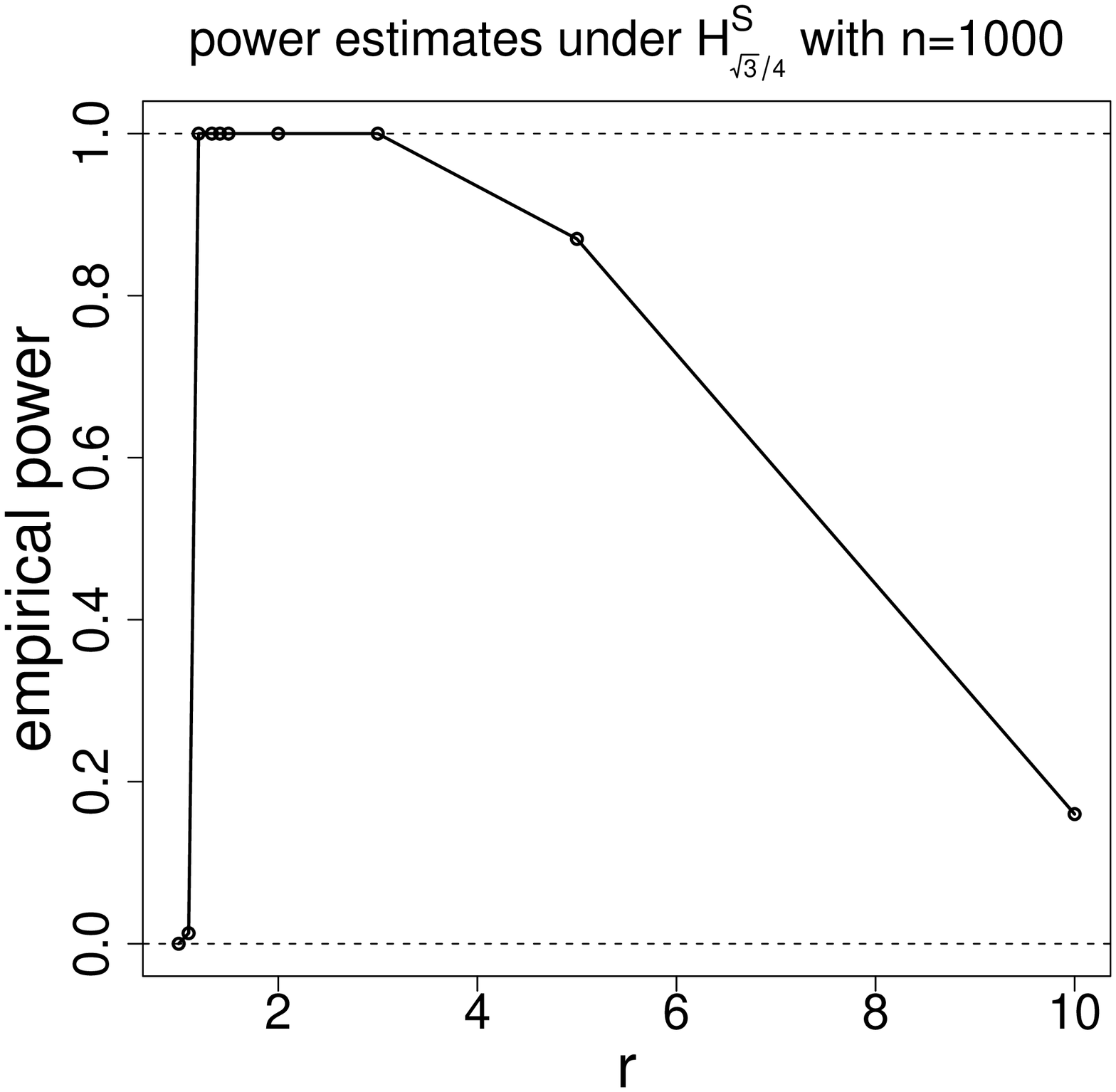}}}
\rotatebox{-90}{ \resizebox{1.7 in}{!}{ \includegraphics{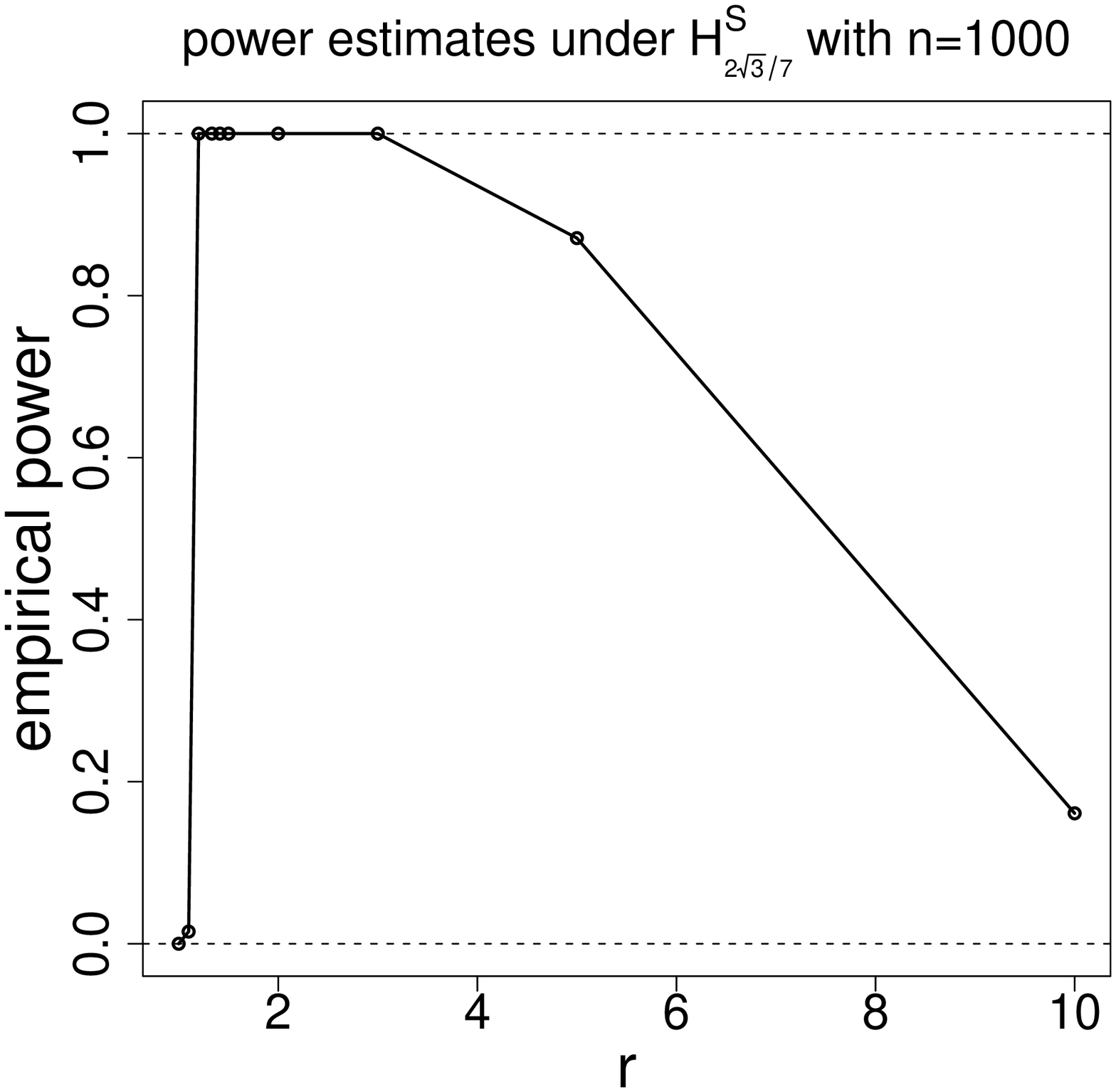}}}
\caption{\label{fig:MT-PE-emp-power-seg}
\textbf{Empirical power for $R_{PE}(r)$ in the multiple triangle case:}
Monte Carlo power estimates for proportional-edge PCDs in the multiple triangle case
using the asymptotic critical value against segregation alternatives
$H^S_{\sqrt{3}/8}$ (left column),
$H^S_{\sqrt{3}/4}$ (middle column),
and
$H^S_{2\,\sqrt{3}/7}$ (right column)
as a function of $r$, for $n=500$ (top) and $n=1000$ (bottom).
}
\end{figure}

In the multiple triangle case,
we generate the $\X$ points uniformly in the support for the segregation alternatives
in the triangles based on the 10 class $\Y$
points given in Figure \ref{fig:deldata}.
We use the parameters $\ve \in \{\sqrt{3}/8,\sqrt{3}/4,2\,\sqrt{3}/7\}$.
We compute the relative density based on the formula given in Corollary \ref{cor:MT-asy-norm}.
The corresponding empirical power estimates as a function of $r$ (using the normal approximation)
are presented in Figure \ref{fig:MT-PE-emp-power-seg} for $n=500$ or 1000.
Observe that the Monte Carlo power estimate
increases as $r$ gets larger and then decreases,
as in the one triangle case.
The empirical power is maximized for $r \in (1.5,2)$ under mild segregation,
and for $r \in (1.5,3)$ under moderate to severe segregation.
Considering the empirical size estimates,
$r \approx 1.5$ is recommended under mild segregation,
while $r \in (2,3)$ seems to be more appropriate (hence recommended for more severe segregation)
since the corresponding test has the desired level with high power.

\subsection{Empirical Power Analysis for Central Similarity PCDs under the Segregation Alternative}
\label{sec:CS-emp-power-seg}
In the one triangle case,
data generation is as in Section \ref{sec:PE-emp-power-seg}.
At each Monte Carlo replicate we compute the relative density
of the central similarity PCDs.
We consider $\tau \in \{0.2,0.4,0.6,\ldots,3.0,3.5,4.0,\ldots,20.0\}$
for the central similarity PCD.
We repeat the simulation procedure $N_{mc}=10000$ times.
Under segregation alternatives with $\ve>0$,
the distribution of $\rho_{_{CS}}(n,\tau)$ is non-degenerate
for all $\tau \in (0,\infty)$ and $\ve \in (0,\sqrt{3}/3)$.

\begin{figure}[]
\centering
\psfrag{kernel density estimate}{ \Huge{\bfseries{kernel density estimate}}}
\psfrag{relative density}{ \Huge{\bfseries{relative density}}}
\rotatebox{-90}{ \resizebox{2. in}{!}{ \includegraphics{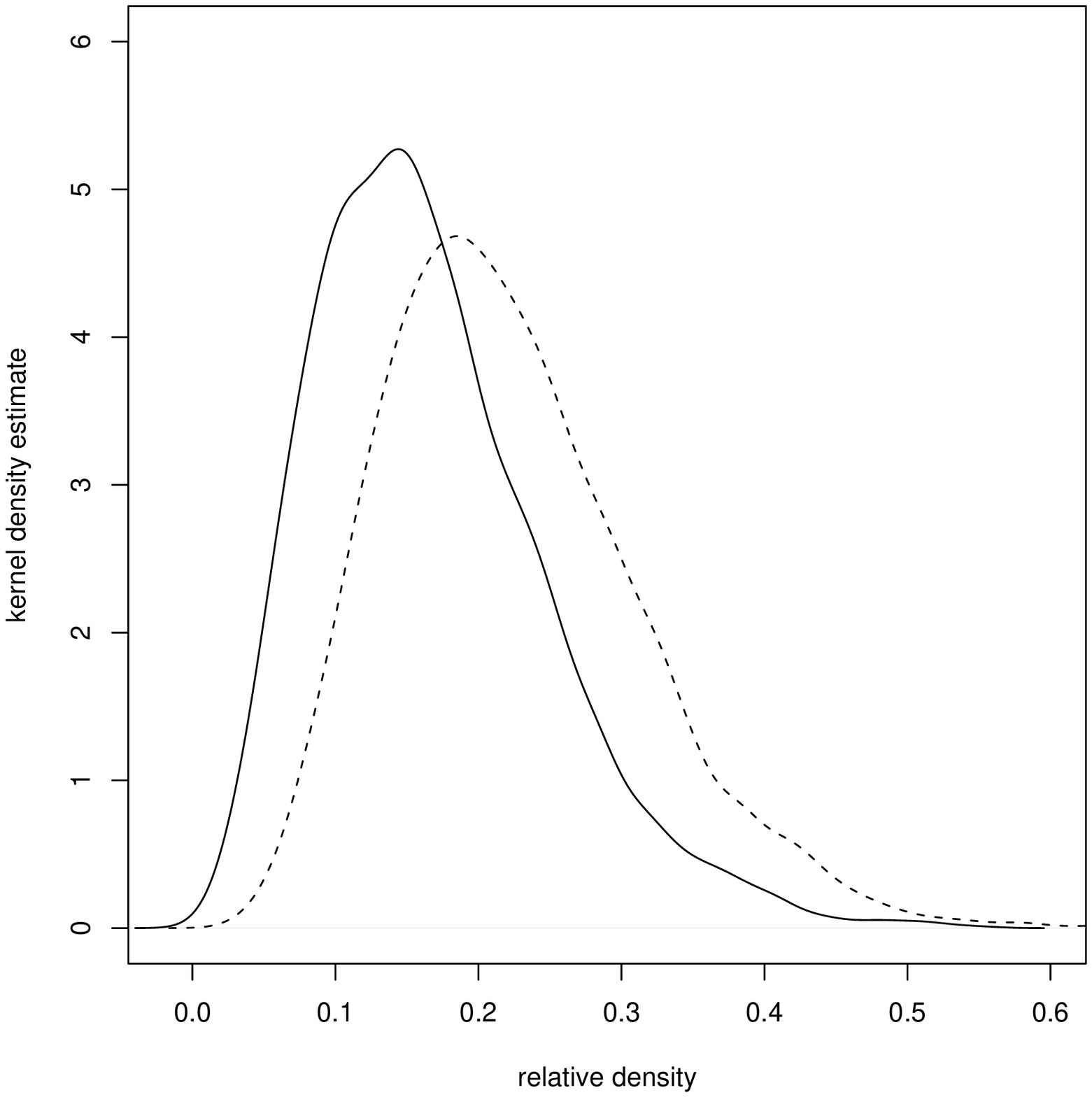}}}
\rotatebox{-90}{ \resizebox{2. in}{!}{ \includegraphics{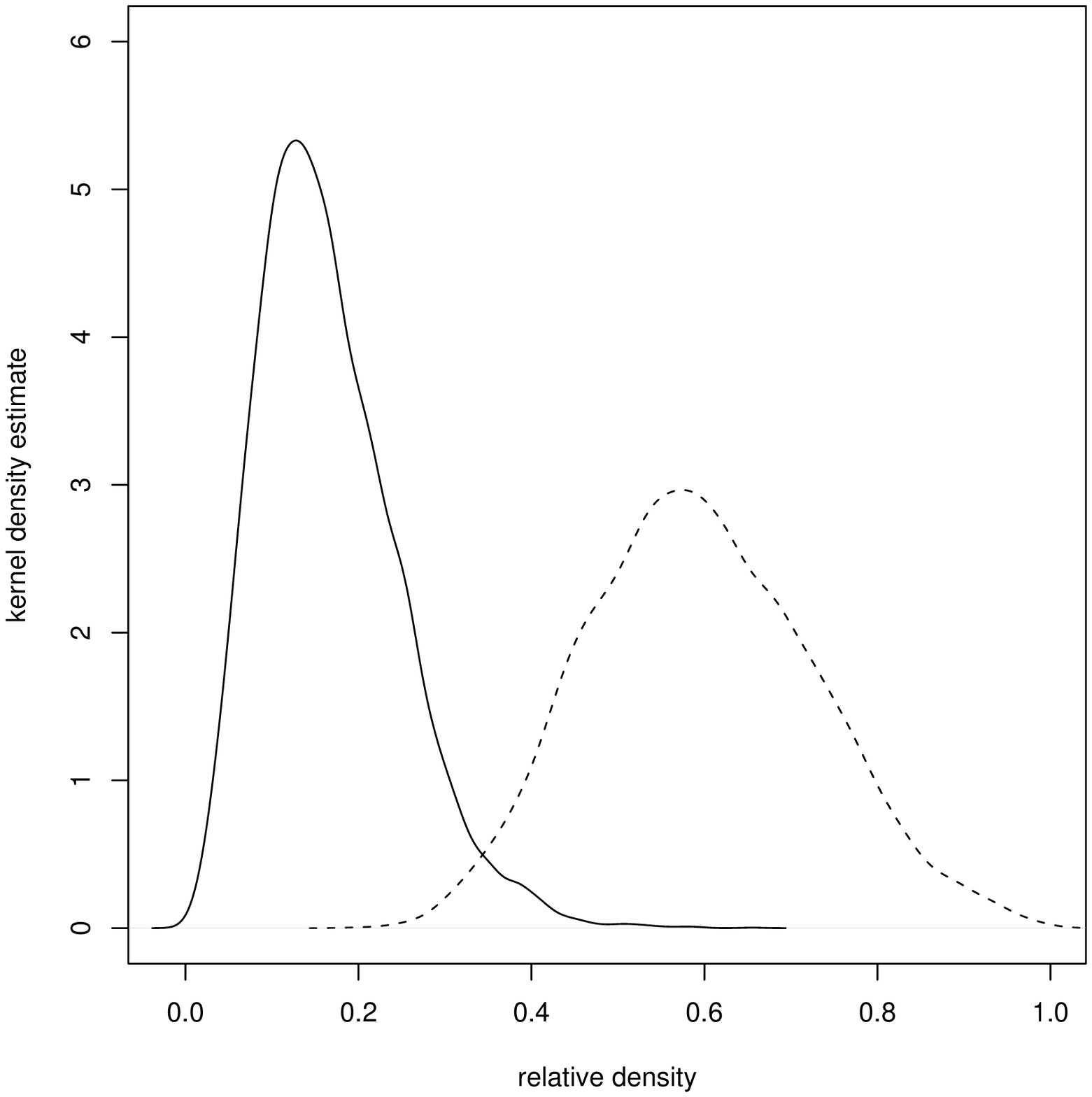}}}
\caption{\label{fig:CSsegsim}
Kernel density estimates of the relative density of central similarity PCD,
$\rho_{_{CS}}(n,\tau)$,
under the null (solid line) and the segregation alternatives (dashed line) with
$H^S_{\sqrt{3}/8}$ (left) and $H^S_{\sqrt{3}/4}$ (right) for $\tau=1$
with $n=10$ based on $N_{mc}=10000$ replicates.
}
\end{figure}

\begin{figure}[ht]
\centering
\psfrag{kernel density estimate}{ \Huge{\bfseries{kernel density estimate}}}
\psfrag{relative density}{ \Huge{\bfseries{relative density}}}
\rotatebox{-90}{ \resizebox{2. in}{!}{ \includegraphics{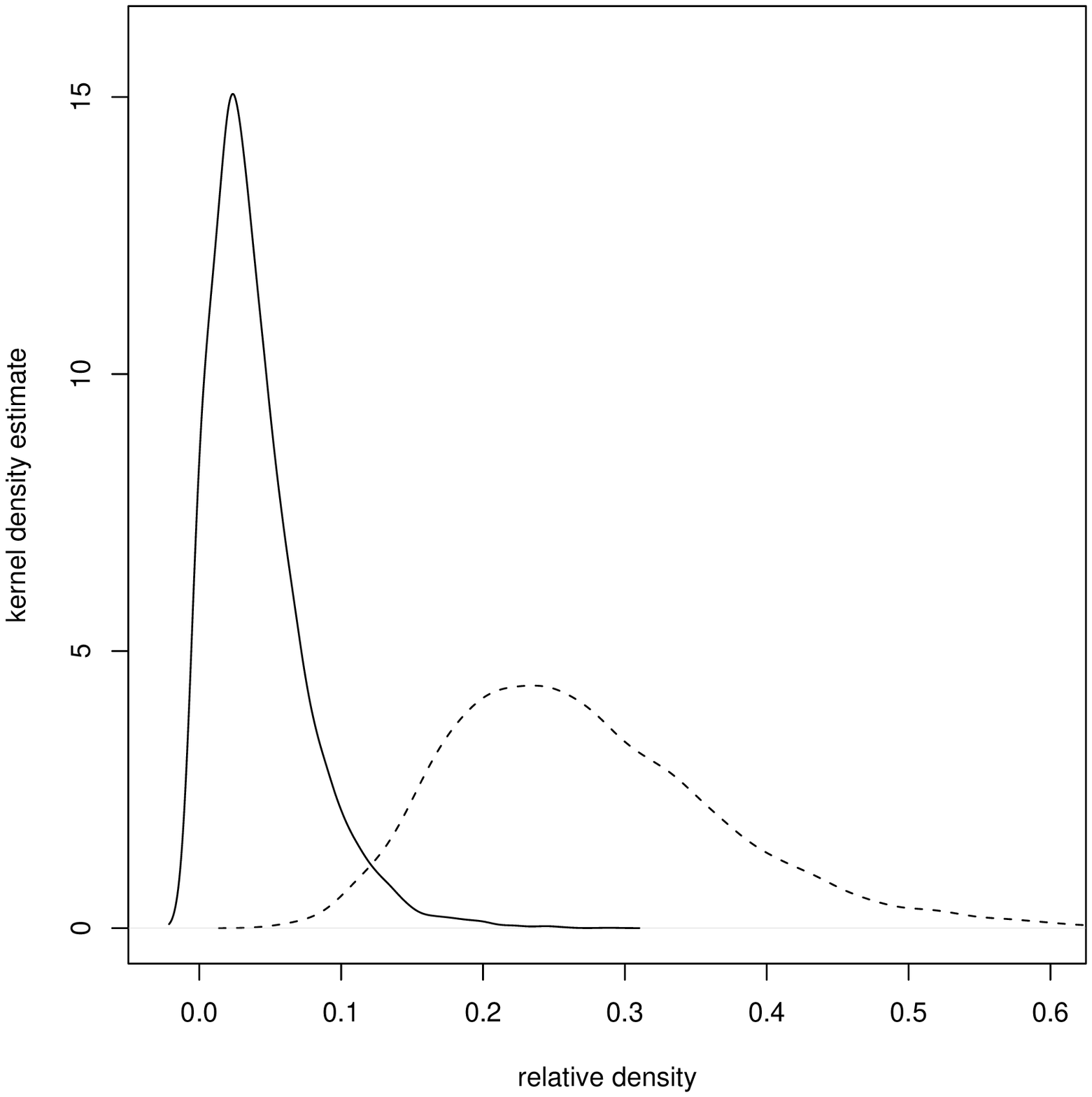}}}
\rotatebox{-90}{ \resizebox{2. in}{!}{ \includegraphics{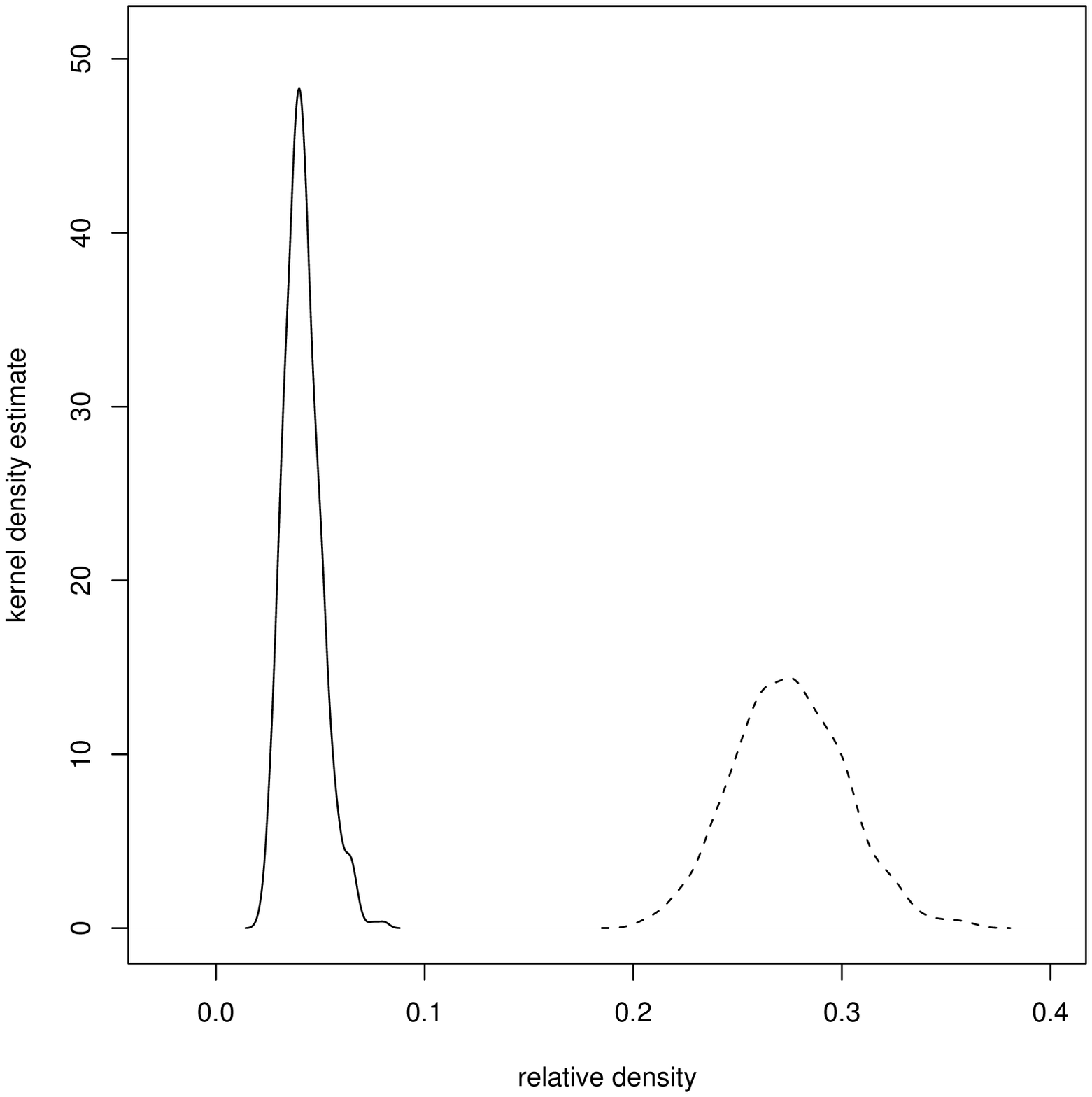}}}
\caption{ \label{fig:CSSegSimPowerPlots}
Depicted are kernel density estimates for $\rho_{_{CS}}(n,0.5)$ for
$n=10$ (left) and $n=100$ (right) under the null (solid line) and segregation alternative $H^S_{\sqrt{3}/4}$ (dashed line).
}
\end{figure}

In the one triangle case,
we plot
the kernel density estimates for the null
case and the segregation alternative with $\ve=\sqrt{3}/8$ and $\ve=\sqrt{3}/4$
with $n=10$ and $N_{mc}=10000$
in Figure \ref{fig:CSsegsim}.
Observe that under
both $H_o$ and alternatives, kernel density estimates
are almost symmetric for $\tau=1$.
Moreover, there is much more separation between
the kernel density estimates of the null and alternatives
for $\ve=\sqrt{3}/4$ compared to $\ve=\sqrt{3}/8$,
implying more power for larger $\ve$ values.
In Figure \ref{fig:CSSegSimPowerPlots},
we present kernel density estimates
for the null case and the segregation alternative
$H^S_{\sqrt{3}/4}$ for $\tau=0.5$, and $n=10$,
$N_{mc}=10000$ (left), $n=100$, $N_{mc}=1000$ (right).
With $n=10$, the null and alternative kernel density functions
for $\rho_{_{CS}}(10,0.5)$ are very similar, implying small power.
With $n=100$,
there is more separation
between null and alternative kernel density functions,
implying higher power.
Notice also that the probability density functions are more skewed for $n=10$,
while approximate normality holds for $n=100$.

\begin{figure}[ht]
\centering
\rotatebox{-90}{ \resizebox{1.7 in}{!}{ \includegraphics{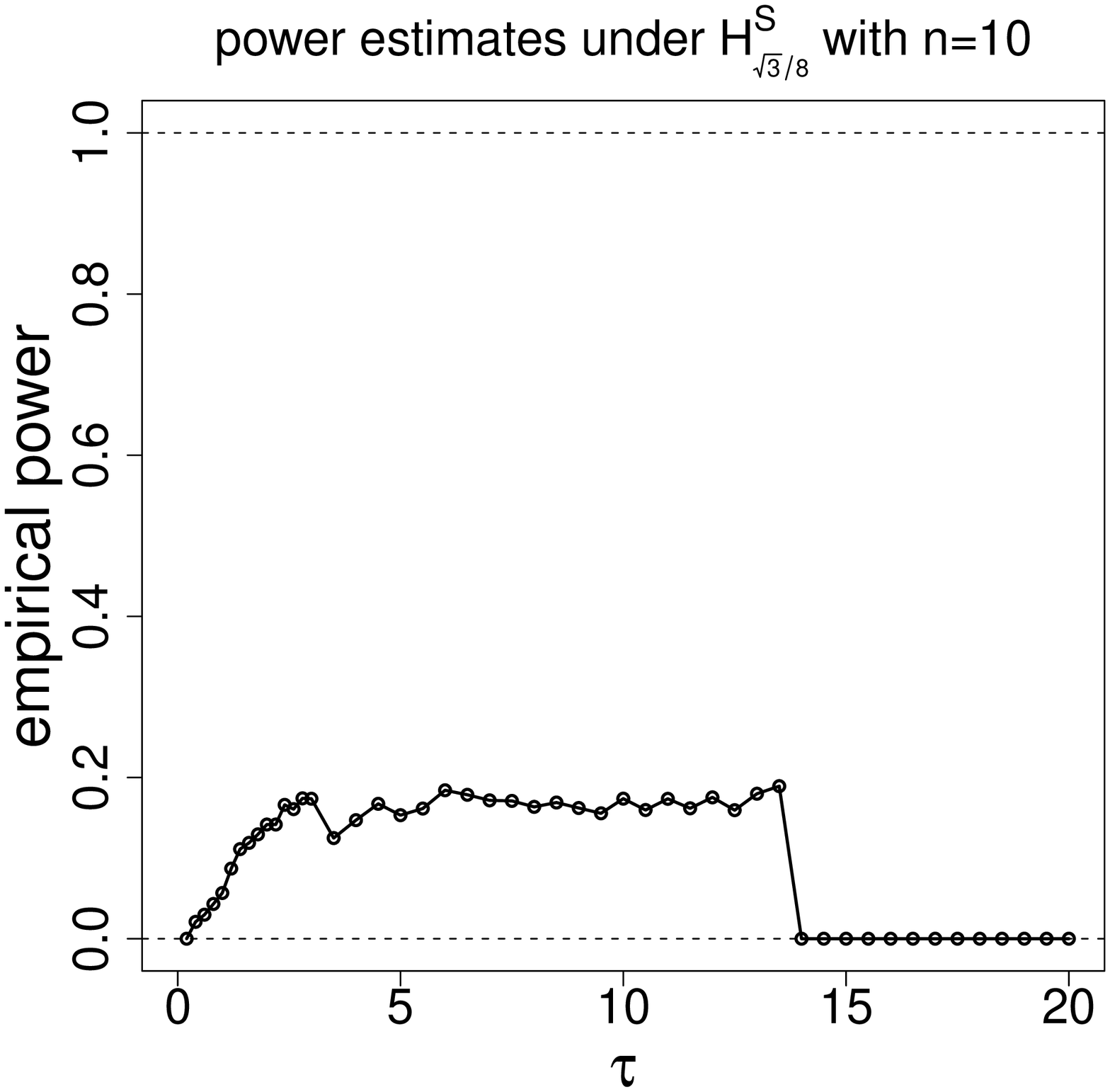}}}
\rotatebox{-90}{ \resizebox{1.7 in}{!}{ \includegraphics{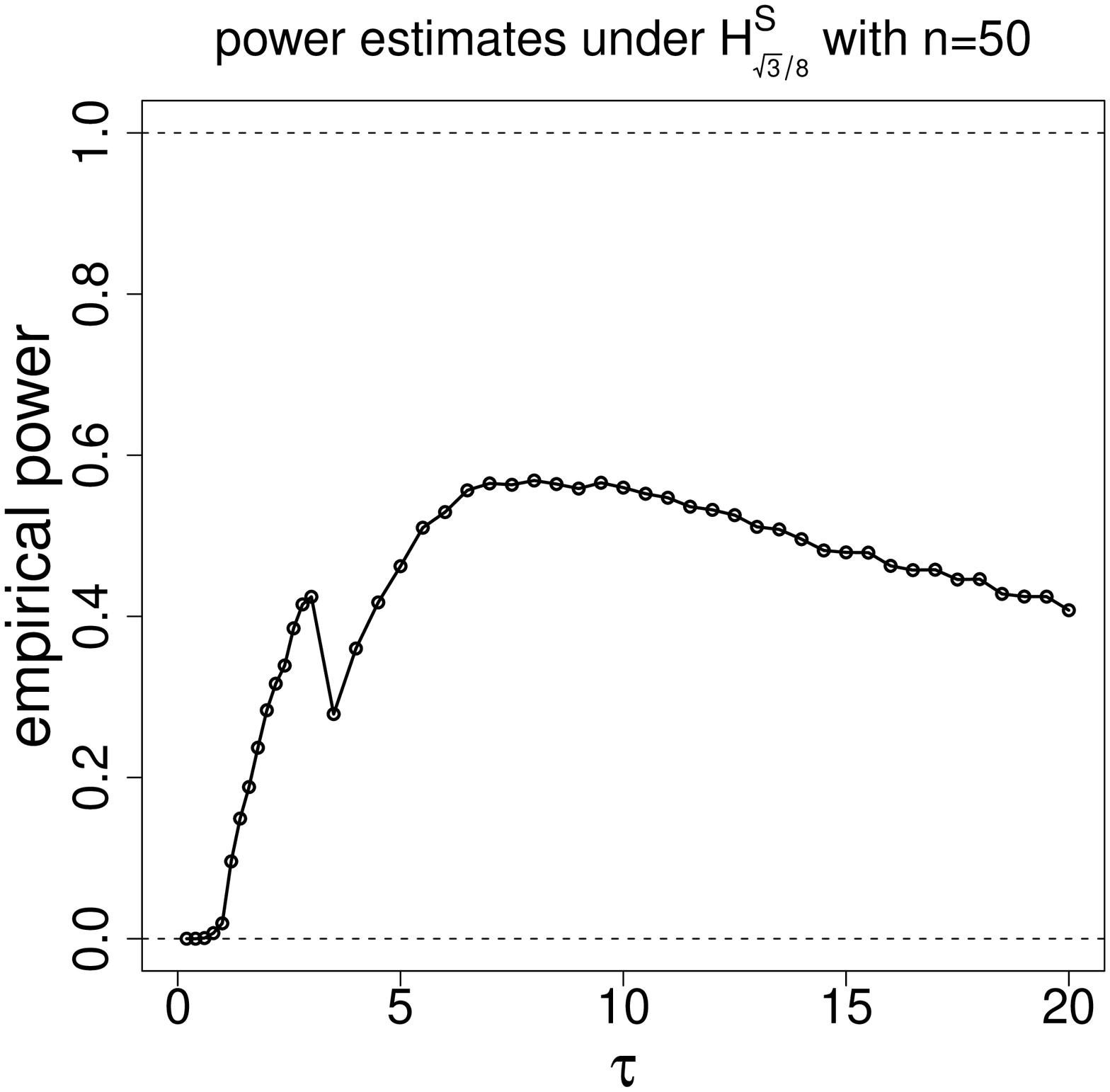}}}
\rotatebox{-90}{ \resizebox{1.7 in}{!}{ \includegraphics{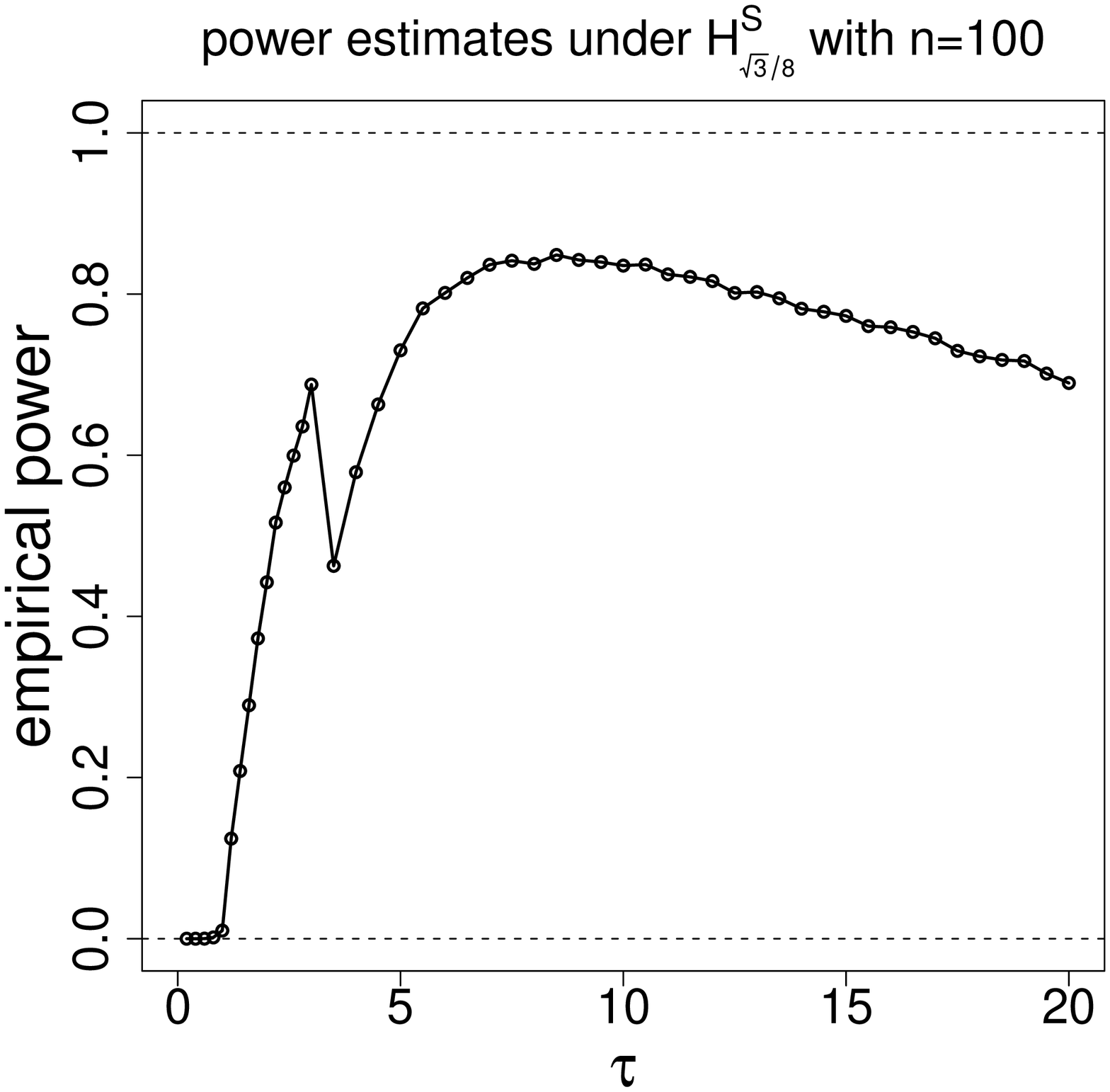}}}
\rotatebox{-90}{ \resizebox{1.7 in}{!}{ \includegraphics{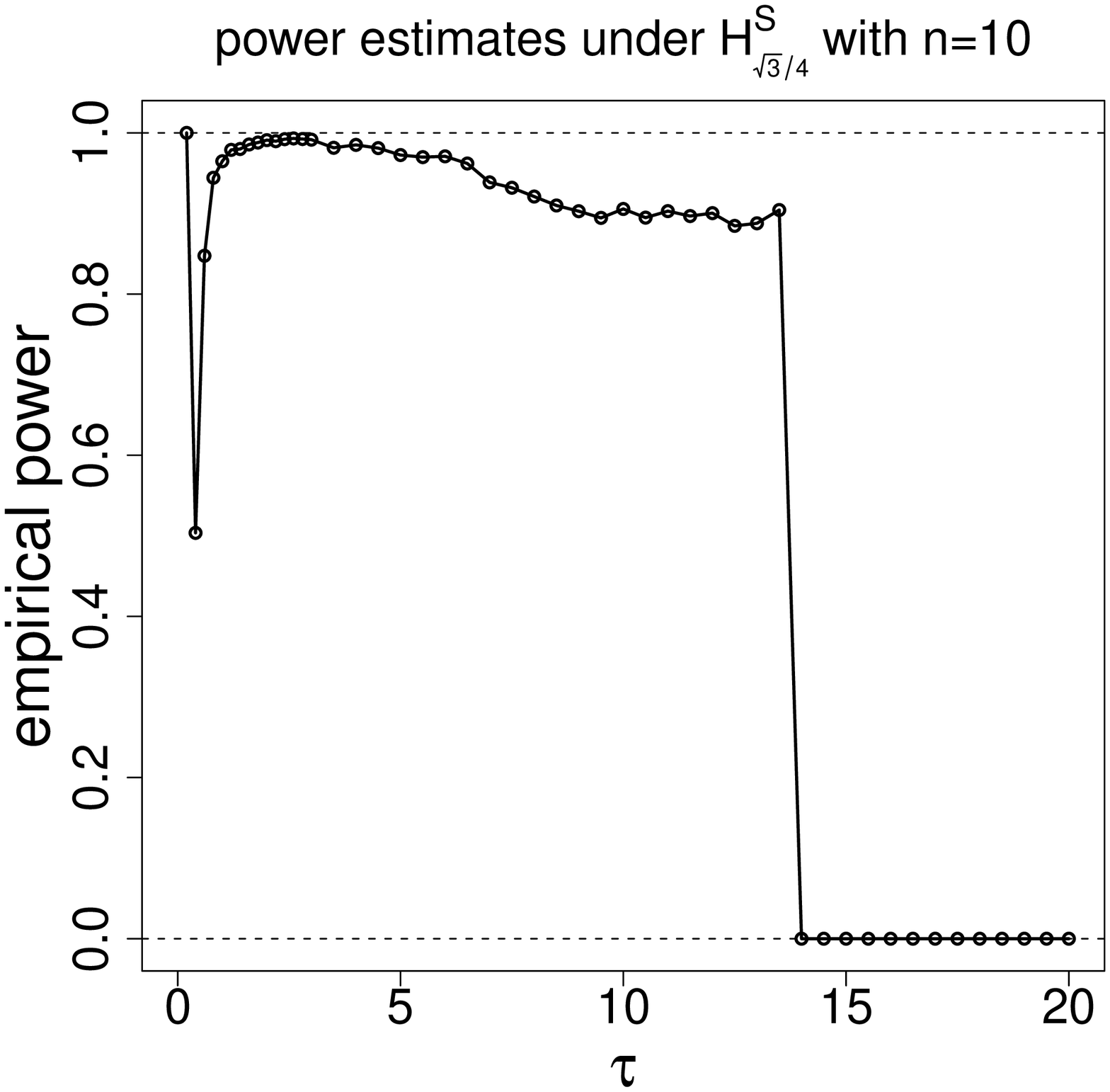}}}
\rotatebox{-90}{ \resizebox{1.7 in}{!}{ \includegraphics{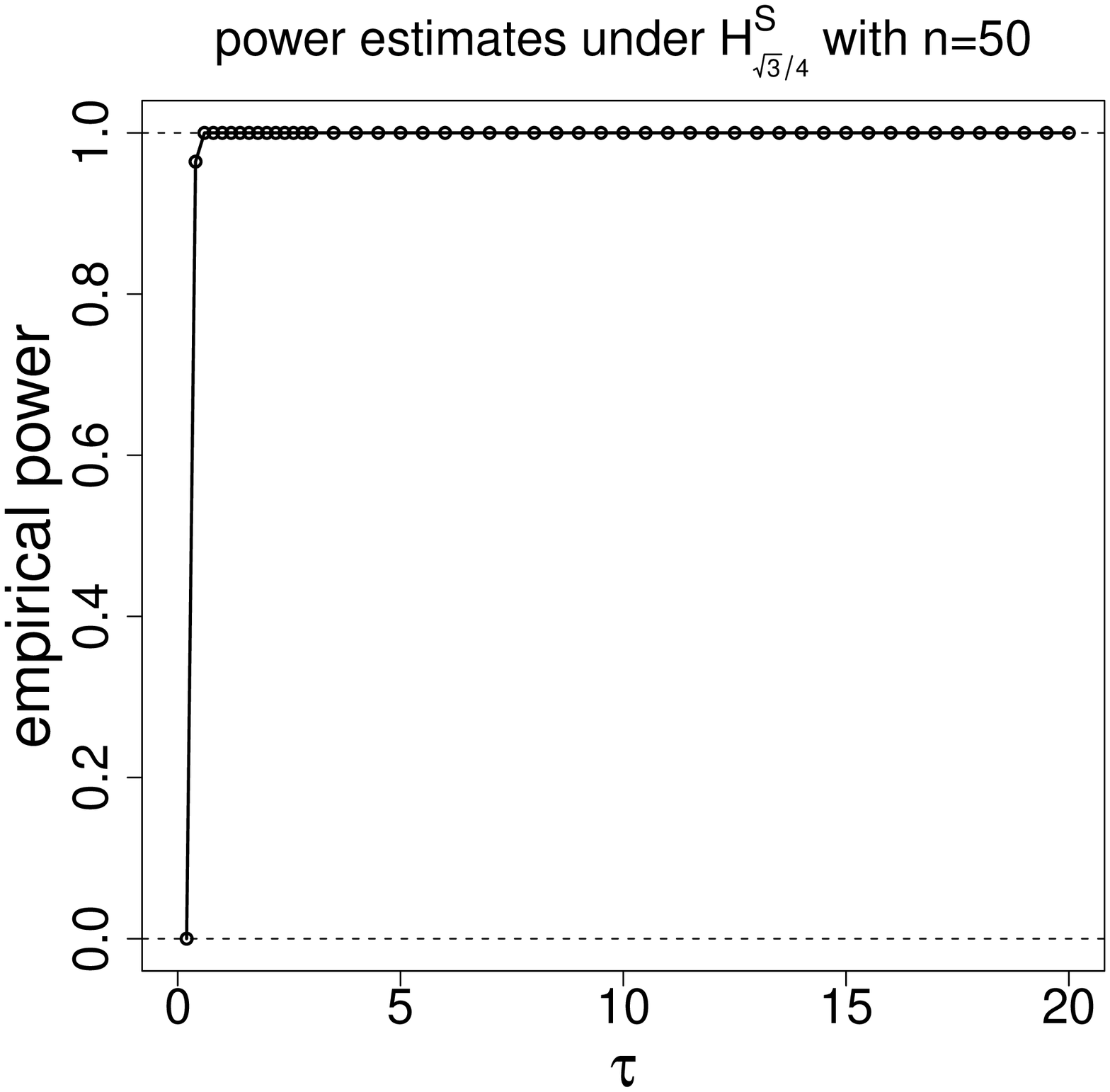}}}
\rotatebox{-90}{ \resizebox{1.7 in}{!}{ \includegraphics{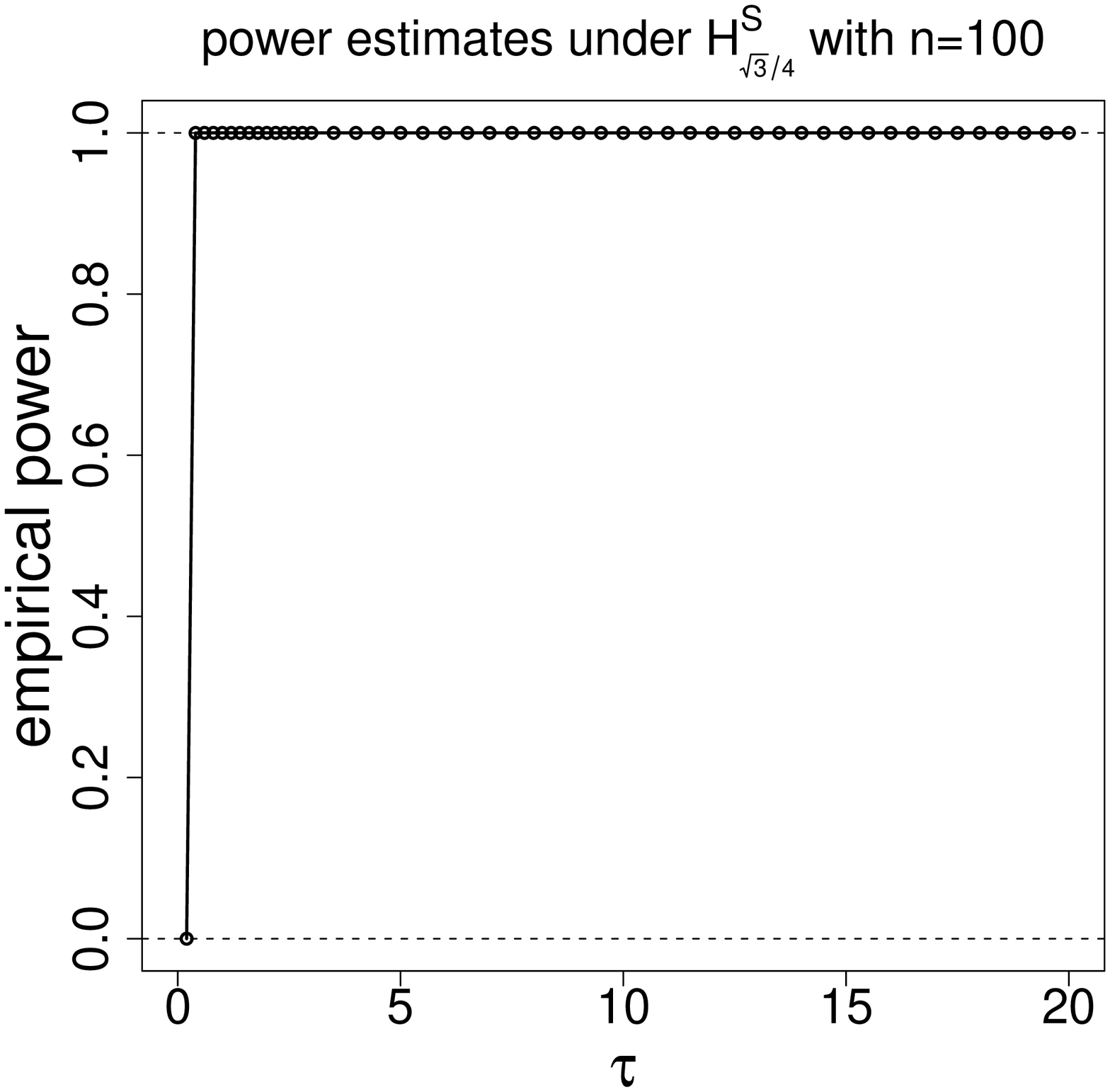}}}
\rotatebox{-90}{ \resizebox{1.7 in}{!}{ \includegraphics{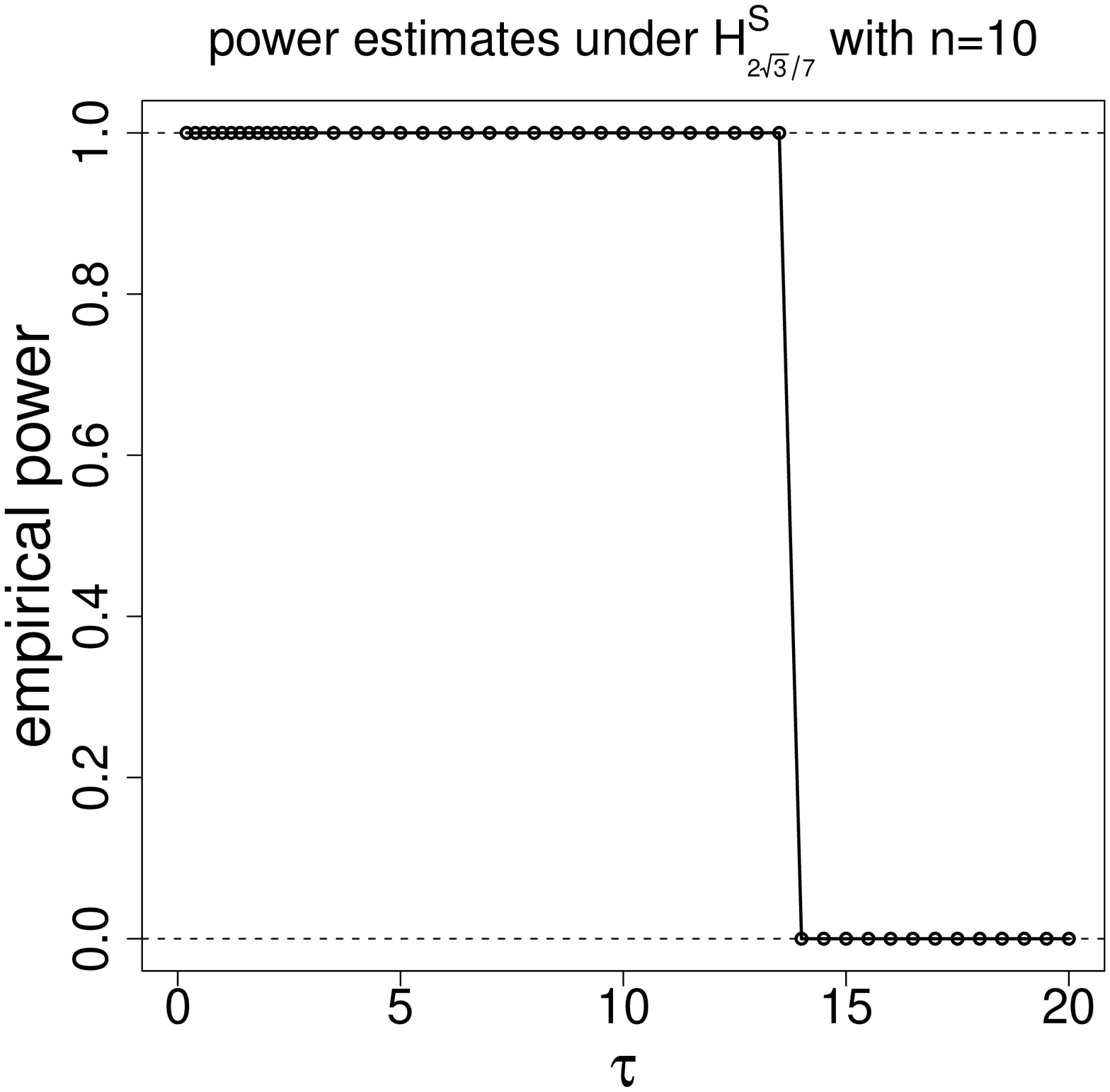}}}
\rotatebox{-90}{ \resizebox{1.7 in}{!}{ \includegraphics{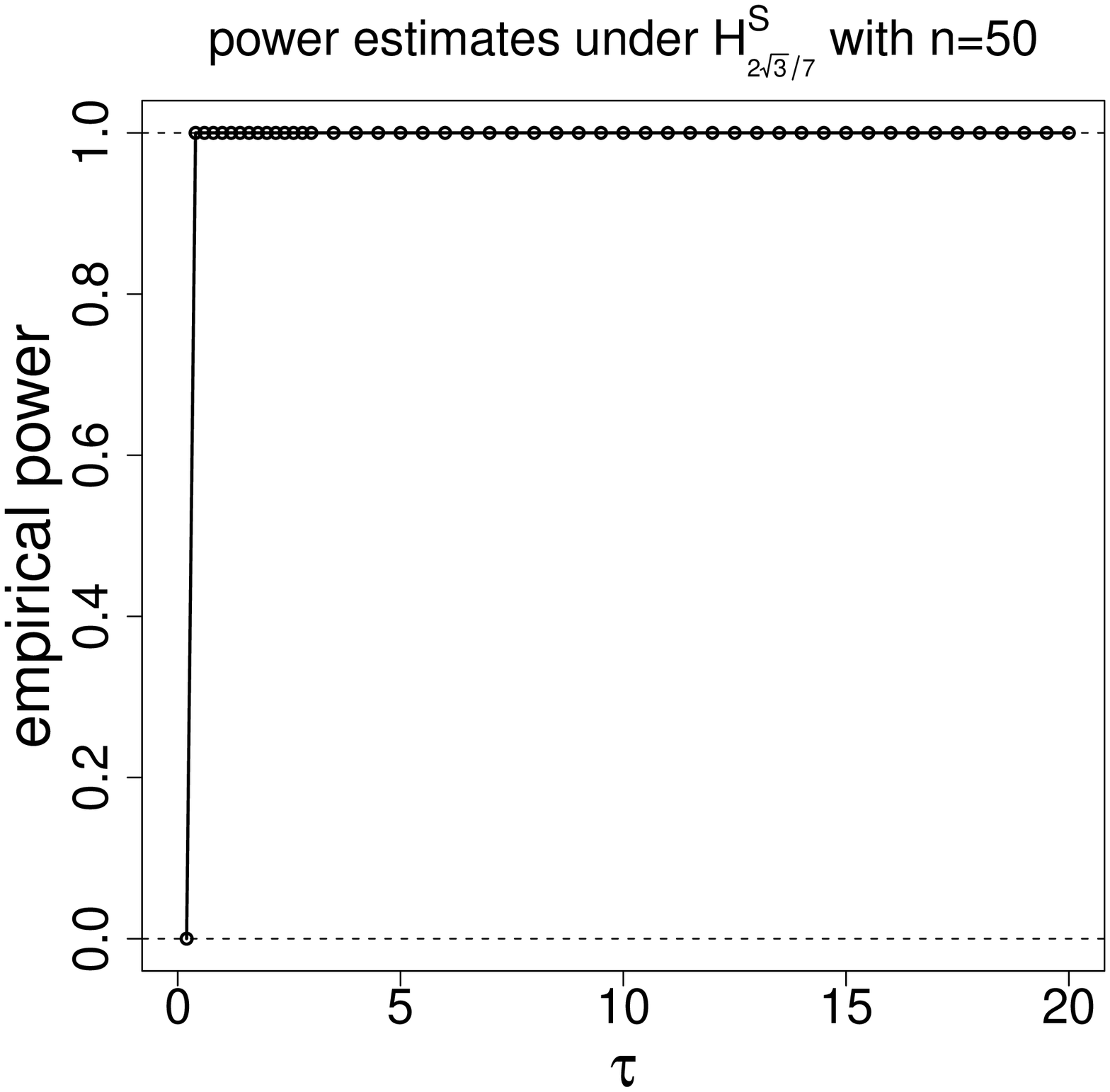}}}
\rotatebox{-90}{ \resizebox{1.7 in}{!}{ \includegraphics{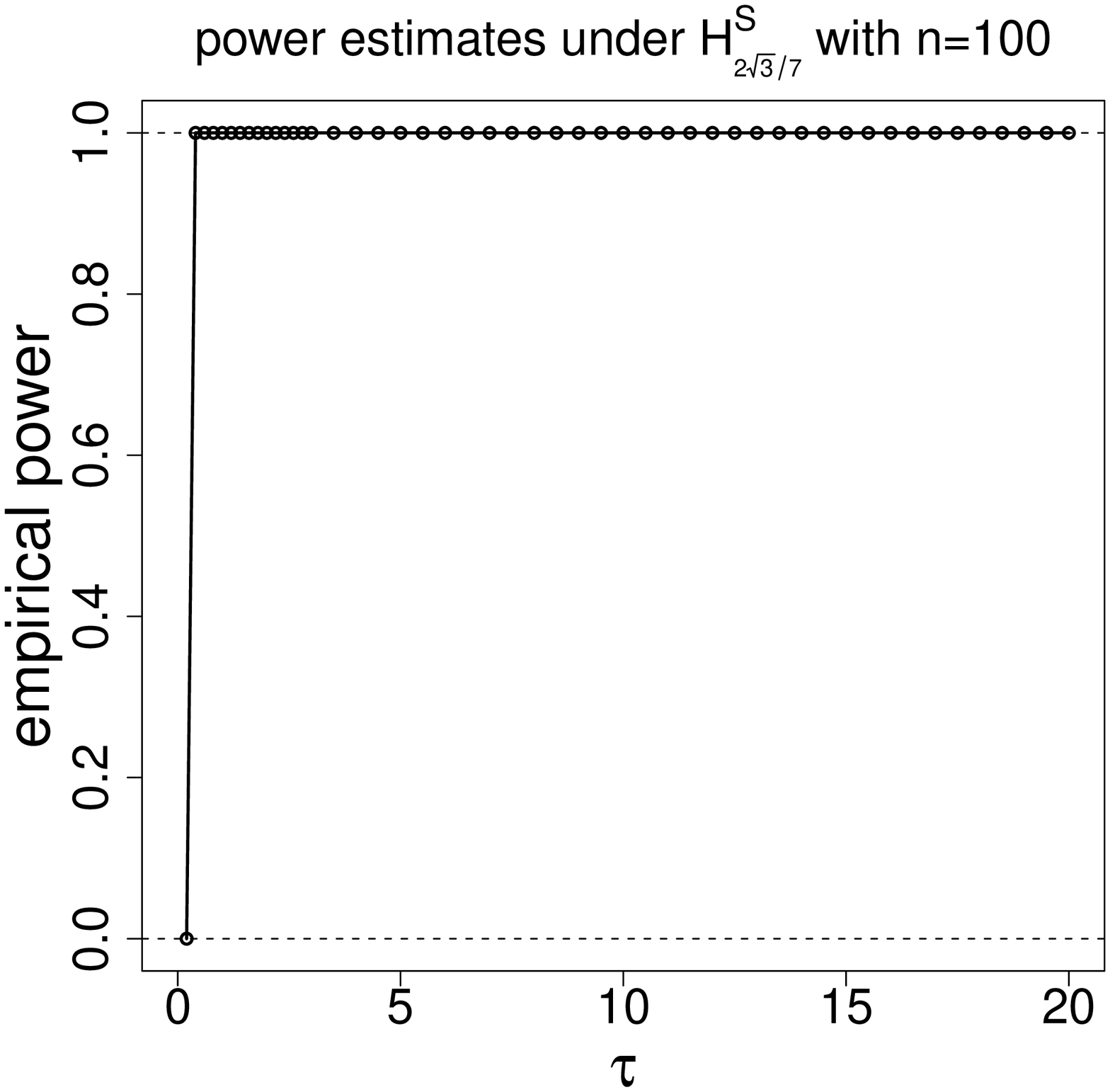}}}
\caption{\label{fig:CS-emp-power-seg}
\textbf{Empirical power for $R_{CS}(\tau)$ in the one triangle case:}
Monte Carlo power estimates for relative density of central similarity PCDs
in the one triangle case
using the asymptotic critical value against segregation alternatives
$H^S_{\sqrt{3}/8}$ (top row),
$H^S_{\sqrt{3}/4}$ (middle row),
and
$H^S_{2\,\sqrt{3}/7}$ (bottom row)
as a function of $\tau$, for $n=10$ (left column), $n=50$ (middle column), and $n=100$ (right column).
}
\end{figure}

We estimate the empirical power as
$\frac{1}{N_{mc}}\sum_{j=1}^{N_{mc}}\I \left(R_{CS}(\tau)(\tau,j) > z_{1-\alpha} \right)$.
In Figure \ref{fig:CS-emp-power-seg},
we present Monte Carlo power estimates for relative density of central similarity PCDs
in the one triangle case
against $H^S_{\sqrt{3}/8}$,
$H^S_{\sqrt{3}/4}$, and $H^S_{2\,\sqrt{3}/7}$ as a function of $\tau$ for $n=10,50,100$.
Notice that Monte Carlo power estimate
increases as $\tau$ gets larger or $n$ gets larger.
Moreover, the more severe the segregation,
the higher the power estimate at each $\tau$.
With $n=10$,
the power estimates are high for $\tau \in (5,14)$
and virtually 0 for $\tau \ge 14$.
With $n=50$ or $100$,
the power values are high for $\tau \ge 1$,
with highest power being attained around $\tau \approx 8$.
However, for $\tau \ge 6$,
the power values are virtually same.
Considering the empirical size estimates,
we recommend $\tau \approx 8$ for mild segregation,
and $\tau \approx 5$ for more severe segregation alternatives.

\begin{figure}[]
\centering
\rotatebox{-90}{ \resizebox{1.7 in}{!}{ \includegraphics{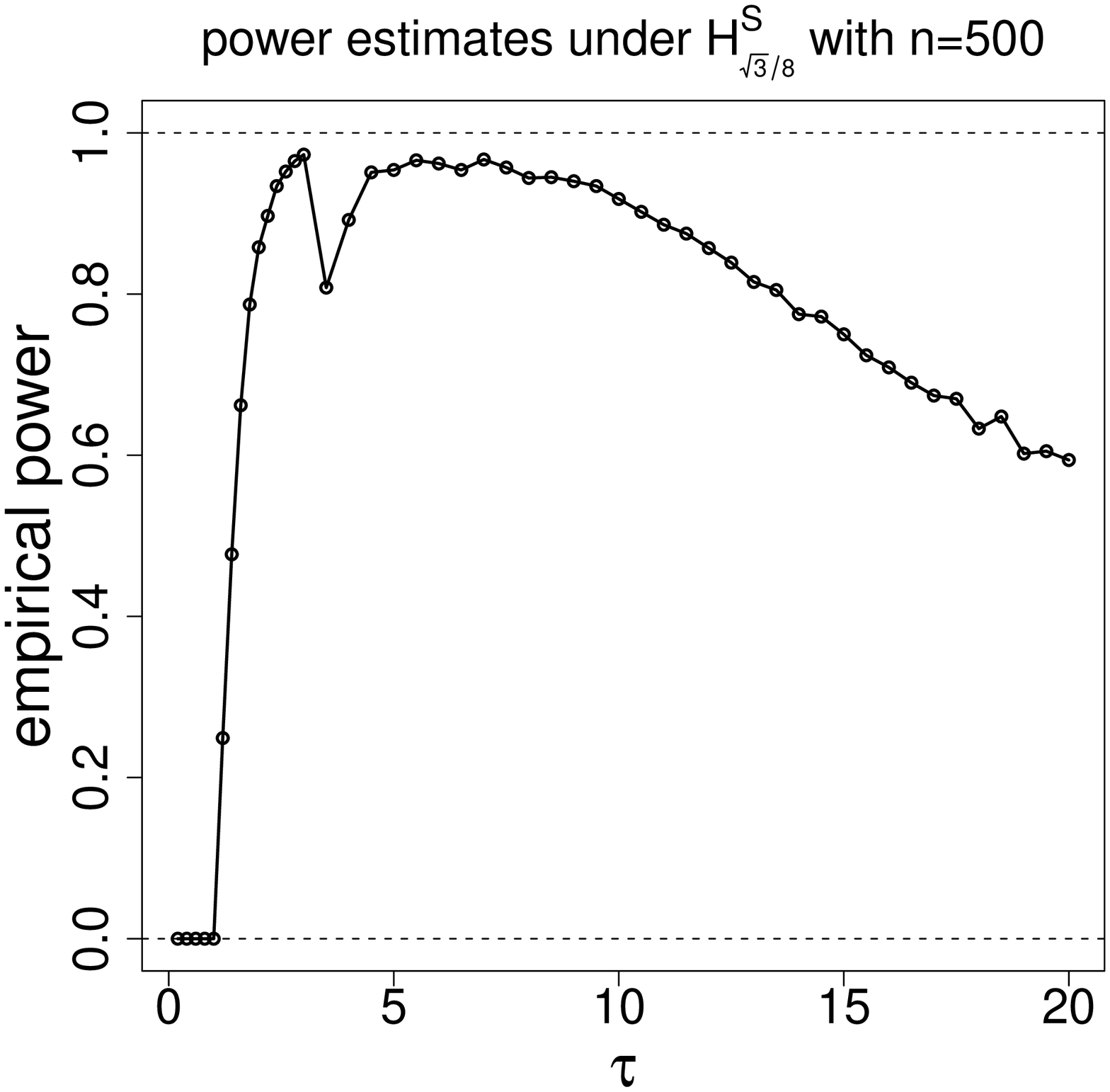}}}
\rotatebox{-90}{ \resizebox{1.7 in}{!}{ \includegraphics{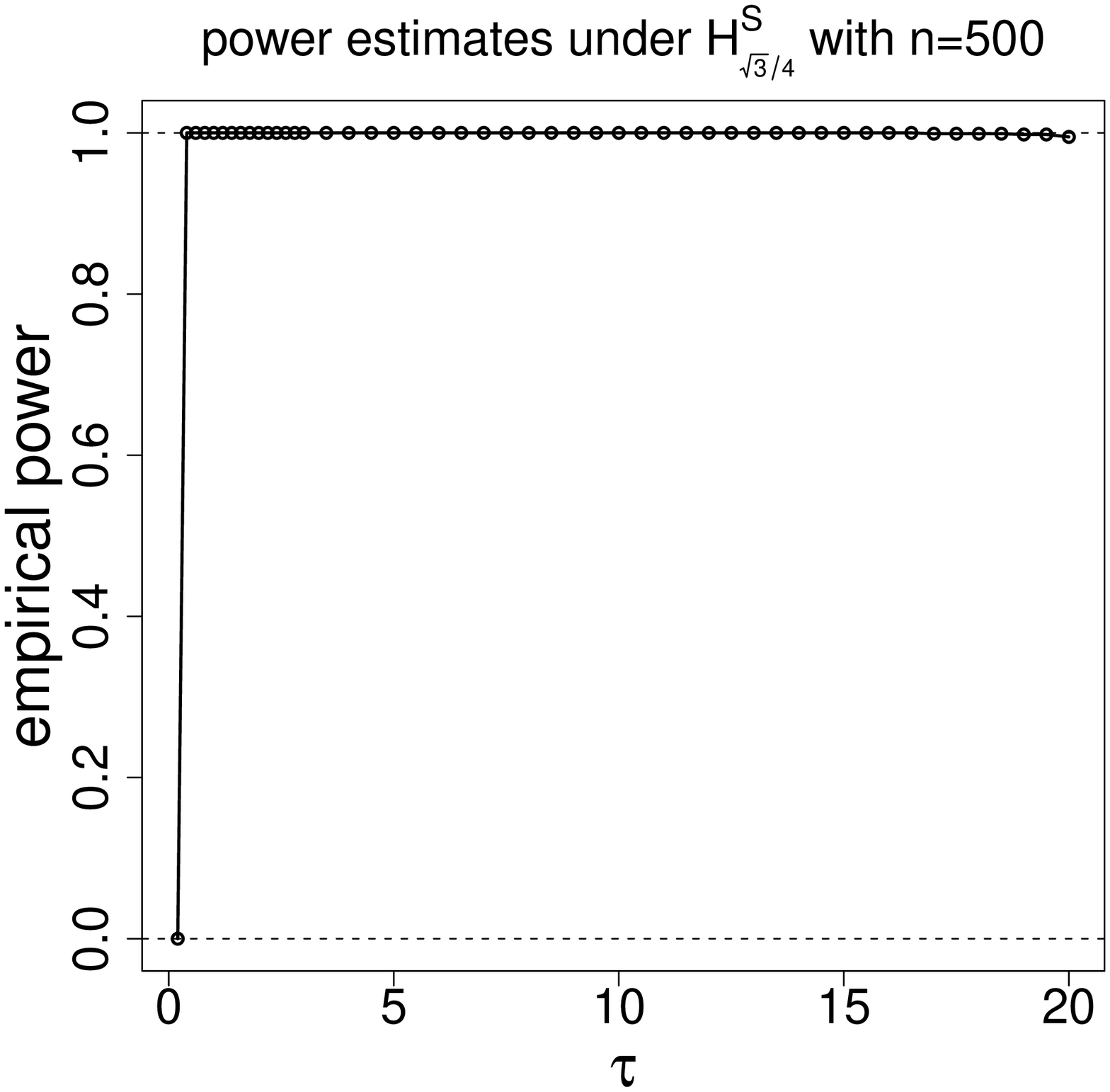}}}
\rotatebox{-90}{ \resizebox{1.7 in}{!}{ \includegraphics{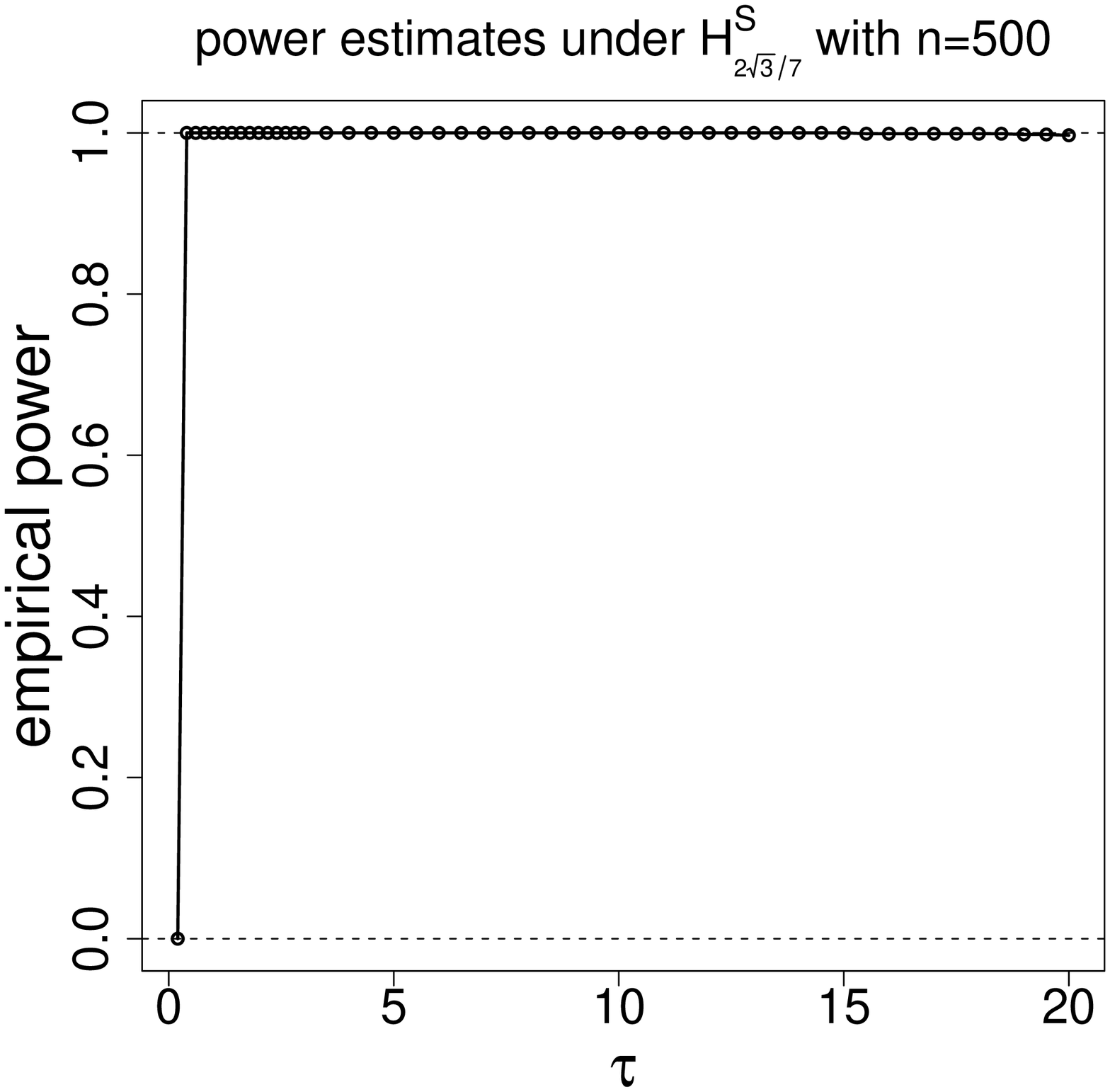}}}
\rotatebox{-90}{ \resizebox{1.7 in}{!}{ \includegraphics{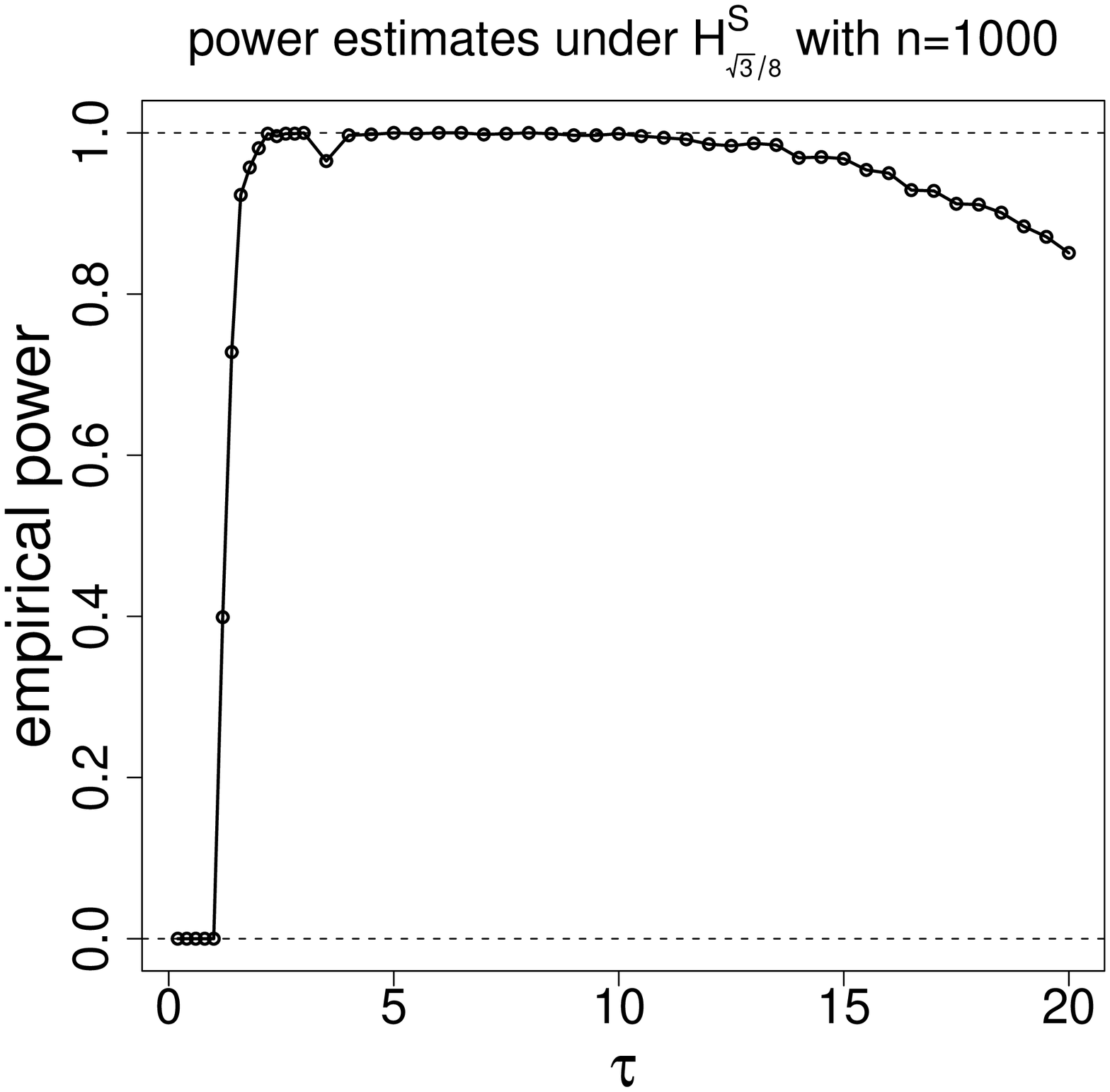}}}
\rotatebox{-90}{ \resizebox{1.7 in}{!}{ \includegraphics{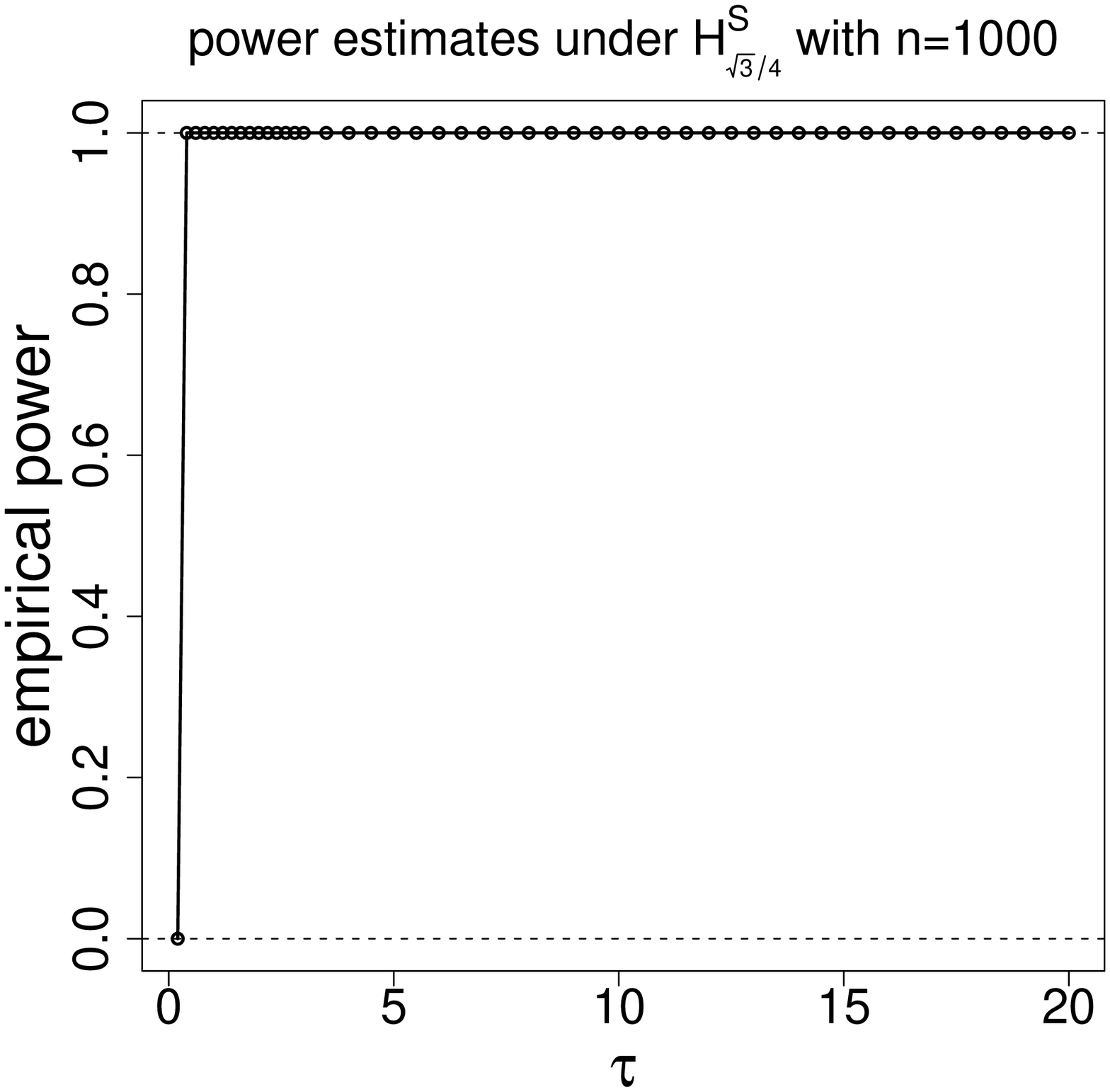}}}
\rotatebox{-90}{ \resizebox{1.7 in}{!}{ \includegraphics{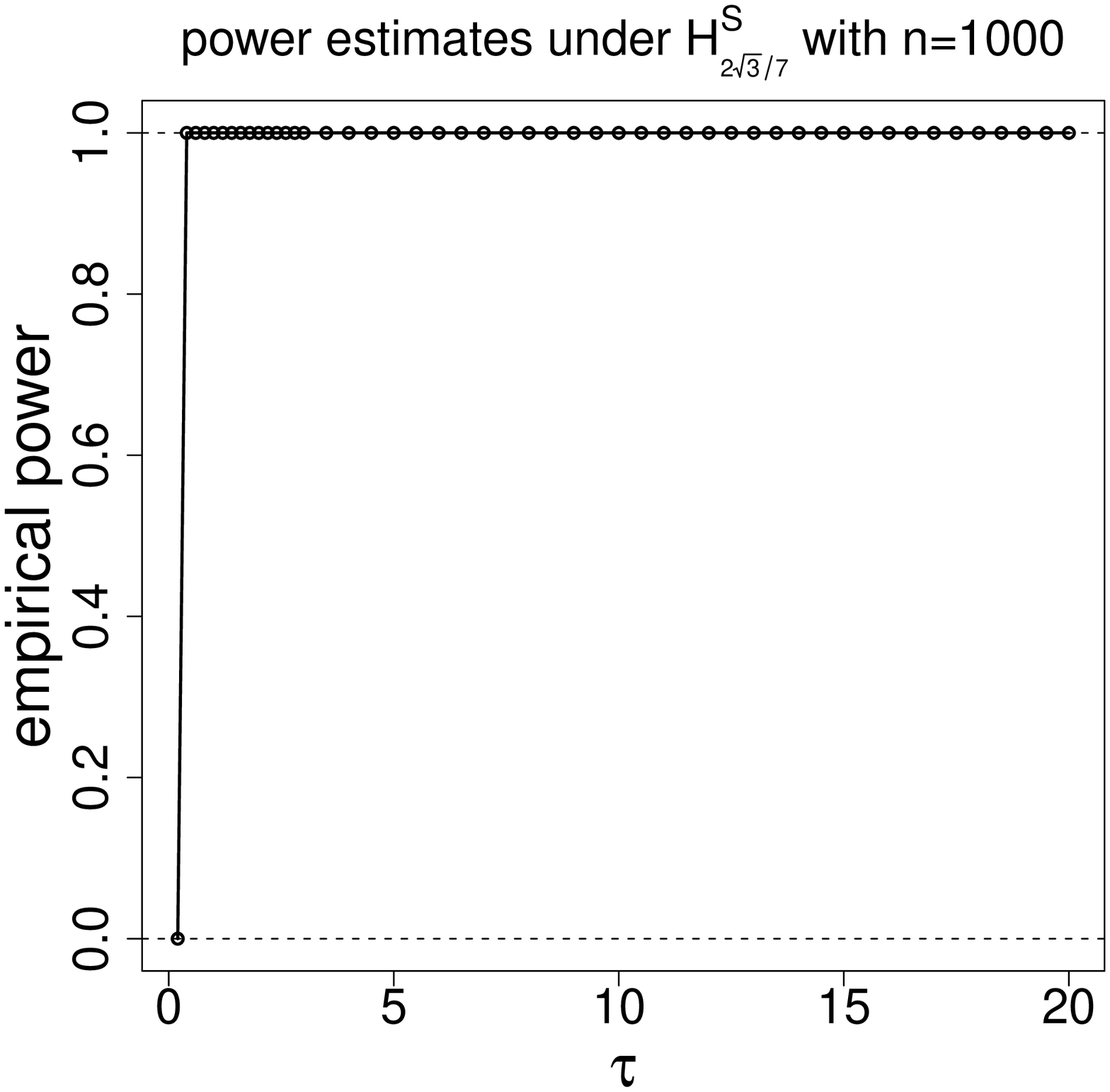}}}
\caption{\label{fig:MT-CS-emp-power-seg}
\textbf{Empirical power for $R_{CS}(\tau)$ in the multiple triangle case:}
Monte Carlo power estimates for central similarity PCDs in the multiple triangle case
using the asymptotic critical value against segregation alternatives
$H^S_{\sqrt{3}/8}$ (left column),
$H^S_{\sqrt{3}/4}$ (middle column),
and
$H^S_{2\,\sqrt{3}/7}$ (right column)
as a function of $\tau$, for $n=500$ (top) and $n=1000$ (bottom).
}
\end{figure}

In the multiple triangle case,
data generation is again as in Section \ref{sec:PE-emp-power-seg}.
We compute the relative density based on the formula given in Corollary \ref{cor:MT-asy-norm}.
The corresponding empirical power estimates as a function of $\tau$ (using the normal approximation)
are presented in
in Figure \ref{fig:MT-CS-emp-power-seg} for $n=500$ and $1000$.
Observe that the Monte Carlo power estimate tends to
increase as $\tau$ gets larger.
Under mild segregation with $\ve = \sqrt{3}/8$,
the empirical power is large for $\tau \ge 2$
with largest being around $\tau \in (4,8)$.
Under moderate to severe segregation,
the empirical power is virtually one for $\tau \ge 0.4$.
Considering the empirical size estimates,
$\tau \approx 7$ seems to be more appropriate (hence recommended for segregation)
since the corresponding test has the desired level with highest power.

\subsection{Empirical Power Analysis for Proportional-Edge PCDs under the Association Alternative}
\label{sec:PE-emp-power-assoc}
In the one triangle case,
at each of $N_{mc}=10000$ Monte Carlo replicates under association $H^A_{\ve}$,
we generate
$X_i \stackrel{iid}{\sim} \U\left(\mathcal T_{\sqrt{3}/3-\ve}\right)$,
for $i=1,2,\ldots,n$ for $n=10,50,100$.
The relative density is computed as in Section \ref{sec:PE-emp-power-seg}.
Unlike the segregation alternatives,
the distribution of $\rho_{_{PE}}(n,r)$ is non-degenerate
for all $\ve \in (0,\sqrt{3}/3)$ and $r \in [1,\infty)$.
We consider $\ve \in \{ 5\, \sqrt{3}/24, \sqrt{3}/12, \sqrt{3}/21 \}$
(which correspond to 18.75 \%, 75 \%, and $4500/49\approx 91.84$ \%
of the triangle being occupied around the $\Y$ points by the $\X$ points, respectively)
for the association alternatives.

\begin{figure}[]
\centering
\psfrag{kernel density estimate}{ \Huge{\bfseries{kernel density estimate}}}
\psfrag{relative density}{ \Huge{\bfseries{relative density}}}
\rotatebox{-90}{ \resizebox{2. in}{!}{ \includegraphics{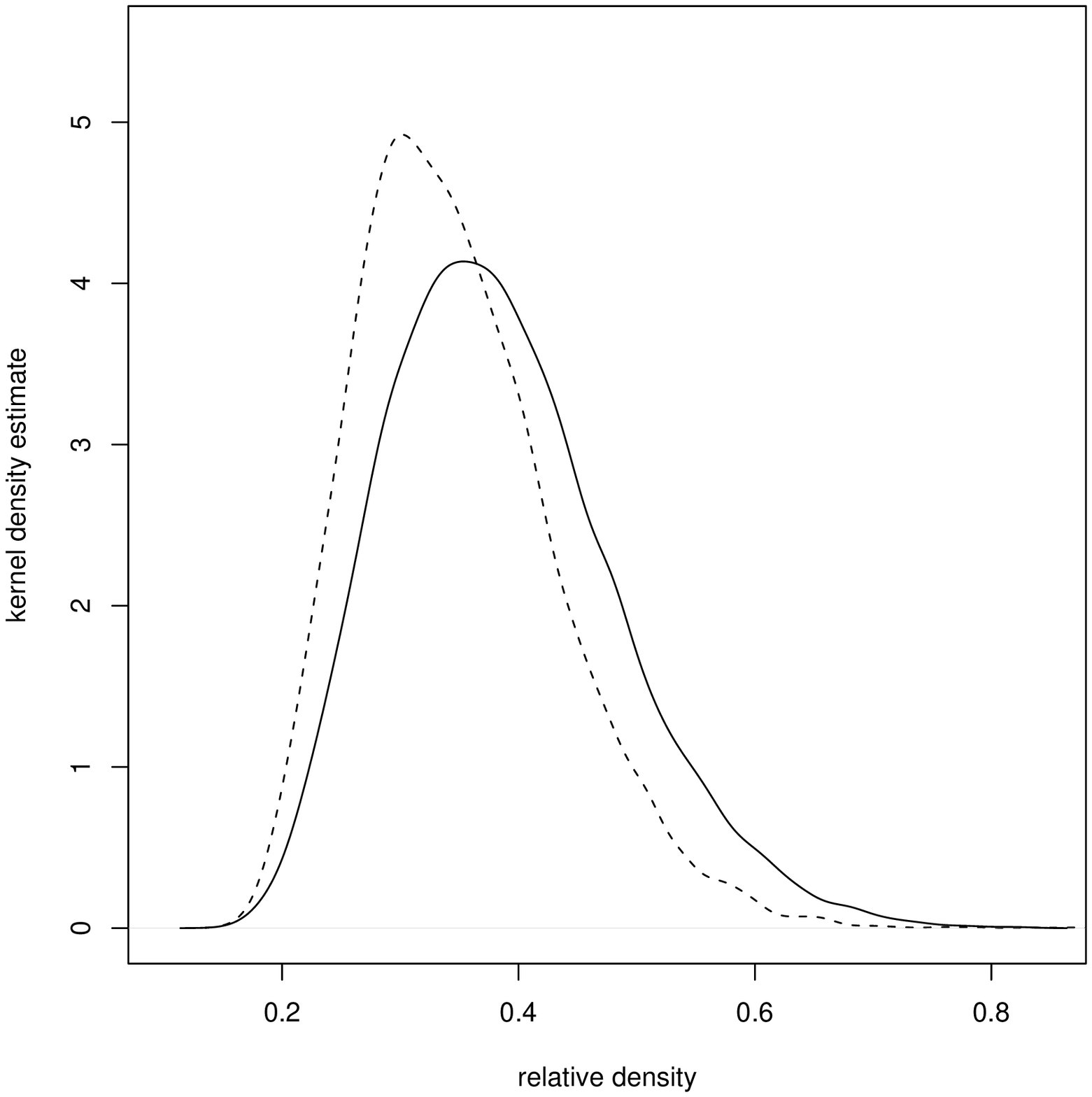}}}
\rotatebox{-90}{ \resizebox{2. in}{!}{ \includegraphics{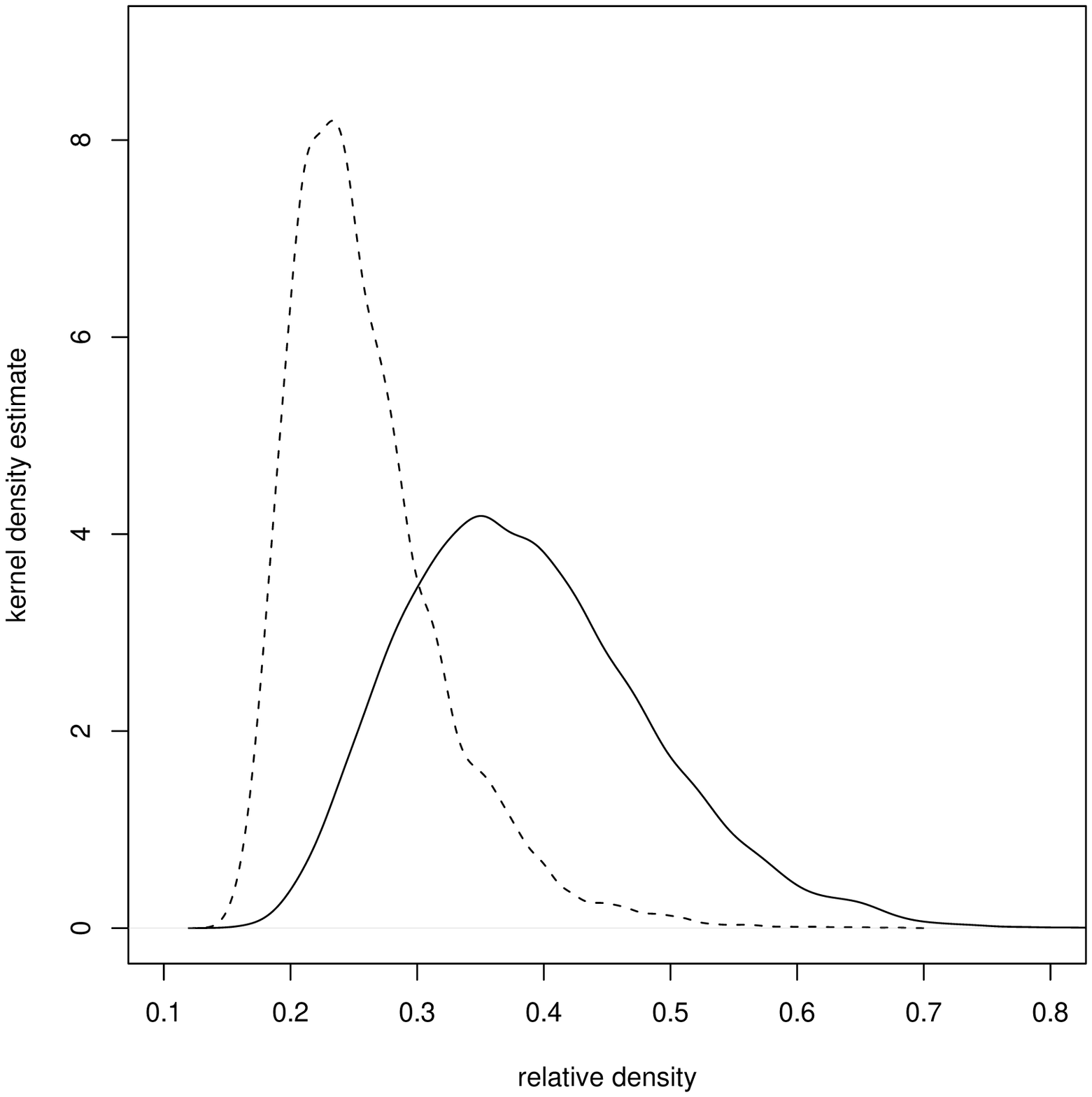}}}
\caption{
\label{fig:aggsim}
Kernel density estimates of the relative density of proportional-edge PCD,
$\rho_{_{PE}}(n,r)$,
under the null (solid line) and the association alternatives (dashed line) with
$H^A_{\sqrt{3}/21}$ (left) and  $H^A_{5\,\sqrt{3}/24}$ (right) for $r=3/2$
with $n=10$ based on $N_{mc}=10000$ replicates.
}
\end{figure}

\begin{figure}[ht]
\centering
\psfrag{kernel density estimate}{ \Huge{\bfseries{kernel density estimate}}}
\psfrag{relative density}{ \Huge{\bfseries{relative density}}}
\rotatebox{-90}{ \resizebox{2. in}{!}{ \includegraphics{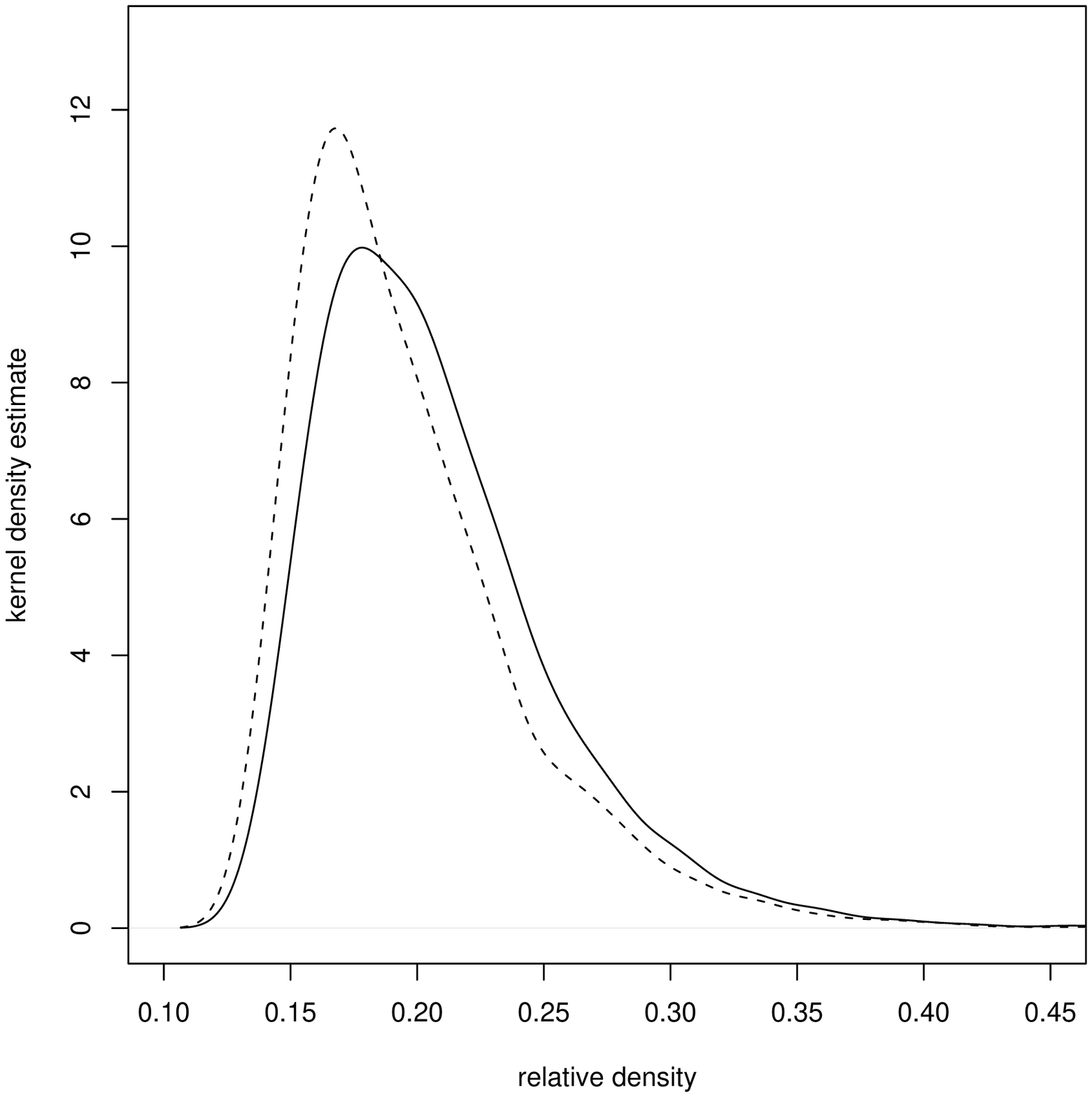}}}
\rotatebox{-90}{ \resizebox{2. in}{!}{ \includegraphics{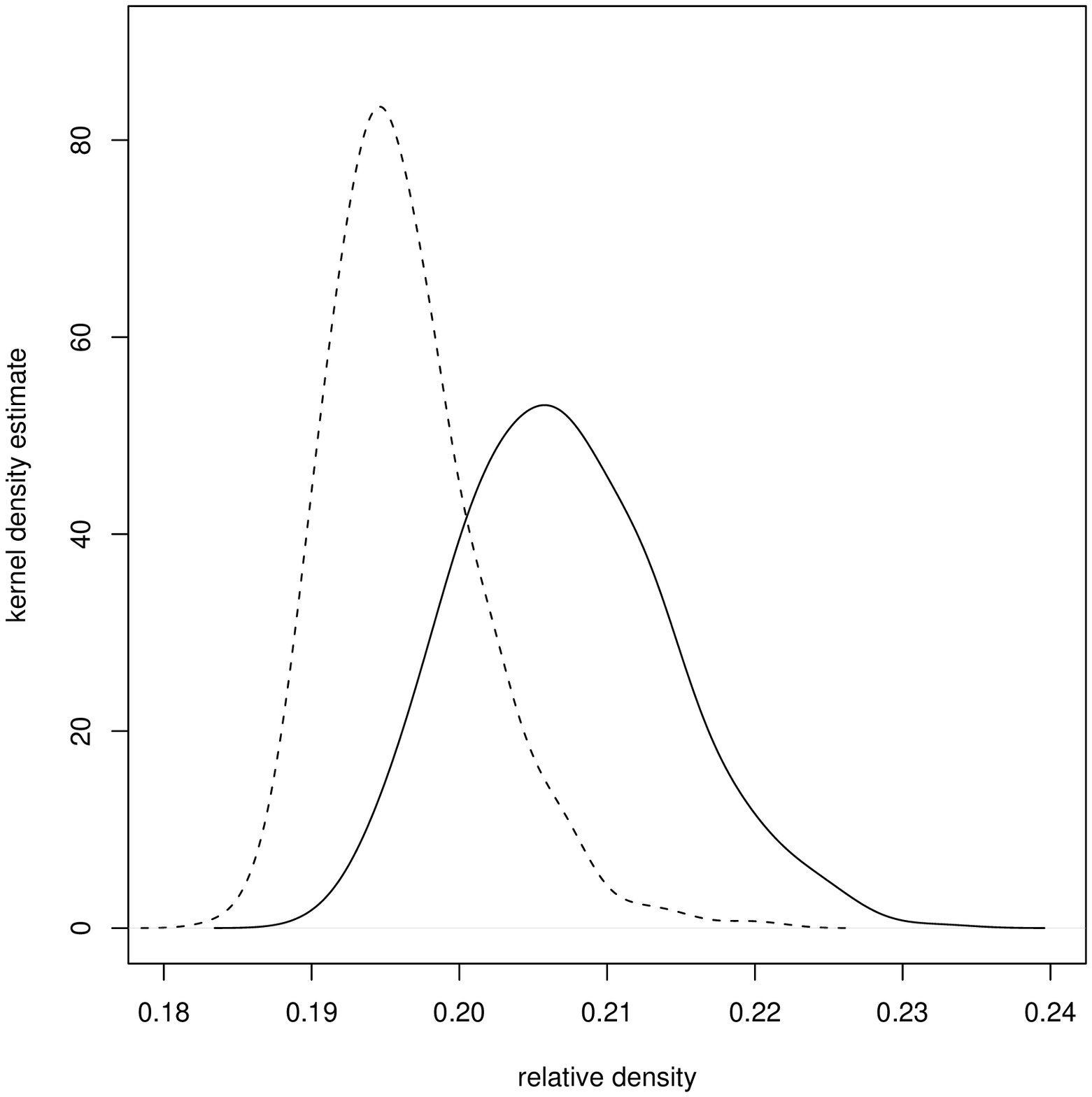}}}
\caption{
\label{fig:AggSimPowerPlots}
Depicted are kernel density estimates for $\rho_{_{PE}}(n,11/10)$ for
$n=10$ (left) and $n=100$ (right) under the null (solid line) and association alternative $H^A_{\sqrt{3}/12}$ (dashed line).
}
\end{figure}

In the one triangle case,
we plot
the kernel density estimates for the null
case and the association alternative with $\ve=\sqrt{3}/21$ and $\ve=5\,\sqrt{3}/24$
with $n=10$ and $N_{mc}=10000$
in Figure \ref{fig:aggsim}.
Observe that under
both $H_o$ and alternatives, kernel density estimates
are almost symmetric for $r=3/2$.
Moreover, there is more separation between
the kernel density estimates of the null and alternatives
for $\ve=5\,\sqrt{3}/24$ compared to $\ve=\sqrt{3}/21$,
implying more power for larger $\ve$ values.
In Figure \ref{fig:AggSimPowerPlots},
we present  kernel density estimates for the null and
the association alternative $H^A_{\sqrt{3}/12}$ for $r=11/10$,
and $n=10$, $N_{mc}=10000$ (left), $n=100$, $N_{mc}=1000$ (right).
With $n=10$, the null and alternative kernel density functions
for $\rho_{10}(11/10)$ are very similar, implying small power.
With $n=100$,
there is more separation
between null and alternative kernel density functions
implying higher power.
Notice also that the probability density functions are more skewed for $n=10$,
while approximate normality holds for $n=100$.

\begin{figure}[ht]
\centering
\rotatebox{-90}{ \resizebox{1.7 in}{!}{ \includegraphics{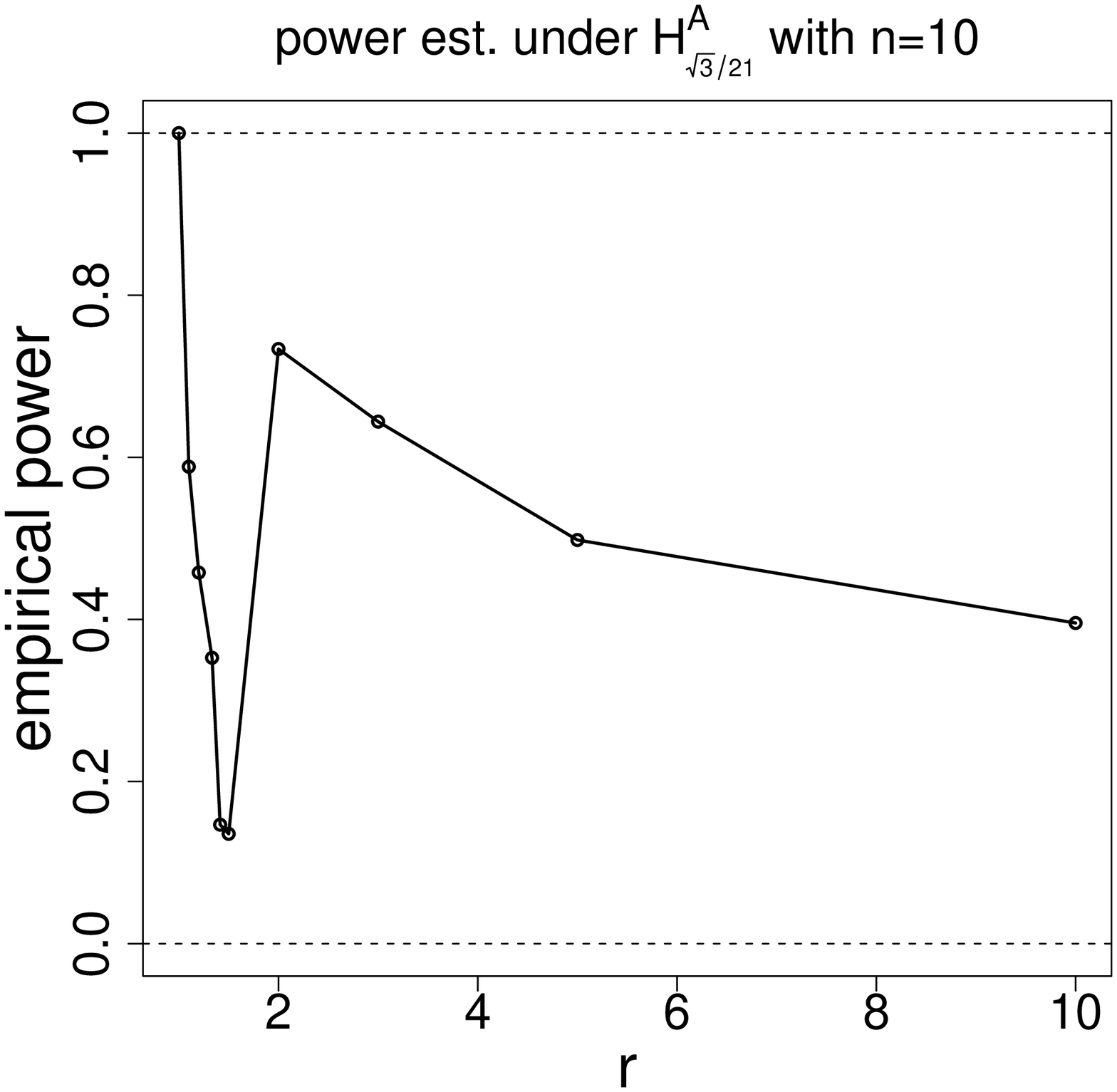}}}
\rotatebox{-90}{ \resizebox{1.7 in}{!}{ \includegraphics{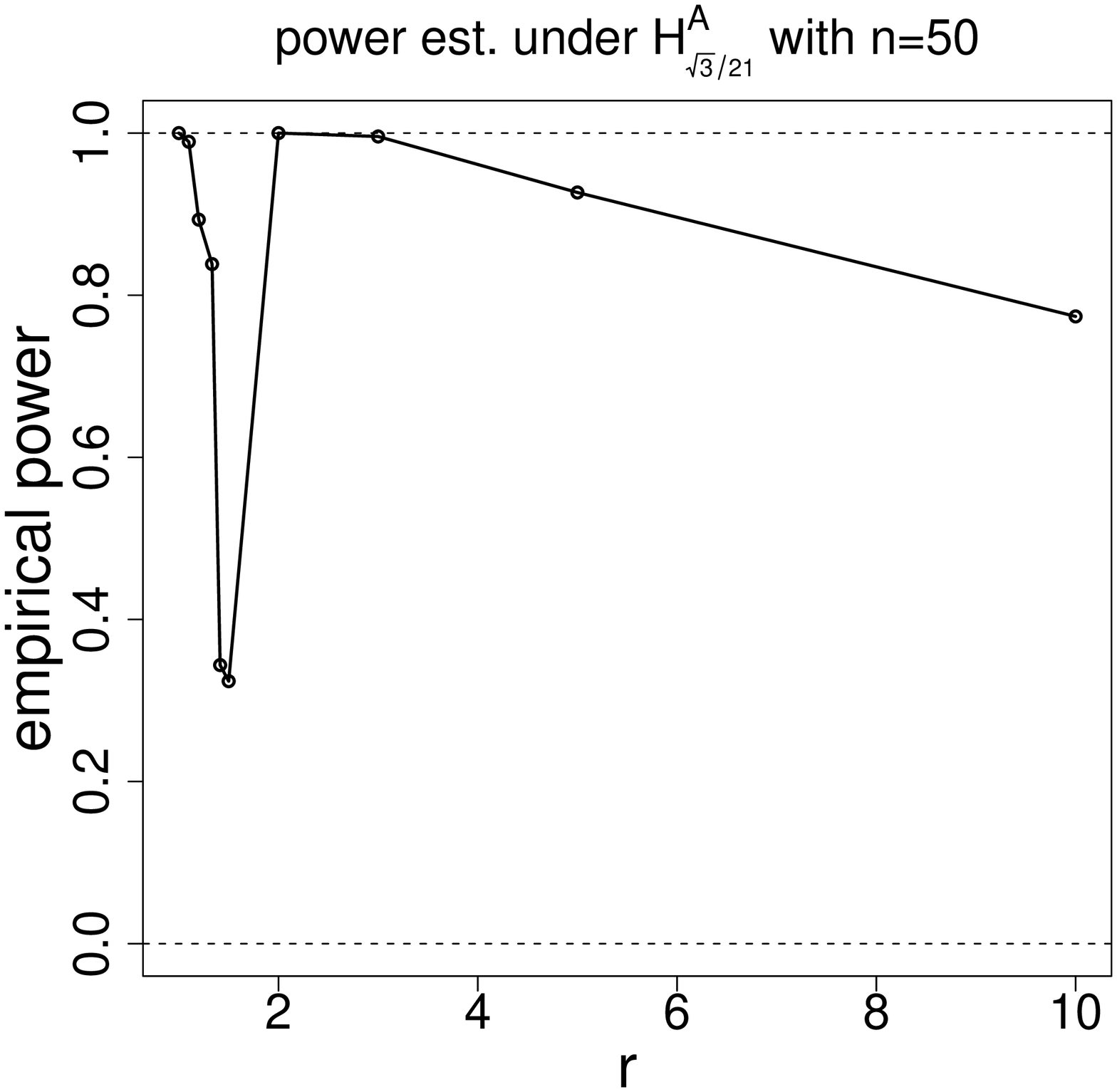}}}
\rotatebox{-90}{ \resizebox{1.7 in}{!}{ \includegraphics{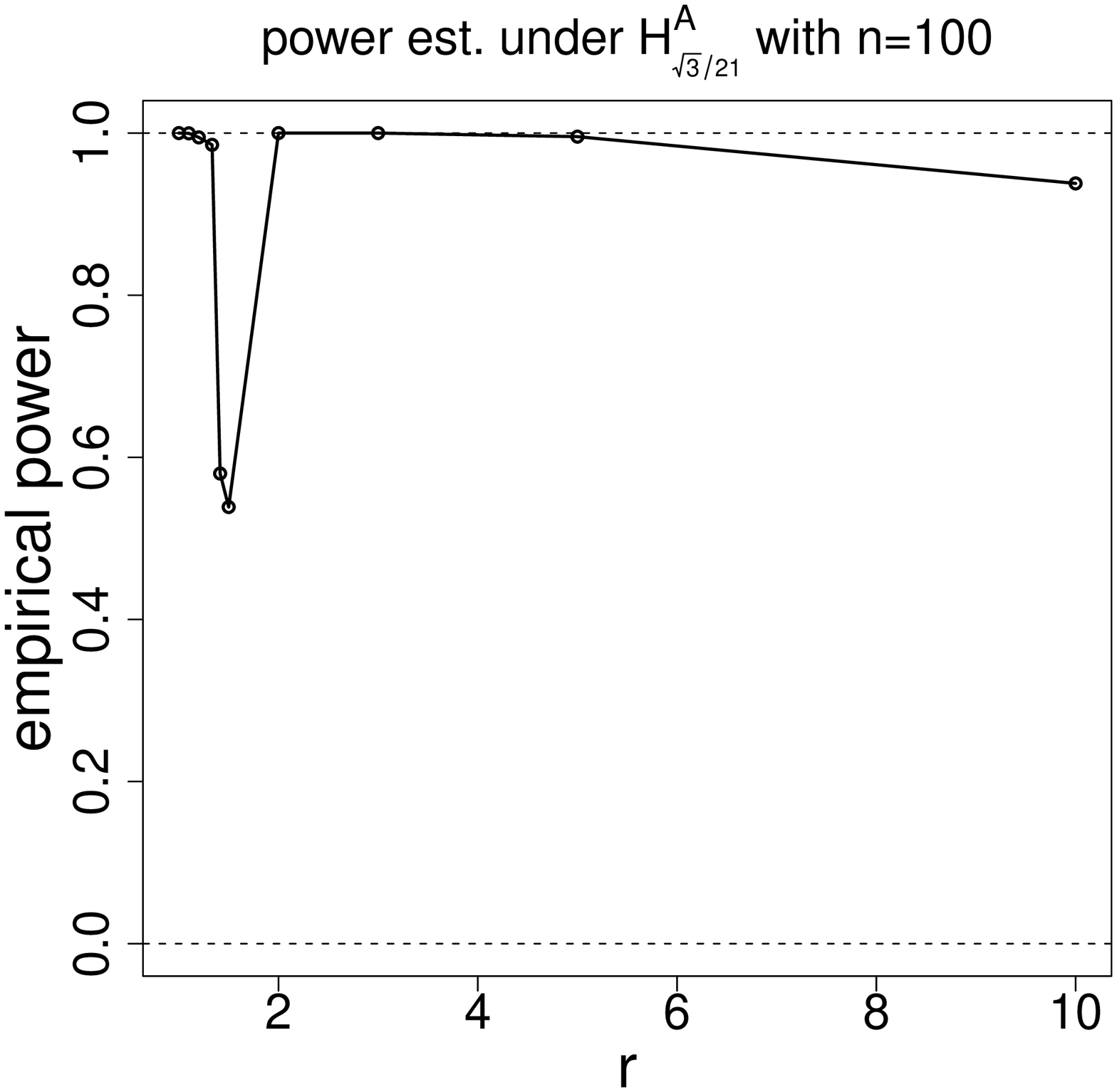}}}
\rotatebox{-90}{ \resizebox{1.7 in}{!}{ \includegraphics{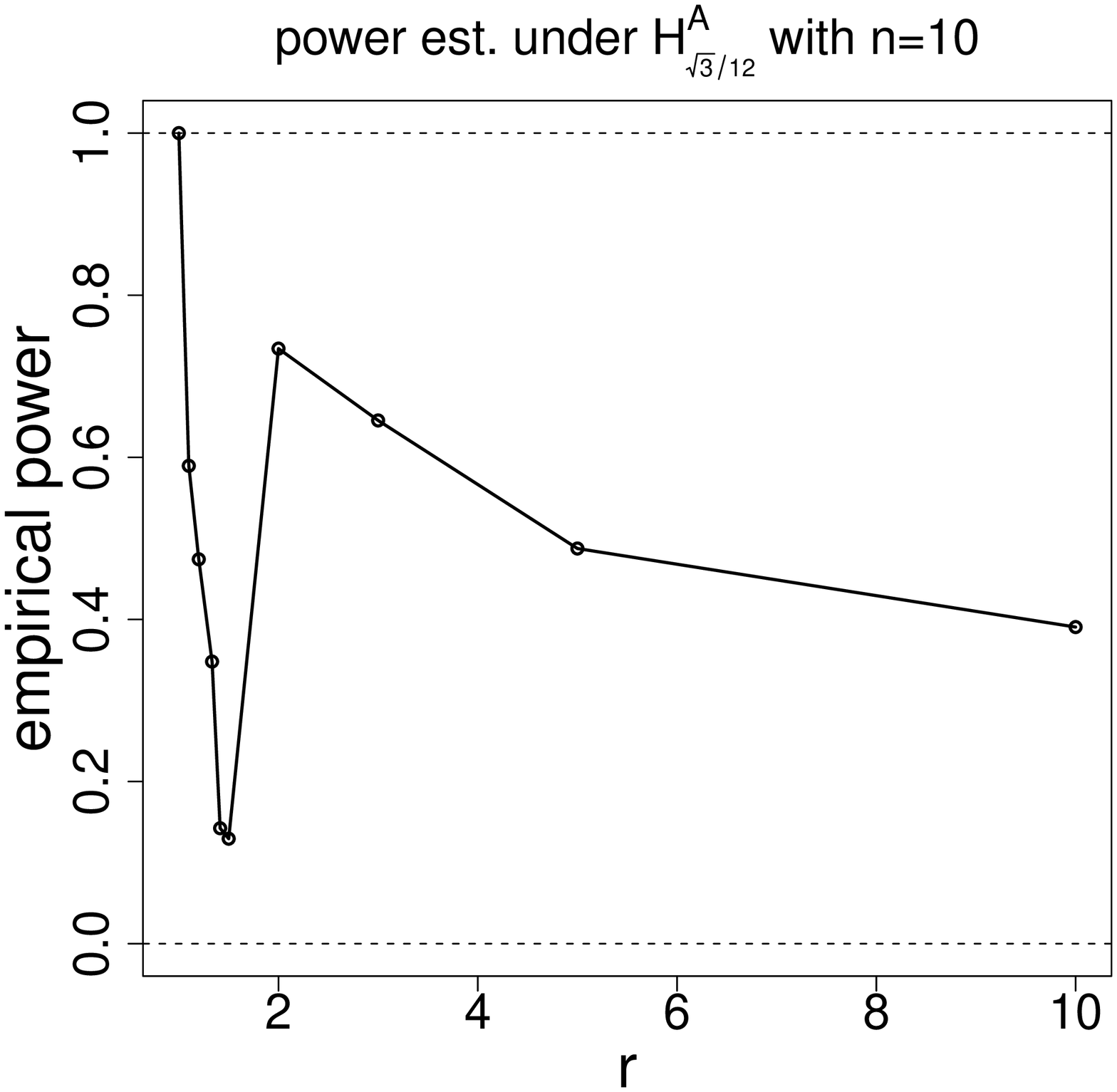}}}
\rotatebox{-90}{ \resizebox{1.7 in}{!}{ \includegraphics{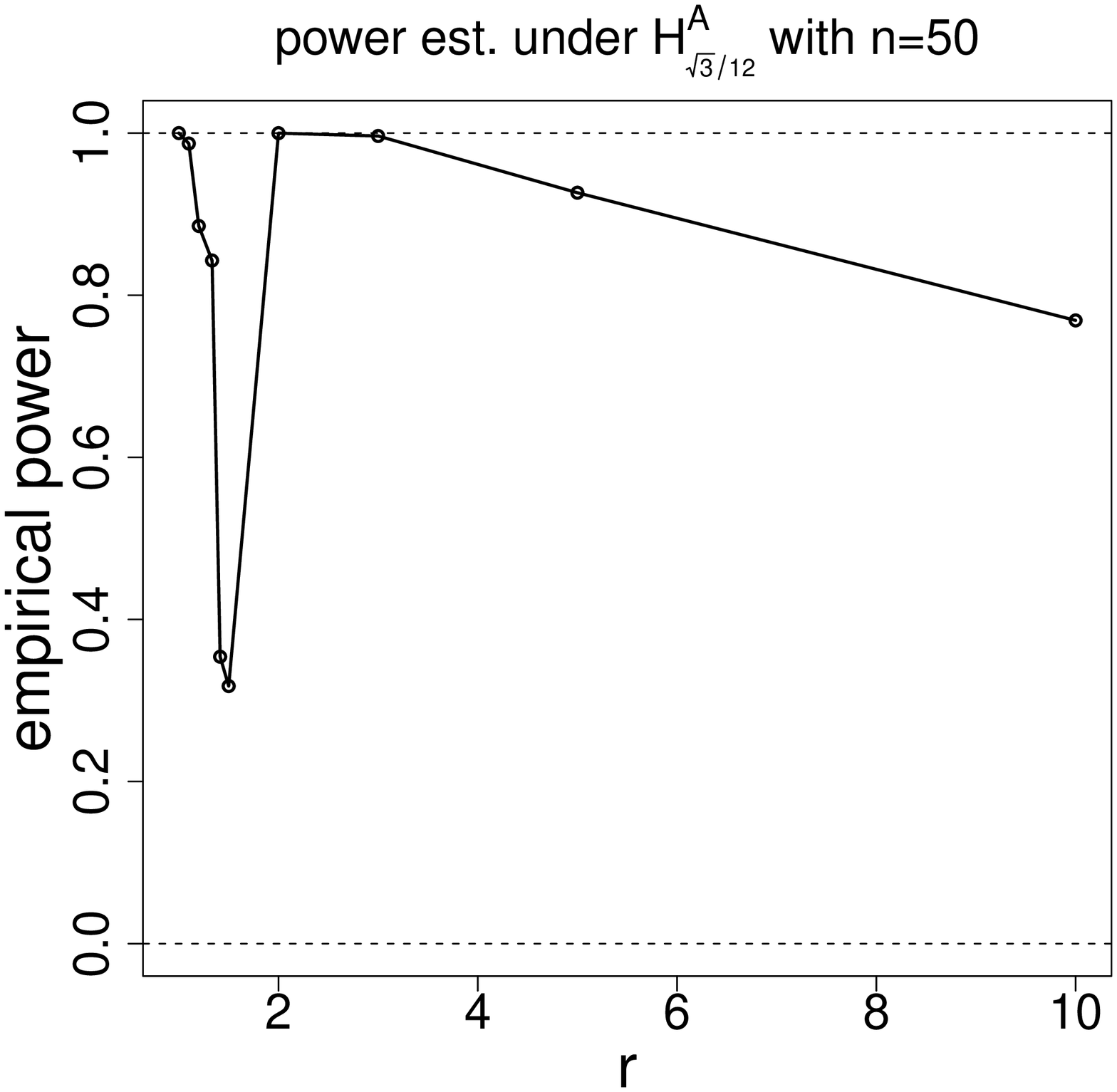}}}
\rotatebox{-90}{ \resizebox{1.7 in}{!}{ \includegraphics{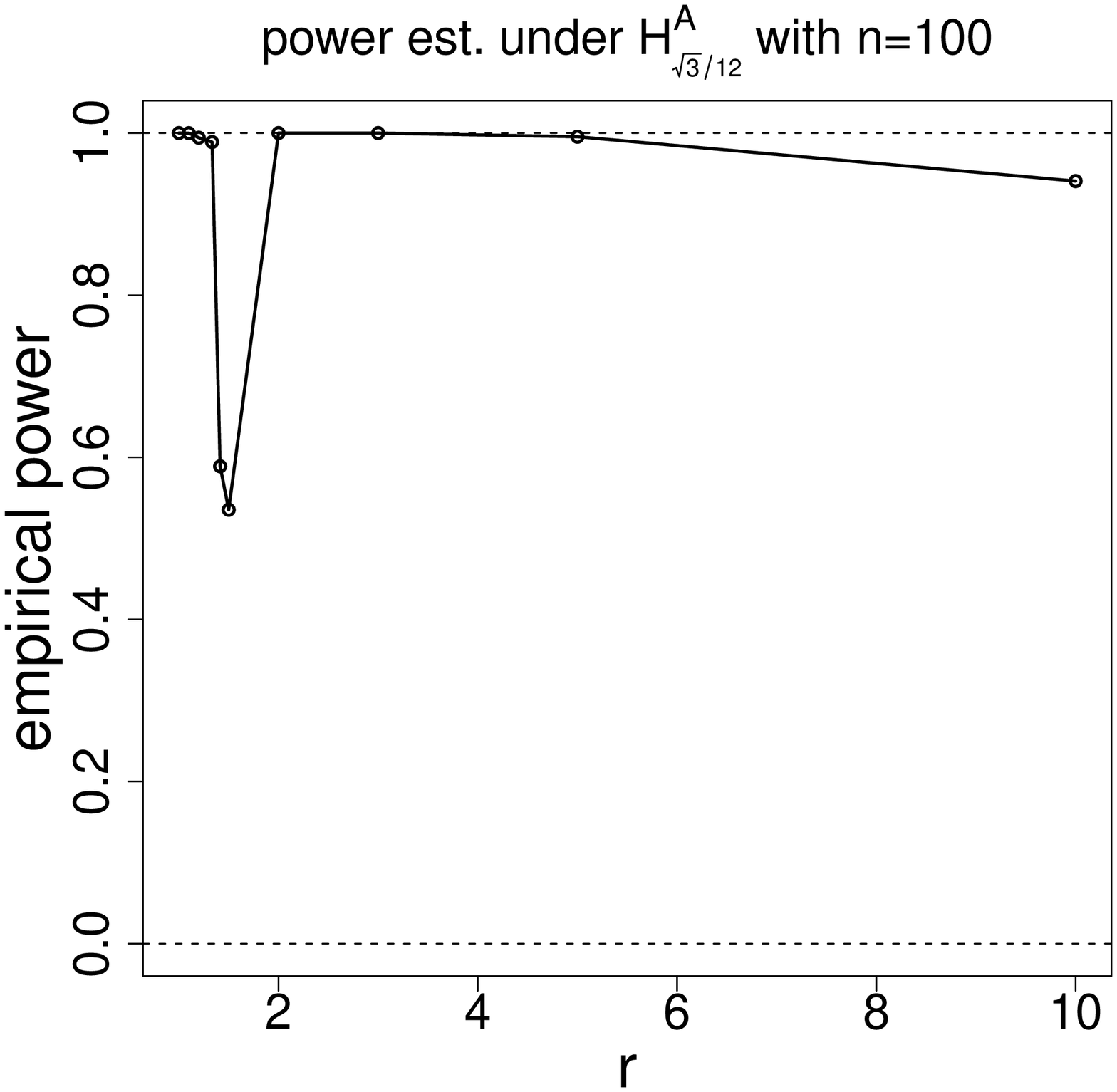}}}
\rotatebox{-90}{ \resizebox{1.7 in}{!}{ \includegraphics{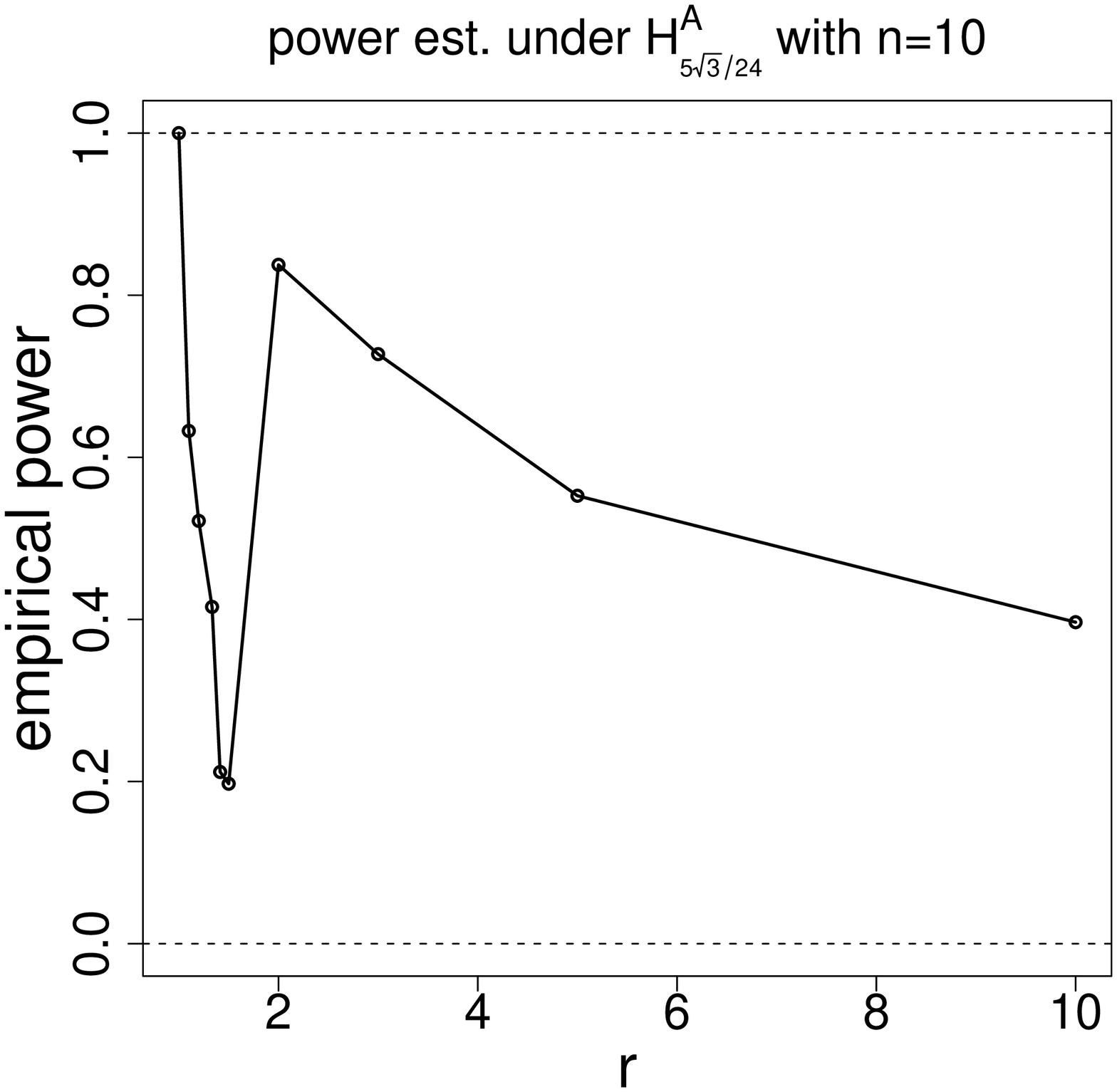}}}
\rotatebox{-90}{ \resizebox{1.7 in}{!}{ \includegraphics{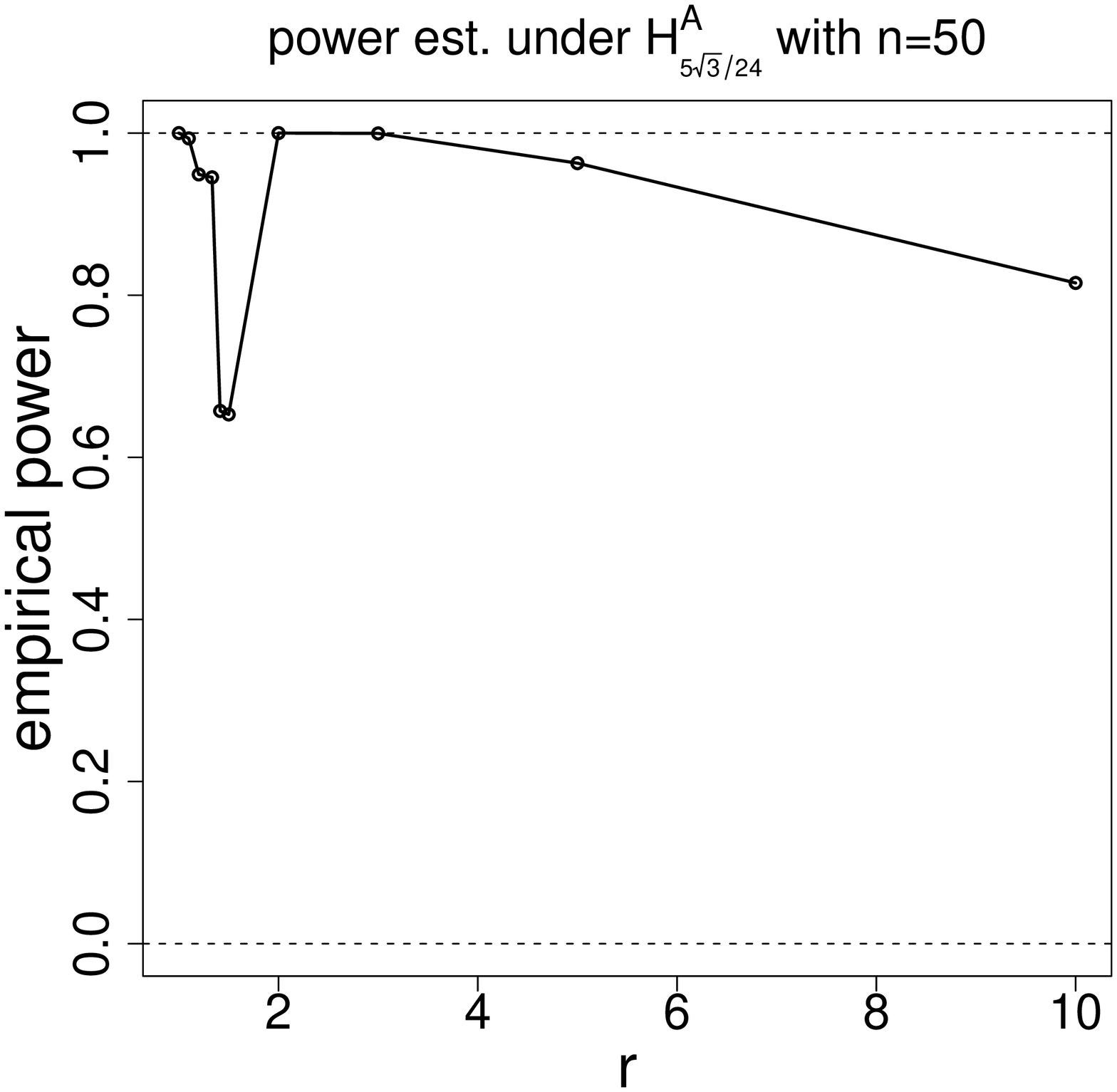}}}
\rotatebox{-90}{ \resizebox{1.7 in}{!}{ \includegraphics{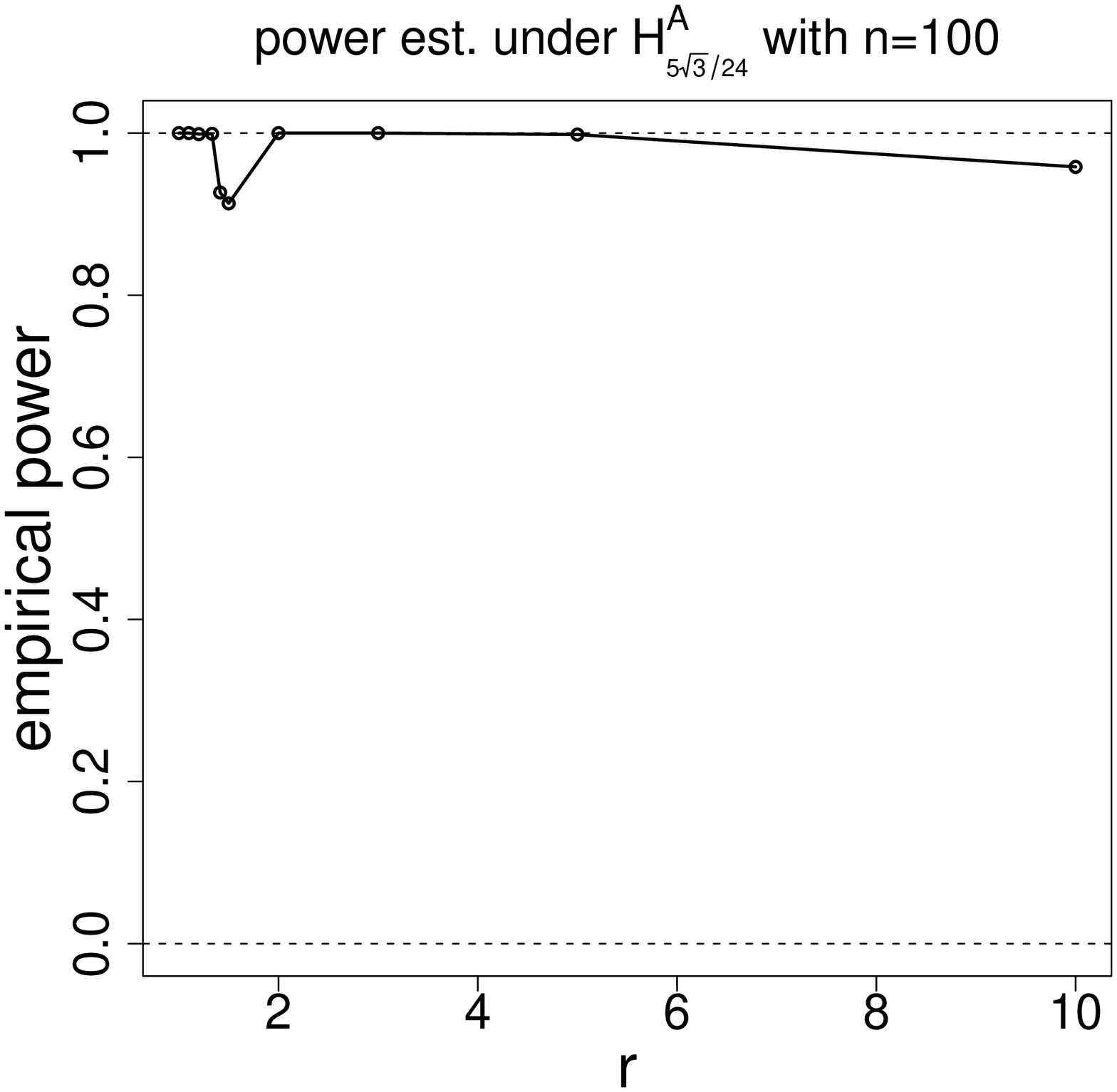}}}
\caption{
\label{fig:PE-emp-power-assoc}
\textbf{Empirical power for $R_{PE}(r)$ in the one triangle case:}
Monte Carlo power estimates for relative density of proportional-edge PCDs
in the one triangle case
using the asymptotic critical value against association alternatives
$H^A_{\sqrt{3}/21}$ (top row),
$H^A_{\sqrt{3}/12}$ (middle row),
 and
$H^A_{5\,\sqrt{3}/24}$ (bottom row)
as a function of $r$, for $n=10$ (left column), $n=50$ (middle column), and $n=100$ (right column).
}
\end{figure}

Under association,
for each $r$ value,
the level $\alpha$ asymptotic critical value is
$\mu_{_{PE}}(r) z_{\alpha} \cdot \sqrt{\nu_{_{PE}}(r)/n}$.
We estimate the empirical power as
$\frac{1}{N_{mc}}\sum_{j=1}^{N_{mc}}\I \left(R_{PE}(r)(r,j) < z_{\alpha} \right)$.
In Figure \ref{fig:PE-emp-power-assoc},
we present Monte Carlo power estimates for relative density of proportional-edge PCDs
in the one triangle case
against $H^A_{5\,\sqrt{3}/24}$,
$H^A_{\sqrt{3}/12}$, and $H^A_{\sqrt{3}/21}$ as a function of $r$ for $n=10,50,100$.
Notice that Monte Carlo power estimate
increases as $r$ gets larger and then decreases,
as in the segregation case.
Because for small $n$ and large $r$, the
critical value is approximately one under $H_o$,
as we get a nearly complete
digraph with high probability.
Moreover, the more severe the association, the higher the power estimate at each $r$.
Highest power is attained for $r \approx 2$,
which is recommended against the association,
as it yields the desired level with high power.

\begin{figure}[]
\centering
\rotatebox{-90}{ \resizebox{1.7 in}{!}{ \includegraphics{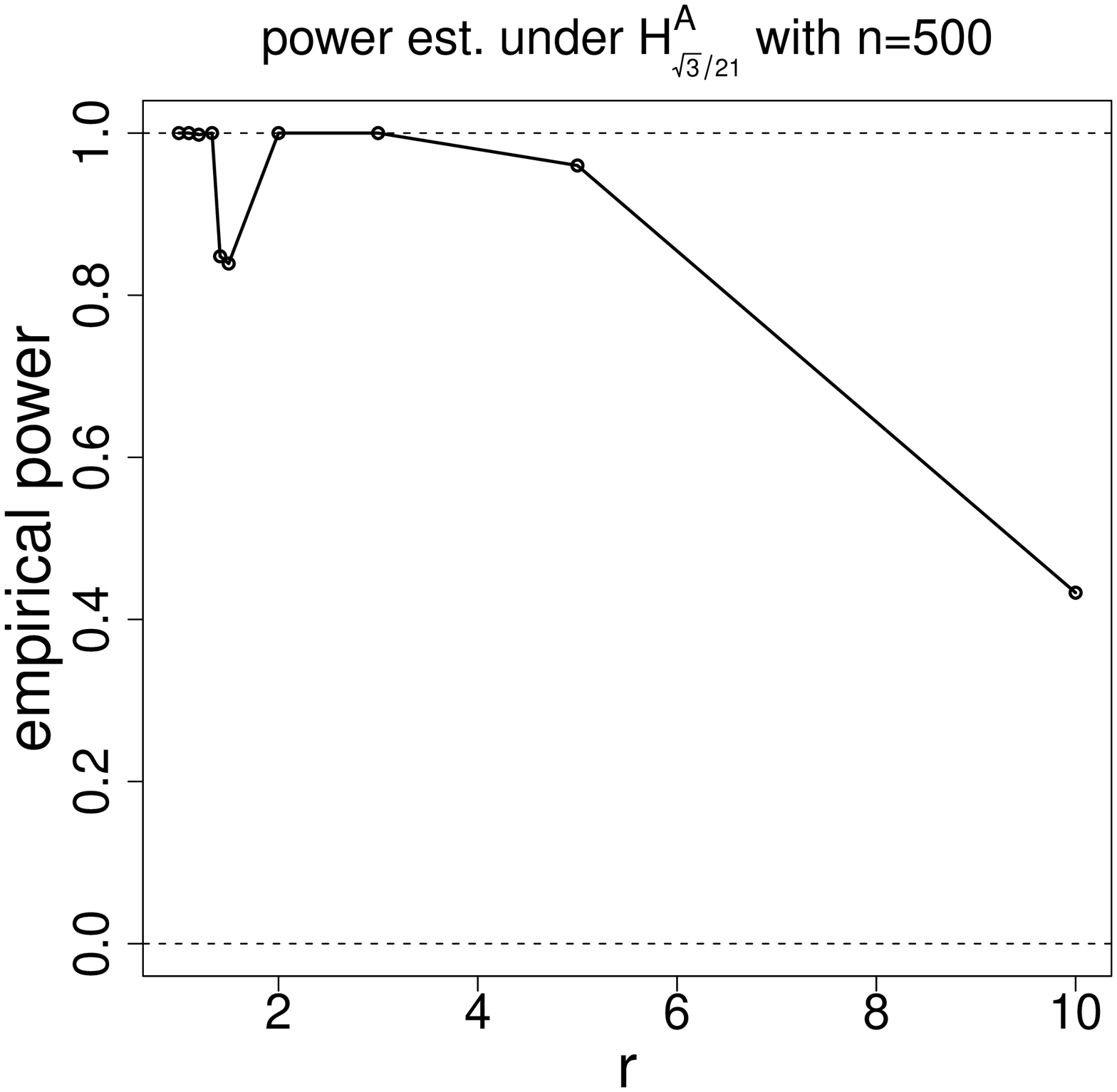}}}
\rotatebox{-90}{ \resizebox{1.7 in}{!}{ \includegraphics{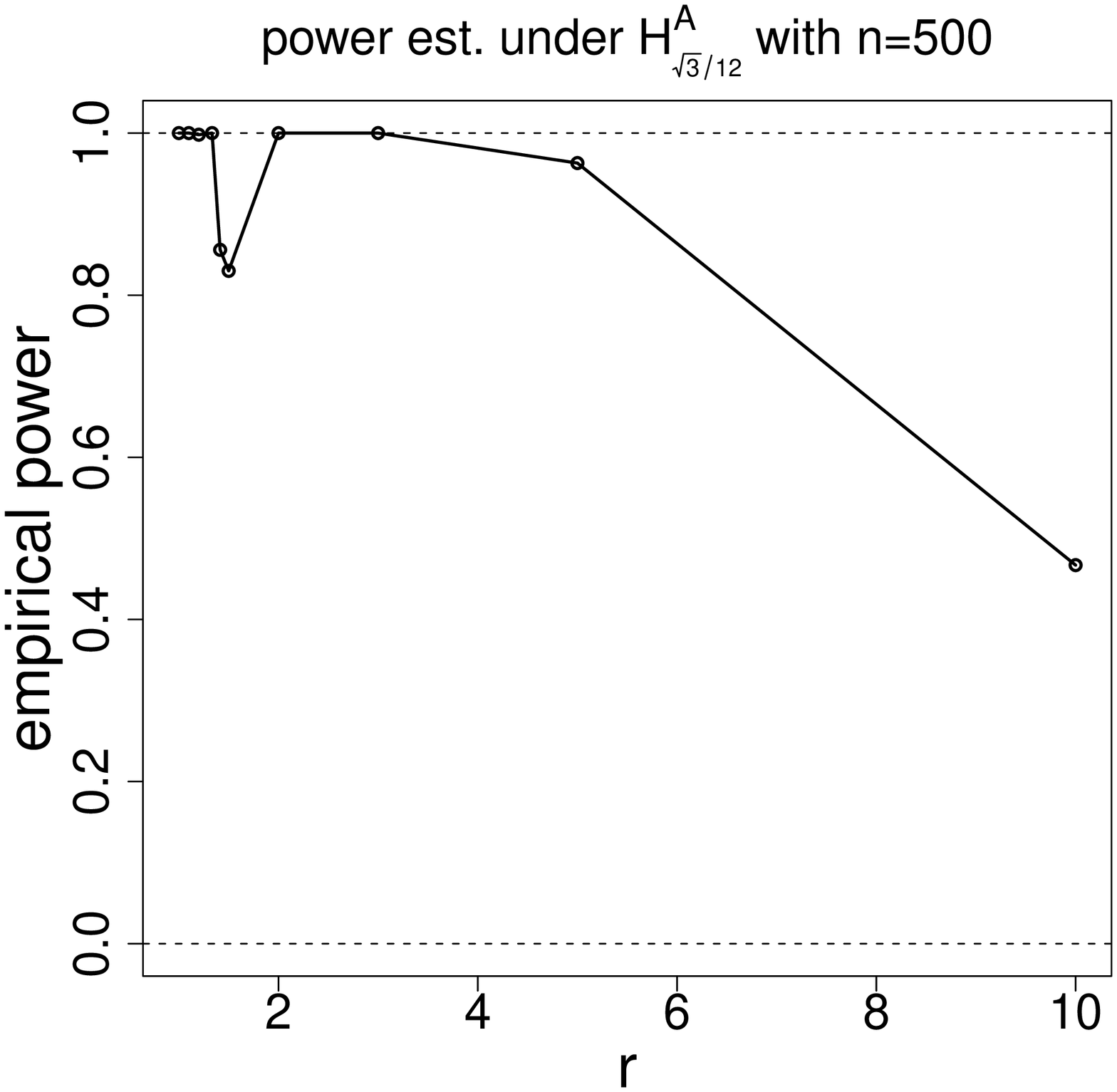}}}
\rotatebox{-90}{ \resizebox{1.7 in}{!}{ \includegraphics{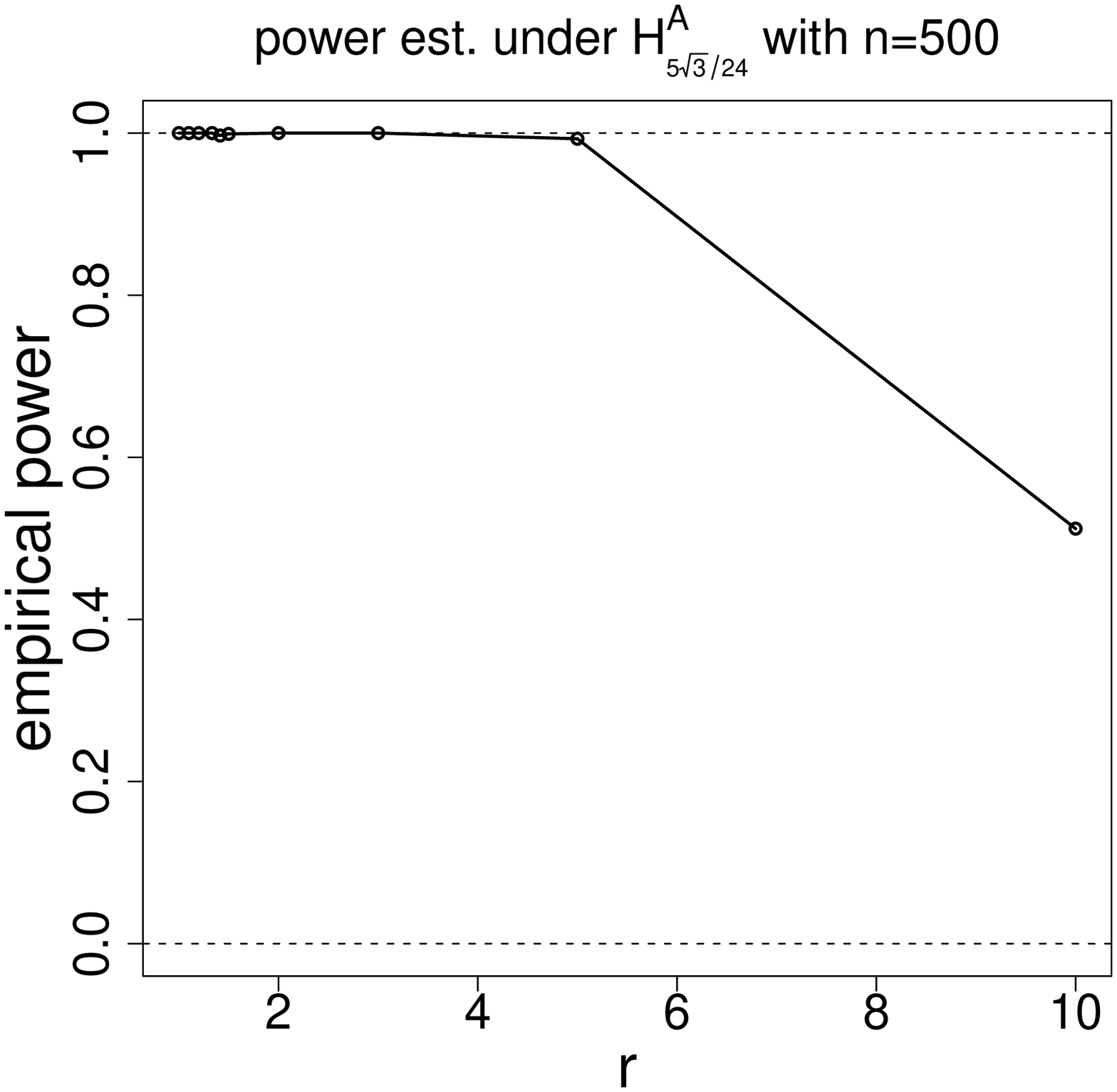}}}
\rotatebox{-90}{ \resizebox{1.7 in}{!}{ \includegraphics{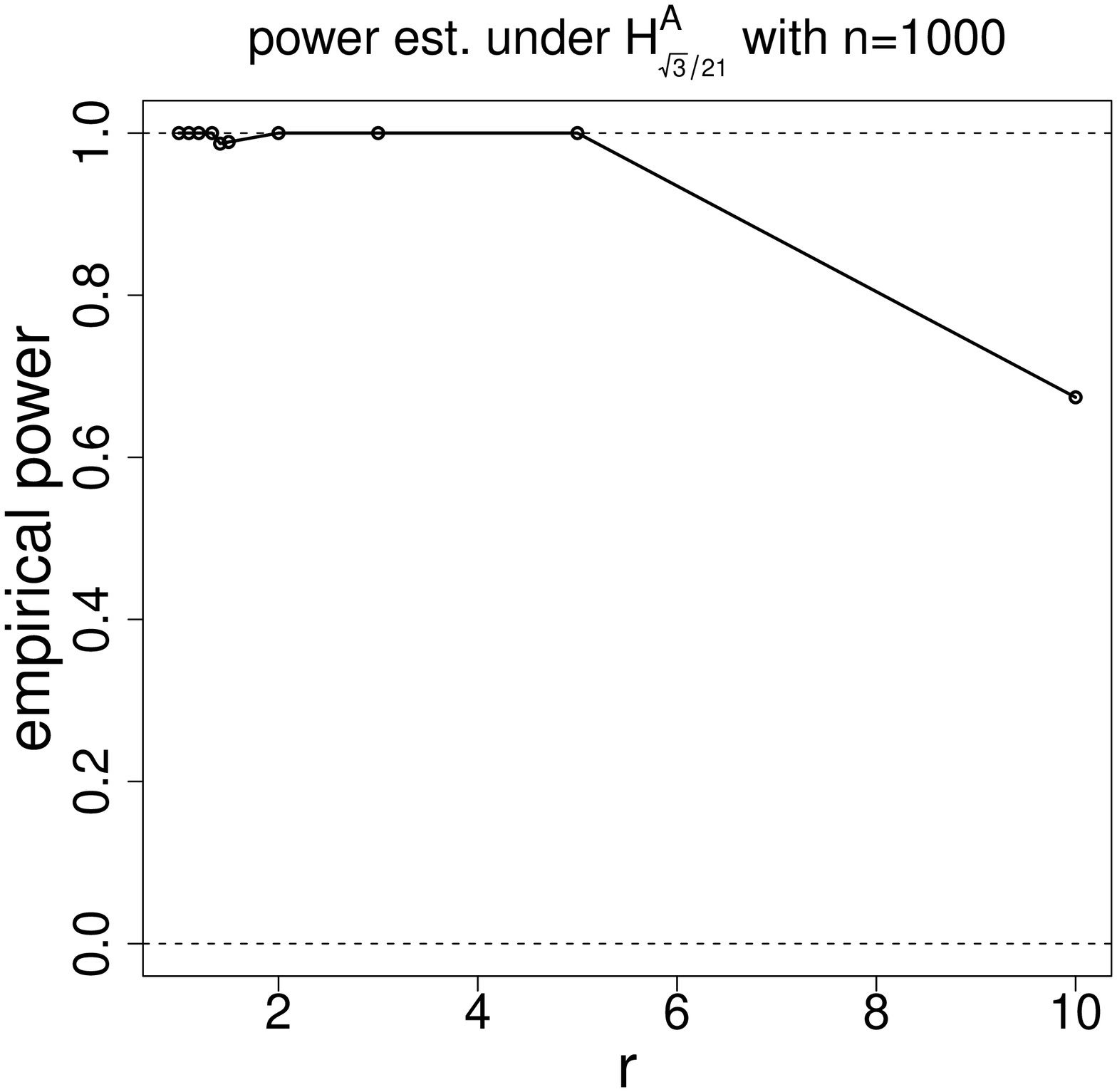}}}
\rotatebox{-90}{ \resizebox{1.7 in}{!}{ \includegraphics{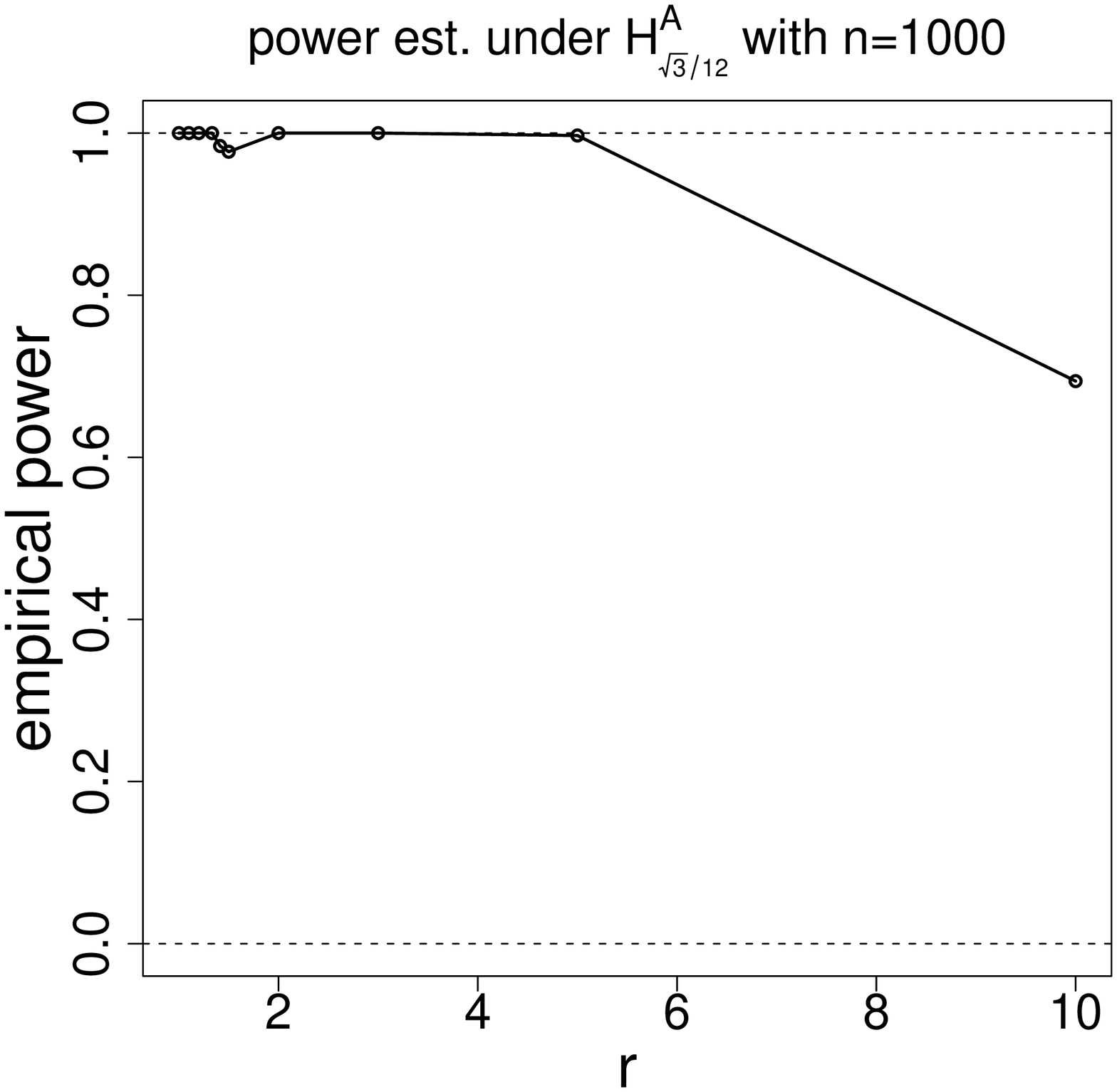}}}
\rotatebox{-90}{ \resizebox{1.7 in}{!}{ \includegraphics{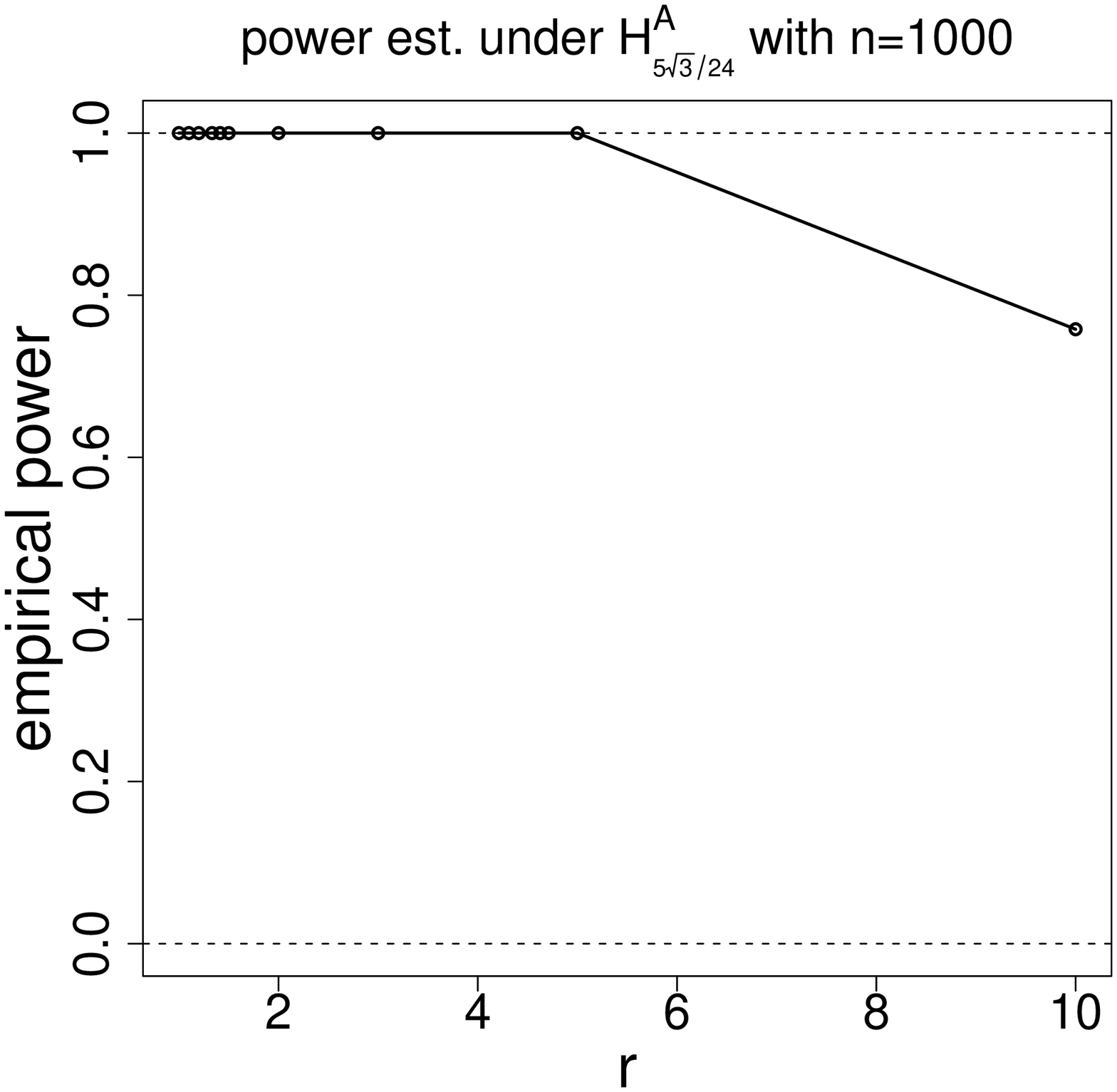}}}
\caption{
\label{fig:MT-PE-emp-power-assoc}
\textbf{Empirical power for $R_{PE}(r)$ in the multiple triangle case:}
Monte Carlo power estimates of the relative density of proportional-edge PCDs
in the multiple triangle case
using the asymptotic critical value against association alternatives
$H^A_{\sqrt{3}/21}$ (left column),
$H^A_{\sqrt{3}/12}$ (middle column),
and
$H^A_{5\,\sqrt{3}/24}$ (right column)
as a function of $r$, for $n=500$ (top) and $n=1000$ (bottom).
}
\end{figure}

In the multiple triangle case,
we generate the $\X$ points uniformly in the support for the association alternatives
in the triangles based on the 10 class $\Y$
points given in Figure \ref{fig:deldata}.
We use the parameters $\ve \in \{5\,\sqrt{3}/24,\sqrt{3}/12,\sqrt{3}/21\}$.
We compute the relative density based on the formula given in Corollary \ref{cor:MT-asy-norm}.
The corresponding empirical power estimates as a function of $r$ (using the normal approximation)
are presented in Figure \ref{fig:MT-PE-emp-power-assoc} for $n=500$ or 1000.
Observe that the Monte Carlo power estimate
decreases as $r$ gets larger unlike the the one triangle case.
The empirical power is large (i.e., close to one) for $r \in (1,5)$.
Considering the empirical size estimates,
we recommend $r \approx 2$
for association alternative
since the corresponding test has the desired level with high power.

\subsection{Empirical Power Analysis for Central Similarity PCDs under the Association Alternative}
\label{sec:CS-emp-power-assoc}
In the one triangle case,
we generate data as in Section \ref{sec:PE-emp-power-assoc}
and compute the relative density as in Section \ref{sec:CS-emp-power-seg}.
The distribution of $\rho_{_{CS}}(n,\tau)$ is non-degenerate
for all $\ve \in (0,\sqrt{3}/6)$ and $\tau \in (0,\infty)$.

\begin{figure}[]
\centering
\psfrag{kernel density estimate}{ \Huge{\bfseries{kernel density estimate}}}
\psfrag{relative density}{ \Huge{\bfseries{relative density}}}
\rotatebox{-90}{ \resizebox{2. in}{!}{ \includegraphics{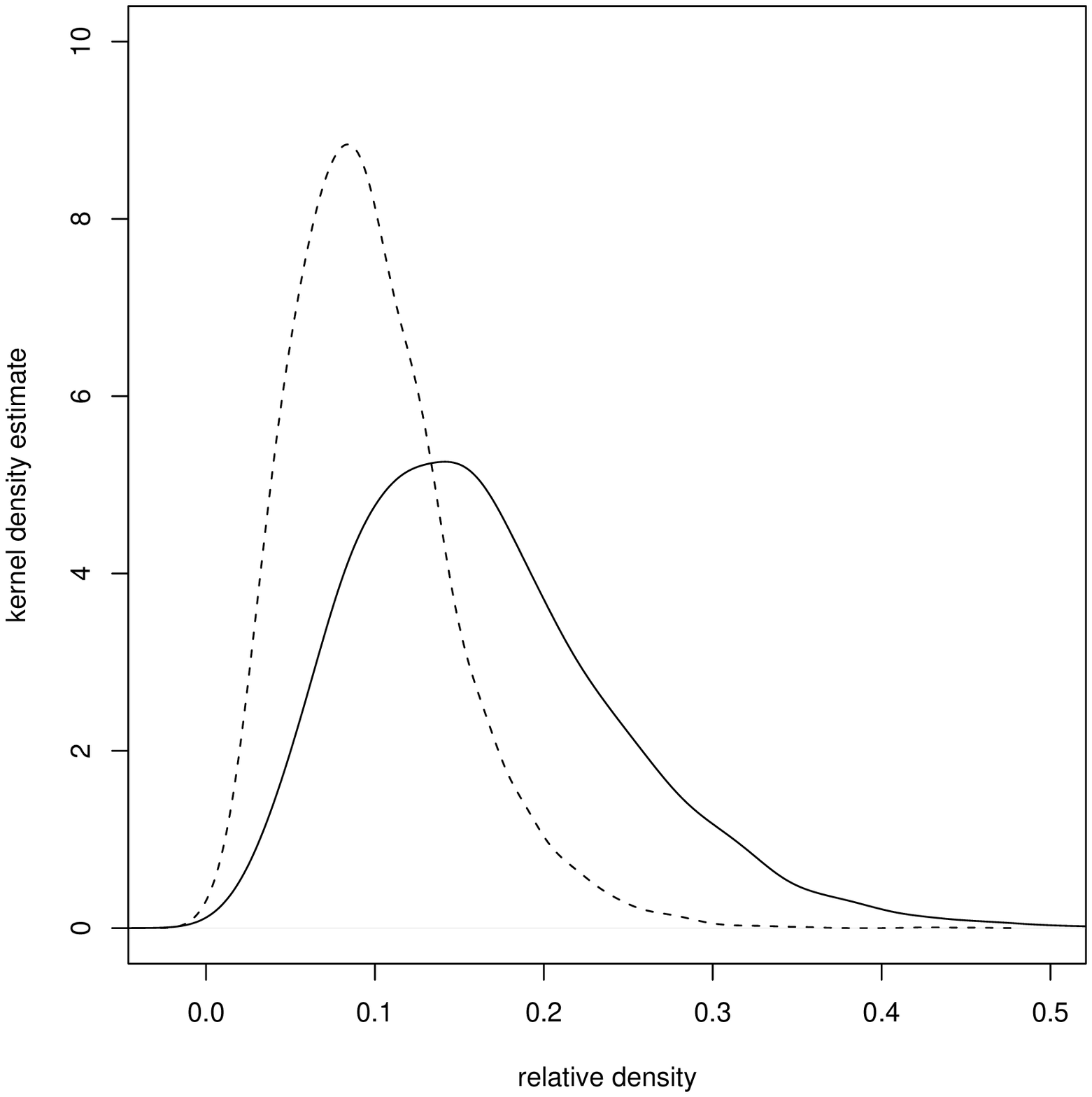}}}
\rotatebox{-90}{ \resizebox{2. in}{!}{ \includegraphics{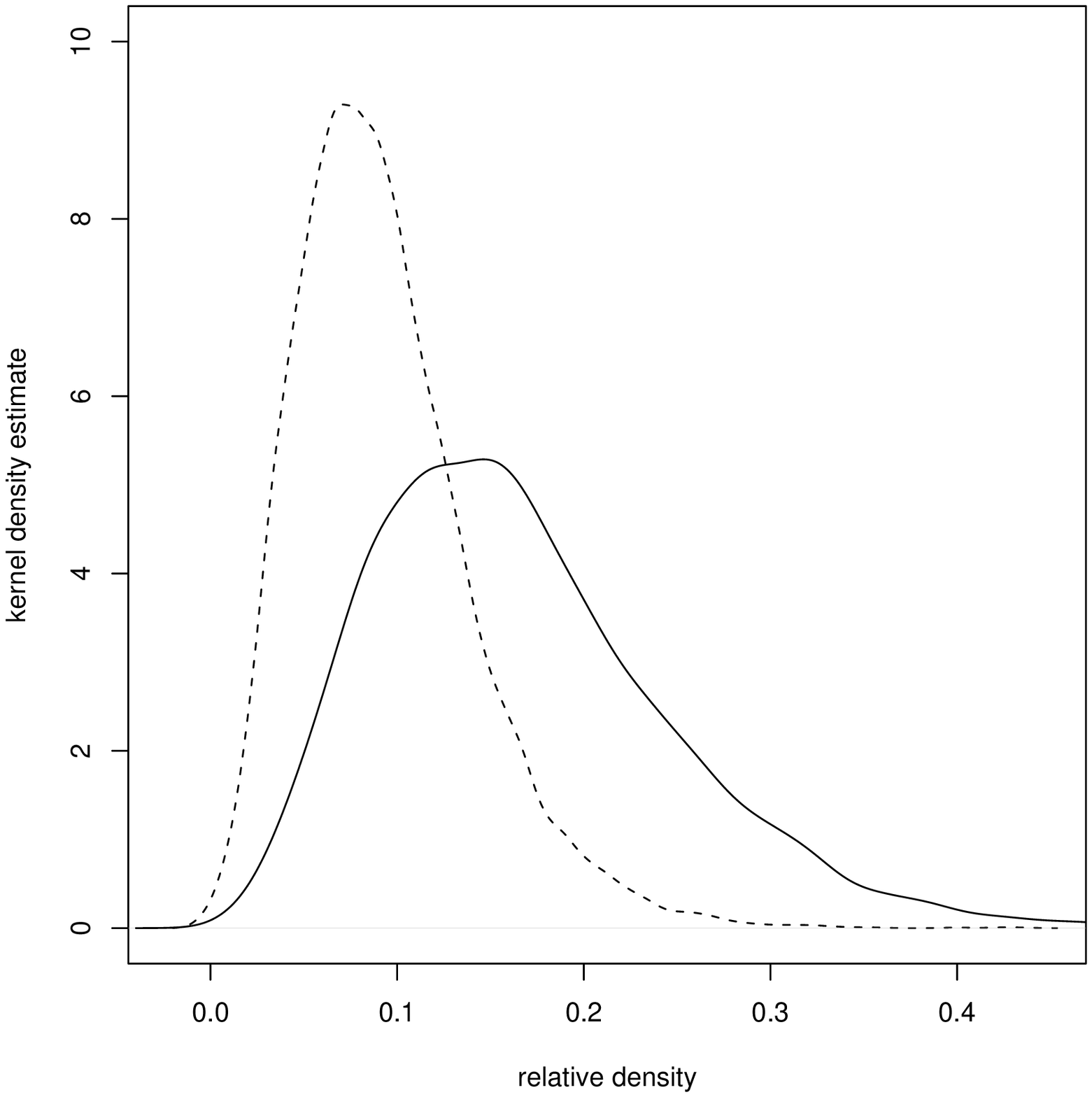}}}
\caption{
\label{fig:CSaggsim}
Kernel density estimates of the relative density of central similarity PCD,
$\rho_{_{CS}}(n,\tau)$,
under the null (solid line) and the association alternatives (dashed line) with
$H^A_{\sqrt{3}/12}$ (left) and  $H^A_{5\,\sqrt{3}/24}$ (right) for $\tau =1$
with $n=10$ based on $N_{mc}=10000$ replicates.
}
\end{figure}

In the one triangle case,
we plot
the kernel density estimates for the null
case and the association alternative with $\ve=\sqrt{3}/12$ and $\ve=5\,\sqrt{3}/24$
with $n=10$ and $N_{mc}=10000$
in Figure \ref{fig:CSaggsim}.
Observe that under
both $H_o$ and alternatives, kernel density estimates
are almost symmetric for $\tau=1$.
However, there is only mild separation between
the kernel density estimates of the null and alternatives
implying small power.

\begin{figure}[ht]
\centering
\psfrag{kernel density estimate}{ \Huge{\bfseries{kernel density estimate}}}
\psfrag{relative density}{ \Huge{\bfseries{relative density}}}
\rotatebox{-90}{ \resizebox{2. in}{!}{ \includegraphics{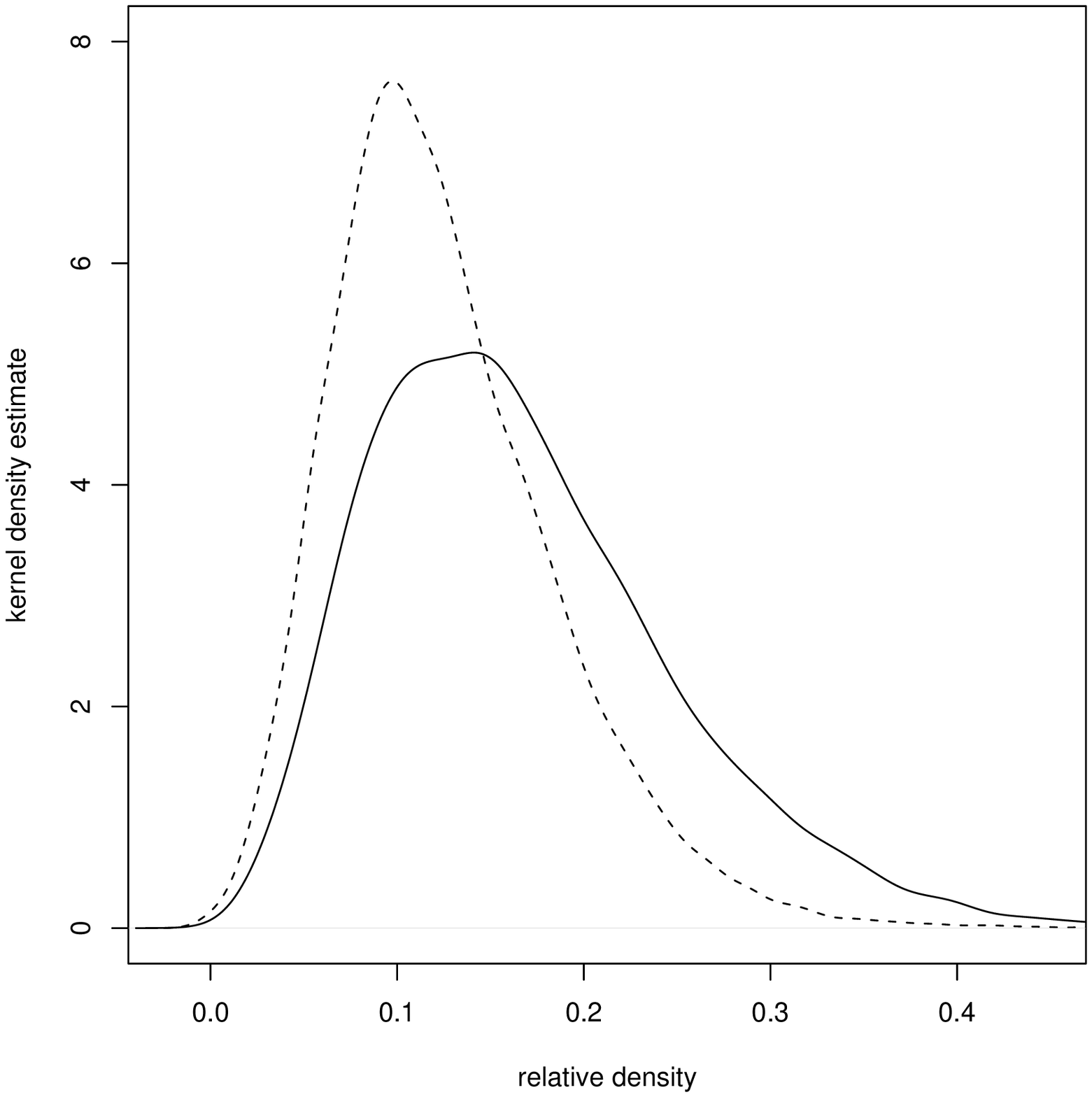}}}
\rotatebox{-90}{ \resizebox{2. in}{!}{ \includegraphics{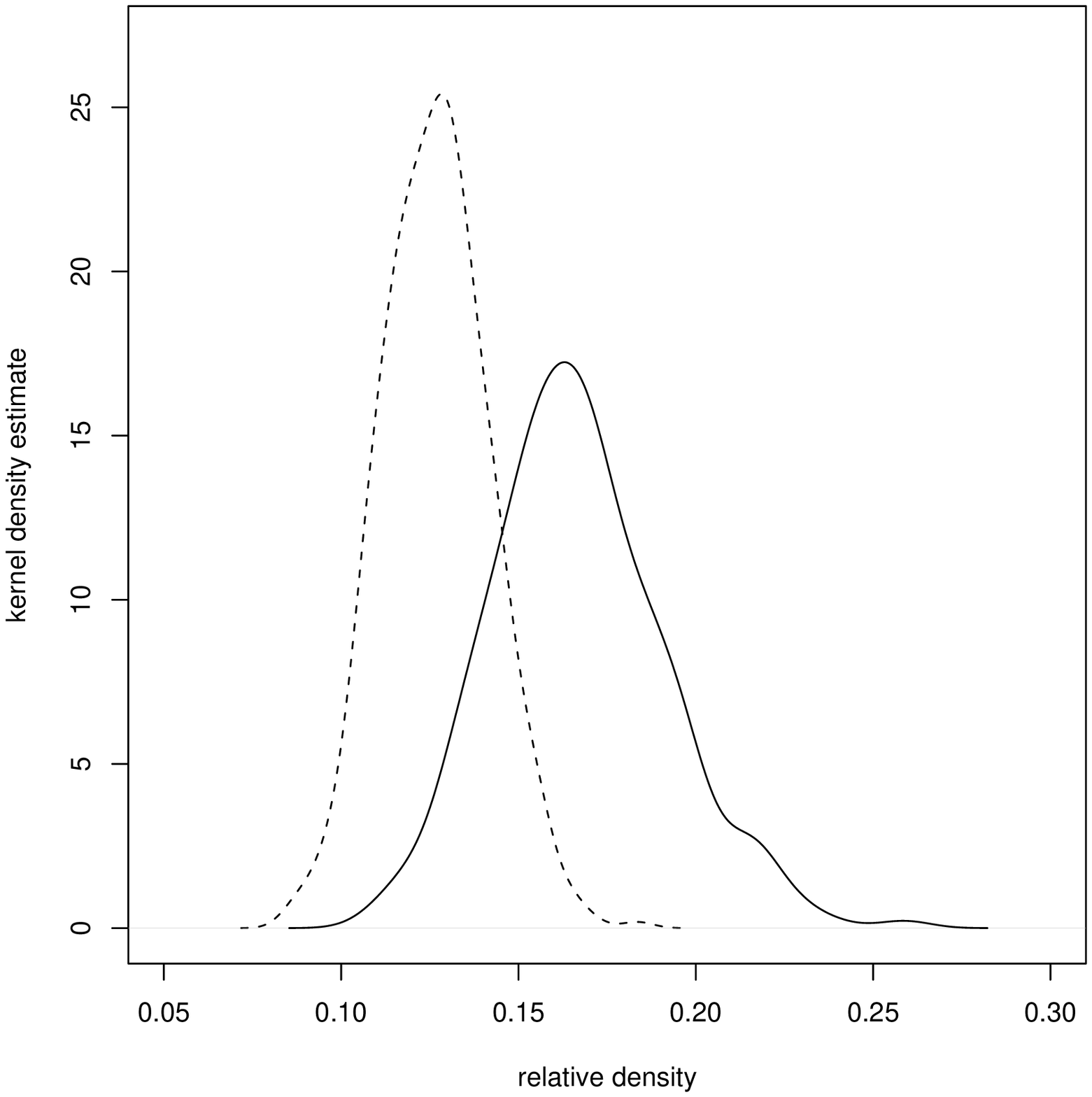}}}
\caption{
\label{fig:CSAggSimPowerPlots}
Depicted are kernel density estimates for $\rho_{_{CS}}(n,1)$ for
$n=10$ (left) and $n=100$ (right) under the null (solid line) and association alternative $H^A_{\sqrt{3}/21}$ (dashed line).
}
\end{figure}

In Figure \ref{fig:CSAggSimPowerPlots}, we present a Monte Carlo investigation
against the association alternative $H^A_{\sqrt{3}/21}$ for $\tau=1$,
and $n=10$, $N_{mc}=10000$ (left), $n=100$, $N_{mc}=1000$ (right).
With $n=10$, the null and alternative kernel density functions
for $\rho_{_{CS}}(n,1)$ are very similar, implying small power.
With $n=100$,
there is more separation
between null and alternative kernel density functions,
implying higher power.

\begin{figure}[ht]
\centering
\rotatebox{-90}{ \resizebox{1.7 in}{!}{ \includegraphics{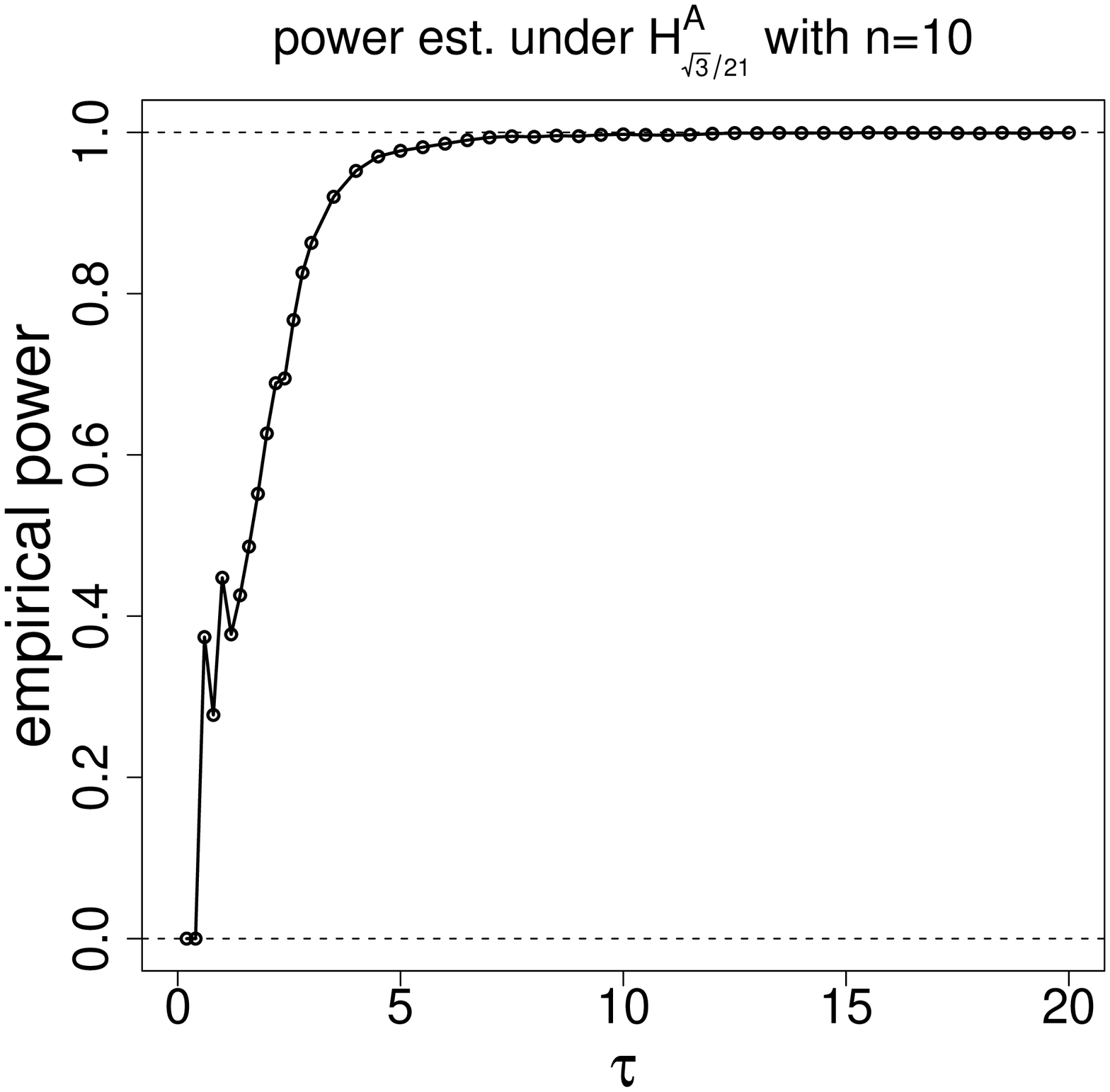}}}
\rotatebox{-90}{ \resizebox{1.7 in}{!}{ \includegraphics{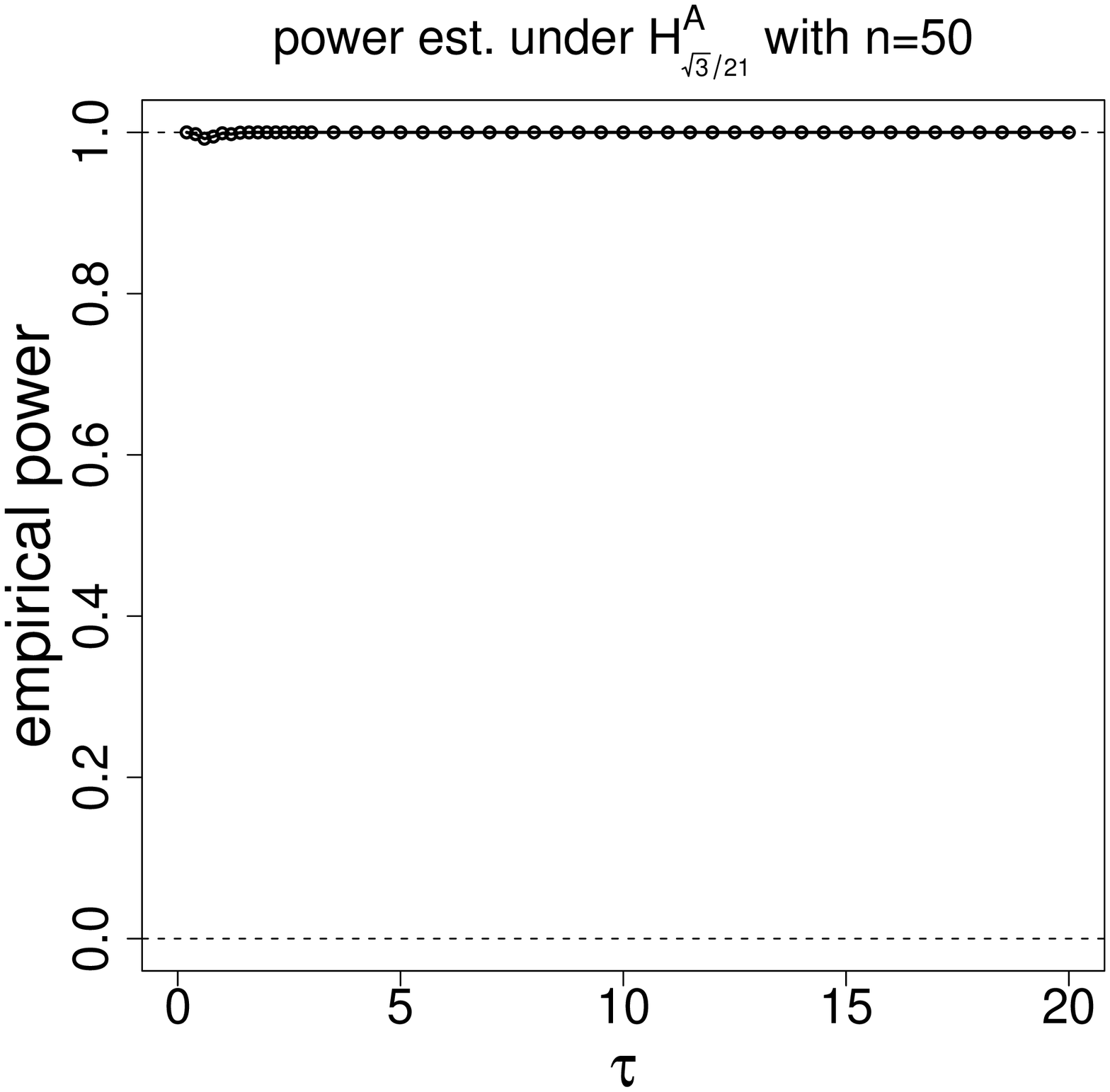}}}
\rotatebox{-90}{ \resizebox{1.7 in}{!}{ \includegraphics{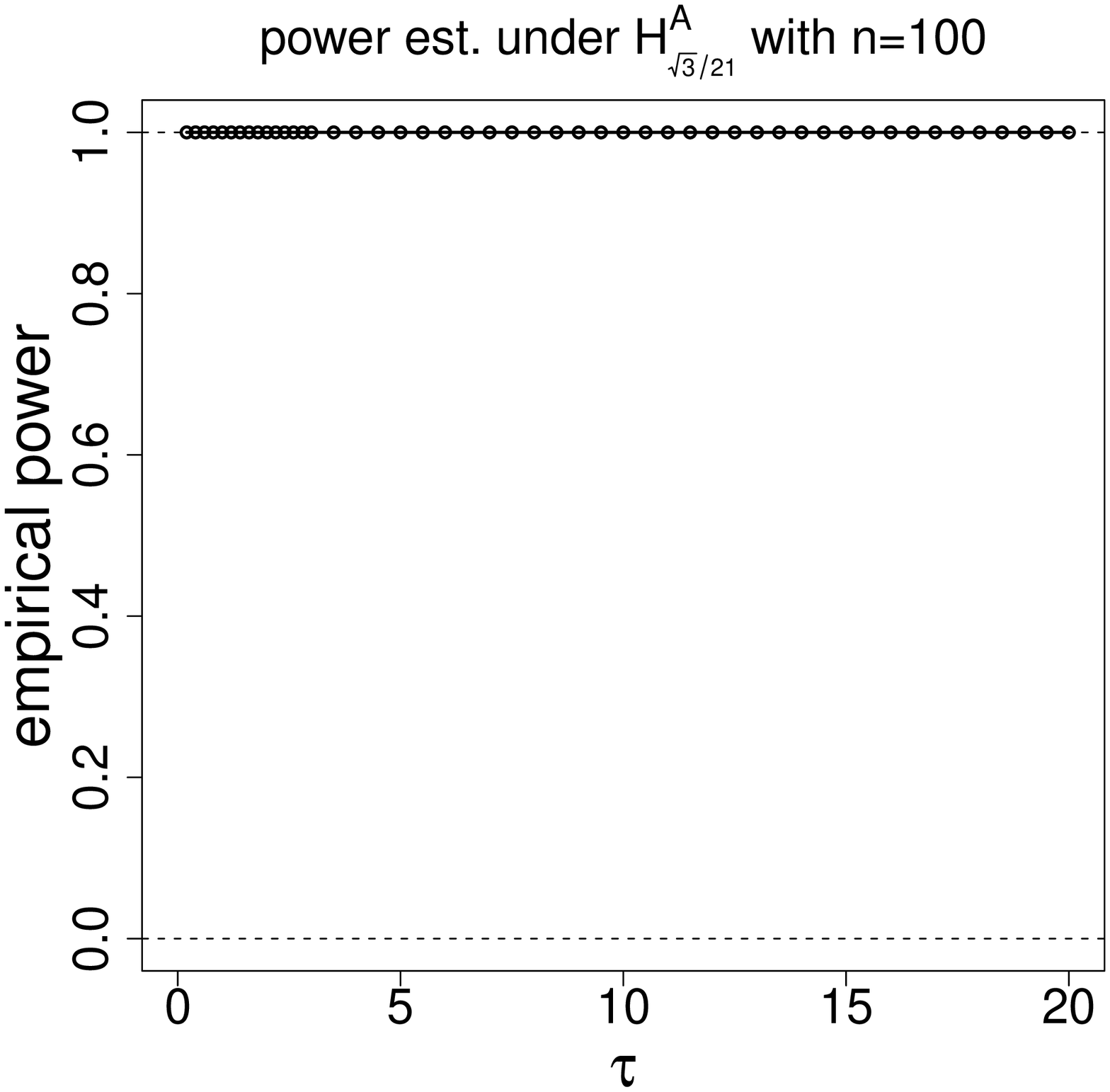}}}
\rotatebox{-90}{ \resizebox{1.7 in}{!}{ \includegraphics{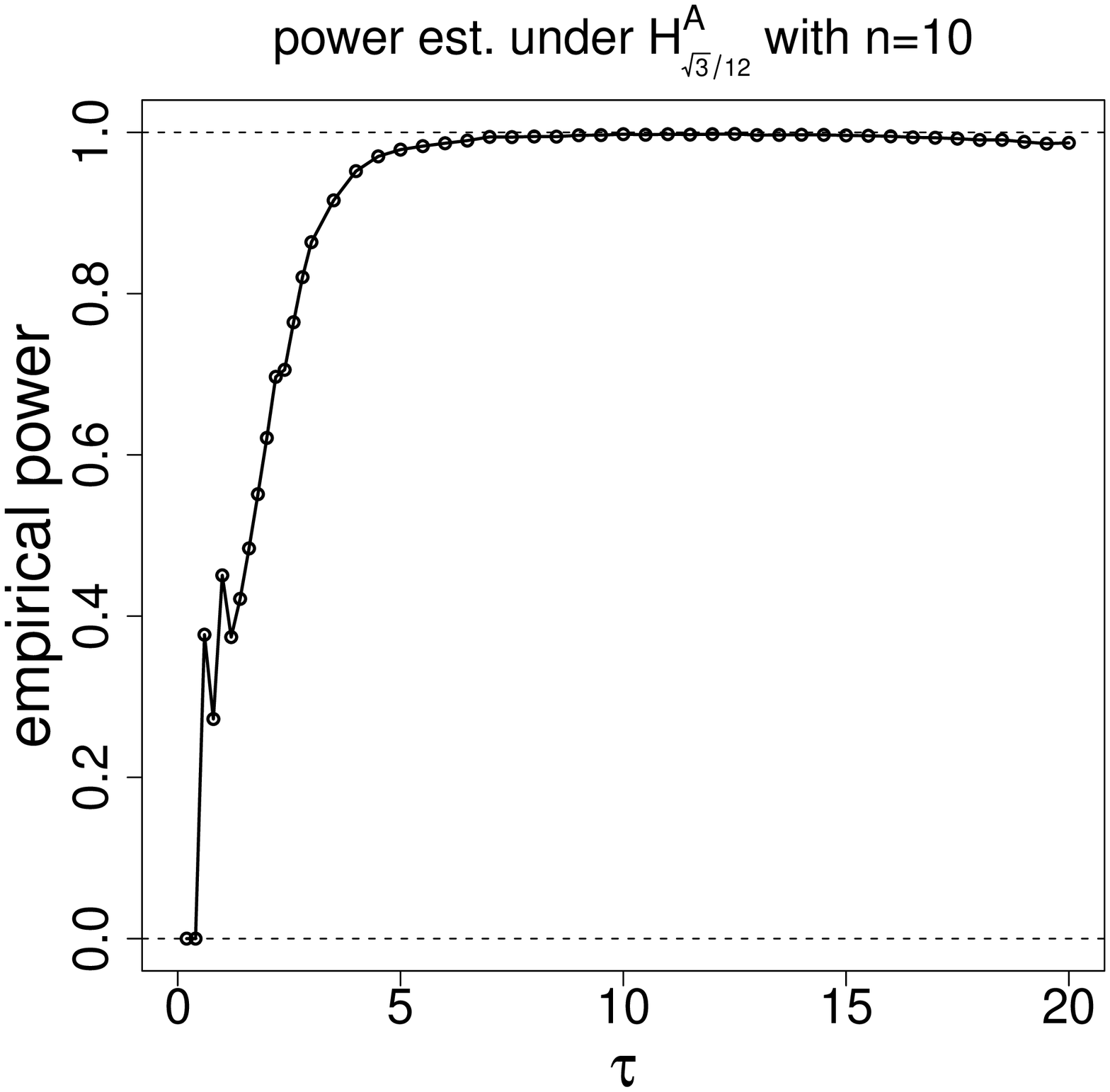}}}
\rotatebox{-90}{ \resizebox{1.7 in}{!}{ \includegraphics{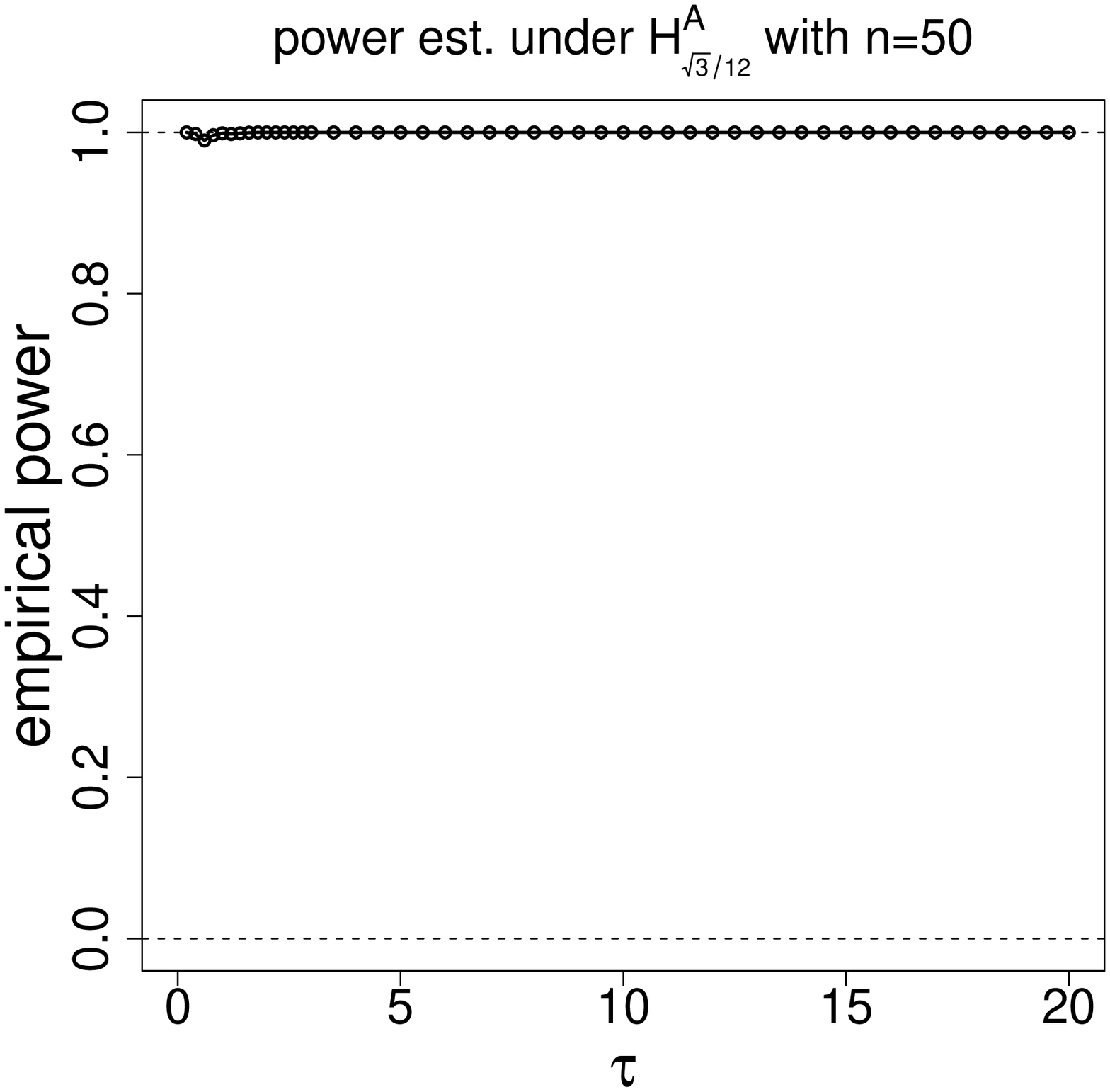}}}
\rotatebox{-90}{ \resizebox{1.7 in}{!}{ \includegraphics{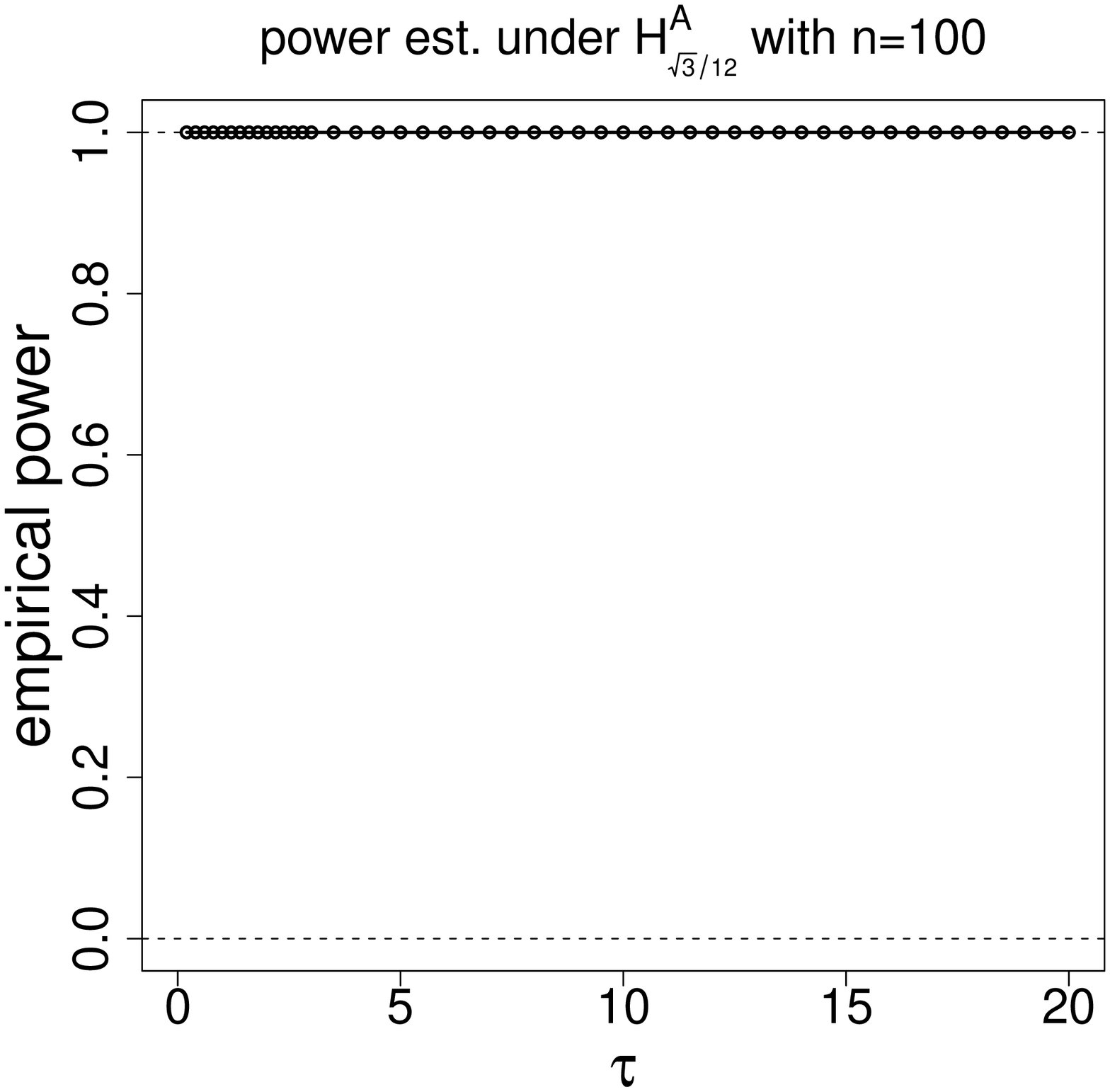}}}
\rotatebox{-90}{ \resizebox{1.7 in}{!}{ \includegraphics{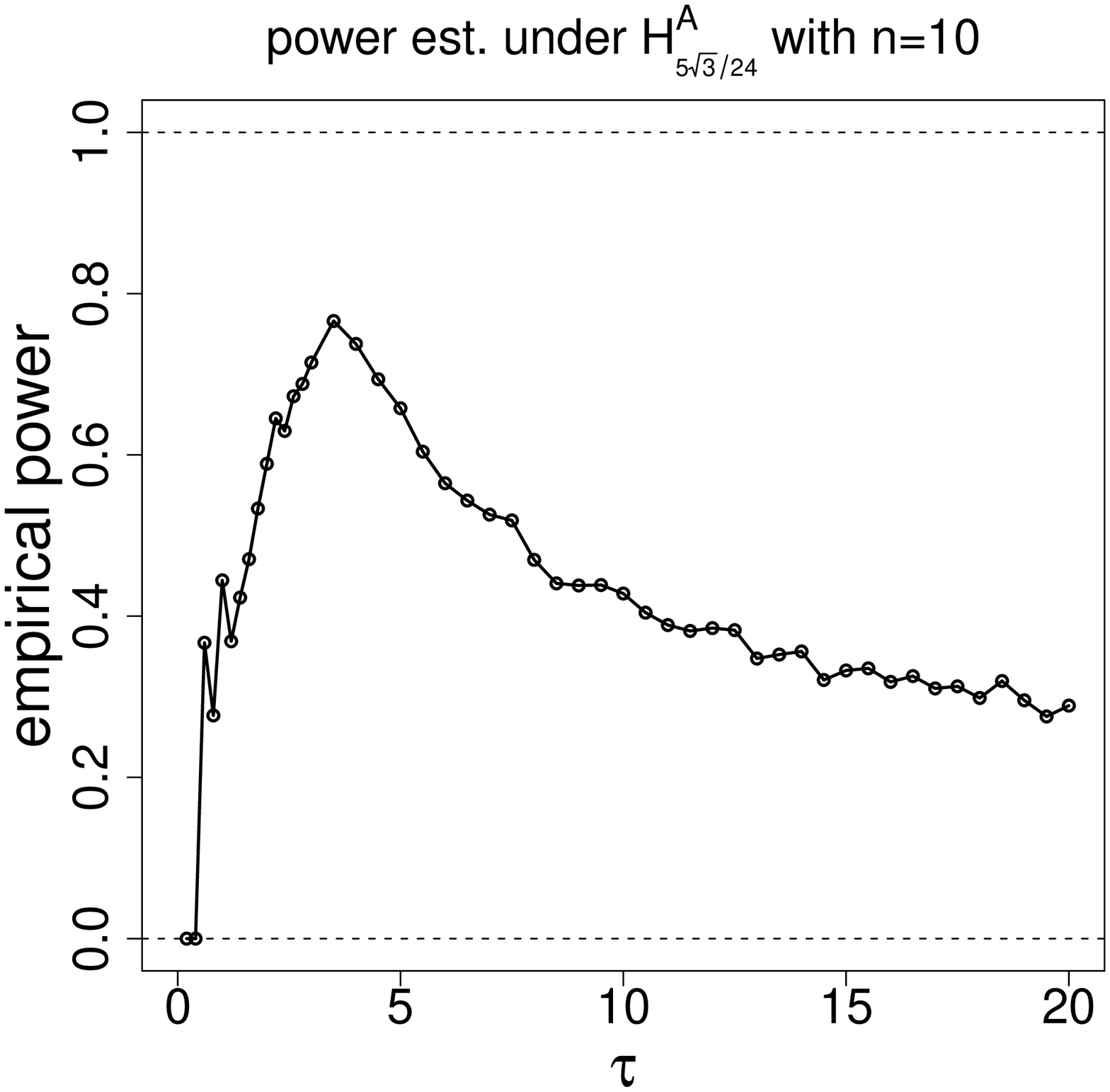}}}
\rotatebox{-90}{ \resizebox{1.7 in}{!}{ \includegraphics{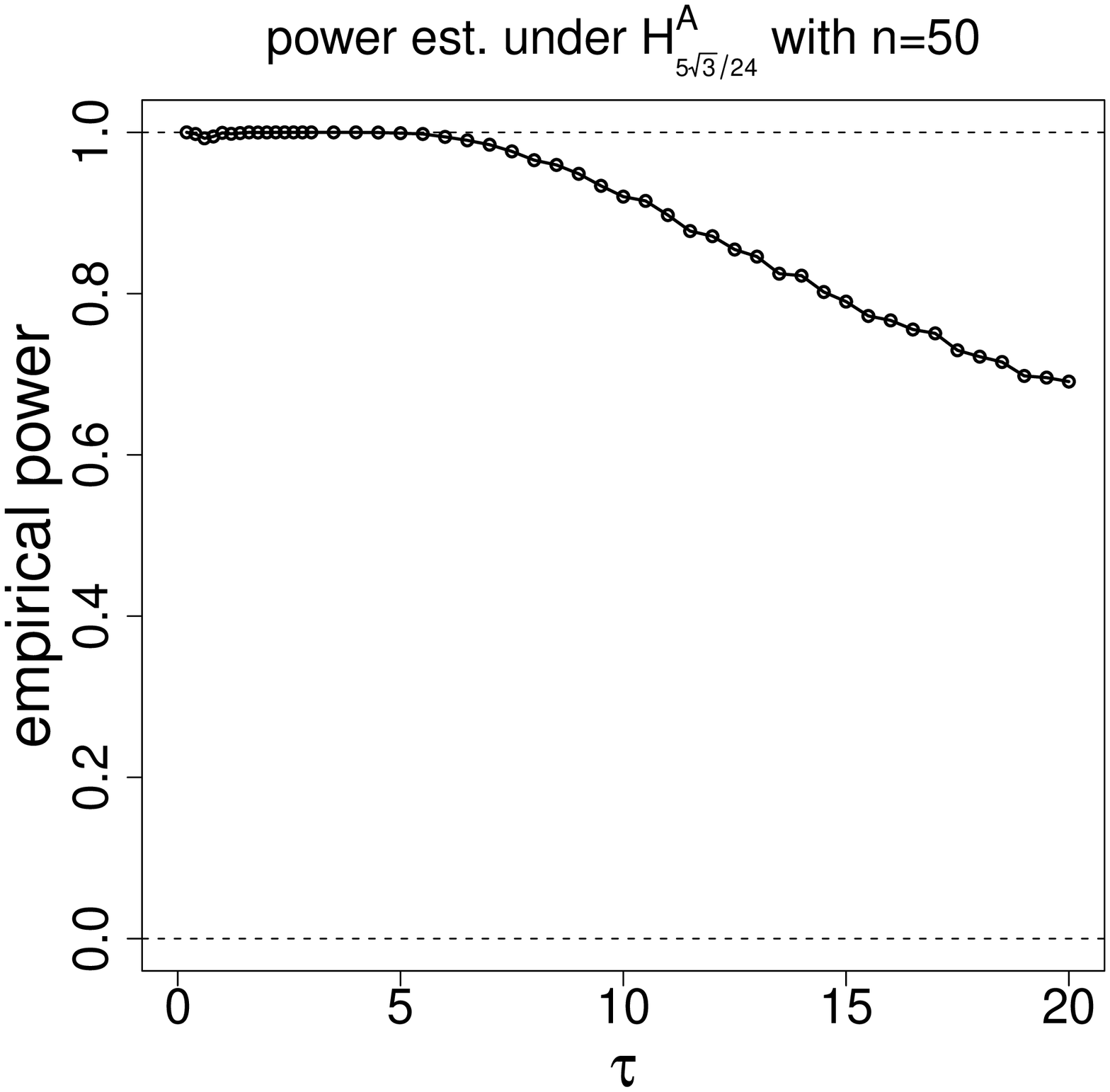}}}
\rotatebox{-90}{ \resizebox{1.7 in}{!}{ \includegraphics{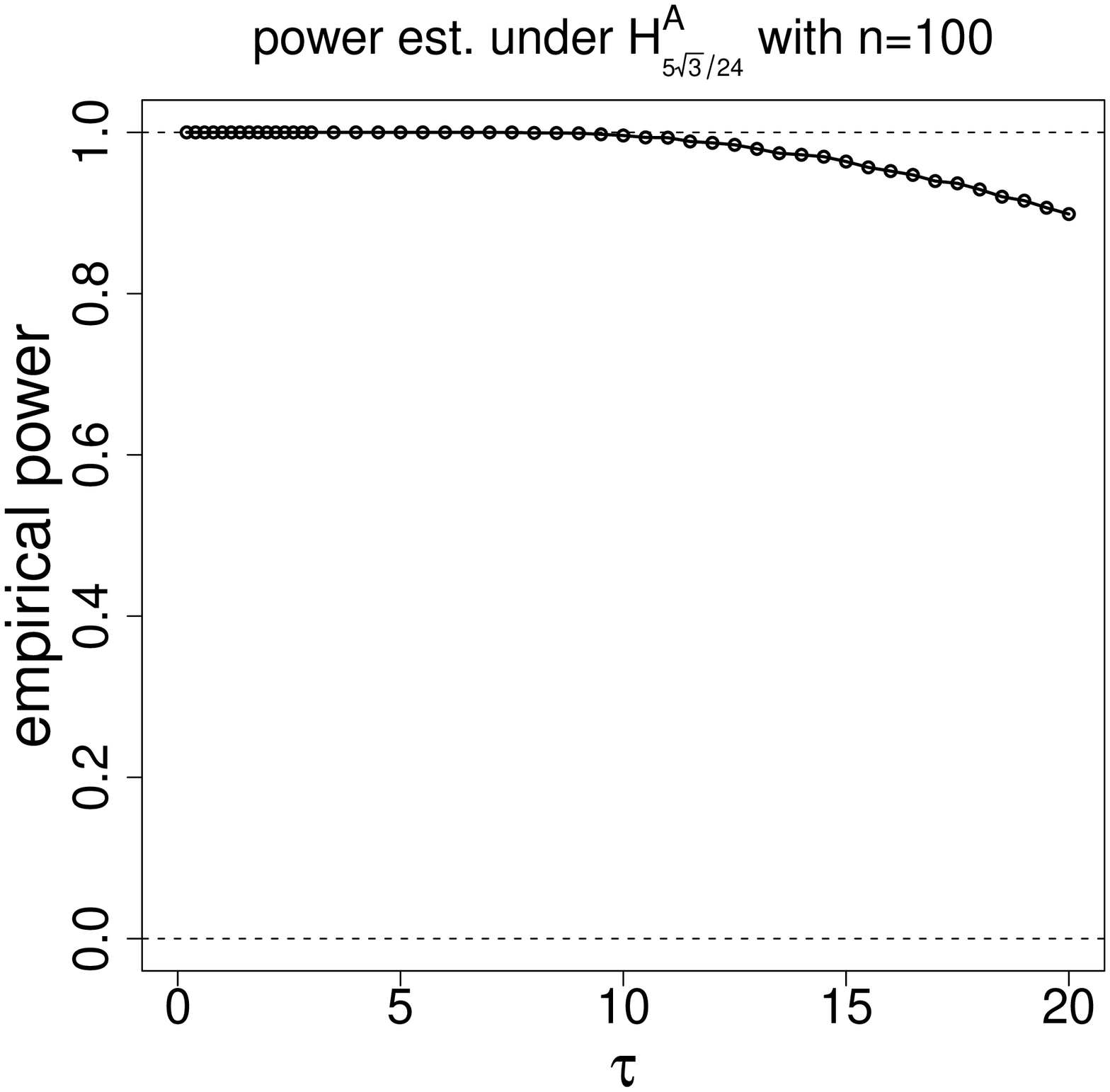}}}
\caption{
\label{fig:CS-emp-power-assoc}
\textbf{Empirical power for $R_{CS}(\tau)$ in the one triangle case:}
Monte Carlo power estimates for relative density of central similarity PCDs
in the one triangle case
using the asymptotic critical value against association alternatives
$H^A_{\sqrt{3}/21}$ (top row),
$H^A_{\sqrt{3}/12}$ (middle row),
 and
$H^A_{5\,\sqrt{3}/24}$ (bottom row)
as a function of $\tau$, for $n=10$ (left column), $n=50$ (middle column), and $n=100$ (right column).
}
\end{figure}

Under association,
we estimate the empirical power as in Section \ref{sec:PE-emp-power-assoc}.
In Figure \ref{fig:CS-emp-power-assoc},
we present Monte Carlo power estimates for relative density of central similarity PCDs
in the one triangle case
against $H^A_{5\,\sqrt{3}/24}$,
$H^A_{\sqrt{3}/12}$, and $H^A_{\sqrt{3}/21}$ as a function of $\tau$ for $n=10,50,100$.
Under mild association and small $n$,
highest power is attained around $\tau \approx 3$,
under mild association with large $n$,
power increases as $\tau$ increases.
For moderate to severe association and large $n$,
power is virtually one for all $\tau$ values considered.
Considering the empirical size performance,
we recommend $\tau \approx 5$,
as it has the desired level and high power.

\begin{figure}[]
\centering
\rotatebox{-90}{ \resizebox{1.7 in}{!}{ \includegraphics{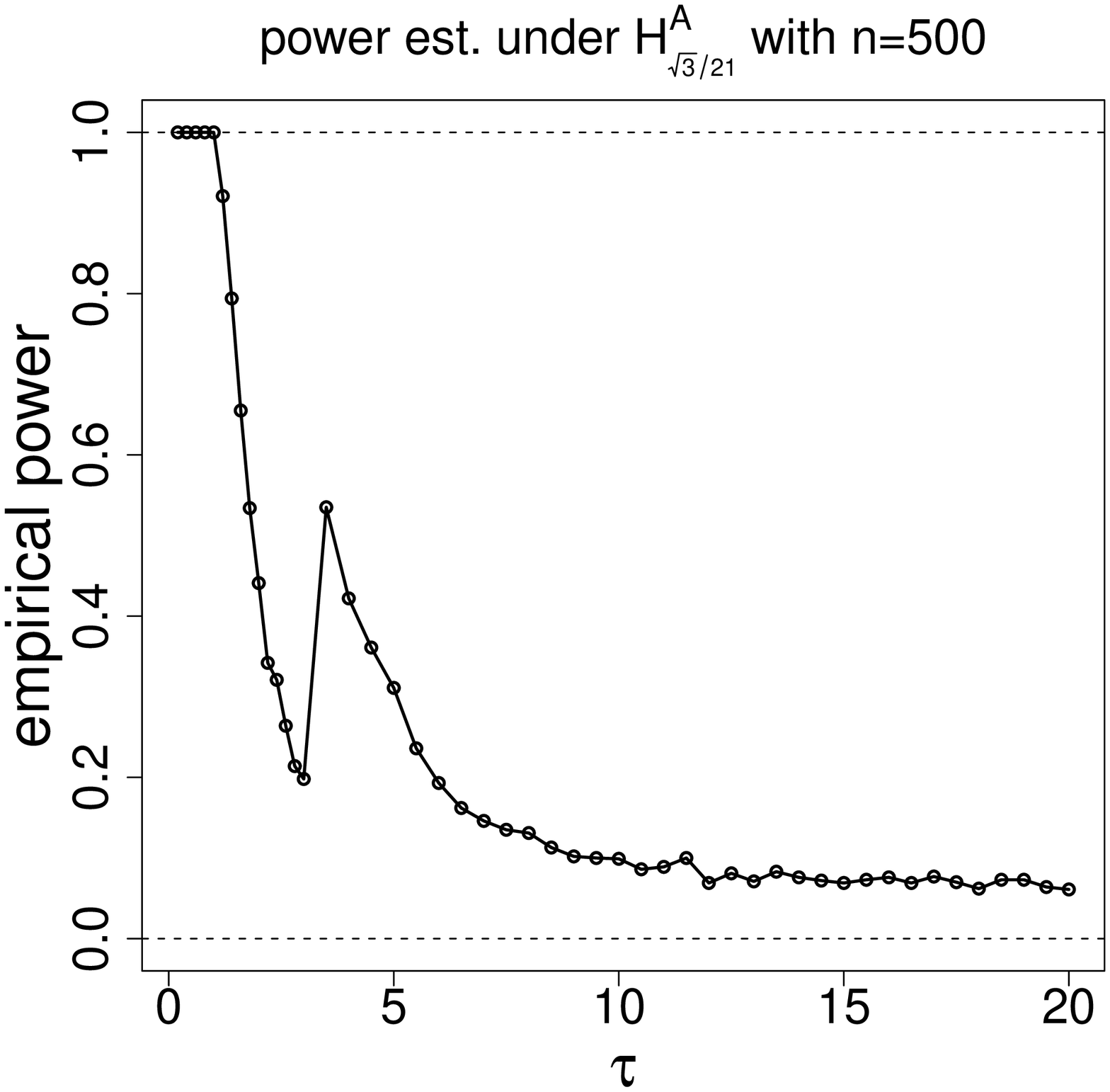}}}
\rotatebox{-90}{ \resizebox{1.7 in}{!}{ \includegraphics{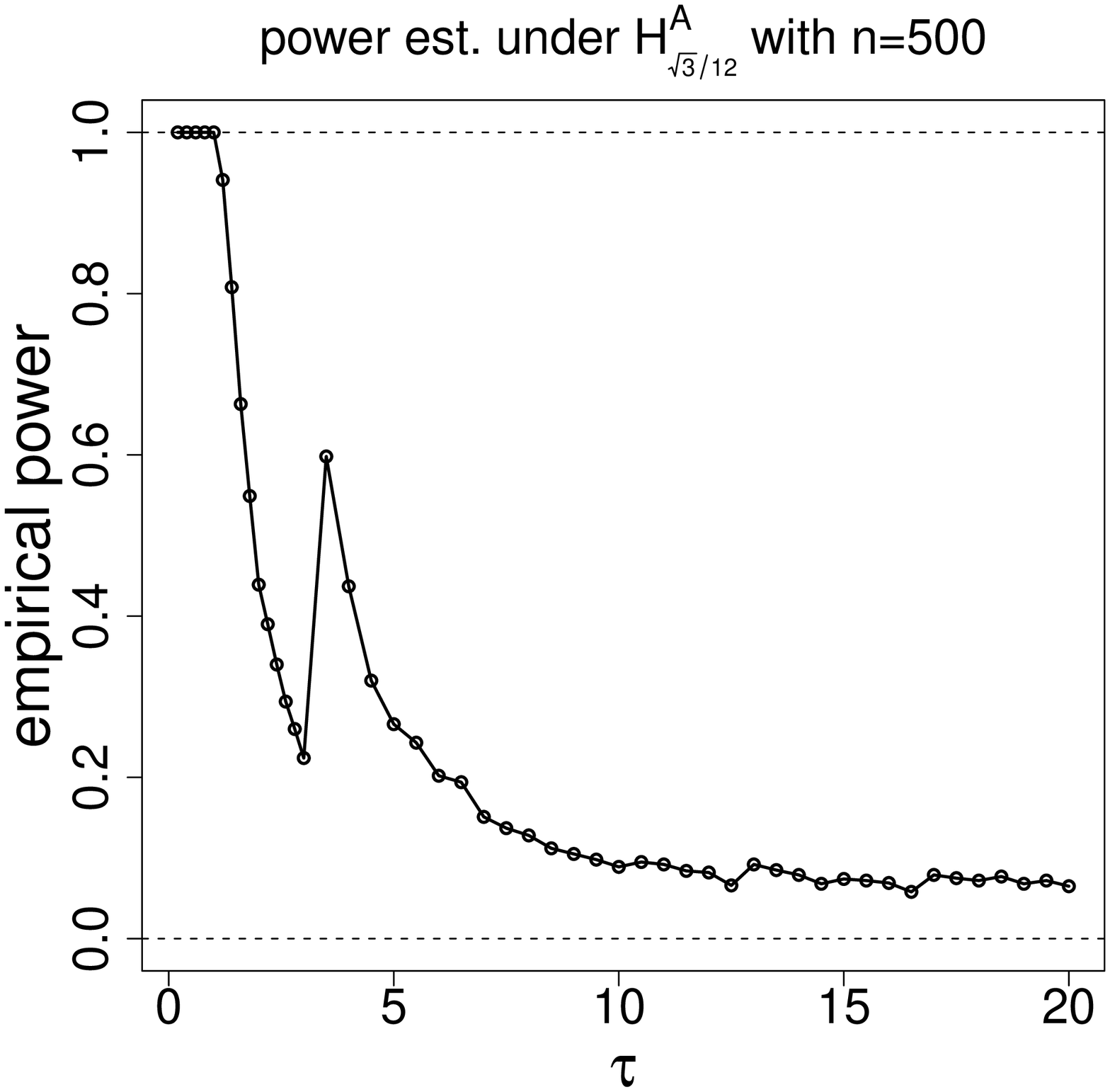}}}
\rotatebox{-90}{ \resizebox{1.7 in}{!}{ \includegraphics{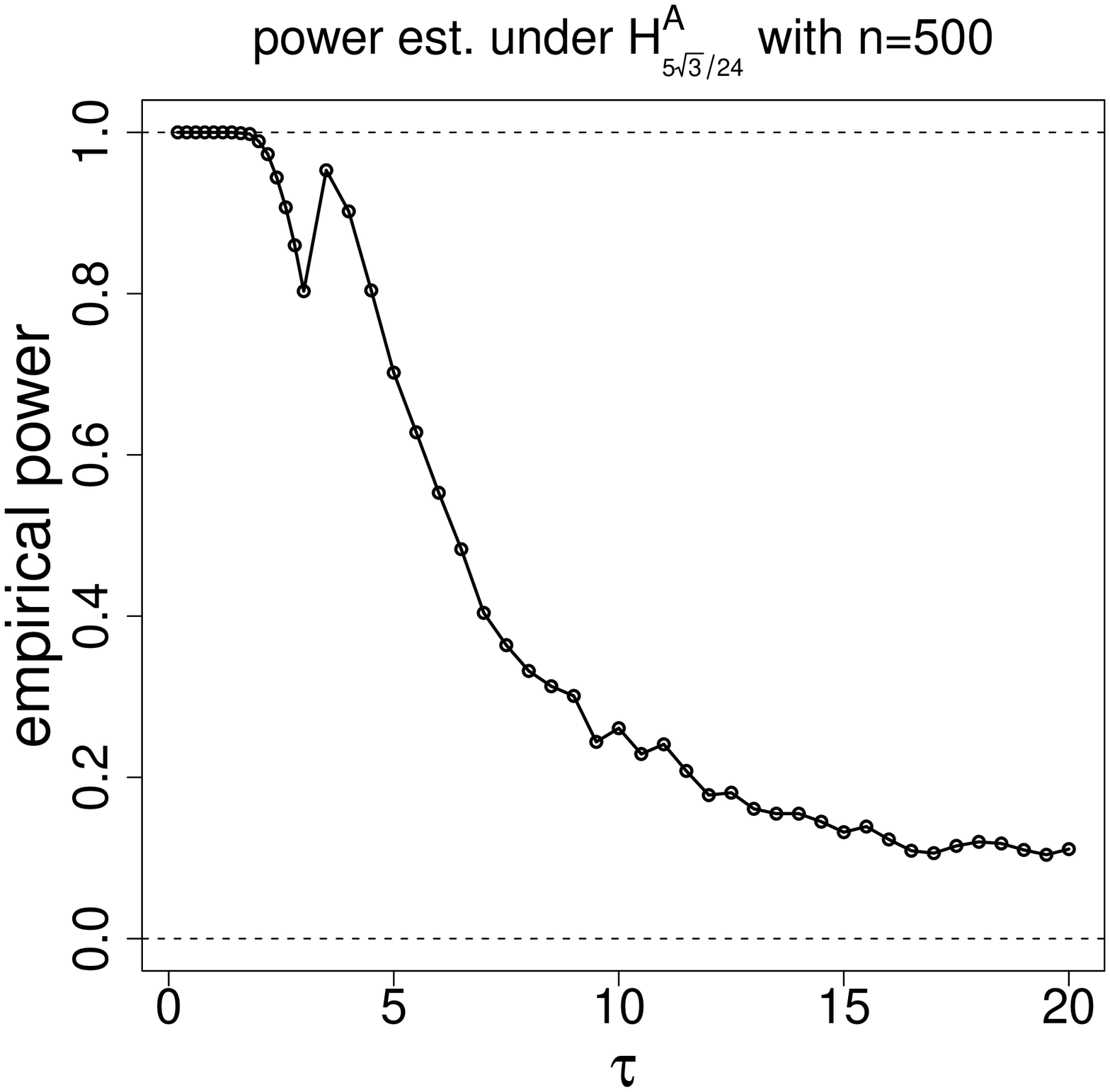}}}
\rotatebox{-90}{ \resizebox{1.7 in}{!}{ \includegraphics{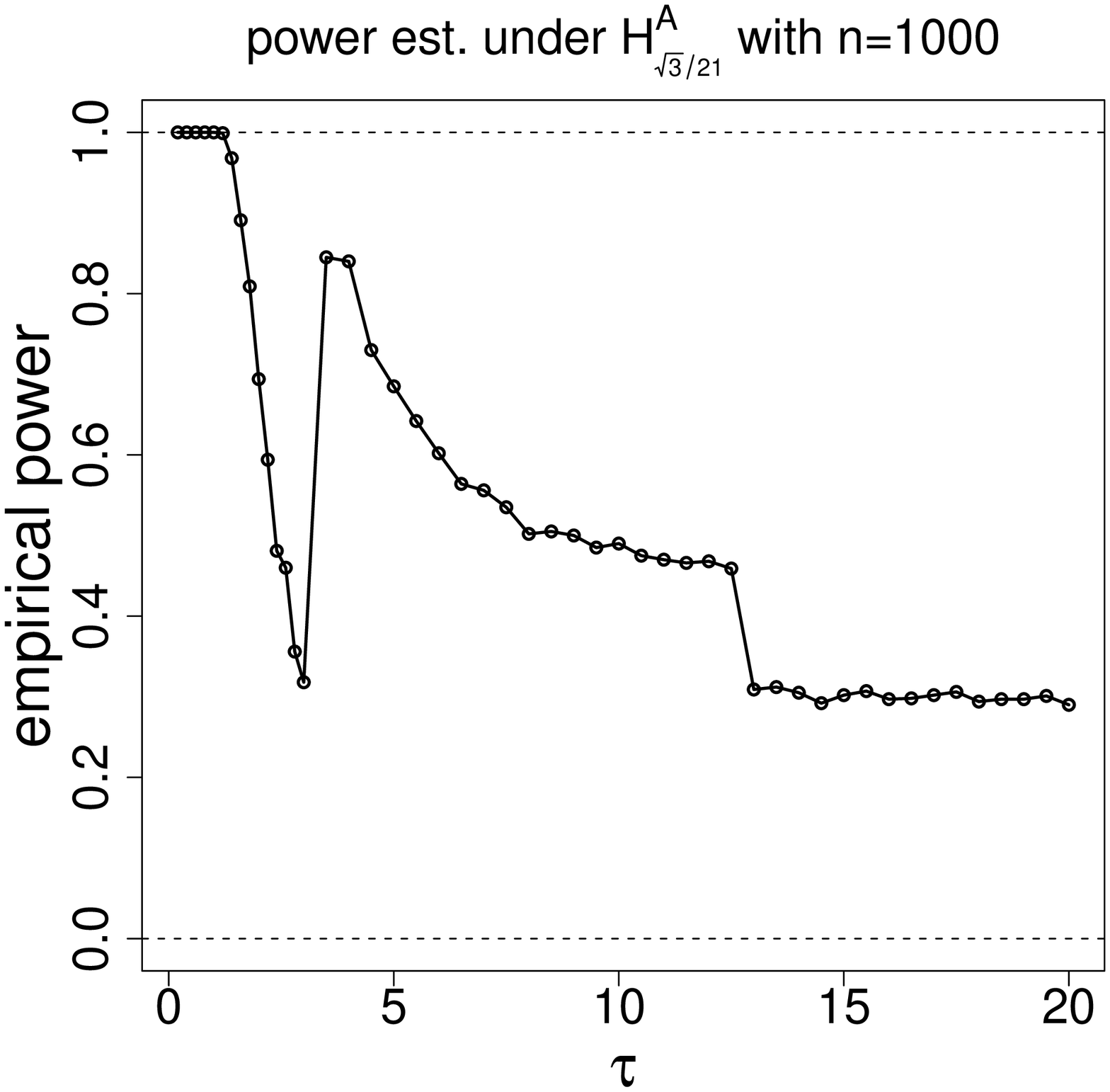}}}
\rotatebox{-90}{ \resizebox{1.7 in}{!}{ \includegraphics{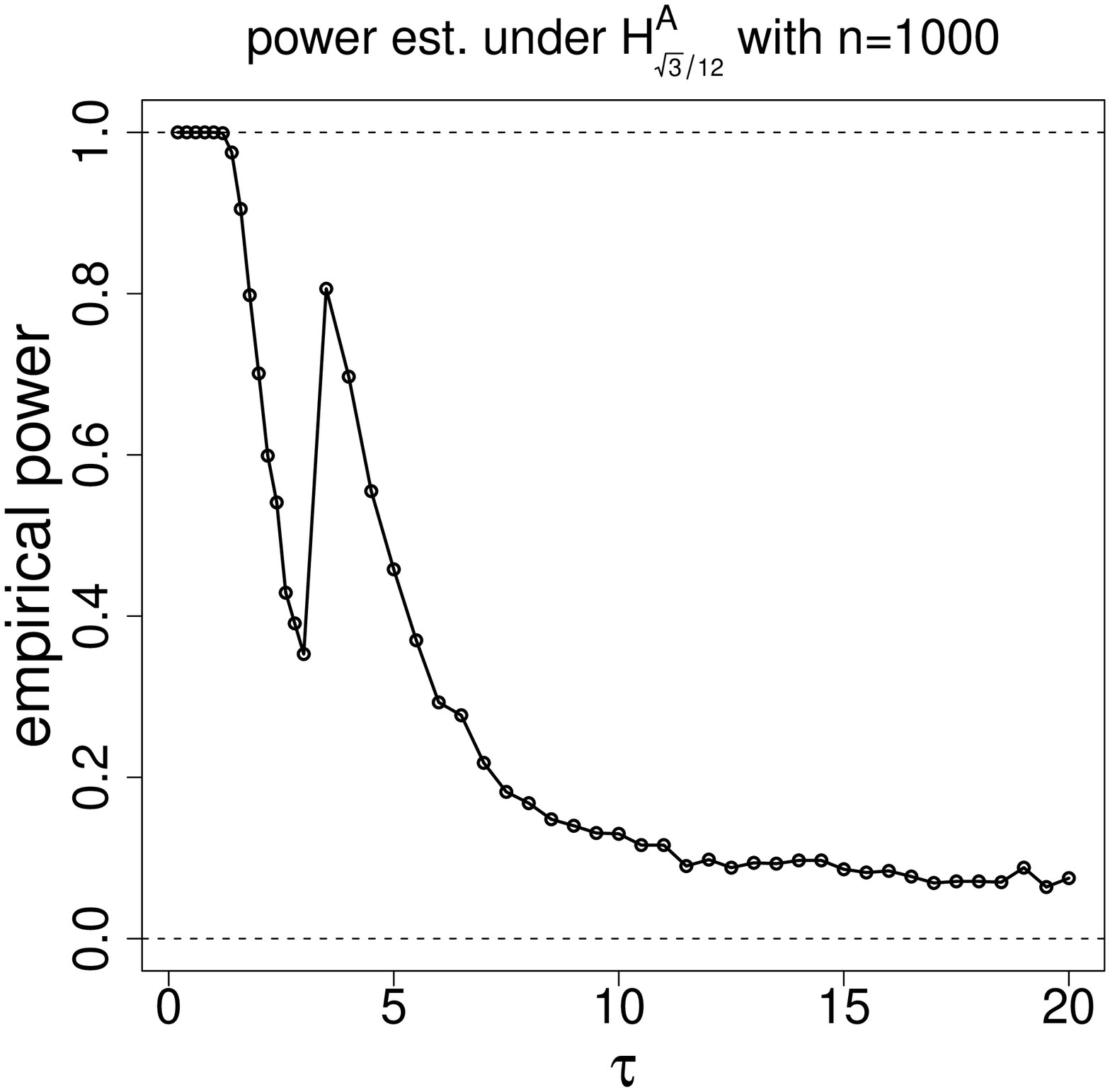}}}
\rotatebox{-90}{ \resizebox{1.7 in}{!}{ \includegraphics{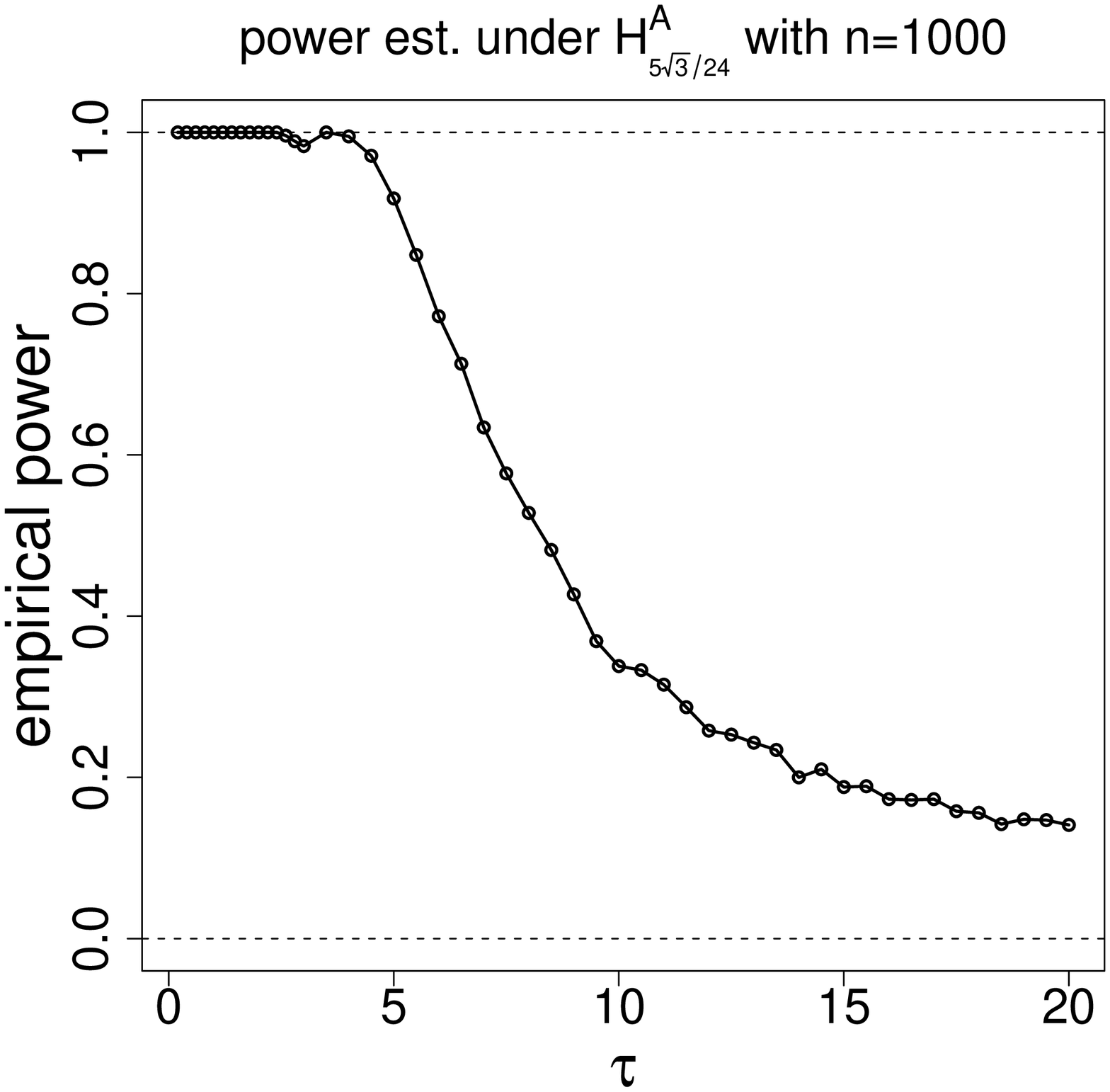}}}
\caption{\label{fig:MT-CS-emp-power-assoc}
\textbf{Empirical power for $R_{CS}(\tau)$ in the multiple triangle case:}
Monte Carlo power estimates of the relative density of central similarity PCDs
in the multiple triangle case
using the asymptotic critical value against association alternatives
$H^A_{\sqrt{3}/21}$ (left column),
$H^A_{\sqrt{3}/12}$ (middle column),
and
$H^A_{5\,\sqrt{3}/24}$ (right column)
as a function of $\tau$, for $n=500$ (top) and $n=1000$ (bottom).
}
\end{figure}

In the multiple triangle case,
we generate data as in Section \ref{sec:PE-emp-power-assoc}.
The corresponding empirical power estimates as a function of $\tau$
are presented in Figure \ref{fig:MT-CS-emp-power-assoc} for $n=500$ or 1000.
Observe that the Monte Carlo power estimate
tends to decrease as $\tau$ gets larger.
The empirical power is maximized for $\tau \le 1$.
Considering the empirical size estimates,
we recommend $\tau \approx 1$ for association,
since the corresponding test has the desired level with high power.

\begin{remark}
\textbf{Empirical Power Comparison for the Two PCD Families:}
In the one triangle case,
under the segregation alternatives,
the power estimates of the central similarity PCDs
tend to be higher than those of the proportional-edge PCDs.
Under mild to moderate association alternatives,
central similarity PCDs have higher power estimates,
while under severe association,
proportional-edge PCD has higher power estimates.
In the multiple triangle case,
under segregation,
central similarity PCDs has higher power estimates,
and
under association,
proportional-edge PCDs has higher power estimates.
$\square$
\end{remark}

\section{Pitman Asymptotic Efficiency}
\label{sec:PAE}

Suppose that the distribution $F $ under consideration may be
indexed by a set $\Theta \subset \R$ and consider
$H_o:\theta=\theta_0 $ versus $H_a: \theta>\theta_0 $.

Pitman asymptotic efficiency or efficacy (PAE)
provides for an investigation of ``local asymptotic power''
--- local around $H_o$.
This involves the limit as $n \rightarrow \infty$ as well as
the limit as $\ve \rightarrow 0$.

Consider the comparison of test sequences $S=\bigl\{ S_n \bigr\}$
satisfying the following conditions in a neighborhood $\theta \in
[\theta_0,\theta_0+\kappa] $ of the null parameter for some
$\kappa>0$.

\noindent
{\bfseries Pitman's Conditions:}
\begin{itemize}
\item[(PC1)] For some functions $\mu_n(\theta)$ and $\sigma_n(\theta)$,
the distribution $F_{\theta}$ of $\bigl[S_n-\mu_n(\theta)\bigr]/\sigma_n(\theta)$
converges to $Z \sim \N(0,1)$ uniformly on $\bigl[ \theta_0,\theta_0+\kappa \bigr] $, i.e.,
$$\sup_{~~~~~~~~~\theta_0\le \theta \le \theta_0+\kappa}\;\;
\sup_{t \in \R} \left| P \left( \frac{S_n-\mu_n(\theta)}{\sigma_n(\theta)}\le t \right)-\Phi(t) \right|\rightarrow 0
\text{ as } n \rightarrow \infty.$$
\item[(PC2)] For $\theta \in [\theta_0,\theta_0+\kappa] $,
$\mu_n(\theta)$ is differentiable with $\mu_n'(\theta_0)>0$,
\item[(PC3)] For $\theta_n=\theta_0 + O\left( n^{-1/2} \right)$,
$\lim_{n\rightarrow \infty}\frac{\mu_n'(\theta_n)}{\mu_n'(\theta_0)}=1$,
\item[(PC4)] For $\theta_n=\theta_0 + O\left( n^{-1/2} \right)$,
$\lim_{n\rightarrow \infty}\frac{\sigma_n(\theta_n)}{\sigma_n(\theta_0)}=1$.
\item[(PC5)] For some constant $c>0$,
$$\lim_{n\rightarrow \infty}\frac{\mu_n'(\theta_0)}{\sqrt{n}\,\sigma_n(\theta_0)}=c, $$
\end{itemize}
Condition (PC1) is equivalent to
\begin{itemize}
\item[(PC1)$^{\prime}$] For some functions $\mu_n(\theta)$ and $\sigma_n(\theta)$,
the distribution $F_{\theta}$ of $\bigl[S_n-\mu_n(\theta_n)\bigr]/\sigma_n(\theta_n)$
converges to a standard normal distribution (see \cite{eeden:1963}).
\end{itemize}
Note that if $\mu_n^{(k)}(\theta_0)>0$ and
$\mu_n^{(l)}(\theta_0)=0$, for all $l=1,\,2,\ldots,k-1$, then
$\mu_n^{\prime}(\theta_0)$ in (PC2), (PC3), and (PC5) can be
replaced by $\mu_n^{(k)}(\theta_0)>0$ and $\mu_n^{\prime}(\theta_n)$
in (PC3) can be replaced by $\mu_n^{(k)}(\theta_n)$
(see \cite{kendall:1979}).

\begin{lemma}
\textbf{(Pitman-Noether)}
\begin{itemize}
\item[(i)]
Let $S=\bigl\{ S_n \bigr\}$ satisfy conditions (PC1)-(PC5).
Consider testing $H_o $ by the critical regions $S_n > u_{\alpha_n}$
with $\alpha_n=P_{\theta_0}\bigl(S_n >u_{\alpha_n}\bigr)\rightarrow
\alpha $ as $n \rightarrow \infty $ where $\alpha \in (0,1)$.
For $\beta \in (0,1-\alpha)$ and $\theta_n=\theta_0 + O\left( n^{-1/2}\right)$,
we have $$\beta_n(\theta_n)=P_{\theta_n}\bigl( T_n >u_{\alpha_n} \bigr)\rightarrow \beta \text{ iff }
c\,\sqrt{n}\bigl( \theta_n-\theta \bigr)\rightarrow \Phi^{-1}(1-\alpha)-\Phi^{-1}(\beta).$$
\item[(ii)]
Let $S=\bigl\{ S_n \bigr\}$ and $Q=\bigl\{ Q_n \bigr\}$ each
satisfy conditions (PC1)-(PC5).
Then the asymptotic relative efficiency
of $S $ relative to $Q $ is given by $ARE(S,Q)=\left(c_S/c_Q
\right)^2$.
\end{itemize}
\end{lemma}

Thus, to evaluate $ARE(S,Q)$ under the conditions (PC1)-(PC5), we
need only calculate the quantities $c_S $ and $c_Q $, where
$$c_S=\lim_{n \rightarrow \infty}\frac{\mu'_{S_n}(\theta_0)}{\sqrt{n}\cdot \sigma_{S_n}(\theta_0)}
\text{ and }
c_Q=\lim_{n \rightarrow \infty}\frac{\mu'_{Q_n}(\theta_0)}{\sqrt{n}\cdot \sigma_{Q_n}(\theta_0)}$$
$PAE(S)=c_S^2$ is called the {\em Pitman Asymptotic Efficiency}
(PAE) of the test based on $S_n $. Using similar notation and
terminology for $Q_n $,
$$ARE(S,Q)=\frac{\PAE(S)}{\PAE(Q)}.$$

A detailed discussion of PAE can be found in \cite{kendall:1979} and \cite{eeden:1963}.

Under segregation or association alternatives,
the PAE of $\rho_{_{PE}}(n,r)$
is given by  $\PAE(r) = \frac{\left( \mu^{(k)}(r,\ve=0) \right)^2}{\nu_{_{PE}}(r)}$
where $k$ is the minimum order of the derivative with respect to $\ve$
for which $\mu^{(k)}(r,\ve=0) \not=0$.
That is, $\mu^{(k)}(r,\ve=0) \not=0$ but $\mu^{(l)}(r,\ve=0)=0$ for $l=1,2,\ldots,k-1$.
Similarly, the PAE of $\rho_{_{CS}}(n,\tau)$ is given by
$\PAE(\tau) = \frac{\left(\mu^{(k)}(\tau,\ve=0) \right)^2}{\nu_{_{CS}}(\tau)}$ where
$k$ is the minimum order of the derivative with respect to
$\ve$ for which $\mu^{(k)}(\tau,\ve=0) \not= 0$.

\subsection{Pitman Asymptotic Efficiency for Proportional-Edge PCDs under the Segregation Alternative}
\label{sec:PAE-PE-seg}
Consider the test sequences $\rho_{_{PE}}(r)=\bigl\{ \rho_{_{PE}}(n,r) \bigr\}$
under segregation alternatives for
sufficiently small $\ve>0$ and $r \in \bigl[ 1,\sqrt{3}/(4\,\ve)\bigr)$.
In the PAE framework above,
the parameters are $\theta=\ve$ and $\theta_0=0$.
Suppose $\mu^S_{_{PE}}(r,\ve)=\E^S_{\ve}[\rho_{_{PE}}(n,r)]$.
For $\ve \in \bigl[0,\sqrt{3}/8 \bigr)$,
$$\mu^S_{_{PE}}(r,\ve)=\sum_{j=1}^5 \varpi_{1,j}(r,\ve)\,\I(r \in \mI_j)$$
with the corresponding intervals
$\mI_1=\bigl[ 1,3/2-\sqrt{3}\,\ve \bigr)$,
$\mI_2=\bigl[3/2-\sqrt{3}\,\ve,3/2 \bigr)$,
$\mI_3=\bigl[3/2,2-4\,\ve/\sqrt{3} \bigr)$,
$\mI_4=\bigl[ 2-4\,\ve/\sqrt{3},2 \bigr)$,
$\mI_5=\bigl[ 2,\sqrt{3}/(2\,\ve) \bigr)$.
See \cite{ceyhan:TR-rel-dens-NPE}
for the explicit form of $\mu^S_{_{PE}}(r,\ve)$ and for derivation.
Notice that as $\ve \rightarrow 0$,
only $\mI_1=\bigl[1,3/2-\sqrt{3}\,\ve \bigr)$,
$\mI_3=\bigl[ 3/2,2-4\,\ve/\sqrt{3}\bigr)$,
$\mI_5=\bigl[ 2,\sqrt{3}/(2\,\ve)\bigr)$ do not vanish,
so we only keep the components of $\mu^S_{_{PE}}(r,\ve)$ on these intervals.

Furthermore,
$\sigma_S^2(n,\ve)=\Var^S_{\ve}(\rho_{_{PE}}(n,r))=\frac{1}{2\,n\,(n-1)}\Var^S_{\ve}[h_{12}]+\frac{(n-2)}{n\,(n-1)}\nu^S_{_{PE}}(r,\ve)$,
with
$\nu^S_{_{PE}}(r,\ve)=\Cov^S_{\ve}[h_{12},h_{13}].$
The explicit forms of
$\Var^S_{\ve}[h_{12}]$ and $\Cov^S_{\ve}[h_{12},h_{13}]$  are not calculated,
since we only need $\lim_{n\rightarrow \infty}\sigma_n^2(\ve=0)=\nu_{_{PE}}(r)$
which is given in Equation \eqref{eqn:PEAsyvar}.

Notice that $\E^S_{\ve}|h_{12}|^3 \le 8 < \infty $ and
$\E^S_{\ve}[h_{12}\,h_{13}]-\E^S_{\ve}[h_{12}]^2=\Cov^S_{\ve}[h_{12},h_{13}]>0$
then by \cite{callaert:1978}
$$\sup\text{}_{t\in \R} \left| P_\ve \left( \sqrt{n}\frac{\bigl( \rho_{_{PE}}(n,r)-\mu_S(r,\ve) \bigr)}
{\sqrt{\nu^S_{_{PE}}(r,\ve)}}\le t \right)-\Phi(t)\right| \le C\,\E^S_{\ve}\left| h_{12} \right|^3\,
\left[ \nu^S_{_{PE}}(r,\ve) \right]^{-\frac{3}{2}}\,n^{-\frac{1}{2}}$$
where $C $ is an absolute constant and $\Phi(\cdot)$ is the standard
normal distribution function.
Then (PC1) follows for each $r \in \bigl[ 1,\sqrt{3}/(2\,\ve) \bigr)$
and
$\ve \in \bigl[0,\sqrt{3}/4 \bigr)$.

Differentiating $\mu^S_{_{PE}}(r,\ve)$ with respect to $\ve$ yields
\begin{multline*}
(\mu^S_{_{PE}})^{\prime}(r,\ve)= \varpi_{1,1}^{\prime}(r,\ve)\,\I\left(r \in \bigl[ 1,3/2-\sqrt{3}\,\ve \bigr)\right)+
\varpi_{1,3}^{\prime}(r,\ve)\,\I\left( r \in [3/2,2-4\,\ve/\sqrt{3}) \right)\\
+\varpi_{1,5}^{\prime}(r,\ve)\,\I\left( r \in \bigl[2,\sqrt{3}/(2\,\ve)\bigr) \right)
\end{multline*}
where
\begin{align*}
\varpi_{1,1}^{\prime}(r,\ve)&=
\frac{2\,\ve\,(144\,\ve^2\,(r^2-1)+36-37\,r^2)}{27\,(2\,\ve-1)^3(2\,\ve+1)^3},\\
\varpi_{1,3}^{\prime}(r,\ve)&=
\Bigl[2\,\sqrt{3}\Bigl( (2\,r-3)\,64\,\ve^3+(7\,r^2+r^4-24\,r+20)\,16\,\sqrt{3}\ve^2+
(r-3)\,48\,\ve+3\,\sqrt{3}\,r^4+96\,\sqrt{3}\,r\\
&-36\,\sqrt{3}-60\,\sqrt{3}\,r^2\Bigr)\ve\Bigr]\Big/\Bigl[9\,(2\,\ve+1)^3(2\,\ve-1)^3r^2 \Bigr],\\
\varpi_{1,5}^{\prime}(r,\ve)&=
\frac{8\,\sqrt{3}\,\ve\,\left(48\,\ve^3+(3\,r^4+3\,r^2-20)\,4\,\sqrt{3}\,\ve^2+36\,\ve+
9\,\sqrt{3}-9\,\sqrt{3}\,r^2 \right)}{27\,r^2(2\,\ve+1)^3(2\,\ve-1)^3}.
\end{align*}

Since $(\mu^S_{_{PE}})^{\prime}(r,\ve=0)=0$,
we need higher order derivatives for (PC2).
A detailed discussion is available in (\cite{kendall:1979}).

Differentiating $(\mu^S_{_{PE}})^{\prime}(r,\ve)$ with respect to $\ve$ yields
\begin{multline*}
(\mu^S_{_{PE}})^{\prime\prime}(r,\ve)=
\varpi_{1,1}^{\prime\prime}(r,\ve)\,\I\left(r \in \bigl[ 1,3/2-\sqrt{3}\,\ve \bigr)\right)+
\varpi_{1,3}^{\prime\prime}(r,\ve)\,\I\left( r \in [3/2,2-4\,\ve/\sqrt{3}) \right)\\
+\varpi_{1,5}^{\prime\prime}(r,\ve)\,\I\left( r \in \bigl[2,\sqrt{3}/(2\,\ve)\bigr) \right)
\end{multline*}
where
\begin{align*}
\varpi_{1,1}^{\prime\prime}(r,\ve)&=-\frac{2\,(r^2-1)\,1728\,\ve^4+(72-77\,r^2)\,4\,\ve^2+36-37\,r^2}{27\,(4\,\ve^2-1)^4},\\
\varpi_{1,3}^{\prime\prime}(r,\ve)&= -2\,\Bigl[(2\,r-3)\,512\,\sqrt{3}\,\ve^5+
(20+r^4+7\,r^2-24\,r)\,576\,\ve^4+(2\,r-3)\,1024\,\sqrt{3}\,\ve^3+(20-108\,r^2\\
&+96\,r+9\,r^4)\,36\,\ve^2+(-3+2\,r)\,96\,\sqrt{3}\,\ve-108+9\,r^4-180\,r^2+288\,r\Bigr]\Big/
\Bigl[9\,r^2(2\,\ve+1)^4(2\,\ve-1)^4\,\Bigr],\\
\varpi_{1,5}^{\prime\prime}(r,\ve)&= -8\,\Bigl[128\,\sqrt{3}\,\ve^5+(-20+3\,r^4+3\,r^2)\,48\,\ve^4+256\,\sqrt{3}\,\ve^3+(-5-12\,r^2+3\,r^4)\,12\,\ve^2+24\,\ve\,\sqrt{3}+9\\
&-9\,r^2\Bigr]\Big/\Bigl[9\,r^2(2\,\ve+1)^4(2\,\ve-1)^4\,\Bigr].
\end{align*}
Thus,
\begin{equation}
\label{eqn:PAE-S-''}
 (\mu^S_{_{PE}})^{\prime\prime}(r,\ve=0)=
\begin{cases}
          -\frac{8}{3}+{\frac{74}{27}}\,r^2 &\text{for} \quad r \in [1,3/2),\\
          -2\,{\frac{(r^2-4\,r+2)(r^2+4\,r-6)}{r^2}} &\text{for} \quad r \in [3/2,2),\\
          -\frac{8\,(1-r^2)}{r^2} &\text{for} \quad r \in [2,\sqrt{3}/(2\,\ve)).
\end{cases}
\end{equation}
Observe that $(\mu^S_{_{PE}})^{\prime\prime}(r,\ve=0)>0$ for all
$r \in \bigl[1,\sqrt{3}/(2\,\ve)\bigr)$,
so (PC2) holds with the second derivative.
(PC3) in the second derivative form follows from
continuity of $(\mu^S_{_{PE}})^{\prime\prime}(r,\ve)$ in $\ve$ and
(PC4) follows from continuity of $\sigma_n^2(r,\ve)$ in $\ve$.

Next,
we find $c^S_{_{PE}}(r)=\lim_{n\rightarrow
\infty}\frac{(\mu^S_{_{PE}})^{\prime\prime}(r,\ve=0)}{\sqrt{n}\,\sigma_n(r,\ve=0)}=
\frac{(\mu^S_{_{PE}})^{\prime\prime}(r,\ve=0)}{\sqrt{\nu_{_{PE}}(r)}}$,
where numerator is given in Equation \eqref{eqn:PAE-S-''} and
denominator is given in Equation \eqref{eqn:PEAsyvar}.
We can easily see that $c^S_{_{PE}}(r)>0$,
since $c^S_{_{PE}}(r)$ is increasing in $r$ and $c^S_{_{PE}}(r=1)>0$.
Then (PC5) follows.
So under segregation
alternatives $H^S_{\ve}$, the PAE of $\rho_{_{PE}}(n,r)$ is given by
$$
\PAE_{PE}^S(r) =\left(c^S_{_{PE}}(r)\right)^2=
   \frac{\left( (\mu^S_{_{PE}})^{\prime\prime}(r,\ve=0) \right)^2}{\nu_{_{PE}}(r)}.
$$

\begin{figure}[h]
\centering
\psfrag{mu}[c]
{\begin{picture}(0,0)
\put(-10.5,-6.5){\makebox(0,0)[l]{\Huge PAE score}}
\end{picture}}
\psfrag{t}{{\Huge expansion parameter}}
\rotatebox{-90}{ \resizebox{2. in}{!}{\includegraphics{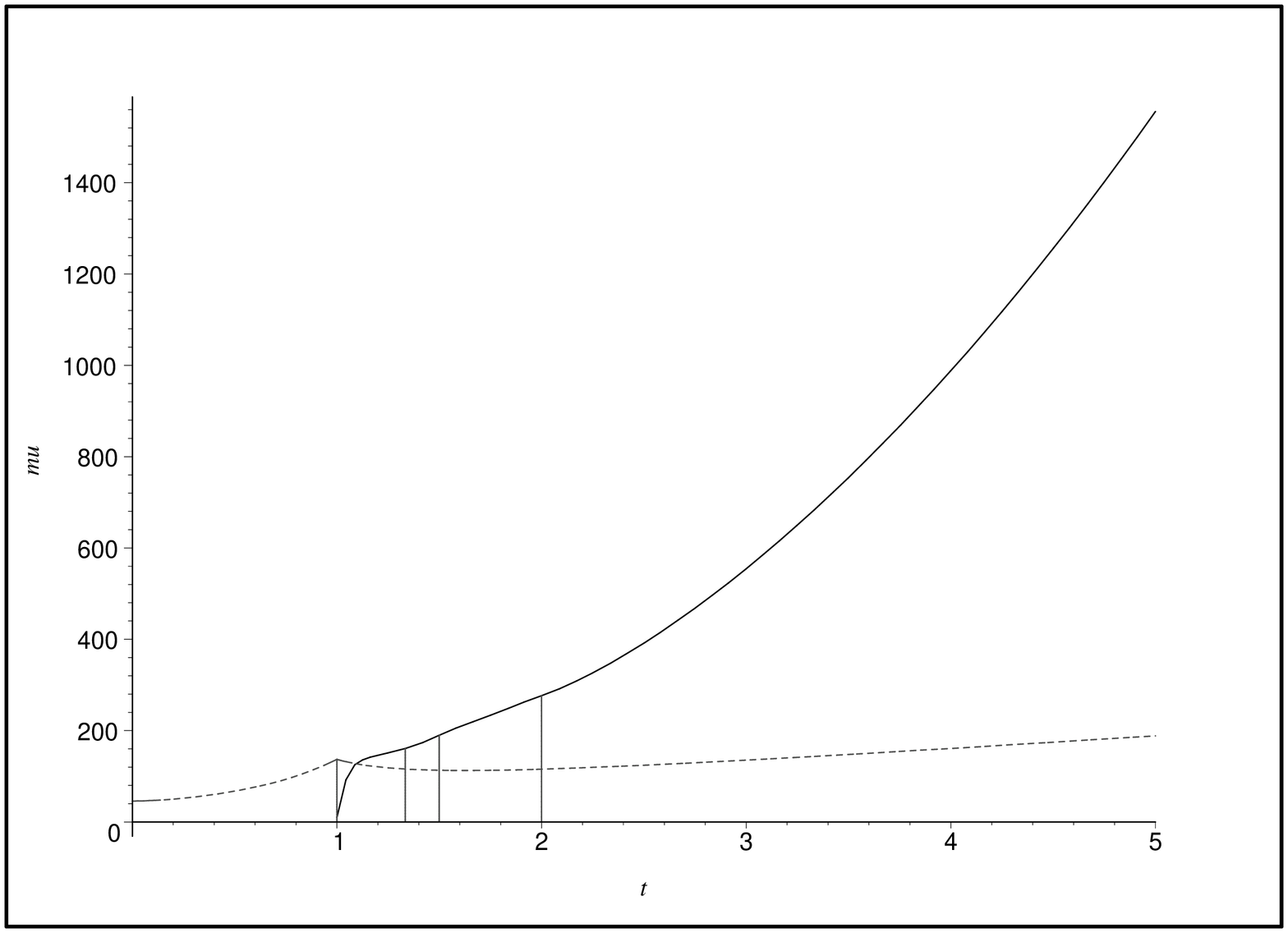}}}
\psfrag{mu}[c]
{\begin{picture}(0,0)
\put(-10.5,-8.5){\makebox(0,0)[l]{\Huge PAE score}}
\end{picture}}
\rotatebox{-90}{ \resizebox{2. in}{!}{\includegraphics{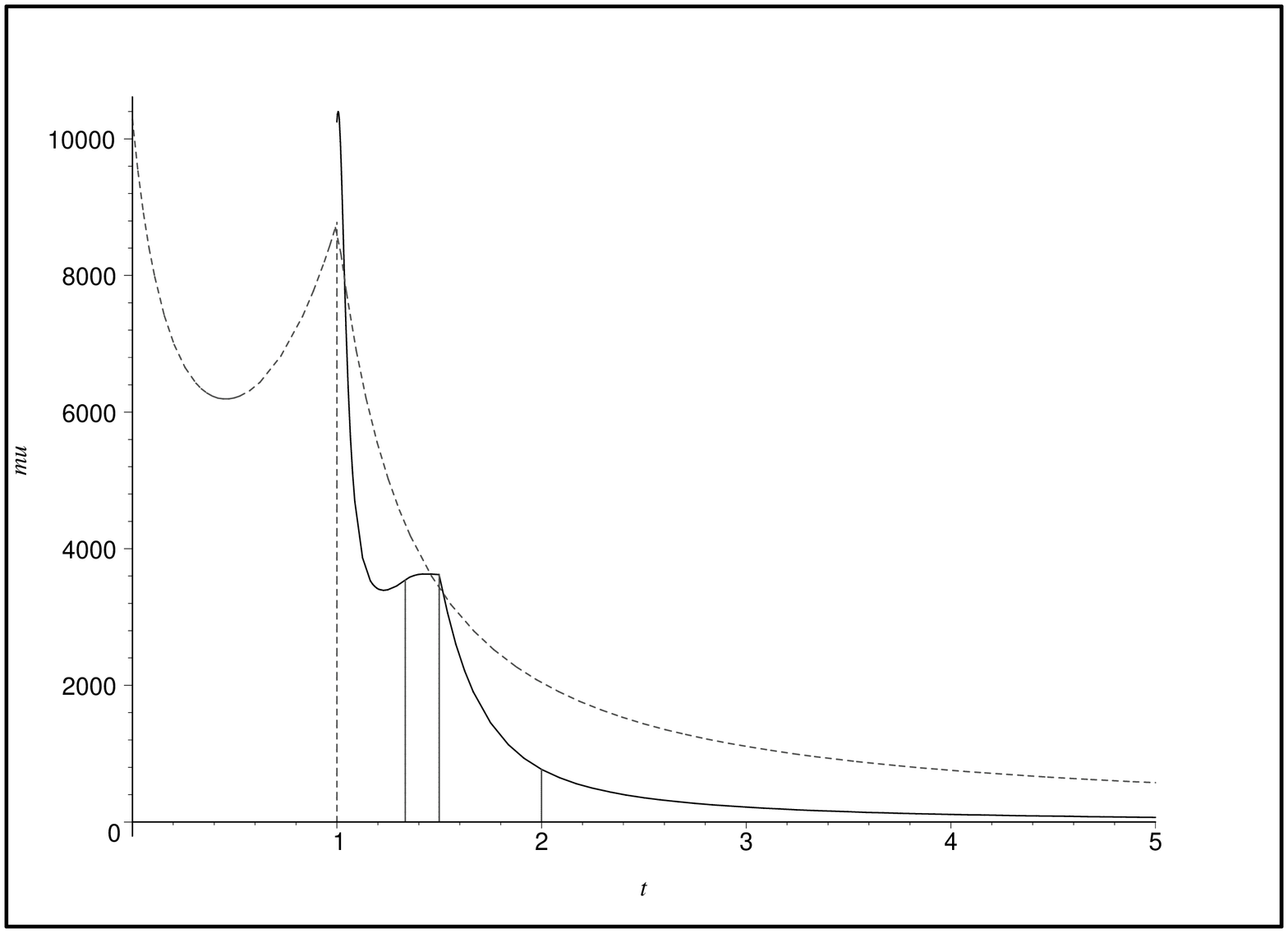}}}
\caption{
\label{fig:PAE-seg-assoc}
Pitman asymptotic efficiency against
segregation (left) and association (right) alternatives
as a function of the expansion parameters
in the one triangle case
for the relative density of proportional-edge PCDs (solid line) and
central similarity PCDs (dashed line).
}
\end{figure}

In Figure \ref{fig:PAE-seg-assoc} (left),
we present the PAE as a function of the expansion parameter for segregation.
Notice that $\PAE_{PE}^S(r=1) = 160/7 \approx 22.8571$,
$\lim_{r \rightarrow \infty} \PAE_{PE}^S(r) = \infty$.
Based on the PAE analysis, we suggest,
for large $n$ and small $\ve$,
choosing $r$ large for testing against segregation.
However, for small and moderate values of $n$,
normal approximation is not appropriate
due to the skewness in the density of $\rho_{_{PE}}(n,r)$.
Therefore, for small $n$, we suggest moderate $r$ values.

\subsection{Pitman Asymptotic Efficiency for Central Similarity PCDs under the Segregation Alternative}
\label{sec:PAE-CS-seg}
Consider the test sequences
$\rho_{_{CS}}(\tau)=\bigl\{ \rho_{_{CS}}(n,\tau)\bigr\}$
for sufficiently small $\ve>0$ and $\tau \in (0,\infty)$.
In the PAE framework above,
the parameters are $\theta=\ve$ and $\theta_0=0$.
Suppose,
$\mu^S_{_{CS}}(\tau,\ve)=\E^S_{\ve}[\rho_{_{CS}}(n,\tau)]$.
For $\ve \in \bigl[ 0,\sqrt{3}/5 \bigr)$,
$$\mu^S_{_{CS}}(\tau,\ve)=\sum_{j=1}^3 \varpi_{1,j}(\tau,\ve)\,\I(\tau \in \mI_j)$$
with the corresponding intervals
$\mI_1=\bigl(0,1-\sqrt{3}\,\ve\bigr)$,
$\mI_2=\bigl[1-\sqrt{3}\,\ve,1\bigr)$,
and
$\mI_3=[1,\infty)$.
See Appendix 2 for the derivation of $\mu^S_{_{CS}}(\tau,\ve)$ for $\tau > 1$
and
Appendix 3 for the explicit form of $\mu^S_{_{CS}}(\tau,\ve)$ for $\tau \in (0,1]$.
Notice that as $\ve \rightarrow 0$,
only $\mI_1$ and $\mI_3$ do not vanish,
so we only keep the components of $\mu^S_{_{CS}}(\tau,\ve)$
on these intervals.

Furthermore,
$\sigma_n^2(\ve)=\Var^S_{\ve}(\rho_{_{CS}}(n,\tau))=\frac{1}{2\,n\,(n-1)}\Var^S_{\ve}[h_{12}]+\frac{(n-2)}{n\,(n-1)}\Cov^S_{\ve}[h_{12},h_{13}].$
The explicit forms of $\Var^S_{\ve}[h_{12}]$ and
$\Cov^S_{\ve}[h_{12},h_{13}]$ are not calculated, since we
only need $\lim_{n\rightarrow \infty}\sigma_n^2(\ve=0)=\nu_{_{CS}}(\tau)$ which is given
in Equation \eqref{eqn:CSAsyvar}.

Notice that $\E^S_{\ve}|h_{12}|^3 \le 8 < \infty $ and
$\E^S_{\ve}[h_{12}\,h_{13}]-\E^S_{\ve}[h_{12}]^2=\Cov^S_{\ve}[h_{12},h_{13}]>0$
then (PC1) follows for each $\tau \in (0,\infty)$ and $\ve \in \bigl[ 0,\sqrt{3}/3 \bigr)$.

Differentiating $\mu^S_{_{CS}}(\tau,\ve)$ with respect to $\ve$ yields
$$
(\mu^S_{_{CS}})^{\prime}(\tau,\ve)= \varpi_{1,1}^{\prime}(\tau,\ve)\,\I\left(\tau \in \bigl(0,1-\sqrt{3}\,\ve \bigr)\right)+
\varpi_{1,3}^{\prime}(\tau,\ve)\,\I\left( \tau \in [1,\infty) \right)
$$
where
$$
\varpi_{1,1}^{\prime}(\tau,\ve)=\frac{8\,\ve\,\tau^2(5\,\ve^2\tau-9\,\ve^2-3\,\tau+3)}{9(1-2\,\ve)^3(2\,\ve+1)^3(1-\tau)}
$$
and
\begin{multline*}
\varpi_{1,3}^{\prime}(\tau,\ve)=
{\frac{8}{9}}\,\Bigl[ \bigl( 2\,\tau^5\ve^2+21\,\tau^4\ve^2+
116\,\tau^3\ve^2+48\,\sqrt{3}\tau^
{2}\ve^3+37\,\tau^2\ve^2+96\,\sqrt{3}\tau\,{
\ve}^3-18\,\tau^3+36\,\sqrt{3}\tau^2\ve-200\,\tau
\,\ve^2+48\,\sqrt{3}\ve^3-\\
45\,\tau^2+72\,\sqrt
{3}\tau\,\ve-132\,\ve^2-36\,\tau+36\,\sqrt{3}\ve-9 \bigr) \tau\,\ve\Bigr]\Big/\Bigl[ \left( 2\,\tau+1 \right)  \left( \tau+2
 \right)  \left( \tau+1 \right)^2 \left( 2\,\ve+1 \right)^3
 \left( 2\,\ve-1 \right)^3\Bigr].
\end{multline*}
hence $(\mu^S_{_{CS}})^{\prime}(\tau,\ve=0)=0$, so we need higher
order derivatives for (PC2).
Differentiating $(\mu^A_{_{CS}})^{\prime}(\tau,\ve)$ with respect to $\ve$,
we get
$$
(\mu^S_{_{CS}})^{\prime\prime}(\tau,\ve)= \varpi_{1,1}^{\prime\prime}(\tau,\ve)\,\I\left(\tau \in \bigl(0,1-\sqrt{3}\,\ve \bigr)\right)+
\varpi_{1,3}^{\prime\prime}(\tau,\ve)\,\I\left( \tau \in [1,\infty) \right)
$$
where
$$
\varpi_{1,1}^{\prime\prime}(\tau,\ve)={\frac{-8\,\tau^2 \left( 20\,\tau\,\ve^4-36\,\ve^4-15\,\tau\,\ve^2+11\,\ve^2-\tau+1
 \right) }{ 3 \left( \tau-1 \right)  \left( 2\,\ve-1 \right)^4 \left( 2\,\ve+1 \right)^4}},
 $$
\begin{multline*}
\varpi_{1,3}^{\prime\prime}(\tau,\ve)= -\frac{8}{3}\,\Bigl[ \bigl( 8\,\tau^5\ve^4+84\,\tau^4\ve^4+
2\,\tau^5\ve^2+464\,\tau^3\ve^4+128\,\sqrt{3}\tau^2\ve^5+21\,\tau^4\ve^2+
148\,\tau^2\ve^4+256\,\sqrt{3}\tau\,\ve^5-4\,\tau^3\ve^2+\\
256\,\sqrt{3}\tau^2\ve^3-800\,
\tau\,\ve^4+128\,\sqrt{3}\ve^5-263\,\tau^2{
\ve}^2+512\,\sqrt{3}\tau\,\ve^3-528\,\ve^4-6
\,\tau^3+24\,\sqrt{3}\tau^2\ve-440\,\tau\,\ve^2
+256\,\sqrt{3}\ve^3-15\,\tau^2+\\
48\,\sqrt{3}\tau\,
\ve-192\,\ve^2-12\,\tau+24\,\sqrt{3}\ve-3 \bigr)
\tau\Bigr]\Big/ \Bigl[ \left( 2\,\tau+1 \right)  \left( \tau+2 \right)  \left( \tau+1
 \right)^2 \left( 2\,\ve+1 \right)^4 \left( 2\,\ve-1
 \right)^4\Bigr].
\end{multline*}
Hence
\begin{equation}
\label{eqn:CS-PAE-S-''}
(\mu^S_{_{CS}})^{\prime\prime}(\tau,\ve=0)=
\begin{cases}
          \frac{8\,\tau^2}{3} &\text{for} \quad \tau \in (0,1),\\
          \frac{8\,\tau}{2+\tau} &\text{for} \quad \tau \in [1,\infty).
\end{cases}
\end{equation}

Observe that $(\mu^S_{_{CS}})^{\prime\prime}(\tau,\ve=0)>0$ for all $\tau \in (0,\infty)$,
so (PC2) holds with the second derivative.
(PC3) in the second derivative form follows from continuity of
$(\mu^S_{_{CS}})^{\prime\prime}(\tau,\ve)$ in $\ve$ and (PC4)
follows from continuity of $\sigma_n^2(\tau,\ve)$ in $\ve$.

Next, we find
$c^S_{_{CS}}(\tau)=
\lim_{n\rightarrow \infty}\frac{(\mu^S_{_{CS}})^{\prime\prime}(\tau,\ve=0)}
{\sqrt{n}\,\sigma_n(\tau,\ve=0)}=\frac{(\mu^S_{_{CS}})^{\prime\prime}(\tau,\ve=0)}{\sqrt{\nu_{_{CS}}(\tau)}}$,
where numerator is given in Equation \eqref{eqn:CS-PAE-S-''} and
denominator is given in Equation \eqref{eqn:CSAsyvar}.
We can easily see that $c^S_{_{CS}}(\tau)>0$,
then (PC5) follows.

So under segregation alternatives $H^S_{\ve}$, the PAE of
$\rho_{_{CS}}(n,\tau)$ is given by
$$
\PAE_{CS}^S(\tau) =\left(c^S_{_{CS}}(\tau)\right)^2=
   \frac{\left( (\mu^S_{_{CS}})^{\prime\prime}(\tau,\ve=0) \right)^2}{\nu_{_{CS}}(\tau)}.
$$

In Figure \ref{fig:PAE-seg-assoc} (left),
we present the PAE as a function of $\tau$ for segregation.
Notice that $\lim_{\tau \rightarrow 0}\PAE_{CS}^S(\tau) = 320/7 \approx 45.7143$,
$\argsup_{\tau \in (0,\infty)}\PAE_{CS}^S(\tau)=1.0$ with
$\PAE_{CS}^S(\tau=1) = 960/7 \approx 137.1429$
and
$\lim_{\tau \rightarrow \infty} \PAE_{CS}^S(\tau) = \infty$.
Moreover a local maximum occurs at $\tau=1$
and a local minimum occurs at $\tau \approx 1.62$ with
PAE score $\approx 112.70$.
Based on the PAE analysis, we suggest, for large $n$ and small $\ve$,
choosing $\tau$ large for testing against segregation.
However, for small and moderate values of $n$,
normal approximation is not appropriate
due to the skewness in the density of $\rho_{_{CS}}(n,\tau)$
for extreme values of $\tau$.
Therefore, for small $n$, we suggest moderate $\tau$ values
(i.e., $\tau \approx 7$ or 8).

Comparing the PAE scores
of the relative density of proportional-edge PCDs
and
central similarity PCDs
under segregation alternatives,
we see that
$\PAE_{PE}^S(t) < \PAE_{CS}^S(t)$
for $ 1 \le t \lesssim 1.093)$;
and
$\PAE_{PE}^S(t) > \PAE_{CS}^S(t)$
for $t \gtrsim 1.093$.
Therefore, under segregation alternative,
overall, relative density of proportional edge PCD
is asymptotically more efficient compared to the central similarity PCD.
Furthermore,
$\PAE_{PE}^S(t)$ tends to $\infty$
as $t \rightarrow \infty$ at rate $O(t^2)$,
while
$\PAE_{CS}^S(t)$ tends to $\infty$
as $t \rightarrow \infty$ at rate $O(t)$.

\subsection{Pitman Asymptotic Efficiency for Proportional-Edge PCDs under the Association Alternative}
\label{sec:PAE-PE-assoc}
Consider the test sequences $\rho_{_{PE}}(r)=\bigl\{ \rho_{_{PE}}(n,r) \bigr\}$ for
sufficiently small $\ve >0$ and $r \in [1,\infty)$.
In the PAE framework above, the parameters are $\theta=\ve$ and $\theta_0=0$.
Suppose,
$\mu^A_{_{PE}}(r,\ve)=\E_\ve[\rho_{_{PE}}(n,r)]$.
For $\ve \in \bigl[0,\left(7\,\sqrt{3}-3\,\sqrt{15}\right)/12 \approx .042 \bigr)$,
$$\mu^A_{_{PE}}(r,\ve)=\sum_{j=1}^6 \varpi_{1,j}(r,\ve)\,\I(r \in \mI_j)$$
with the corresponding intervals
$\mI_1=\bigl[1,\left( 1+2\,\sqrt{3}\,\ve \right)/\left( 1-\sqrt{3}\,\ve \right)\bigr)$,
$\mI_2=\bigl[\left( 1+2\,\sqrt{3}\,\ve \right)/\left( 1-\sqrt{3}\,\ve \right),\\
4\,\left( 1-\sqrt{3}\,\ve \right)/3 \bigr)$,
$\mI_3=\bigl[4\,\left( 1-\sqrt{3}\,\ve \right)/3,4\,\left( 1+2\,\sqrt{3}\,\ve \right)/3 \bigr)$,
$\mI_4=\bigl[ 4\,\left( 1+2\,\sqrt{3}\,\ve \right)/3,3/(2\,\left( 1-\sqrt{3}\,\ve \right))\bigr)$,
$\mI_5=\bigl[3/(2\,\left( 1-\sqrt{3}\,\ve \right)),2 \bigr)$ and
$\mI_6=[2,\infty)$.
Notice that as $\ve \rightarrow 0$,
only $\mI_j$ for $j=2,4,5,6$ do not vanish,
so we only keep the components of $\mu^A_{_{PE}}(r,\ve)$ on these intervals.
See \cite{ceyhan:TR-rel-dens-NPE}
for the explicit form of $\mu^A_{_{PE}}(r,\ve)$ and for derivation.

Furthermore,
$\sigma_n^2(\ve)=\Var^A_{\ve}(\rho_{_{PE}}(n,r))=
\frac{1}{2\,n\,(n-1)}\Var^A_{\ve}[h_{12}]+\frac{(n-2)}{n\,(n-1)} \, \Cov^A_{\ve}[h_{12},h_{13}] $
whose explicit form is not calculated,
since we only need $\lim _{n \rightarrow \infty}\sqrt{n}\, \sigma_n(\ve=0)=\nu_{_{PE}}(r)$
which is given in  Equation \eqref{eqn:PEAsyvar}.

(PC1) follows for each $r \in [1,\infty)$ and $\ve \in \bigl[0,\sqrt{3}/3 \bigr)$ as in the segregation case.

Differentiating $\mu^A_{_{PE}}(r,\ve)$ with respect to $\ve$,
we get
\begin{multline*}
(\mu^A_{_{PE}})^{\prime}(r,\ve)=\varpi_{1,2}^{\prime}(r,\ve)\,\I(r \in [1,4/3))+\varpi_{1,4}^{\prime}(r,\ve)\,\I(r \in [4/3,3/2))\\
+\varpi_{1,5}^{\prime}(r,\ve)\,\I(r \in [3/2,2))+\varpi_{1,6}^{\prime}(r,\ve)\,\I(r \in [2,\infty))
\end{multline*}
where
\begin{align*}
\varpi_{1,2}^{\prime}(r,\ve)&=
-2\,\Bigl[\sqrt{3}\Bigl(-1152\,r^4\ve^3+720\,\sqrt{3}\,r^4\ve^2-288\,r^4\,\ve+
11\,\sqrt{3}\,r^4+2592\,\sqrt{3}\,r^2\ve^2-10368\,\sqrt{3}\,r\ve^2\\
&+432\,\sqrt{3}\,r^2+6480\,\sqrt{3}\ve^2-864\,\sqrt{3}\,r+432\,\sqrt{3}\Bigr)\ve\Bigr]\Big/
\Bigl[\left( -6\,\ve+\sqrt{3} \right)^3\left( 6\,\ve+\sqrt{3} \right)^3r^2\Bigr],\\
\varpi_{1,4}^{\prime}(r,\ve)&=
-2\,\Bigl[\sqrt{3}\Bigl(-1152\,r^4\ve^3+720\,\sqrt{3}\,r^4\ve^2-288\,r^4\,\ve+
11\,\sqrt{3}\,r^4-1296\,\sqrt{3}\,r^2\ve^2+108\,\sqrt{3}\,r^2\\
&-2160\,\sqrt{3}\ve^2-144\,\sqrt{3}\Bigr)\ve\Bigr]\Big/\Bigl[\left( -6\,\ve+\sqrt{3} \right)^3\left( 6\,\ve+\sqrt{3} \right)^3r^2\Bigr],\\
\varpi_{1,5}^{\prime}(r,\ve)&=\frac{2\,\ve\,(3\,r^4-72\,r^2-240\,\ve^2+192\,r-124)}{r^2(12\,\ve^2-1)^3},\\
\varpi_{1,6}^{\prime}(r,\ve)&=-\frac{40\,\ve}{r^2(12\,\ve^2-1)^2}.
\end{align*}
Hence $(\mu^A_{_{PE}})^{\prime}(r,\ve=0)=0$,
so we differentiate $(\mu^A_{_{PE}})^{\prime}(r,\ve)$ with respect to $\ve$ and get
\begin{multline*}
(\mu^A_{_{PE}})^{\prime\prime}(r,\ve)=
\varpi_{1,2}^{\prime\prime}(r,\ve)\,\I(r \in [1,4/3))+\varpi_{1,4}^{\prime\prime}(r,\ve)\,\I(r \in [4/3,3/2))\\
+\varpi_{1,5}^{\prime\prime}(r,\ve)\,\I(r \in [3/2,2))+\varpi_{1,6}^{\prime\prime}(r,\ve)\,\I(r \in [2,\infty))
\end{multline*}
where
\begin{align*}
\varpi_{1,2}^{\prime\prime}(r,\ve)&=
-6\,\Bigl[\sqrt{3}\Bigl(-27648\,r^4\ve^5+25920\,\sqrt{3}\,r^4\ve^4-18432\,r^4\ve^3+
2820\,\sqrt{3}\,r^4\ve^2+93312\,\sqrt{3}\,r^2\ve^4\\
&-576\,r^4\,\ve-373248\,\sqrt{3}\,r\ve^4+11\,\sqrt{3}\,r^4+33696\,\sqrt{3}\,r^2\ve^2+
233280\,\sqrt{3}\ve^4-82944\,\sqrt{3}\,r\ve^2\\
&+432\,\sqrt{3}\,r^2+45360\,\sqrt{3}\ve^2-864\,\sqrt{3}\,r+432\,\sqrt{3}\Bigr)\Bigr]\Big/
\Bigl[\left( 6\,\ve+\sqrt{3} \right)^4\,\left( -6\,\ve+\sqrt{3} \right)^4r^2\Bigr],\\
\varpi_{1,4}^{\prime\prime}(r,\ve)&=
-6\,\Bigl[\sqrt{3}\Bigl(-27648\,r^4\ve^5+25920\,\sqrt{3}\,r^4\ve^4-18432\,r^4\ve^3+
2820\,\sqrt{3}\,r^4\ve^2-46656\,\sqrt{3}\,r^2\ve^4\\
&-576\,r^4\,\ve+11\,\sqrt{3}\,r^4+2592\,\sqrt{3}\,r^2\ve^2-77760\,\sqrt{3}\ve^4+
108\,\sqrt{3}\,r^2-15120\,\sqrt{3}\ve^2-144\,\sqrt{3}\Bigr)\Bigr]\\
&\Big/\Bigl[\left( 6\,\ve+\sqrt{3} \right)^4\,\left( -6\,\ve+\sqrt{3} \right)^4r^2\Bigr],\\
\varpi_{1,5}^{\prime\prime}(r,\ve)&=-\frac{2\,(180\,r^4\ve^2+3\,r^4-4320\,r^2\ve^2-8640\,\ve^4+
11520\,r\ve^2-72\,r^2-8160\,\ve^2+192\,r-124)}{r^2(12\,\ve^2-1)^4},\\
\varpi_{1,6}^{\prime\prime}(r,\ve)&=\frac{40\,(36\,\ve^2+1)}{r^2(12\,\ve^2-1)^3}.
\end{align*}
Thus,
\begin{equation}
\label{eqn:PAE-A-''}
(\mu^A_{_{PE}})^{\prime\prime}(r,\ve=0)=
\begin{cases}
           -\frac{22}{9}\,r^2+192\,r^{-1}-96\,r^{-2}-96 &\text{for} \quad r \in [1,4/3),\\
           -\frac{22}{9}\,r^2+32\,r^{-2}-24 &\text{for} \quad r \in [4/3,3/2),\\
           -6\,r^2-384\,r^{-1}+248\,r^{-2}+144 &\text{for} \quad r \in [3/2,2),\\
            -40\,r^{-2} &\text{for} \quad r \in [2,\infty).
\end{cases}
\end{equation}
Note that $(\mu^A_{_{PE}})^{\prime\prime}(r,\ve=0) < 0$ for all $r \in [1,\infty)$,
so (PC2) follows with the second derivative.
(PC3) and (PC4) follow from continuity of $(\mu^A_{_{PE}})^{\prime\prime}(r,\ve)$
and $\sigma_n^2(r,\ve)$ in $\ve$.

Next, we find $c^A_{_{PE}}(r)=\lim_{n\rightarrow \infty}
\frac{(\mu^A_{_{PE}})^{\prime\prime}(r,\ve=0)}{\sqrt{n}\,\sigma_n(r,\ve=0)}=
\frac{\mu_A^{\prime\prime}(r,0)}{\sqrt{\nu_{_{PE}}(r)}}$,
by substituting the numerator from Equation \eqref{eqn:PAE-A-''} and
denominator from Equation \eqref{eqn:PEAsyvar}.
We can easily see that $c^A_{_{PE}}(r)<0$,
for all $r \ge 1 $.
Then (PC5) holds, so under
association alternatives $H^A_{\ve}$, the PAE of $\rho_{_{PE}}(n,r)$ is
$$
\PAE_{PE}^A(r) = \left(c^A_{_{PE}}(r)\right)^2=
   \frac{\bigl( (\mu^A_{_{PE}})^{\prime\prime}(r,\ve=0) \bigr)^2}{\nu_{_{PE}}(r)}.
$$

In Figure \ref{fig:PAE-seg-assoc} (right),
we present the PAE as a function of $r$ for association.
Notice that $\PAE_{PE}^A(r=1) = 174240/17 \approx 10249.41$,
$\lim_{r \rightarrow \infty} \PAE_{PE}^A(r) = 0$,
$\argsup_{r \in [1,\infty)} \PAE_{PE}^A(r) \approx 1.01$
with supremum $\approx 10399.77$.
$\PAE_{PE}^A(r)$ has also a local
supremum at $r_l\approx 1.44$ with local supremum $\approx 3630.89$.
Based on the Pitman asymptotic efficiency analysis,
we suggest,
for large $n$ and small $\ve$, choosing $r$ small for testing against association.
However, for small and moderate values of $n$,
normal approximation is not appropriate due to the skewness
in the density of $\rho_{_{PE}}(n,r)$.
Therefore, for small $n$,
we suggest moderate $r$ values.

\subsection{Pitman Asymptotic Efficiency for Central Similarity PCDs under the Association Alternative}
\label{sec:PAE-CS-assoc}

Consider the test sequences
$\rho_{_{CS}}(\tau)=\bigl\{ \rho_{_{CS}}(n,\tau) \bigr\}$
for sufficiently small $\ve >0$ and $\tau \in (0,\infty)$.
In the PAE framework above, the parameters are $\theta=\ve$ and $\theta_0=0$.
Suppose,
$\mu^A_{_{CS}}(r,\ve)=\E^A_{\ve}[\rho_{_{CS}}(n,\tau)]$.
For $\ve \in [0,\sqrt{3}/21)$,
$$\mu^A_{_{CS}}(\tau,\ve)=\sum_{j=1}^7 \varpi_{1,j}(\tau,\ve)\,\I(\tau \in \mI_j)$$
with the corresponding intervals
$\mI_1=\Bigl[0,\frac{3\,\sqrt{3}\,\ve}{2\,\left(1-\sqrt{3}\,\ve \right)}\Bigr)$,
$\mI_2=\Bigl[\frac{3\,\sqrt{3}\,\ve}{2\,\left(1-\sqrt{3}\,\ve\right)},\frac{2\,\sqrt{3}\,\ve}{1-2\,\sqrt{3}\,\ve}\Bigr)$,
$\mI_3=\Bigl[\frac{2\,\sqrt{3}\,\ve}{1-2\,\sqrt{3}\,\ve},\frac{3\,\sqrt{3}\,\ve}{1-\sqrt{3}\,\ve}\Bigr)$,
$\mI_4=\Bigl[\frac{3\,\sqrt{3}\,\ve}{1-\sqrt{3}\,\ve},\frac{3\,\sqrt{3}\,\ve}{1-4\,\sqrt{3}\,\ve}\Bigr)$,
$\mI_5=\Bigl[\frac{3\,\sqrt{3}\,\ve}{1-4\,\sqrt{3}\,\ve},\frac{6\,\sqrt{3}\,\ve}{1-\sqrt{3}\,\ve}\Bigr)$,
$\mI_6=\Bigl[\frac{6\,\sqrt{3}\,\ve}{1-\sqrt{3}\,\ve},
\frac{6\,\sqrt{3}\,\ve}{1-4\,\sqrt{3}\,\ve}\Bigr)$,
$\mI_7=\Bigl[\frac{6\,\sqrt{3}\,\ve}{1-4\,\sqrt{3}\,\ve},1\Bigr)$,
and
$\mI_8=[1,\infty)$.
Notice that as $\ve \rightarrow 0$, only the intervals $\mI_7$ and $\mI_8$ do not vanish,
so we only keep the component of $\mu^A_{_{CS}}(\tau,\ve)$ on these intervals.
See Section \ref{sec:arc-prob-agg-CS} for the explicit form of $\mu(\tau,\ve)$.

Furthermore,
$\sigma_n^2(\ve)=\Var^A_{\ve}(\rho_{_{CS}}(n,\tau))=
\frac{1}{2\,n\,(n-1)}\Var^A_{\ve}[h_{12}]+\frac{(n-2)}{n\,(n-1)} \, \Cov^A_{\ve}[h_{12},h_{13}]$
whose explicit form is not calculated,
since we only need $\lim_{n \rightarrow \infty}\sqrt{n}\, \sigma_n(\ve=0)=\nu_{_{CS}}(\tau)$ which
is given Equation \eqref{eqn:CSAsyvar}.

(PC1) follows for each $\tau \in (0,\infty)$ and $\ve \in \bigl[0,\sqrt{3}/3 \bigr)$ as in the segregation case.

Differentiating $\mu^A_{_{CS}}(\tau,\ve)$ with respect to $\ve$,
we get
$$
(\mu^A_{_{CS}})^{\prime}(\tau,\ve)=
\varpi_{1,7}^{\prime}(\tau,\ve)\,\I\left(r \in \bigl(0,1 \bigr)\right)+
\varpi_{1,8}^{\prime}(\tau,\ve)\,\I\left( r \in [1,\infty \bigr) \right)
$$
where
\begin{multline*}
\varpi_{1,7}^{\prime}(\tau,\ve) =
-72\,\Bigl[\sqrt{3} \bigl( -360\,\tau^4\ve^3+198\,\sqrt{3}\tau^4\ve^2-900\,\tau^3\ve^3-90\,{
\tau}^4\ve+495\,\sqrt{3}\tau^3\ve^2-360\,\tau^{
2}\ve^3+4\,\sqrt{3}\tau^4-225\,\tau^3\ve+\\
126\,\sqrt{3}\tau^2\ve^2+10\,\sqrt{3}\tau^3-90\,\tau^{
2}\ve-387\,\sqrt{3}\tau\,\ve^2+10\,\sqrt{3}\tau^2-
126\,\sqrt{3}\ve^2 \bigr) \ve \Bigr]\Big/\Bigl[ \left( 2\,\tau+1
\right)  \left( \tau+2 \right)  \left( -6\,\ve+\sqrt{3}
\right)^3 \left( 6\,\ve+\sqrt{3} \right)^3\Bigr],
\end{multline*}
and
$$
\varpi_{1,8}^{\prime}(\tau,\ve) = -24\,{\frac{\ve\, \left( 69\,\tau\,\ve^2+18\,\ve^2-6\,\tau-2 \right) }{ \left( \tau+2 \right)  \left( 2\,\tau+1
 \right)  \left( 12\,\ve^2-1 \right)^3}}.
$$

Hence $(\mu^A_{_{CS}})^{\prime}(\tau,\ve=0)=0$,
so we differentiate $(\mu^A_{_{CS}})^{\prime}(\tau,\ve)$
with respect to $\ve$ and get
\begin{multline*}
\varpi_{1,7}^{\prime\prime}(\tau,\ve) =
-216\,\Bigl[\sqrt{3} \bigl( -8640\,\tau^4\ve^5+7128\,\sqrt{3}\tau^4\ve^4-21600\,\tau^3\ve^5-5760
\,\tau^4\ve^3+17820\,\sqrt{3}\tau^3\ve^4-
8640\,\tau^2\ve^5+834\,\sqrt{3}\tau^4\ve^2-\\
14400\,\tau^3\ve^3+4536\,\sqrt{3}\tau^2\ve^{4
}-180\,\tau^4\ve+2085\,\sqrt{3}\tau^3\ve^2-5760
\,\tau^2\ve^3-13932\,\sqrt{3}\tau\,\ve^4+4\,
\sqrt{3}\tau^4-450\,\tau^3\ve+978\,\sqrt{3}\tau^2\ve^2-\\
4536\,\sqrt{3}\ve^4+10\,\sqrt{3}\tau^3-
180\,\tau^2\ve-1161\,\sqrt{3}\tau\,\ve^2+10\,\sqrt
{3}\tau^2-378\,\sqrt{3}\ve^2 \bigr) \Bigr]\Big/\Bigl[ \left( 2\,\tau+1
 \right)  \left( \tau+2 \right)  \left( -6\,\ve+\sqrt{3}
 \right)^4 \left( 6\,\ve+\sqrt{3} \right)^4\Bigr],
\end{multline*}
and
$$
\varpi_{1,8}^{\prime\prime}(\tau,\ve) = 24\,{\frac{2484\,\tau\,\ve^4+648\,\ve^4-153\,\tau\,\ve^2-66\,\ve^2-6\,\tau-2}{ \left( \tau+2 \right)
 \left( 2\,\tau+1 \right)  \left( 12\,\ve^2-1 \right)^4}}.
$$

Thus
\begin{equation}
\label{eqn:CS-PAE-A-1}
 (\mu^A_{_{CS}})^{\prime\prime}(\tau,\ve=0)=
\begin{cases}
         {\frac{-16\,\tau^2 \left( 2\,\tau^2+5\,\tau+5 \right) }{ \left( 2\,\tau+1 \right)  \left( \tau+2 \right) }} &\text{for} \quad r \in (0,1),\\
         {\frac{-48\,(3\,\tau+1)}{ \left( 2\,\tau+1 \right)  \left( \tau+2 \right) }} &\text{for} \quad r \in [1,\infty).
\end{cases}
\end{equation}
Note that $(\mu^A_{_{CS}})^{\prime\prime}(\tau,\ve=0)<0$ for all $\tau \in (0,\infty)$,
so (PC2) follows with the second derivative.
(PC3) and (PC4) follow from continuity of
$\mu^{\prime\prime}(\tau,\ve)$ and
$\sigma_n^2(\tau,\ve)$ in $\ve$.

Next, we find
$c^A_{_{CS}}(\tau)=\lim_{n\rightarrow
\infty}\frac{(\mu^A_{_{CS}})^{\prime\prime}(\tau,\ve=0)}{\sqrt{n}\,\sigma_n(\tau,\ve=0)}=
\frac{\mu_A^{\prime\prime}(\tau,0)}{\sqrt{\nu_{_{CS}}(\tau)}}$,
by substituting the numerator from Equation \eqref{eqn:CS-PAE-A-1} and
denominator from Equation \eqref{eqn:CSAsyvar}.
We can easily see that
$c^A_{_{CS}}(\tau)<0$, for all $\tau \in (0,\infty)$.
Then (PC5) holds, so under association alternatives $H^A_{\ve}$,
the PAE of $\rho_{_{CS}}(n,\tau)$ is
$$
\PAE_{CS}^A(\tau) = \left(c^A_{_{CS}}(\tau)\right)^2=
   \frac{\left( (\mu^A_{_{CS}})^{\prime\prime}(\tau,\ve=0) \right)^2}{\nu_{_{CS}}(\tau)}.
$$

In Figure \ref{fig:PAE-seg-assoc} (right), we present the PAE as a
function of $\tau$ for association.
Notice that $\lim_{\tau \rightarrow 0} \PAE_{CS}^A(\tau) = 72000/7 \approx 10285.71$
which is also the global maximum,
$\PAE_{CS}^A(\tau=1) = 61440/7 \approx 8777.14$
which is also a local maximum.
Moreover,
a local minimum of $\PAE_{CS}^A(\tau)$ occurs at $\tau \approx .45$
with PAE score being equal to $\approx 6191.67$.
Based on the Pitman asymptotic efficiency analysis,
we suggest, for large $n$ and small $\ve$,
choosing $\tau$ small for testing against association.
However, for small and moderate values of $n$,
normal approximation is not
appropriate due to the skewness in the density of $\rho_{_{CS}}(n,\tau)$.
Therefore, for small $n$, we suggest $\tau \approx 1$.

Comparing the PAE scores
of the relative density of proportional-edge PCDs
and
central similarity PCDs
under association alternatives,
we see that $\PAE_{PE}^A(t) < \PAE_{CS}^A(t)$
for $1 \le t \lesssim 1.4564)$
and for $t \gtrsim 1.5192)$;
and
$\PAE_{PE}^A(t) > \PAE_{CS}^A(t)$
for $1.4564 \lesssim t \lesssim 1.5192$.
Under association,
relative density of central similarity PCD is asymptotically more
efficient compared to the proportional-edge PCD.
Furthermore,
$\PAE_{PE}^A(t)$ goes to 0
as $t \rightarrow \infty$ at rate $O(t^{-2})$,
while
$\PAE_{CS}^A(t)$ goes to $\infty$
as $t \rightarrow \infty$ at rate $O(t^{-1})$.

\begin{remark}
\textbf{Hodges-Lehmann Asymptotic Efficiency:}
PAE analysis is local (around $\ve=0$) and for arbitrarily large $n$.
The comparison would hold in general
provided that $\mu(r,\ve)$ is convex in $\ve$ for all $\ve \in \bigl[ 0,\sqrt{3}/3 \bigr)$.
As an alternative, we fix an $\ve$ under segregation alternative and then compare the
asymptotic behavior of $\rho_{_{PE}}(n,r)$ with Hodges-Lehmann asymptotic
efficiency in (\cite{ceyhan:TR-rel-dens-NPE}).

Hodges-Lehmann asymptotic efficiency (HLAE) of
$\rho_{_{CS}}(n,\tau)$  is given by
$$
\HLAE(\tau,\ve):=\frac{\bigl(
\mu(\tau,\ve)-\mu_{_{CS}}(\tau)
\bigr)^2}{\nu(\tau,\ve)}.
$$
Unlike PAE, HLAE does only involve $n \rightarrow \infty$  at a
fixed $\ve > 0$. Hence HLAE requires the mean and,
especially, the asymptotic variance of $\rho_{_{CS}}(n,\tau)$ \emph{ under a
fixed alternative}.  So, one can  investigate HLAE for specific
values of $\ve$, if not for all $\ve \in \left(0,\sqrt{3}/3 \right)$.
$\square$
\end{remark}

\begin{remark}
The asymptotic power function allows investigation as
a function of the expansion parameter, $n$, and $\ve$ using the asymptotic
critical value and an appeal to normality.
The asymptotic power functions of $\rho_{_{PE}}(n,r)$ under the alternatives
is investigated in (\cite{ceyhan:TR-rel-dens-NPE}).

Under a specific
segregation alternative $H^S_{\ve}$,
the asymptotic power function of $\rho_{_{CS}}(n,\tau)$ is given by
\begin{eqnarray*}
\Pi_S(n,\tau,\ve) = 1-\Phi \left(\frac{z_{(1-\alpha)} \cdot
\sqrt{\nu_{_{CS}}(\tau)}+\sqrt{n}\,(\mu_{_{CS}}(\tau)-\mu^S_{_{CS}}(\tau,\ve))}{\sqrt{\nu^S_{_{CS}}(\tau,\ve)}}
\right).
\end{eqnarray*}

Under $H^A_{\ve}$,
we have
\begin{eqnarray*}
\Pi_A(n,\tau,\ve) = \Phi
\left(\frac{z_{\alpha}\,\sqrt{\nu_{_{CS}}(\tau)}+\sqrt{n}\cdot(\mu_{_{CS}}(\tau)-\mu^A_{_{CS}}(\tau,\ve))}{\sqrt{\nu^A_{_{CS}}(\tau,\ve)}}
\right). \square
\end{eqnarray*}
\end{remark}

\subsection{Pitman Asymptotic Efficiency Analysis in the Multiple Triangle Case}
\label{sec:PAE-multi-triangle}
For $J_m>1$ (i.e., $m>3$),
in addition to the expansion parameter,
PAE analysis depends on
the number of triangles as well as the relative sizes of the triangles
(i.e., on $\Y_m$).
So the optimal expansion parameter values with respect to the PAE
criteria in the multiple triangle case might be different
than that of the one triangle case.

Given the values of $J_m$ and $\mathcal W$,
under segregation alternative $H^S_{\ve}$,
the PAE for the relative density of proportional-edge PCDs
is given by
\begin{equation}
\label{eqn:PAE-PE-multi-tri}
\PAE_{_{PE}}^S(m,r) = \frac{\bigl( (\widetilde \mu^S_{_{PE}})^{\prime\prime}(m,r,\ve=0)\bigr)^2}{\widetilde \nu_{_{PE}}(m,r)}=
   \frac{\left( (\mu^S_{_{PE}})^{\prime\prime}(r,\ve=0)\,\sum_{j=1}^{J_m}w_j^2 \right)^2}
   {\nu_{_{PE}}(r) \,\sum_{j=1}^{J_m}w_j^3 +4\,\mu_{_{PE}}(r)^2\left(\sum_{j=1}^{J_m}w_j^3-\left(\sum_{j=1}^{J_m}w_j^2 \right)^2\right)}.
\end{equation}
PAE score for the relative density of proportional-edge PCDs
under the association alternative is similar.

Similarly,
the PAE for the relative density of central similarity PCDs under segregation alternative $H^S_{\ve}$
is given by
\begin{equation}
\label{eqn:PAE-CS-multi-tri}
\PAE_{_{CS}}^S(m,\tau) = \frac{\bigl( (\widetilde \mu^S_{_{CS}})^{\prime\prime}(m,\tau,\ve=0)\bigr)^2}{\widetilde \nu_{_{CS}}(m,\tau)}=
   \frac{\left( (\mu^S_{_{CS}})^{\prime\prime}(\tau,\ve=0)\,\sum_{j=1}^{J_m}w_j^2 \right)^2}
   {\nu_{_{CS}}(\tau) \,\sum_{j=1}^{J_m}w_j^3 +4\,\mu_{_{CS}}(\tau)^2\left(\sum_{j=1}^{J_m}w_j^3-\left(\sum_{j=1}^{J_m}w_j^2 \right)^2\right)}.
\end{equation}
PAE score for the relative density of central similarity PCDs
under the association alternative is similar.

\begin{figure}[h]
\centering
\psfrag{mu}{{\Huge PAE score}}
\psfrag{t}{{\Huge expansion parameter}}
\rotatebox{-90}{ \resizebox{2. in}{!}{\includegraphics{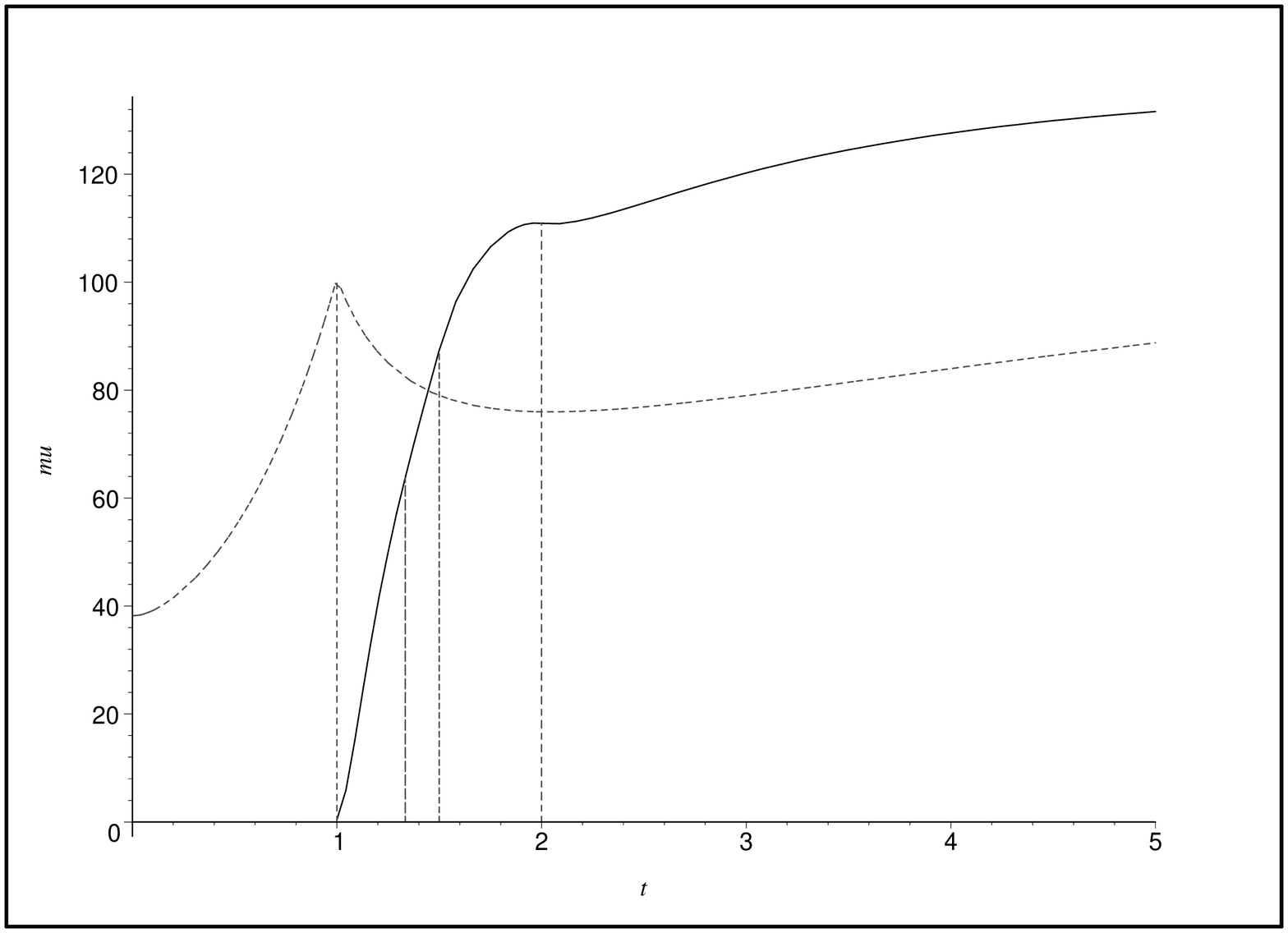}}}
\psfrag{mu}[c]
{\begin{picture}(0,0)
\put(-10.5,-0.5){\makebox(0,0)[l]{\Huge PAE score}}
\end{picture}}
\rotatebox{-90}{ \resizebox{2. in}{!}{\includegraphics{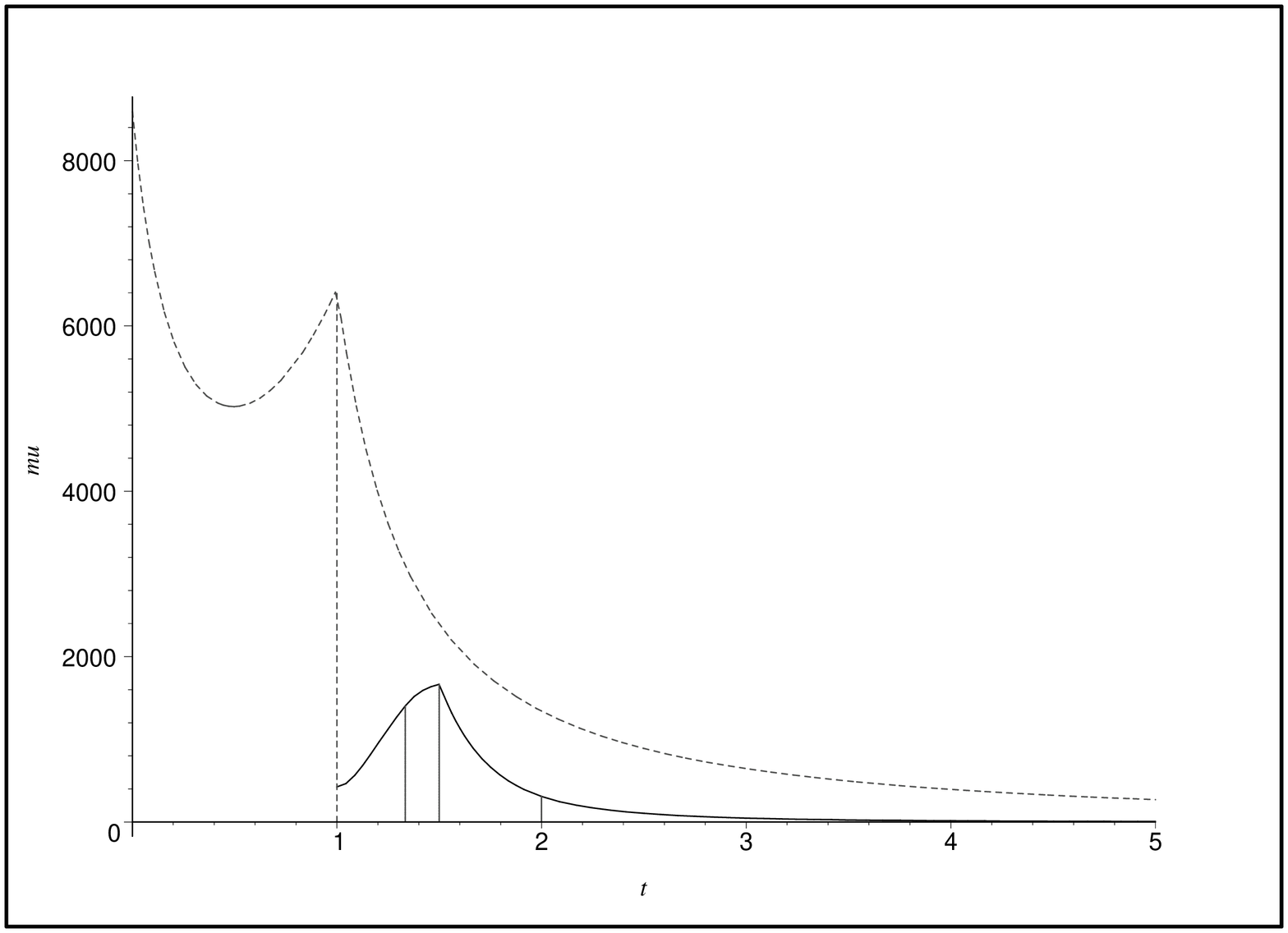}}}
\caption{
\label{fig:MT-PAE-seg-assoc}
Pitman asymptotic efficiency against
segregation (left) and association (right) alternatives
as a function of expansion parameters
in the multiple triangle case with $J_m=13$
for the relative density of proportional-edge PCDs (solid line) and
central similarity PCDs (dashed line).
Notice that vertical axes are differently scaled.
}
\end{figure}

In Figure \ref{fig:MT-PAE-seg-assoc} (left),
we present the PAE scores as a function
the expansion parameter under segregation alternative
conditional on the realization of $\Y_m$ given in Figure \ref{fig:deldata}.
Notice that,
unlike the one triangle case,
$\PAE_{_{PE}}^S(m,r)$ is bounded with $\lim_{r \rightarrow \infty}
\PAE_{_{PE}}^S(m,r) =
\frac{8\,\sum_{j=1}^{J_m}w_j^2}{256\,\left(\sum_{j=1}^{J_m}w_j^3-\left(\sum_{j=1}^{J_m}w_j^2 \right)^2\right)}
\approx 139.34$.
Some values of interest are
$\PAE_{_{PE}}^S(m,r=1) \approx .39$, and a local maximum value
of $\approx 110.97$ is attained at the
$\argsup_{r \in [1,2]}\PAE_{_{PE}}^S(m,r) \approx 1.97$.
On the other hand,
the PAE curve for the central similarity PCDs in the multiple triangle case
is similar that in the one triangle case
(See Figure \ref{fig:PAE-seg-assoc} (left)).
But unlike the one triangle case,
$\PAE_{_{CS}}^S(m,\tau)$ is bounded with $\lim_{\tau \rightarrow \infty}
\PAE_{_{CS}}^S(m,\tau) \approx 139.34$.
Some values of note are
$\lim_{\tau \rightarrow 0}\PAE_{_{CS}}^S(m,\tau) \approx 38.20$,
and a local maximum of $\approx 100.77$ is attained at
$\argsup_{\tau \in (0,2)}\PAE_{_{CS}}^S(m,\tau) =1$;
and
a local minimum of $\approx 75.97$ is attained at
$\arginf_{\tau \in (1,3)}\PAE_{_{CS}}^S(m,\tau) \approx 2.04$.
Based on the PAE analysis of the relative density of proportional-edge PCDs,
under segregation alternative
larger $r$ values have larger asymptotic relative efficiency.
However, due to the skewness of the pdf of $\rho_{_{PE}}(m,r)$,
moderate $r$ values ($r$ around 1.5 or 2) are recommended.
As for the central similarity PCDs,
larger $\tau$ values have larger asymptotic relative efficiency.
However, due to the skewness of the pdf of $\rho_{_{CS}}(m,\tau)$,
moderate $\tau$ values ($\tau$ around 1) are recommended.

Comparing the PAE scores for proportional-edge and central similarity PCDs
under the segregation alternative,
we see that
for $1 \le t \lesssim 1.45$
asymptotic relative efficiency of relative density of central similarity PCDs is larger
since $\PAE_{_{CS}}^S(m,t) > \PAE_{_{PE}}^S(m,t)$,
and
for $t \gtrsim 1.45$
asymptotic relative efficiency of relative density of proportional-edge PCDs is larger
since $\PAE_{_{CS}}^S(m,t) < \PAE_{_{PE}}^S(m,t)$.
Therefore, proportional-edge PCD tends to be more asymptotically efficient
compared to the central similarity PCD under segregation.

In Figure \ref{fig:MT-PAE-seg-assoc} (right),
we present the PAE scores as a function
the expansion parameter under association alternative
conditional on the realization of $\Y_m$ given in Figure \ref{fig:deldata}.
Notice that,
as in the one triangle case,
$\PAE_{_{PE}}^A(m,r)$ tends to 0 as $r \rightarrow \infty$.
Some values of interest are
$\PAE_{_{PE}}^A(m,r=1) \approx 422.96$, and a global maximum value
of $\approx 1855.97$ is attained at $r=1.5$
On the other hand,
the PAE curve for the central similarity PCDs in the multiple triangle case
is similar to the one in the one triangle case
(See Figure \ref{fig:PAE-seg-assoc} (left)).
Note that
$\lim_{\tau \rightarrow 0} \PAE_{_{CS}}^A(m,\tau) \approx 8593.97$;
a local maximum value of $\approx 6449.54$ is attained at $\tau=1$;
and
a local minimum value of $\approx 5024.22$ is attained at $\tau \approx 0.49$.
Moreover,
$\lim_{\tau \rightarrow \infty} \PAE_{_{CS}}^A(m,\tau) = 0$
at rate $O(\tau^{-2})$.
Based on the PAE analysis for relative density of proportional-edge PCDs,
smaller $\tau$ values tend to have larger asymptotic relative efficiency.
However,
we suggest,
for large $n$ and small $\ve$,
choosing moderate $\tau$ for testing
against association
due to the skewness of the density of $\rho_{_{CS}}(n,\tau)$
for very small $\tau$ values.

Comparing the PAE scores for proportional-edge and central similarity PCDs
under the association alternative,
we see that
for $t \ge 1$
asymptotic relative efficiency of relative density of central similarity PCDs is larger
since $\PAE_{_{CS}}^A(m,t) > \PAE_{_{PE}}^A(m,t)$.
Therefore, central similarity PCD tends to be more asymptotically efficient
compared to the proportional-edge PCD under association.

\begin{remark}
\textbf{Empirical Power Comparison versus PAE Comparison for the Two PCD Families:}
Notice that the finite sample performance
(based on the Monte Carlo simulations)
and the asymptotic efficiency
(based on PAE scores)
seem to give conflicting results.
The reason for this is two fold:
(i) in the Monte Carlo simulations,
we only have a finite number of observations,
and the asymptotic normality of the relative density of the PCDs
require smaller sample sizes for moderate values of the expansion parameters,
and
(ii) PAE is designed for infinitesimal deviations from the null hypothesis
(i.e., as close as possible to the null case),
while in our simulations we use mild to severe but fixed levels of deviations.
Hence, if we had extremely large samples,
the results of our finite sample and asymptotic comparisons would agree
under extremely mild segregation or association.

Furthermore,
when the PAE scores are compared at the optimal expansion parameters,
the comparison results agree with that of the Monte Carlo simulation results.
In particular,
recall that in the one triangle case,
the optimal parameters for proportional-edge PCDs were 1.5 and 2
and for central similarity PCDs, they were 8 and 5
against mild segregation and association, respectively.
Under segregation, central similarity PCD is asymptotically more efficient,
while under association proportional-edge PCD is asymptotically more efficient
at these optimal parameters.
This agrees with the conclusion of empirical power comparison.
In the multiple triangle case,
the optimal parameters for proportional-edge PCDs were 1.5 and 2
and for central similarity PCDs, they were 7 and 1
against mild segregation and association, respectively.
Under both alternatives,
central similarity PCD is asymptotically more efficient.
In this case,
only the segregation results are in agreement.
The power estimates under association were virtually same at these optimal values
with both PCD families.
$\square$
\end{remark}

An extension of proportional-edge proximity regions
and central similarity proximity regions
to higher dimensions
(hence the corresponding PCDs to data in higher dimensions)
are provided in
\cite{ceyhan:arc-density-PE} and \cite{ceyhan:arc-density-CS}, respectively.

\section{Correction for $\X$ Points Outside the Convex Hull of $\Y_m$}
\label{sec:conv-hull-correction}
Our null hypothesis in Equation \eqref{eqn:null-pattern-mult-tri} is rather restrictive,
in the sense that,
it might not be realistic to assume the support of $\X$
being $C_H(\Y_m)$ in practice.
Up to now, our inference was restricted to the $C_H(\Y_m)$.
However, crucial information from the data (hence power) might be lost,
since a substantial proportion of $\X$ points, denoted $\pi_{\out}$,
might fall outside the $C_H(\Y_m)$.
A correction is suggested in (\cite{ECarXivGamforSpat:2009,ceyhan:dom-num-NPE-Spat2010})
to mitigate the effect of $\pi_{\out}$ (or restriction to the $C_H(\Y_m)$)
on the use of the domination number for the proportional-edge PCDs.
We propose a similar correction for the points outside the $C_H(\Y_m)$
for the relative density in this article.

Along this line,
\cite{ECarXivGamforSpat:2009,ceyhan:dom-num-NPE-Spat2010}
estimated the $\pi_{\out}$ values
for independently generated $\X_n$ and $\Y_m$ as random samples from  $\U((0,1)\times(0,1))$.
The considered values were $n=100,200,\ldots,900,1000$, $2000,\ldots,9000$, $10000$
for each of $m=10,20,\ldots,50$.
The procedure is repeated $N_{mc}=1000$ times for each $n,m$ combination.
Let $\widehat \pi_{\out}$ be the estimate
of the proportion of $\X$ points outside the $C_H(\Y_m)$
which is obtained by averaging the $\pi_{\out}$ values
(over $n$) for each $m,n$ combination.
The simulation results suggested that
$\widehat \pi_{\out} \approx 1.7932/m+1.2229/\sqrt{m}$ (\cite{ceyhan:dom-num-NPE-Spat2010}).
Notice that as $m \rightarrow \infty$,
$\widehat \pi_{\out} \rightarrow 0$.

Based on the Monte Carlo simulation results,
we propose a coefficient to adjust for the proportion of $\X$ points outside $C_H(\Y_m)$,
namely,
\begin{equation}
\label{eqn:conv-hull-correct}
C_{ch}:=\signum(p_{\out}-\E[\widehat \pi_{\out}]) \times (p_{\out}-\E[\widehat \pi_{\out}])^2
\end{equation}
where $\signum(p_{\out}-\E[\widehat \pi_{\out}])$ is the sign of the difference $p_{\out}-\E[\widehat \pi_{\out}]$
and
$p_{\out}$ is the observed and $\E[\widehat \pi_{\out}] \approx 1.7932/m+1.2229/\sqrt{m}$
is the expected proportion of $\X$ points outside $C_H(\Y_m)$.
For the test statistics in
Section \ref{sec:test-stat-analysis},
we suggest
\begin{equation}
\label{eqn:test-stat-adj}
\widetilde R^{ch}_{PE}(r):= \widetilde R_{PE}(r)+ |\widetilde R_{PE}(r)| \cdot C_{ch}
\text{ and }
\widetilde R^{ch}_{CS}(r):= \widetilde R_{CS}(\tau)+ |\widetilde R_{CS}(\tau)| \cdot C_{ch}
\end{equation}
Note that this (convex hull) adjustment slightly affects the empirical size estimates under CSR
of $\X$ and $\Y$ points in the same rectangular supports,
since $p_{\out}$ and $\E[\widehat \pi_{\out}]$ values would be very similar.
On the other hand,
under segregation alternatives,
we expect $\widetilde R^{ch}_{PE}(r)$ value and $p_{\out}-\E[\widehat \pi_{\out}]$ to be positive,
so the convex hull correction increases the value of $\widetilde R_{PE}(r)$
in favor of the right-sided alternative (i.e., segregation).
Under association alternatives,
we expect $\widetilde R^{ch}_{PE}(r)$ value and $p_{\out}-\E[\widehat \pi_{\out}]$ to be negative,
so the convex hull correction decreases the value of $\widetilde R_{PE}(r)$
in favor of the left-sided alternative (i.e., association).


\section{Example Data Set}
\label{sec:example}
We illustrate the method on an ecological data set
(namely, swamp tree data of \cite{dixon:EncycEnv2002}).
\cite{good:1982} considered the spatial patterns of tree species
along the Savannah River, South Carolina, U.S.A.
From this data, \cite{dixon:EncycEnv2002} used a single 50m $\times$ 200m rectangular plot
(denoted as the $(0,200)\times(0,50)$ rectangle)
to illustrate his nearest neighbor contingency table (NNCT) methods.
All live or dead trees with 4.5 cm or more dbh (diameter at breast height)
were recorded together with their species.
Hence it is an example of a realization of a marked multi-variate point pattern.
The plot contains 13 different tree species,
four of which comprising over 90 \% of the 734 tree stems.
See \cite{ceyhan:class2009} for more detail on the data.

\begin{figure}[ht]
\centering
\rotatebox{-90}{ \resizebox{3. in}{!}{\includegraphics{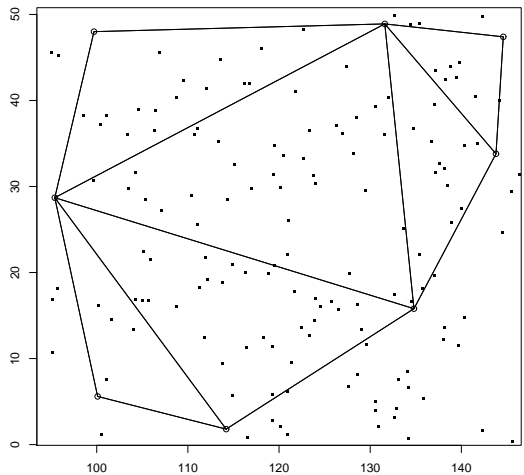} }}
 \caption{
\label{fig:SwampTrees}
The scatter plot of the locations of bald cypresses (circles $\circ$)
and other trees (black squares $\blacksquare$)
in the swamp tree data.
The Delaunay triangulation is based on the locations of the bald cypresses.}
\end{figure}

In this article,
we only consider the middle
50m $\times$ 55m rectangular plot from the original study area
(i.e., the subset $(95,150)\times(0,50)$ of the 50m $\times$ 200m rectangular plot)
and
investigate the spatial interaction of all other tree species (i.e., other than bald cypress trees)
with bald cypresses
(i.e., bald cypresses are taken to be the $\Y$ points,
while all other trees are taken to be the $\X$ points;
hence Delaunay triangulation is based on the locations of bald cypresses).
The study area contains 8 bald cypress trees and 156 other trees.
See also Figure \ref{fig:SwampTrees} which is suggestive of segregation
of other trees from bald cypresses.

\begin{figure}[]
\centering
\rotatebox{-90}{ \resizebox{2.1 in}{!}{ \includegraphics{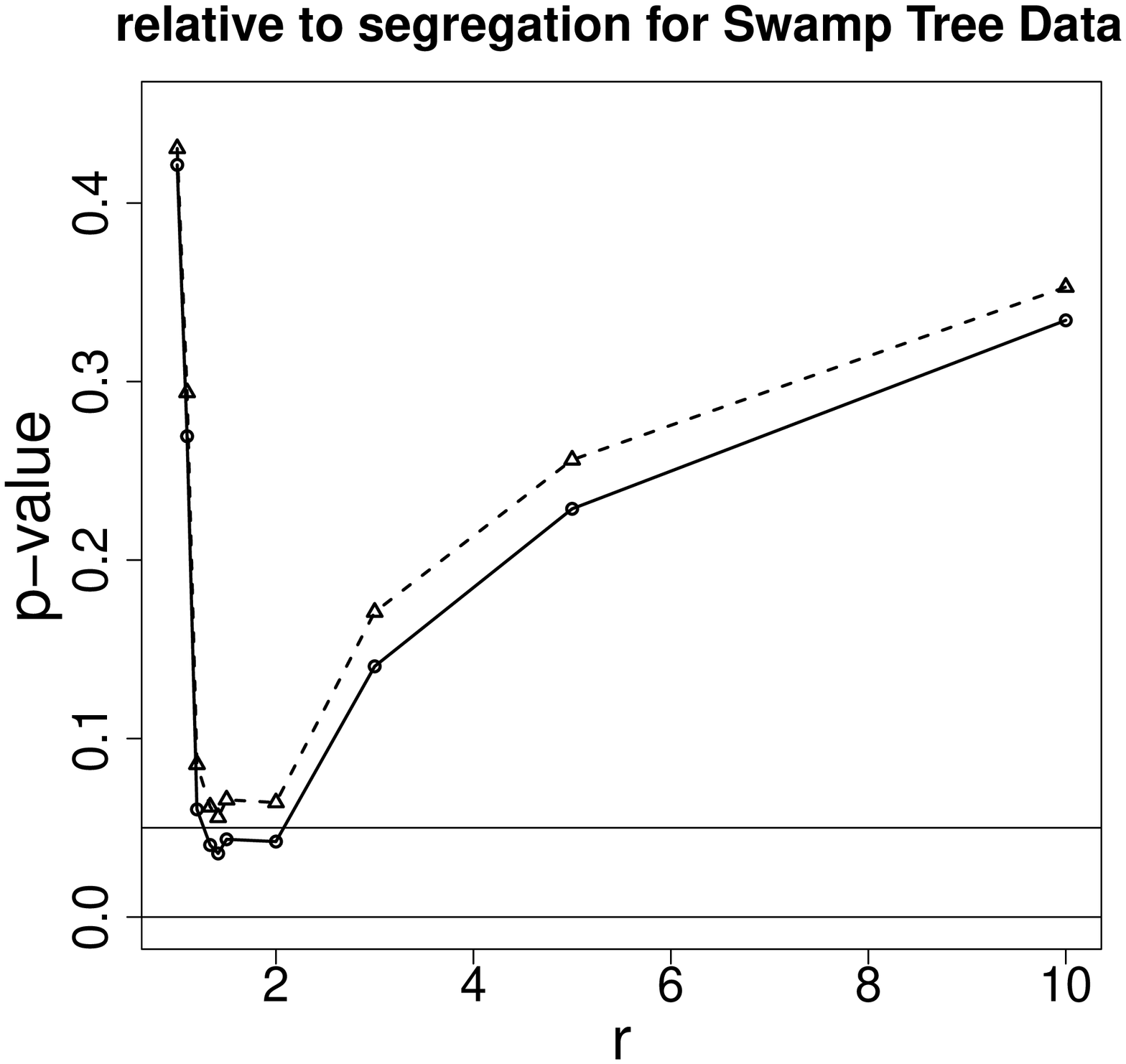}}}
\rotatebox{-90}{ \resizebox{2.1 in}{!}{ \includegraphics{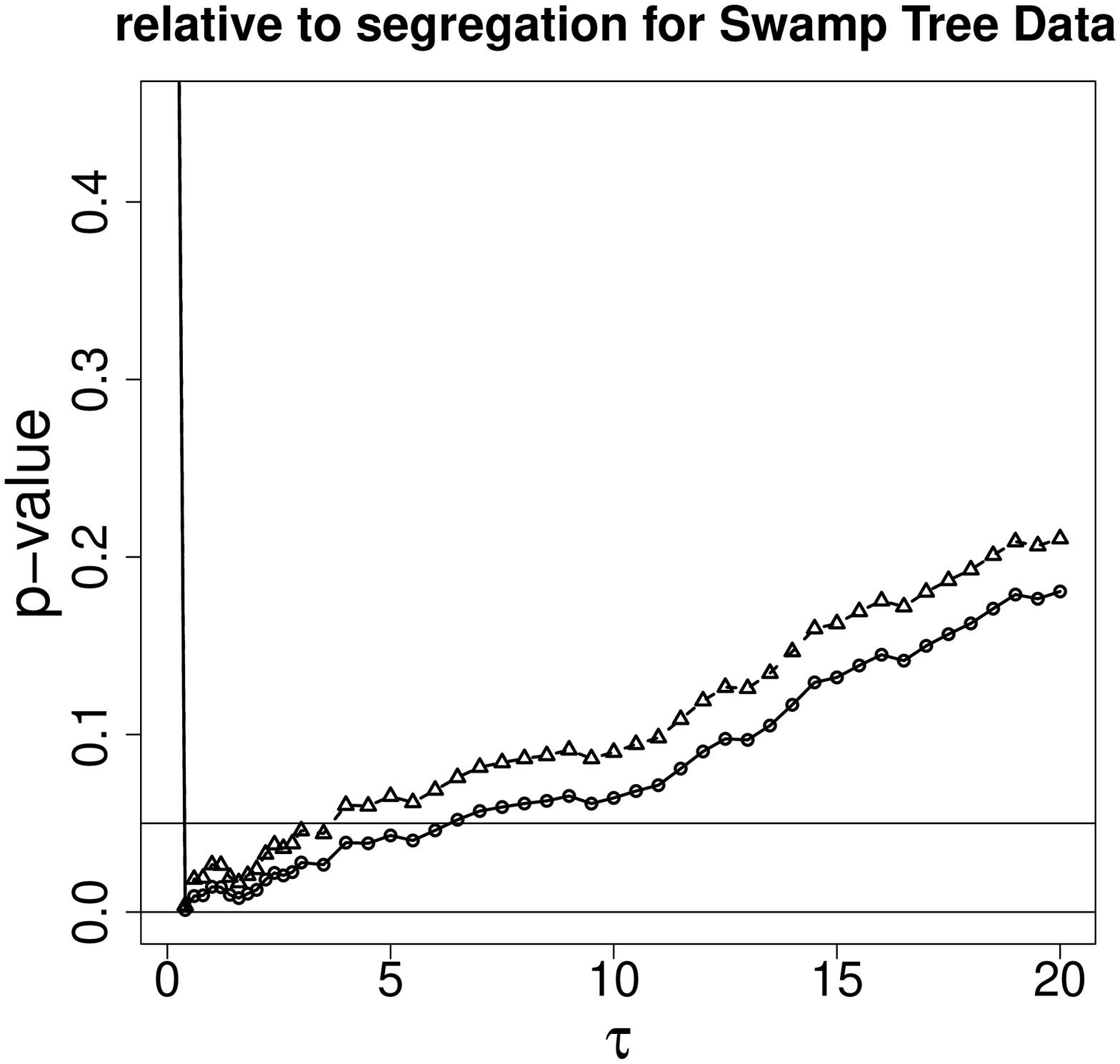}}}
\caption{
\label{fig:p-values-swamp-tree}
The $p$-values based on proportional-edge PCDs (left)
and
central similarity PCDs (right)
with convex hull corrected test statistics (circles connected with solid lines)
and
uncorrected test statistics (triangles connected with dashed lines).
The horizontal lines are at 0 and 0.05 values.
Notice that the horizontal axes are differently scaled.
}
\end{figure}

For this data,
we find that 108 other trees are inside and 48 are outside of
the convex hull of bald cypresses.
Hence
the proportion of other trees outside the convex hull of bald cypresses is $p_{\out}=0.3077$
and the expected proportion is $\pi_{\out}=0.6515$.
Hence the convex hull correction decreases the magnitude of the raw test statistics.
We calculate the standardized test statistics,
$R_{PE}(r)$, for $r=1, 11/10, 6/5, 4/3, \sqrt{2},3/2, 2, 3, 5, 10$ values
and, $R_{CS}(\tau)$,
for $\tau= 0.2,0.4,0.6.\ldots,3.0,3.5,4.0,\ldots,20.0$ values
and the corresponding convex hull corrected versions.
The $p$-values based on the normal approximation are presented
in Figure \ref{fig:p-values-swamp-tree}.
Observe that with $R_{PE}(r)$,
the convex hull corrected version is not significant
(for both the right- and the left-sided alternatives)
at 0.05 level at any of the $r$ values considered
(only significant at 0.10 level at $r$ between 1.4 and 2.0 for the right-sided alternative),
while the uncorrected version is significant (at 0.05 level)
for $r$ values between 1.4 and 2.0.
On the other hand,
with $R_{CS}(\tau)$,
the convex hull corrected version is significant (for the right-sided alternatives)
at 0.05 level at $\tau$ values between 0.2 and 4.0,
while the uncorrected version is significant (at 0.05 level)
for $\tau$ values between 0.2 and 7.
Hence, there is significant evidence for segregation of other trees from bald cypresses.

\begin{figure}[]
\centering
\rotatebox{-90}{ \resizebox{2.1 in}{!}{ \includegraphics{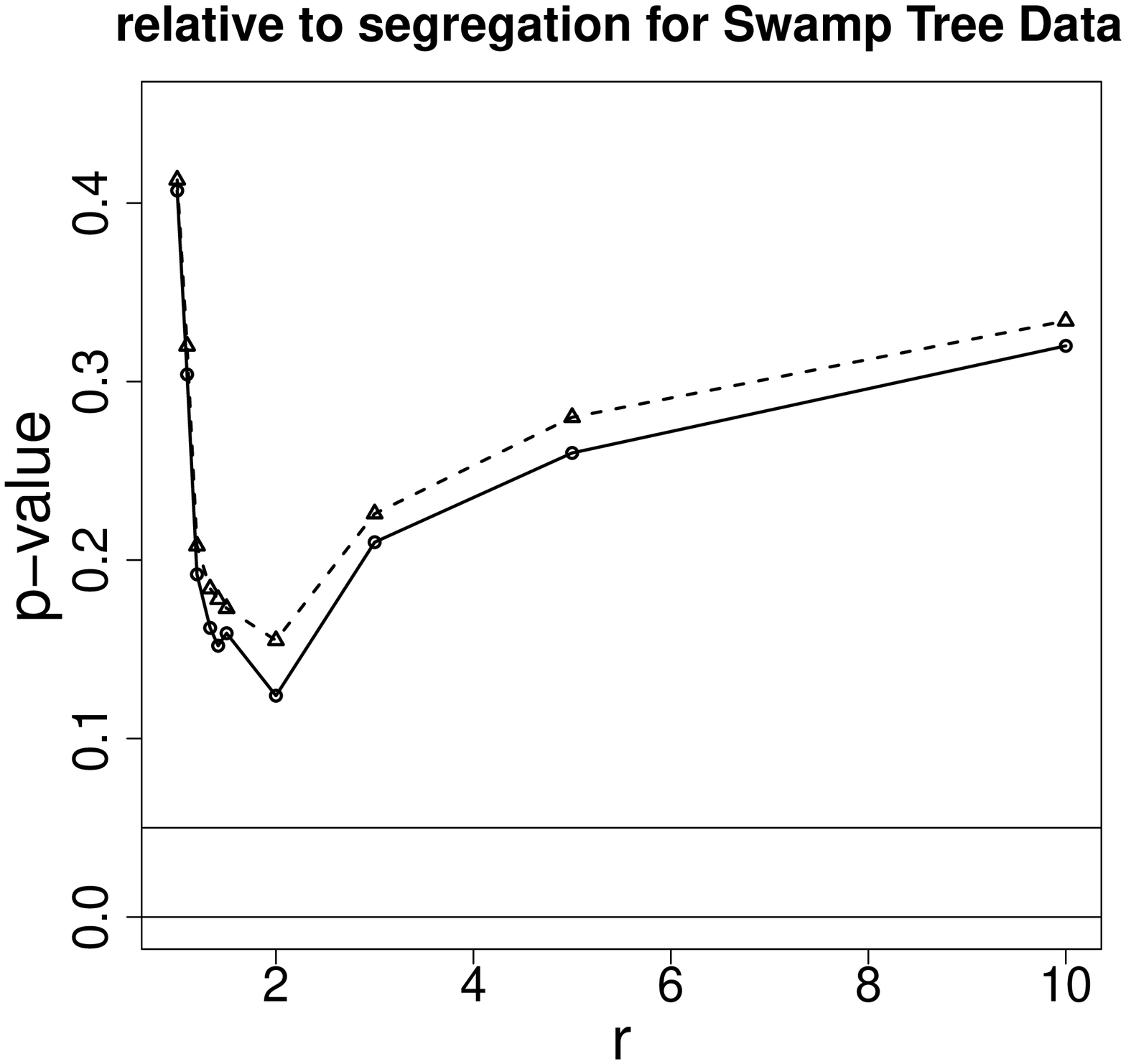}}}
\rotatebox{-90}{ \resizebox{2.1 in}{!}{ \includegraphics{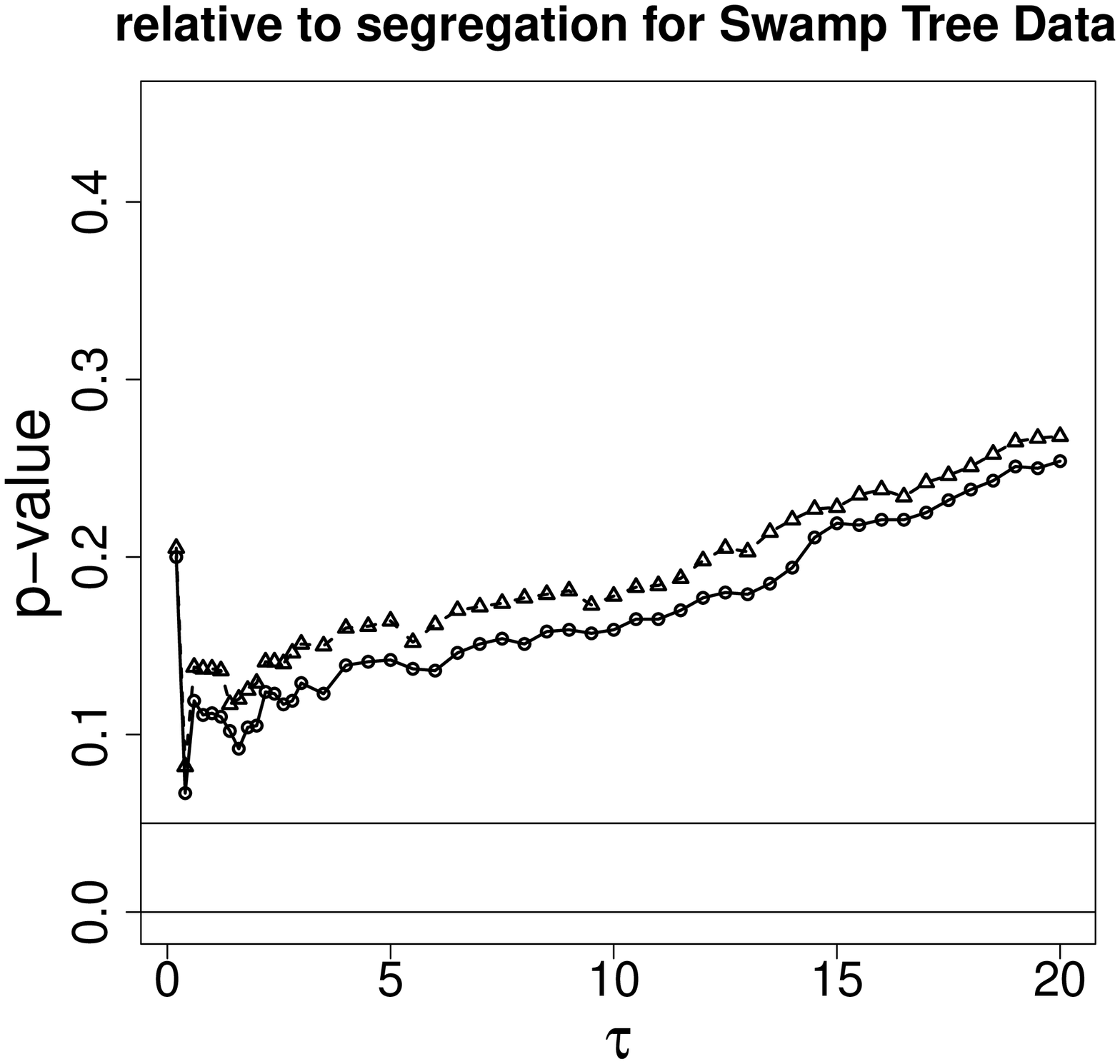}}}
\caption{
\label{fig:MC-p-values-swamp-tree}
The $p$-values estimated by Monte Carlo randomization for
proportional-edge PCDs (left)
and
central similarity PCDs (right)
with convex hull corrected test statistics (circles connected with solid lines)
and
uncorrected test statistics (triangles connected with dashed lines).
The horizontal lines are at 0 and 0.05 values.
Notice that the horizontal axes are differently scaled.
}
\end{figure}

We also perform a Monte Carlo randomization test as follows.
First we calculate the standardized relative density values,
denoted $\widetilde R_{PE}^{obs}(r)$ and $\widetilde R_{CS}^{obs}(\tau)$ and for the current data set,
so they are observed test statistics.
Then we randomly assign 8 of the trees as ``bald cypresses" (without replacement)
and the remaining trees as ``the other trees",
then calculate the test statistics (standardized relative density scores)
for the other trees within the convex hull of the bald cypresses.
We repeat this procedure 999 times.
Combining the observed $\widetilde R_{PE}^{obs}(r)$ and $\widetilde R_{CS}^{obs}(\tau)$ values
with these Monte Carlo randomization test statistic values,
we obtain 1000 values.
We sort these test statistics values and determine the ranks of the
$\widetilde R_{PE}^{obs}(r)$ and $\widetilde R_{CS}^{obs}(\tau)$ values
within the respective Monte Carlo randomized test statistic values.
These ranks divided by 1000 (or 1000 minus the rank divided by 1000)
will yield the estimated $p$-values for the left-sided alternative (or the right-sided alternative).
Here we also apply the convex hull correction as in Equation \eqref{eqn:test-stat-adj}
by determining the proportion of other trees outside the convex hull of bald cypresses.
Then we determine the estimated $p$-values for these
convex hull corrected test statistic values as before.
See Figure \ref{fig:MC-p-values-swamp-tree} for the $p$-values
based on the Monte Carlo randomization tests.
Observe that among the Monte Carlo randomized test statistics,
none are significant at .05 level,
but $\widetilde R_{CS}(\tau)$ yields significant results at .10 level for some of the small $\tau$ values.
Notice the discrepancy between the significance in the original test (with the asymptotic normality)
and the Monte Carlo randomization results.

We also analyze the same data in a $2 \times 2$ NNCT
with Dixon's overall test of segregation (\cite{dixon:NNCTEco2002}).
See Table  \ref{tab:NNCT-swamp}
for the corresponding NNCT and the percentages
(observe that the row sum for live trees is 157 instead of 156
due to ties in nearest neighbor (NN) distances).
The cell percentages are relative to the row sums (i.e., number of other or bald cypress trees)
and marginal percentages are relative to the overall sum.
Notice that the table is not suggestive of segregation.
Dixon's overall test statistic is $C_D=0.9735$  ($p=0.6146$)
and Ceyhan's test is $C_N=0.1825$  ($p=0.6692$),
both of which are suggestive of no significant deviation from CSR independence.
So, NNCT-analysis and our relative density approach seem to yield conflicting results
about the spatial interaction of other trees with bald cypresses.
However, NNCT and our relative density approach answer different questions.
More specifically,
NNCT-tests in this example tests the spatial interaction
between the two tree groups,
while the relative density approach only tests the spatial interaction of other trees with bald cypresses,
but not vice versa.
Furthermore,
this situation is an example where relative density is more appropriate
since there is much more other trees compared to bald cypresses.
On the other hand,
the NNCT tests are more appropriate in the cases where the relative abundance
of the two species are similar and cell sizes are larger than 5
(\cite{dixon:NNCTEco2002} and \cite{ceyhan:corrected}).

\begin{table}
\centering
\begin{tabular}{cc}

\begin{tabular}{cc|cc|c}
\multicolumn{2}{c}{}& \multicolumn{2}{c}{NN}& \\
\multicolumn{2}{c}{}&    O.T. &  B.C.   &   sum  \\
\hline
& O.T. &    151  &   6    &   157  \\
\raisebox{1.5ex}[0pt]{base}
& B.C. &    8 &  0    &   8  \\
\hline
&sum     &    159   & 6           &  736  \\
\end{tabular}
&
\begin{tabular}{cc|cc|c}
\multicolumn{2}{c}{}& \multicolumn{2}{c}{NN}& \\
\multicolumn{2}{c}{}&    O.T. &  B.C.   &    \\
\hline
& O.T. &   96 \%  &   4 \%    &   95 \%  \\
& B.C. &    100 \% &  0 \%    &   5 \%  \\
\hline
&     &    96 \%   & 4 \%             &  100 \%  \\
\end{tabular}

\end{tabular}
\caption{ \label{tab:NNCT-swamp}
The NNCT for swamp tree data (left) and the corresponding percentages (right).
O.T. stands for ``other trees" and
B.C. for ``bald cypresses".
}
\end{table}

To find out the level of interaction between the tree species
at different scales (i.e., distances between the trees),
we also present the second-order
analysis of the swamp tree data (\cite{diggle:2003})
using the functions (or some modified version of them) provided
in spatstat package in R (\cite{baddeley:2005}).
We use Ripley's bivariate $L$-functions
which are modified versions of his $K$-functions.
For a rectangular region to remove the bias in estimating $K(t)$,
it is recommended to use $t$ values up to 1/4 of the smaller side length of the rectangle.
So we take the values $t \in [0,12.5]$ in our analysis,
since the rectangular region is $50 \times 55$ m.

\begin{figure}
\centering
\rotatebox{-90}{ \resizebox{2.5 in}{!}{\includegraphics{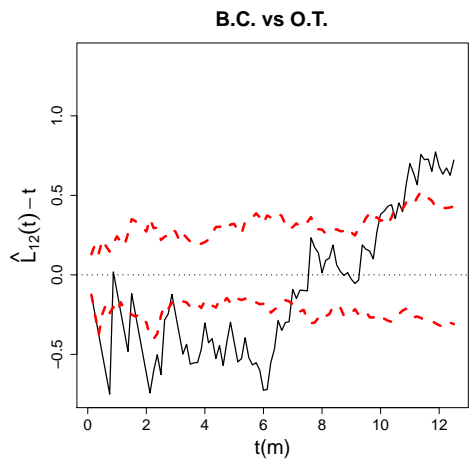} }}
\caption{
\label{fig:swamp-second-order-multi}
Ripley's bivariate $L$-function $\widehat{L}_{12}(t)-t$ for the part of the swamp tree data
we considered.
Wide dashed lines are the upper and lower (pointwise) 95 \% confidence bounds for the
functions based on Monte Carlo simulations under the CSR independence pattern.
B.C. = bald cypresses and O.T. = other trees.}
\end{figure}

Ripley's bivariate $L$-function $L_{ij}(t)$ is symmetric in $i$ and $j$ in theory,
that is, $L_{ij}(t)=L_{ji}(t)$ for all $i,j$.
In practice although edge corrections will render it slightly asymmetric,
i.e., $\widehat{L}_{ij}(t)\not=\widehat{L}_{ji}(t)$ for $i \not= j$.
The corresponding estimates are pretty close in our example,
so we only present one bivariate.
Ripley's bivariate $L$-function for the
bald cypresses and other trees are plotted in Figure \ref{fig:swamp-second-order-multi},
which suggests that
bald cypresses and other trees are significantly segregated for distances
about 0.5 to 7 meters,
and do not significantly deviate from CSR for distances from 7 to 10 meters.

\section{Discussion}
\label{sec:discussion}
In this article,
we consider the asymptotic distribution of the
relative density of two proximity catch digraphs (PCDs),
namely, proportional-edge PCDs and central similarity PCDs
for testing bivariate spatial point patterns of segregation and association
against complete spatial randomness (CSR).
To our knowledge the PCD-based methods are the only graph theoretic
tools for testing spatial point patterns in literature
(\cite{ceyhan:dom-num-NPE-SPL}, \cite{ceyhan:arc-density-PE},
\cite{ceyhan:arc-density-CS}, and \cite{ceyhan:dom-num-NPE-Spat2010}).

We first extend the expansion parameter of the central similarity PCD
which was introduced in \cite{ceyhan:CS-JSM-2003} and \cite{ceyhan:arc-density-CS}
to values higher than one.
We demonstrate that the relative density of the PCDs
can be expressed as $U$-statistic of order 2 (in estimating the arc probability)
and thereby prove the asymptotic normality of the relative density of the PCDs.
For finite samples,
we assess the empirical size and power of the relative density of the PCDs by extensive Monte Carlo simulations.
For the proportional-edge PCDs,
the optimal expansion parameters (in terms of appropriate empirical size and high power)
are about 1.5 under mild segregation and values in $(2,3)$ under moderate to severe segregation;
and about 2 under association.
On the other hand,
for central similarity PCDs,
the optimal parameters are about 7 under segregation,
and about 1 under association.
Furthermore, we have shown that
relative density of central similarity PCDs has better empirical size performance;
and also,
it has higher power against the segregation alternatives.
On the other hand,
relative density of proportional-edge PCDs
has higher power against the association alternatives.

We also compare the asymptotic relative efficiency of the
relative densities of the two PCD families.
Based on Pitman asymptotic efficiency,
we have shown that
in general the relative density of proportional-edge PCDs is
asymptotically more efficient under segregation,
while relative density of central similarity PCDs is more efficient
under association.
However, for the above optimal expansion parameter values
(optimal with respect to empirical size and power),
the asymptotic efficiency and empirical power
analysis yields the same ordering in terms of performance.

Let the two samples of sizes $n$ and $m$ be from classes $\X$ and $\Y$, respectively,
with $\X$ points being used as the vertices of the PCDs and
$\Y$ points being used in the construction of Delaunay triangulation.
The null hypothesis is assumed to be CSR of $\X$ points,
i.e., the uniformness of $\X$ points in the convex hull of $\Y$ points, $C_H(\Y_m)$.
Although we have two classes here, the null pattern is not the CSR independence,
since for finite $m$,
we condition on $m$ and the locations of the $\Y$ points
(assumed to have no more than three co-circular points)
are irrelevant.
That is, the $\Y$ points can result from any pattern
that results in a unique Delaunay triangulation.
The relative density of the two PCD families lend themselves for spatial pattern testing conveniently,
because of the geometry invariance property for uniform data on Delaunay triangles.

For the relative density approach to be appropriate,
the size of $\X$ points (i.e., $n$) should be much larger compared to size of $\Y$ points (i.e., $m$).
This implies that $n$ tends to infinity while $m$ is assumed to be fixed.
That is, the imbalance in the relative abundance of the two classes
should be large for our method to be appropriate.
Such an imbalance usually confounds the
results of other spatial interaction tests.
Furthermore, by construction our method uses only the $\X$ points in $C_H(\Y_m)$
which might cause substantial data (hence information) loss.
To mitigate this, we propose a correction for the proportion of $\X$ points
outside $C_H(\Y_m)$,
because the pattern inside $C_H(\Y_m)$ might not be the same as the pattern
outside $C_H(\Y_m)$.
We suggest a two-stage analysis with our relative density approach:
(i) analysis for $C_H(\Y_m)$,
which provides inference restricted to $\X$ points in $C_H(\Y_m)$,
(ii) overall analysis with convex hull correction
(i.e., for all $\X$ points with respect to $\Y_m$).
We recommend the use of normal approximation if $n \approx 10\,m$ or more,
although Monte Carlo simulations suggest smaller $n$ might also work fine.

There are many possible types of parameterizations for the alternatives.
The particular parametrization of the alternatives in Equation \eqref{eqn:eps-alt-1}
is chosen so that the distribution of the relative density under the alternatives
would be geometry invariant also
(i.e., independent of the geometry of the support triangles).
The more natural alternatives (i.e.,
the alternatives that are more likely to be found in practice)
can be similar to or might be approximated by our parametrization.
Because under a segregation alternative,
the $\X$ points will tend to be further away from $\Y$ points
and under an association alternative $\X$ points will tend to cluster around the $\Y$ points.
Such patterns can be detected by the test
statistics based on the relative density,
since under segregation (whether it is parametrized as in Section \ref{sec:alternatives}
or not) we expect them to be larger,
and under association (regardless of the parametrization)
they tend to be smaller.

\section*{Acknowledgments}
Supported by TUBITAK Kariyer Project Grant 107T647.
Most of the Monte Carlo simulations presented in this article
were executed at Ko\c{c} University High Performance Computing Laboratory.


\section*{APPENDIX}

\section*{Appendix 1: Derivation of $\mu_{_{CS}}(\tau)$ and $\nu_{_{CS}}(\tau)$ for $\tau > 1$}
\begin{figure} []
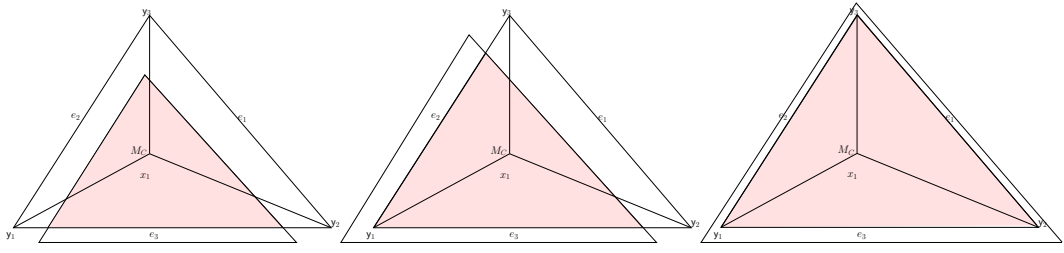

   \centering
   \scalebox{.2}{\input{n_cscase1.pstex_t}}
   \scalebox{.2}{\input{n_cscase2.pstex_t}}
   \scalebox{.2}{\input{n_cscase3.pstex_t}}
   \caption{The prototypes of the three cases of $\NCS(x_1,\tau)$ for $x_1 \in T(\y_1,M_3,M_{C})$ with $\tau=2.5$.}
\label{fig:NCS-Cases-1}
\end{figure}

Let $M_C$ be the center of mass of the standard equilateral triangle $T_e$.
By symmetry
$\mu_{_{CS}}(\tau)=P\bigl(X_2 \in N_{CS}(X_1,\tau)\bigr)=6\,P\bigl(X_2 \in N_{CS}(X_1,\tau),\, X_1 \in T(\y_1,M_3,M_C)\bigr)$.
To calculate this mean,
we need to find the possible types of $N_{CS}(x_1,\tau)$ for $\tau > 1$.
There are three cases regarding $N_{CS}(x_1,\tau)$.
See Figure \ref{fig:NCS-Cases-1} for the prototypes of these three cases of $N_{CS}(x_1,\tau)$ for $x_1=(u_1,v_1) \in \TY$.

\begin{figure} []
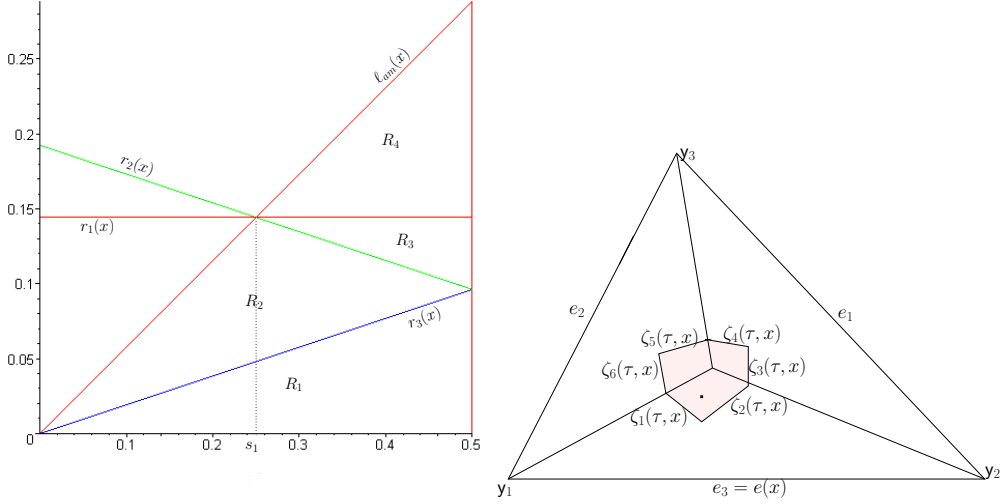

   \centering
   \scalebox{.3}{\input{ncsregionstgt1.pstex_t}}
   \scalebox{.3}{\input{N_CSGam4.pstex_t}}
   \caption{
\label{fig:regions-for-NCS-G1}
The regions corresponding to the prototypes of the four cases (left)
and the prototype of the one case of $\G_1^{^{CS}}(x_1,\tau)$ for $x_1 \in T(\y_1,M_3,M_{C})$ (right) with $\tau=2.5$. }
\end{figure}

Each case $j$, corresponds to region $R_j$ in Figure \ref{fig:regions-for-NCS-G1},
and the bounding lines are $r_6(x)$, $\ell_{am}(x)$, and $r_7(x)$ 
with
$\ell_{am}(x) =x/\sqrt{3},\;\;\; r_6(x)={\frac{\sqrt{3}u_1}{1+2\,\tau}},\;\;\; r_7(x)=-{\frac{\sqrt{3} \left( -1+u_1 \right) }{1+2\,\tau}}$ and
$s_2=\frac{3}{2\, \left( 2+\tau \right)}.$

The explicit forms of $R_j$, $j=1,2,3$ are as follows:
\begin{align*}
R_1&=\bigl\{(x,y)\in  [0,1/2]\times [0,r_6(x)]\bigr\},\\
R_2&=\bigl\{(x,y)\in [0,s_2] \times [r_6(x),\ell_{am}(x)]\cup [s_2,1/2] \times [r_6(x),r_7(x)]\bigr\},\\
R_3&=\bigl\{(x,y)\in [s_2,1/2] \times [r_7(x),\ell_{am}(x)]\bigr\}.
\end{align*}

$$P\bigl(X_2 \in N_{CS}(X_1,\tau),\; X_1 \in T_s \bigr)=
\sum_{j=1}^3 P\bigl( X_2 \in N_{CS}(X_1,\tau),\; X_1 \in R_j \bigr).$$
For $x_1 \in R_1$,
\begin{eqnarray*}
P\bigl(X_2 \in N_{CS}(X_1,\tau),\; X_1 \in R_1\bigr)=
\int_0^{1/2}\int_0^{r_6(x)} \frac{A\left( N_{CS}(X_1,\tau) \right)}{A(\TY)^2}dydx
=  \frac{1}{12\,\left( 1+2\,\tau \right)},
\end{eqnarray*}
where $A\left( N_{CS}(X_1,\tau) \right)=\frac{1}{\sqrt{3}} \left( 1+2\,\tau \right)^2v_1^2$.

For $x_1 \in R_2$,
$$
P\bigl( X_2 \in N_{CS}(X_1,\tau),\; X_1 \in R_2\bigr) =
\left(\int_{0}^{s_2}\int_{r_6(x)}^{\ell_{am}(x)} +
\int_{s_2}^{1/2}\int_{r_6(x)}^{r_7(x)} \right)
 \frac{A\left( N_{CS}(X_1,\tau) \right)^2}{A(\TY)^3}dydx =
{\frac{\tau-1}{2 \left( 1+2\,\tau \right)  \left( 2+\tau\right) }}.$$
where $A\left( N_{CS}(X_1,\tau) \right)=\frac{\sqrt{3}}{12}\, \left( v_1+2\,\tau\,v_1+\sqrt{3}u_1 \right)^2$.

For $x_1 \in R_3$,
\begin{eqnarray*}
P\bigl(X_2 \in N_{CS}(X_1,\tau),\; X_1 \in R_3\bigr)=
\int_{s_2}^{1/2}\int_{r_7(x)}^{\ell_{am}(x)} \frac{A\left( N_{CS}(X_1,\tau) \right)}{A(\TY)^2}dydx
= {\frac{ \left( \tau-1 \right)^2}{3 \left( 1+2\,\tau \right) \left( 2+\tau \right) }}.
\end{eqnarray*}
where $A\left( N_{CS}(X_1,\tau) \right)=\frac{\sqrt{3}}{4}$.

So $P\bigl( X_2 \in N_{CS}(X_1,\tau)\bigr)=
{\frac{\tau\, \left( 4\,\tau-1 \right) }{2\,\left( 1+2\,\tau \right)  \left( 2+\tau \right) }}.$

Next, we find the asymptotic variance term.
Let
\begin{gather*}
P^{2N}_{_{CS}}:=P\bigl(\{X_2,X_3\} \subset N_{CS}(X_1,\tau)\bigr),\;\;\;P^{2G}_{_{CS}}:=
P\bigl(\{X_2,X_3\} \subset \G_1^{^{CS}}(X_1,\tau)\bigr)\; \; \text{ and } \;\;\\
P^M_{_{CS}}:=P\bigl(X_2 \in N_{CS}(X_1,\tau), X_3 \in \G_1^{^{CS}}(X_1,\tau)\bigr).
\end{gather*}
where $\G_1^{^{CS}}(x,\tau)$ is the \emph{$\G_1$-region} of $x$
based on $N_{CS}(\cdot,\tau)$ and defined as
$\G_1^{^{CS}}(x,\tau):=\{y \in \TY:\;x \subset N_{CS}(y,\tau)\}$.
See \cite{ceyhan:TR-dom-num-NPE-spatial} for more detail.

Then $\Cov[h_{12},h_{13}]=\E[h_{12}\,h_{13}]-\E[h_{12}]\E[h_{13}]$ where
\begin{eqnarray*}
\E[h_{12}\,h_{13}] &=&P\bigl(\{X_2,X_3\} \subset N_{CS}(X_1,\tau)\bigr)+2\,P\bigl(X_2 \in N_{CS}(X_1,\tau), X_3 \in \G_1^{^{CS}}(X_1,\tau)\bigr)\\
& & +P\bigl(\{X_2,X_3\} \subset \G_1^{^{CS}}(X_1,\tau)\bigr) = P^{2N}_{_{CS}}+2\,P^M_{_{CS}}+P^{2G}_{_{CS}}.
\end{eqnarray*}
Hence $\nu_{_{CS}}(\tau)=\Cov[h_{12},h_{13}]  = \bigl(P^{2N}_{_{CS}}+2\,P^M_{_{CS}}+P^{2G}_{_{CS}}\bigr)-[2\,\mu_{_{CS}}(\tau)]^2.$

To find the covariance,
we need to find the possible types of $\G_1^{^{CS}}(x_1,\tau)$
and $N_{CS}(x_1,\tau)$ for $\tau \in (1,\infty)$.
There are three cases regarding $N_{CS}(x_1,\tau)$ and
one case for $\G_1^{^{CS}}(x_1,\tau)$.
See Figure \ref{fig:regions-for-NCS-G1} for the prototype of this one case of
$\G_1^{^{CS}}(x_1,\tau)$ for $x_1=(u_1,v_1) \in \TY$,
the explicit forms of $\zeta_j(\tau,x)$ are
\begin{align*}
\zeta_1(\tau,x)&=\frac{\left( \sqrt{3}\,v_1+3\,u_1-3\,x \right)}{\sqrt{3}\,(1+2\,\tau)},\;\; \zeta_2(\tau,x)=-\frac{\left(-\sqrt{3}\,v_1+3\,u_1-3\,x\right)}{\sqrt{3}\,(1+2\,\tau)},\\
\zeta_3(\tau,x)&=\frac{\left(3\,u_1+3\,\tau-3\,\tau\,x-3\,x-\sqrt{3}\,v_1\right)}{\sqrt{3}\,(-1+\tau)},\;\;\;  \zeta_4(\tau,x)=-\frac{-\sqrt{3}\,\tau+\sqrt{3}\,\tau\,x-2\,v_1}{2+\tau},\\
\zeta_5(\tau,x)&=\frac{\sqrt{3}\,\tau\,x+2\,v_1}{2+\tau},\;\;\; \zeta_6(\tau,x)= \frac{\left(-3\,x-3\,\tau\,x+3\,u_1+\sqrt{3}\,v_1 \right)}{\sqrt{3}\,(1-\tau)},\;\;\;\zeta_7(\tau,x)=\frac{v_1}{1-\tau}.
\end{align*}

By symmetry,  $P^{2N}_{_{CS}}=6\,P\bigl( \{X_2,X_3\} \subset N_{CS}(X_1,\tau),\; X_1 \in T_s \bigr),$
and
$$P\bigl(\{X_2,X_3\} \subset N_{CS}(X_1,\tau),\; X_1 \in T_s \bigr)=
\sum_{j=1}^3 P\bigl(\{X_2,X_3\} \subset N_{CS}(X_1,\tau),\; X_1 \in R_j \bigr).$$
The limits of integration are as in
$P\bigl(X_2 \in N_{CS}(X_1,\tau),\, X_1 \in T_s\bigr)$ with the integrand being
$\frac{A(N_{CS}(x_1,\tau))^2}{A(\TY)^3}$.
Hence,
$P^{2N}_{_{CS}}={\frac{10\,\tau^2-9\,\tau+2}{ 5\,\left( 2\,\tau+1 \right) \left( \tau+2 \right) }}.$

Next, by symmetry, $P^{2G}_{_{CS}}=6\,P\bigl(\{X_2,X_3\} \subset \G_1^{^{CS}}(X_1,\tau),\; X_1 \in T_s \bigr),$
and
$$P\bigl( \{X_2,X_3\} \subset \G_1^{^{CS}}(X_1,\tau),\; X_1 \in T_s \bigr)=
\int_0^{1/2}\int_0^{\ell_{am}(x)} \frac{A\left( \G_1^{^{CS}}(x_1,\tau) \right)^2}{A(\TY)^3}dydx
 =  {\frac{\tau^2 \left( 10\,\tau^2-5\,\tau+1 \right) }{15\,\left( 2\,\tau+1 \right)^2 \left( \tau+2 \right)^2}},$$
where
$A\left( \G_1^{^{CS}}(x_1,\tau) \right)=
{\frac{\sqrt{3} \left( 3\,u_1+\sqrt{3}v_1+\tau-1-3\,v_1^2-3\,u_1^2 \right) \tau}
{2\,\left( 2\,\tau+1 \right)  \left( \tau+2 \right) }}$.

So $P^{2G}_{_{CS}}=
{\frac{2\,\tau^2 \left( 10\,\tau^2-5\,\tau+1 \right) }{ 5\,\left( 2\,\tau+1 \right)^2 \left( \tau+2 \right)^2}}.$

Furthermore, by symmetry,
$P^M_{_{CS}}=6\,P\bigl( X_2 \in N_{CS}(X_1,\tau),\;X_3\in \G_1^{^{CS}}(X_1,\tau),\; X_1 \in T_s\bigr),$
and
{\small
$$
P\bigl( X_2 \in N_{CS}(X_1,\tau),\;X_3\in \G_1^{^{CS}}(X_1,\tau),\; X_1 \in T_s \bigr)=
\sum_{j=1}^3 P\left( X_2 \in N_{CS}(X_1,\tau),\;X_3\in \G_1^{^{CS}}(X_1,\tau),\; X_1 \in R_j \right).
$$
}
where $ P\bigl( X_2 \in N_{CS}(X_1,\tau),\;X_3\in \G_1^{^{CS}}(X_1,\tau),\; X_1 \in R_j \bigr)$
can be calculated with the same region of integration with integrand being replaced by
$\frac{A(N_{CS}(x_1,\tau))\,A\left( \G_1^{^{CS}}(x_1,\tau) \right)}{A(\TY)^3}$.

Then
$P^M_{_{CS}}={\frac{\tau^2 \left( 10-54\,\tau-99\,\tau^2+388\,\tau^3+1062\,\tau^4+720\,\tau^5+160\,\tau^6 \right) }
{10\,\left( 2\,\tau+1 \right)^4 \left( \tau+2 \right)^4}}.$

Hence
$$\E[h_{12}\,h_{13}]={\frac{2\,(160\,\tau^8-265\,\tau^3-135\,\tau^4+751\,\tau^5-47\,\tau^2+1373\,\tau^6+804\,\tau^7+24\,\tau+8)}
{ 5\,\left( 2\,\tau+1 \right)^4 \left( \tau+2 \right)^4}}.$$

Therefore,
$$\nu_{_{CS}}(\tau)=
{\frac{168\,\tau^7+886\,\tau^6+1122\,\tau^5+45\,\tau^4-470\,\tau^3-114\,\tau^2+48\,\tau+16}
{ 5\,\left( 2\,\tau+1 \right)^4 \left( \tau+2 \right)^4}}.$$

For $\tau=\infty$,
it is trivial to see that $\nu_{_{CS}}(\tau)=0$.

\section*{Appendix 2: Derivation of $\mu^S_{_{CS}}(\tau,\ve)$ for $\tau>1$}
\label{sec:derivation-patau-eps}
\begin{figure} [ht]
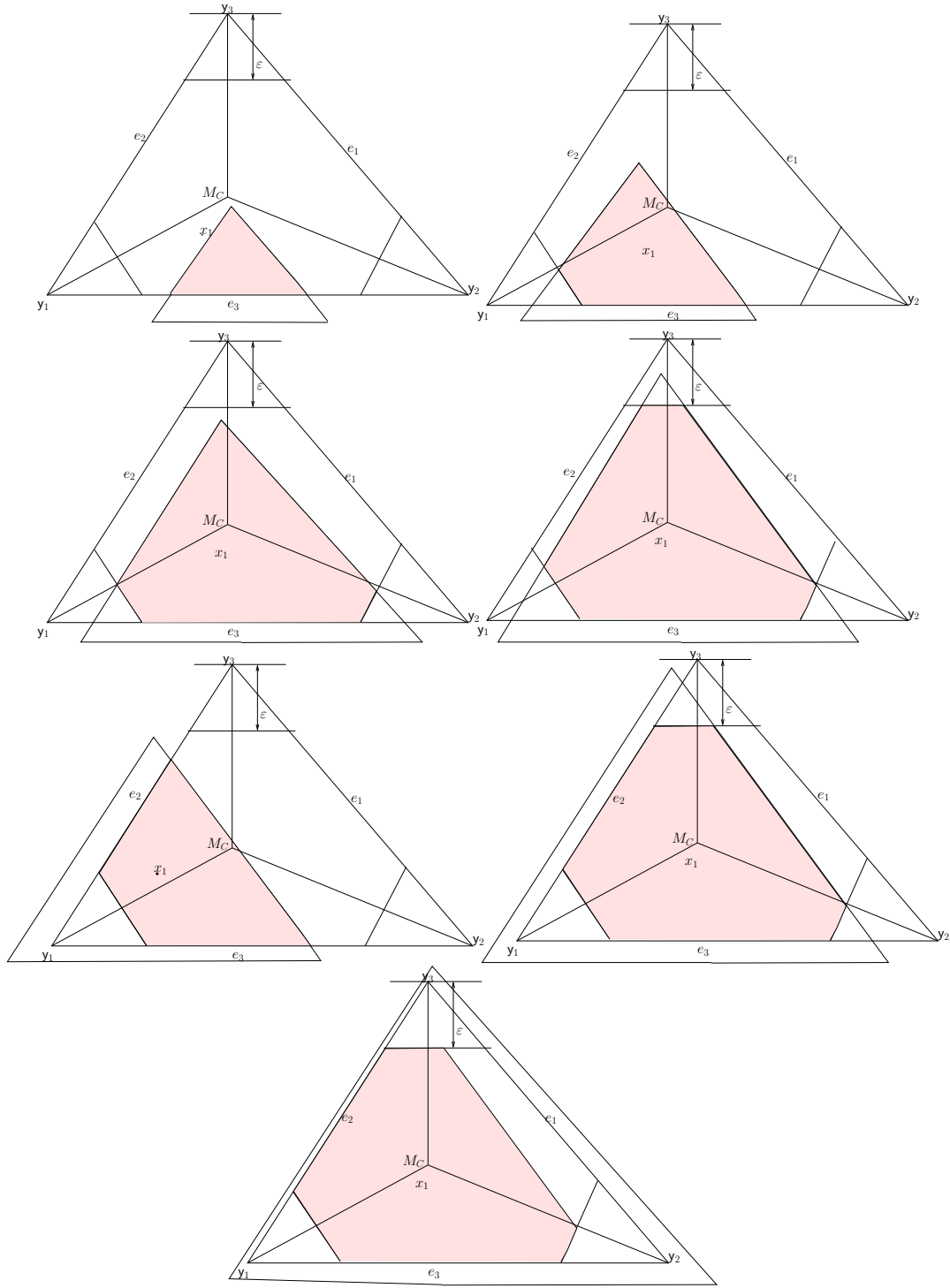

    \centering
   \scalebox{.3}{\input{n_cseps5.pstex_t}}
   \scalebox{.3}{\input{n_cseps6.pstex_t}}
   \scalebox{.3}{\input{n_cseps7.pstex_t}}
   \scalebox{.3}{\input{n_cseps8.pstex_t}}
   \scalebox{.3}{\input{n_cseps9.pstex_t}}
   \scalebox{.3}{\input{n_cseps10.pstex_t}}
   \scalebox{.3}{\input{n_cseps11.pstex_t}}
    \caption{The seven possible cases of $N_{CS}(\cdot,\tau)(x)$ for four distinct $x \in R_{CM}(e_3)$ (shaded regions). }
\label{fig:NCS-Eps1-7}
\end{figure}

We pick the interval $\ve \in \bigl[ 0,\sqrt{3}/5 \bigr)$
for a demonstrative example in derivation of $\mu^S_{_{CS}}(\tau,\ve)$.
For $\tau \in \bigl[ 1-\sqrt{3}\,\ve,1 \bigr)$,
there are seven cases to consider for the form of $N_{CS}(\cdot,\tau)(x_1,\ve)$.
See Figure \ref{fig:NCS-Eps1-7} for the prototypes of these seven cases of $N_{CS}(\cdot,\tau)(x_1,\ve)$.

\begin{figure} [ht]
    \centering
   \scalebox{.6}{\input{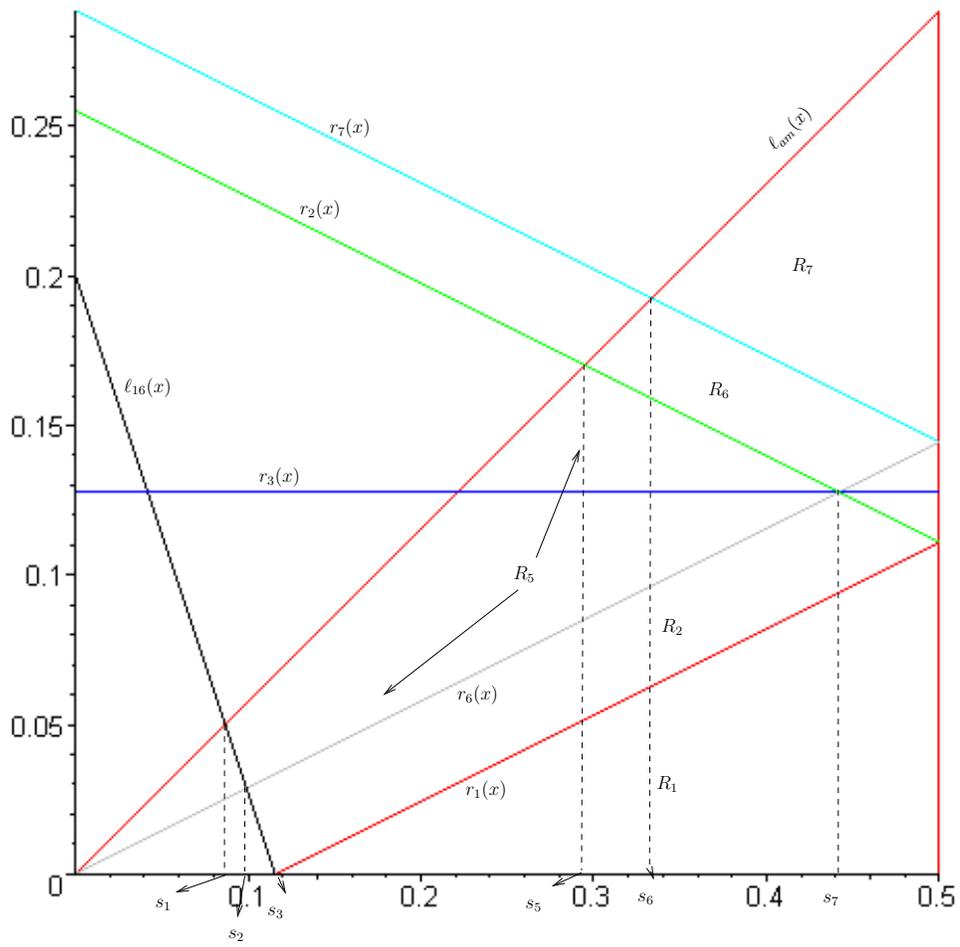}}
    \caption{The regions corresponding to the prototypes of the seven cases shown in Figure \ref{fig:NCS-Eps1-7}.}
    \label{regions-for-NCS-eps}
\end{figure}

Each case $j$, corresponds to the region $R_j$ in Figure \ref{regions-for-NCS-eps},
where
$
q_1(x)=-\sqrt{3}\,x+2\,\ve,\;\;\;
q_2(x)=\frac{-\sqrt{3}\,x+2\,\ve}{1-4\,\tau},\;\;\;
q_3(x)=\frac{\sqrt{3}\,x+2\,\ve-\sqrt{3}}{1-4\,\tau},\;\;\;
q_4(x)=\frac{\sqrt{3}-2\,\ve}{2\,(1+2\,\tau)}$
and
$s_1=\sqrt{3}\,\ve/2,\;\;\;
s_2=2\,\sqrt{3}\,\ve/3,\;\;\;
s_3=\frac{\sqrt{3}\,(\sqrt{3}-2\,\ve)}{2\,(1+2\,\tau)}.$

The explicit forms of $R_j$, $j=1,2,\ldots,7$ are as follows:
\begin{align*}
R_1&=\bigl\{(x,y)\in  [s_3,1/2]\times [0,r_1(x)]\bigr\},\\
R_2&=\bigl\{(x,y)\in  [s_2,s_3]\times [\ell_{16}(x),r_6(x)] \cup [s_3,s_7]\times [r_1(x),r_6(x)] \cup [s_7,1/2]\times [r_1(x),r_2(x)]\bigr\},\\
R_3&=\bigl\{(x,y)\in  [s_7,1/2]\times [r_2(x),r_3(x)]\bigr\},\\
R_4&=\bigl\{(x,y)\in  [s_7,1/2]\times [r_3(x),r_6(x)]\bigr\},\\
R_5&=\bigl\{(x,y)\in  [s_1,s_2]\times [\ell_{16}(x),\ell_{am}(x)] \cup [s_2,s_5]\times [r_6(x),\ell_{am}(x)] \cup [s_5,s_7]\times [r_6(x),r_2(x)]\bigr\},\\
R_6&=\bigl\{(x,y)\in  [s_5,s_6]\times [r_2(x),\ell_{am}(x)] \cup [s_6,s_7]\times [r_2(x),r_7(x)] \cup [s_7,1/2]\times [r_6(x),r_7(x)]\bigr\},\\
R_7&=\bigl\{(x,y)\in [s_6,1/2] \times [r_7(x),\ell_{am}(x)]\bigr\}.
\end{align*}

By symmetry,
$P\bigl(X_2 \in N_{CS}(X_1,\tau,\ve)\bigr)=6\,P\bigl( X_2 \in N_{CS}(X_1,\tau,\ve),\; X_1 \in T_s \bigr).$

For $x_1 \in R_1$,
{\small
$$
P\bigl(X_2 \in N_{CS}(X_1,\tau,\ve),\; X_1 \in R_1\bigr)=
\int_{s_3}^{1/2}\int_{0}^{r_1(x)}\frac{A( N_{CS}(X_1,\tau,\ve) )^2}{A(\TY)^2}dydx =
{\frac{9+256\,\ve^4-48\,\ve\,\sqrt{3}+288\,\ve^2-256\,\ve^3\sqrt{3}}
{ 108\,\left( 1+2\,\tau \right)  \left( 4\,\ve^2-1 \right)^2}},
$$
}
where $A\left( N_{CS}(X_1,\tau,\ve) \right)=\frac{1}{\sqrt{3}} \left( 1+2\,\tau \right)^2v_1^2 $.

For $x_1 \in R_2$,
{\small
\begin{multline*}
P\bigl(X_2 \in N_{CS}(X_1,\tau,\ve),\; X_1 \in R_2\bigr)=
\left(\int_{s_2}^{s_3}\int_{\ell_{16}(x)}^{r_6(x)} + \int_{s_3}^{s_7}\int_{r_1(x)}^{r_6(x)}
+\int_{s_7}^{1/2}\int_{r_1(x)}^{r_2(x)} \right)
\frac{A( N_{CS}(X_1,\tau,\ve) )^2}{A(\TY)^2}dydx =\\
\Big[4\,\ve\, \Big(3\,\sqrt{3}-27\,\ve+9\,\sqrt{3}\tau-16\,\ve^3+24\,\ve^2\sqrt{3}-24\,\ve^3\tau^3-
48\,\ve^3\tau-54\,\ve^3 \tau^2-27\,\ve\,\tau^3-81\,\tau\,\ve-81\,\ve\,\tau^2+\\
9\,\sqrt{3}\tau^2+3\,\sqrt{3}\tau^3+72\,\ve^2\sqrt{3}\tau+72\,\ve^2\sqrt{3}\tau^2+
24\,\ve^2\sqrt{3}\tau^3 \Big) \Big]\Big/\Big[27\,\left( 1+\tau \right)^3 \left( 1+2\,\tau \right)  \left( 4\,\ve^2-1 \right)^2\Big],
\end{multline*}
}
where
$A\left( N_{CS}(X_1,\tau,\ve) \right)=
\frac{\sqrt{3}}{12} \left( -2\,\ve+v_1+2\,\tau\,v_1+\sqrt{3}u_1 \right)
\left(2\,\ve-\sqrt{3}u_1+3\,v_1+6\,\tau\,v_1 \right)$.

For $x_1 \in R_3$,
{\small
$$
P\bigl(X_2 \in N_{CS}(X_1,\tau,\ve),\; X_1 \in R_3\bigr)=
\int_{s_7}^{1/2}\int_{r_2(x)}^{r_3(x)} \frac{A( N_{CS}(X_1,\tau,\ve) )^2}{A(\TY)^2}dydx =
{\frac{2\,\ve^2 \left(9+18\,\ve^2-16\,\ve\,\sqrt{3} \right) }
{27\,\left( 1+2\,\tau \right)  \left( 4\,\ve^2-1 \right)^2}},
$$
}
where
{\small
$$
A\left( N_{CS}(X_1,\tau,\ve) \right)=
\frac{\sqrt{3}}{12} \left( 4\,\ve\,\sqrt{3}-3+4\,\sqrt{3}\tau\,v_1+6\,u_1+2\,\sqrt{3}v_1-8\,\ve^2-16\,
\ve\,\tau\,v_1+2\,v_1^2-8\,\ve\,v_1+8\,\tau\,v_1^2+8\,\tau^2v_1^2-6\,u_1^2 \right).
$$
}

For $x_1 \in R_4$,
{\small
$$
P\bigl(X_2 \in N_{CS}(X_1,\tau,\ve),\; X_1 \in R_4\bigr)=
\int_{s_7}^{1/2}\int_{r_3(x)}^{r_6(x)} \frac{A( N_{CS}(X_1,\tau,\ve) )^2}{A(\TY)^2}dydx =
{\frac{2\,\ve^2 \left(9-8\,\ve\,\sqrt{3}-8\,\ve^2 \right) }
{ 27\,\left( 1+2\,\tau \right)  \left( 4\,\ve^2-1 \right)^2}},
$$
}
where
{\small
$$
A\left( N_{CS}(X_1,\tau,\ve) \right)=
\frac{\sqrt{3}}{6}\left( 4\,\ve\,\sqrt{3}-3+6\,\sqrt{3}\tau\,v_1+3\,u_1+3\,\sqrt{3}v_1-
6\,\ve^2-16\,\ve\,\tau\,v_1-v_1^2-8\,\ve\,v_1-4\,\tau\,v_1^2-4\,\tau^2v_1^2-3\,u_1^2 \right).
$$
}

For $x_1 \in R_5$,
{\small
\begin{multline*}
P\bigl(X_2 \in N_{CS}(X_1,\tau,\ve),\; X_1 \in R_5\bigr)=
\left(\int_{s_1}^{s_2}\int_{\ell_{16}(x)}^{\ell_{am}(x)} + \int_{s_2}^{s_5}\int_{r_6(x)}^{\ell_{am}(x)}
+\int_{s_5}^{s_7}\int_{r_6(x)}^{r_2(x)} \right)
\frac{A( N_{CS}(X_1,\tau,\ve) )^2}{A(\TY)^2}dydx =\\
\Big[27\,\tau^4-4\,\ve^4\tau^7-46\,\ve^4 \tau^6-130\,\ve^4\tau^5-72\,\ve\,\sqrt{3}\tau^4+144\,\ve^2\tau^4-
74\,\ve^4 \tau^4+288\,\ve^2\tau^3+54\,\tau^3+142\,\ve^4\tau^3-144\,\ve\,\sqrt{3}\tau^3+112\,\tau^2 \ve^4-\\
54\,\tau+144\,\tau\,\ve\,\sqrt{3}-288\,\ve^2\tau-27+72\,\ve\,\sqrt{3}-144\,\ve^2\Big]
\Big/
\Big[ 54\,\left( 4\,\ve^2-1 \right)^2 \left( 1+2\,\tau \right)  \left( 2+\tau \right)  \left( 1+\tau \right)^3\Big],
\end{multline*}
}
where
$A\left( N_{CS}(X_1,\tau,\ve) \right)=
\frac{\sqrt{3}}{12}\left( -2\,\ve+v_1+2\,\tau\,v_1+\sqrt{3}u_1 \right)
\left( \sqrt{3}u_1+v_1+2\,\tau\,v_1+2\,\ve \right)$.

For $x_1 \in R_6$,
{\small
\begin{multline*}
P\bigl(X_2 \in N_{CS}(X_1,\tau,\ve),\; X_1 \in R_6\bigr)=
\left( \int_{s_5}^{s_6}\int_{r_2(x)}^{\ell_{am}(x)} + \int_{s_6}^{s_7}\int_{r_2(x)}^{r_7(x)} +
\int_{s_7}^{1/2}\int_{r_6(x)}^{r_7(x)} \right)
\frac{A( N_{CS}(X_1,\tau,\ve) )^2}{A(\TY)^2}dydx =\\
\Big[ 4\,\left( 10\,\ve^3\tau+9\,\sqrt{3}\tau-27\,\tau\,\ve-8\,\ve^2\sqrt{3}\tau-10\,\ve^3-
9\,\sqrt{3}+27\,\ve+8\,\ve^2\sqrt{3} \right)\ve\Big]
\Big/
\Big[27\,\left( 2+\tau \right)  \left( 1+2\,\tau \right)  \left( 4\,\ve^2-1 \right)^2\Big],
\end{multline*}
}
where
{\small
\begin{multline*}
A\left( N_{CS}(X_1,\tau,\ve) \right)=
\frac{\sqrt{3}}{12} \Big(4\,\sqrt{3}v_1-8\,\ve\,v_1+8\,\sqrt{3}\tau\,v_1-16\,\ve\,\tau\,v_1+12\,u_1-
8\,\ve\,\sqrt{3}u_1-6+8\,\ve\,\sqrt{3}-12\,\ve^2-v_1^2-4\,\tau\,v_1^2-\\
2\,\sqrt{3}u_1\,v_1-4\,\sqrt{3}u_1\,\tau\,v_1-3\,u_1^2-4\,\tau^2v_1^2 \Big).
\end{multline*}
}

For $x_1 \in R_7$,
{\small
$$
P\bigl(X_2 \in N_{CS}(X_1,\tau,\ve),\; X_1 \in R_7\bigr)=
\int_{s_7}^{1/2}\int_{r_3(x)}^{r_6(x)} \frac{A( N_{CS}(X_1,\tau,\ve) )^2}{A(\TY)^2}dydx =
{\frac{1-2\,\tau+\tau^2}{ 3\,\left( 1-4\,\ve^2 \right)\left( 1+2\,\tau \right)  \left( 2+\tau \right) }},
$$
}
where
$A\left( N_{CS}(X_1,\tau,\ve) \right)=
\frac{\sqrt{3}}{4} \left( 2\,\ve+1 \right)  \left( 1- 2\,\ve \right)$.

So
{\small
\begin{multline*}
P\bigl( X_2 \in N_{CS}(X_1,\tau,\ve)\bigr)=
\Big[\Big(36\,\tau^3-8\,\ve^4\tau^5-84\,\ve^4\tau^4-144\,\ve^2\tau^3-176\,\ve^4\tau^3-
192\,\tau^2\sqrt{3}\ve^3+63\,\tau^2+140\,\tau^2\ve^4-144\,\ve^2\tau^2-\\
384\,\sqrt{3}\tau\,\ve^3+18\,\tau+512\,\ve^4\tau+
144\,\ve^2\tau-9-192\,\ve^3\sqrt{3}+144\,\ve^2+240\,\ve^4 \Big) \tau\Big]
\Big/
\Big[ 18\,\left( 1+\tau \right)^2 \left( 2+\tau \right)  \left( 1+2\,\tau \right)  \left( 4\,\ve^2-1 \right)^2\Big].
\end{multline*}
}

\section*{Appendix 3: The Mean $\mu^S_{_{CS}}(\tau,\ve)$ under the Segregation and Association Alternatives}
\label{sec:pi_a-tau-epsilon}
Derivation of $\mu^S_{_{CS}}(\tau,\ve)$
involves detailed geometric calculations  and partitioning of the
space of $(\tau,\ve,x_1)$ for $\tau \in [0,1]$,
$\ve \in \bigl[ 0,\sqrt{3}/3 \bigr)$, and
$x_1 \in T_s=T(\y_1,M_3,M_C)$.
A demonstrative calculation is given in Appendix 2 for
$\ve \in [ 0,\sqrt{3}/5 )$ and $\tau \in [1-\sqrt{3}\,\ve,1]$.

\subsubsection*{$\mu^S_{_{CS}}(\tau,\ve)$ under the Segregation Alternatives}
\label{sec:arc-prob-seg-CS}
Under segregation, we compute $\mu^S_{_{CS}}(\tau,\ve)$ explicitly.
For $\ve \in [ 0,\sqrt{3}/5 )$,
$$\mu^S_{_{CS}}(\tau,\ve)=\sum_{j=1}^2 \varpi_{1,j}(\tau,\ve)\,\I(\tau \in \mI_j)$$
where
\begin{align*}
\varpi_{1,1}(\tau,\ve)&=\frac{(20\,\ve^4\tau-36\,\ve^4-3\,\tau+3)\tau^2}{18\,(1-\tau)(2\,\ve+1)^2(2\,\ve-1)^2},\\
\varpi_{1,2}(\tau,\ve)&=-\Bigl[\Bigl(80\,\ve^4\tau^3+16\,\ve^4\tau^2-108\,\ve^4\tau-144\,\ve^4+192\,\sqrt{3}\ve^3+288\,\ve^2\tau+64\,\sqrt{3}\,\ve\,\tau^2+4\,\tau^3\\
&-288\,\ve^2-128\,\sqrt{3}\,\ve\,\tau-60\,\tau^2+64\,\sqrt{3}\,\ve+45\,\tau-16 \Bigr)\tau \Bigr]
\Big/
\bigl[ 18\,(2\,\tau+1)^2(2\,\ve-1)^2(2\,\ve+1)^2 \bigr],
\end{align*}
with the corresponding intervals
$\mI_1=\bigl[1,1-\sqrt{3}\,\ve\bigr)$
and
$\mI_2=\bigl[1-\sqrt{3}\,\ve,1\bigr)$.

For $\ve \in \bigl[ \sqrt{3}/5,\sqrt{3}/4 \bigr)$,
$$\mu^S_{_{CS}}(\tau,\ve)=\sum_{j=1}^3 \varpi_{2,j}(\tau,\ve)\,\I(\tau \in \mI_j)$$
where $\varpi_{2,j}(\tau,\ve)=\varpi_{1,j}(\tau,\ve)$ for $j=1,2$, and for $j=3$,
{\small
\begin{multline*}
\varpi_{2,3}(\tau,\ve)=
\Bigl[\Bigl(-96\,\ve^4\tau^4-1584\,\ve^4\tau^3-3688\,\ve^4\tau^2+512\,\sqrt{3}\ve^3\tau^3-
3636\,\ve^4\tau+2304\,\sqrt{3}\ve^3\tau^2 -1152\,\ve^4+3840\,\sqrt{3}\ve^3\tau-576\,\ve^2\tau^2+\\
128\,\sqrt{3}\,\ve\,\tau^3+8\,\tau^4+1536\,\sqrt{3}\ve^3-4320\,\ve^2\tau-320\,\sqrt{3}\,\ve\,\tau^2
-124\,\tau^3-2304\,\ve^2+640\,\sqrt{3}\,\ve\,\tau+150\,\tau^2+512\,\sqrt{3}\,\ve-\\
77\,\tau-128\Bigr)\tau\Bigr]\Big/\Bigl[18\,(2\,\tau+1)^2(2\,\ve-1)^2(2\,\ve+1)^2(1-2\,\tau)\Bigr],
\end{multline*}
}
with the corresponding intervals
$\mI_1=\bigl[1,1-\sqrt{3}\,\ve\bigr)$,
$\mI_2=\bigl[1-\sqrt{3}\,\ve,\sqrt{3}/(2\,\ve)-3/2\bigr)$,
and
$\mI_3=\bigl[\sqrt{3}/(2\,\ve)-3/2,1\bigr)$.

For $\ve \in \bigl[\sqrt{3}/4,2\,\sqrt{3}/7 \bigr)$,
$\mu^S_{_{CS}}(\tau,\ve)=\sum_{j=1}^5 \varpi_{3,j}(\tau,\ve)\,\I(\tau \in \mI_j)$
where
{\small
\begin{multline*}
\varpi_{3,1}(\tau,\ve)=
\Bigl[(984\,\ve^4\tau^3-3452\,\ve^4\tau^2-1024\,\sqrt{3}\ve^3\tau^3+4992\,\ve^4\tau+3584\,\sqrt{3}\ve^3\tau^2+
1152\,\ve^2\tau^3-2268\,\ve^4 -5120\,\sqrt{3}\ve^3\tau-4032\,\ve^2\tau^2-\\
192\,\sqrt{3}\,\ve\,\tau^3+2304\,\sqrt{3}\ve^3+5760\,\ve^2\tau+672\,\sqrt{3}\,\ve\,\tau^2+42\,\tau^3
-2592\,\ve^2-960\,\sqrt{3}\,\ve\,\tau-141\,\tau^2+432\,\sqrt{3}\,\ve+192\,\tau-84)\tau^2\Bigr]
\Big/\\
\Bigl[32\,(\tau-1)^2\left( 3\,\ve-\sqrt{3} \right)^4(2\,\tau-1)\Bigr],
\end{multline*}
}
{\small
\begin{multline*}
\varpi_{3,2}(\tau,\ve)=
\Bigl[(3936\,\ve^4\tau^6-9872\,\ve^4\tau^5-4096\,\sqrt{3}\ve^3\tau^6+7144\,\ve^4\tau^4+
10240\,\sqrt{3}\ve^3\tau^5+4608\,\ve^2\tau^6+7444\,\ve^4\tau^3-7168\,\sqrt{3}\ve^3\tau^4-\\
11520\,\ve^2\tau^5-768\,\sqrt{3}\,\ve\,\tau^6-4368\,\ve^4\tau^2-8064\,\sqrt{3}\ve^3\tau^3+
7488\,\ve^2\tau^4+1792\,\sqrt{3}\,\ve\,\tau^5+136\,\tau^6-1836\,\ve^4\tau+5056\,\sqrt{3}\ve^3\tau^2+\\
10656\,\ve^2\tau^3-768\,\sqrt{3}\,\ve\,\tau^4-220\,\tau^5-144\,\ve^4+1536\,\sqrt{3}\ve^3\tau
-7200\,\ve^2\tau^2-2464\,\sqrt{3}\,\ve\,\tau^3-154\,\tau^4+192\,\sqrt{3}\ve^3-1152\,\ve^2\tau+\\
1664\,\sqrt{3}\,\ve\,\tau^2+771\,\tau^3-288\,\ve^2+48\,\sqrt{3}\,\ve\,\tau-464\,\tau^2+64\,\sqrt{3}\,\ve+
28\,\tau-16)\tau\Bigr]\Big/\Bigl[32\,\left( 3\,\ve-\sqrt{3} \right)^4(2\,\tau-1)(2\,\tau+1)^2(\tau-1)^2\Bigr],
\end{multline*}
}
{\small
\begin{multline*}
\varpi_{3,3}(\tau,\ve)=
\Bigl[(2096\,\ve^4\tau^5-3376\,\ve^4\tau^4-2048\,\sqrt{3}\ve^3\tau^5+2204\,\ve^4\tau^3+
3840\,\sqrt{3}\ve^3\tau^4+2304\,\ve^2\tau^5+3864\,\ve^4\tau^2-2304\,\sqrt{3}\ve^3\tau^3-\\
4608\,\ve^2\tau^4-384\,\sqrt{3}\,\ve\,\tau^5-1332\,\ve^4\tau-4864\,\sqrt{3}\ve^3\tau^2+
2016\,\ve^2\tau^3+704\,\sqrt{3}\,\ve\,\tau^4+68\,\tau^5-1152\,\ve^4+1536\,\sqrt{3}\ve^3\tau+
6912\,\ve^2\tau^2-32\,\sqrt{3}\,\ve\,\tau^3-\\
76\,\tau^4+1536\,\sqrt{3}\ve^3-1728\,\ve^2\tau-1440\,\sqrt{3}\,\ve\,\tau^2-115\,\tau^3-
2304\,\ve^2+208\,\sqrt{3}\,\ve\,\tau+328\,\tau^2+512\,\sqrt{3}\,\ve+4\,\tau-128)\tau\Bigr]\Big/\Bigl[32\,(\tau-1)^2\\
\left( -3\,\ve+\sqrt{3} \right)^4(2\,\tau+1)^2\Bigr],
\end{multline*}
}
{\small
\begin{multline*}
\varpi_{3,4}(\tau,\ve)=
\Bigl[(2064\,\ve^4\tau^6-1696\,\ve^4\tau^5-2048\,\sqrt{3}\ve^3\tau^6-3292\,\ve^4\tau^4+
1920\,\sqrt{3}\ve^3\tau^5+2304\,\ve^2\tau^6+788\,\ve^4\tau^3+2816\,\sqrt{3}\ve^3\tau^4-\\
2304\,\ve^2\tau^5-384\,\sqrt{3}\,\ve\,\tau^6-3948\,\ve^4\tau^2-2528\,\sqrt{3}\ve^3\tau^3-
3168\,\ve^2\tau^4+320\,\sqrt{3}\,\ve\,\tau^5+68\,\tau^6-5940\,\ve^4\tau+3872\,\sqrt{3}\ve^3\tau^2+4896\,\ve^2\tau^3+\\
672\,\sqrt{3}\,\ve\,\tau^4-8\,\tau^5-1800\,\ve^4+7392\,\sqrt{3}\ve^3\tau-
3600\,\ve^2\tau^2-1088\,\sqrt{3}\,\ve\,\tau^3-191\,\tau^4+2400\,\sqrt{3}\ve^3-
10080\,\ve^2\tau+304\,\sqrt{3}\,\ve\,\tau^2+213\,\tau^3-3600\,\ve^2\\
+1968\,\sqrt{3}\,\ve\,\tau+44\,\tau^2+800\,\sqrt{3}\,\ve-412\,\tau-200)\tau\Bigr]\Big/
\Bigl[32\,(2\,\tau+1)^2\left( -3\,\ve+\sqrt{3} \right)^4(\tau+1)(\tau-1)^2\Bigr],
\end{multline*}
}
{\small
\begin{multline*}
\varpi_{3,5}(\tau,\ve)=
\Bigl[(1032\,\ve^4\tau^5+3280\,\ve^4\tau^4-1024\,\sqrt{3}\ve^3\tau^5+2186\,\ve^4\tau^3-
3136\,\sqrt{3}\ve^3\tau^4+1152\,\ve^2\tau^5-1806\,\ve^4\tau^2-1920\,\sqrt{3}\ve^3\tau^3+\\
3456\,\ve^2\tau^4-192\,\sqrt{3}\,\ve\,\tau^5-2376\,\ve^4\tau+2384\,\sqrt{3}\ve^3\tau^2+
2448\,\ve^2\tau^3-576\,\sqrt{3}\,\ve\,\tau^4+36\,\tau^5-648\,\ve^4+3024\,\sqrt{3}\ve^3\tau-
2736\,\ve^2\tau^2-\\
624\,\sqrt{3}\,\ve\,\tau^3+108\,\tau^4+864\,\sqrt{3}\ve^3-4104\,\ve^2\tau
+264\,\sqrt{3}\,\ve\,\tau^2+207\,\tau^3-1296\,\ve^2+768\,\sqrt{3}\,\ve\,\tau+
54\,\tau^2+288\,\sqrt{3}\,\ve-144\,\tau-72)\tau\Bigr]\Big/\\
\Bigl[16\,(2\,\tau+1)^2\left( 3\,\ve-\sqrt{3} \right)^4(\tau+1)(\tau+2)\Bigr],
\end{multline*}
}
with the corresponding intervals
$\mI_1=\bigl[1,1-\sqrt{3}\,\ve\bigr)$,
$\mI_2=\bigl[1-\sqrt{3}\,\ve,\sqrt{3}/(2\,\ve)-3/2\bigr)$,
$\mI_3=\bigl[\sqrt{3}/(2\,\ve)-3/2,\sqrt{3}/\ve-3\bigr)$,
$\mI_4=\bigl[\sqrt{3}/\ve-3,4\,\left(1-\sqrt{3}\,\ve \right)\bigr)$,
and
$\mI_5=\bigl[4\,\left(1-\sqrt{3}\,\ve \right),1\bigr)$.

For $\ve \in \bigl[2\,\sqrt{3}/7,\sqrt{3}/3\bigr)$,
$$\mu^S_{_{CS}}(\tau,\ve)=\sum_{j=1}^2 \varpi_{4,j}(\tau,\ve)\,\I(\tau \in \mI_j)$$
where $\varpi_{4,j}(\tau,\ve)=\varpi_{3,j}(\tau,\ve)$ for $j=1,\ldots,5$, and for $j=6,7,8$
{\small
\begin{multline*}
\varpi_{4,6}(\tau,\ve)=
\Bigl[-254\,\ve^4\tau^4-1998\,\ve^4\tau^3+256\,\sqrt{3}\ve^3\tau^4-4752\,\ve^4\tau^2+
2160\,\sqrt{3}\ve^3\tau^3-288\,\ve^2\tau^4-4320\,\ve^4\tau+5328\,\sqrt{3}\ve^3\tau^2-\\
2592\,\ve^2\tau^3+48\,\sqrt{3}\,\ve\,\tau^4-1296\,\ve^4+5184\,\sqrt{3}\ve^3\tau-
6552\,\ve^2\tau^2+456\,\sqrt{3}\,\ve\,\tau^3-9\,\tau^4+1728\,\sqrt{3}\ve^3-6912\,\ve^2\tau+
1152\,\sqrt{3}\,\ve\,\tau^2-90\,\tau^3-\\
2592\,\ve^2+1344\,\sqrt{3}\,\ve\,\tau-216\,\tau^2+576\,\sqrt{3}\,\ve-288\,\tau-144\Bigr]\Big/
\Bigl[16\,(\tau+2)\left( -3\,\ve+\sqrt{3} \right)^4(\tau+1)\Bigr],
\end{multline*}
}
{\small
\begin{multline*}
\varpi_{4,7}(\tau,\ve)=
\Bigl[-256\,\ve^4\tau^4-1536\,\ve^4\tau^3+256\,\sqrt{3}\ve^3\tau^4-2160\,\ve^4\tau^2+
1664\,\sqrt{3}\ve^3\tau^3-288\,\ve^2\tau^4-2160\,\ve^4\tau+2304\,\sqrt{3}\ve^3\tau^2-\\
2016\,\ve^2\tau^3+48\,\sqrt{3}\,\ve\,\tau^4-1296\,\ve^4+2592\,\sqrt{3}\ve^3\tau-
2664\,\ve^2\tau^2+360\,\sqrt{3}\,\ve\,\tau^3-9\,\tau^4+1728\,\sqrt{3}\ve^3-
3456\,\ve^2\tau+432\,\sqrt{3}\,\ve\,\tau^2-72\,\tau^3\\
-2592\,\ve^2+672\,\sqrt{3}\,\ve\,\tau-72\,\tau^2+576\,\sqrt{3}\,\ve-144\,\tau-144\Bigr]\Big/
\Bigl[16\,\left( -3\,\ve+\sqrt{3} \right)^4(\tau+1)\tau\Bigr],
\end{multline*}
}
$$\varpi_{4,8}(\tau,\ve)=1,$$
with the corresponding intervals
$\mI_1=\bigl[1,1-\sqrt{3}\,\ve\bigr)$,
$\mI_2=\bigl[1-\sqrt{3}\,\ve,\sqrt{3}/(2\,\ve)-3/2\bigr)$,
$\mI_3=\bigl[\sqrt{3}/(2\,\ve)-3/2,\sqrt{3}/\ve-3\bigr)$,
$\mI_4=\bigl[\sqrt{3}/\ve-3,4\,\left(1-\sqrt{3}\,\ve \right)\bigr)$,
$\mI_5=\bigl[4\,\left(1-\sqrt{3}\,\ve \right),
\frac{\sqrt{3}\,\left(1-\sqrt{3}\,\ve \right)}{4\,\ve-\sqrt{3}}\bigr)$,
$\mI_6=\Bigl[\frac{\sqrt{3}\,\left( 1-\sqrt{3}\,\ve \right)}{4\,\ve-\sqrt{3}},2\,(\sqrt{3}/\ve-3)\Bigr)$,
$\mI_7=\Bigl[2\,\left( \sqrt{3}/\ve-3\right),\frac{2\,\sqrt{3}\,(1-3\,\ve)}{4\,\ve-\sqrt{3}}\Bigr)$,
and
$\mI_8=\Bigl[\frac{2\,\sqrt{3}\,(1-3\,\ve)}{4\,\ve-\sqrt{3}},1\Bigr)$.

\subsubsection*{$\mu^A_{_{CS}}(\tau,\ve)$ under the Association Alternatives}
\label{sec:arc-prob-agg-CS}
Under association, we compute
$\mu^A_{_{CS}}(\tau,\ve)$ explicitly. For $\ve \in
\bigl[0,\sqrt{3}/21 \approx .0825\bigr)$,
$\mu^A_{_{CS}}(\tau,\ve)=\sum_{j=1}^7
\varpi_{1,j}(\tau,\ve)\,\I(\tau \in \mI_j)$ where
{\small
\begin{multline*}
\varpi_{1,1}(\tau,\ve)=
\Bigl[(-63936\,\tau^6\ve^4+20736\,\sqrt{3}\tau^6\ve^3-145728\,\tau^5\ve^4-6912\,\tau^6\ve^2+
46848\,\sqrt{3}\tau^5\ve^3+181872\,\tau^4\ve^4+192\,\sqrt{3}\tau^6\,\ve-14976\,\tau^5\ve^2-\\
60480\,\sqrt{3}\tau^4\ve^3+346896\,\tau^3\ve^4+36\,\tau^6+
256\,\sqrt{3}\tau^5\,\ve+22464\,\tau^4\ve^2-107712\,\sqrt{3}\tau^3\ve^3-
296640\,\tau^2\ve^4+128\,\tau^5-1200\,\sqrt{3}\tau^4\,\ve+\\
28512\,\tau^3\ve^2+93696\,\sqrt{3}\tau^2\ve^3-228528\,\tau\,\ve^4+63\,\tau^4+
1056\,\sqrt{3}\tau^3\,\ve-27360\,\tau^2\ve^2+74304\,\sqrt{3}\tau\,\ve^3+81648\,\ve^4-
726\,\tau^3-208\,\sqrt{3}\tau^2\,\ve-\\
25056\,\tau\,\ve^2-25920\,\sqrt{3}\ve^3+445\,\tau^2+768\,\sqrt{3}\tau\,\ve+
7776\,\ve^2+108\,\tau-108)\tau^2\Bigr]\Big/\Bigl[18\,(2\,\tau+1)(2\,\tau-1)
\left( -6\,\ve+\sqrt{3} \right)^2\\
\left( 6\,\ve+\sqrt{3} \right)^2(\tau+2)^2(\tau-1)^2\Bigr],
\end{multline*}
}
{\small
\begin{multline*}
\varpi_{1,2}(\tau,\ve)=
\Bigl[(-62784\,\tau^7\ve^4+19200\,\sqrt{3}\tau^7\ve^3-199872\,\tau^6\ve^4-4608\,\tau^7\ve^2+
56832\,\sqrt{3}\tau^6\ve^3+34992\,\tau^5\ve^4 -320\,\sqrt{3}\tau^7\,\ve-9216\,\tau^6\ve^2-\\
28224\,\sqrt{3}\tau^5\ve^3+476640\,\tau^4\ve^4+164\,\tau^7-1600\,\sqrt{3}\tau^6\,\ve+22464\,\tau^5\ve^2
-151104\,\sqrt{3}\tau^4\ve^3+82368\,\tau^3\ve^4+484\,\tau^6-1872\,\sqrt{3}\tau^5\,\ve+31392\,\tau^4\ve^2+\\
6528\,\sqrt{3}\tau^3\ve^3-405216\,\tau^2\ve^4-15\,\tau^5+3872\,\sqrt{3}\tau^4\,\ve-31392\,\tau^3\ve^2+
153792\,\sqrt{3}\tau^2\ve^3-252720\,\tau\,\ve^4-1214\,\tau^4+3280\,\sqrt{3}\tau^3\,\ve
-47520\,\tau^2\ve^2+\\
67392\,\sqrt{3}\tau\,\ve^3-46656\,\ve^4-13\,\tau^3-768\,\sqrt{3}\tau^2\,\ve-
7776\,\tau\,\ve^2+324\,\tau^2+108\,\tau)\tau\Bigr]\Big/
\Bigl[18\,(\tau-1)^2(2\,\tau+1)^2(\tau+2)^2\left( -6\,\ve+\sqrt{3} \right)^2\left( 6\,\ve+\sqrt{3} \right)^2\Bigr],
\end{multline*}
}
{\small
\begin{multline*}
\varpi_{1,3}(\tau,\ve)=
\Bigl[(-62784\,\tau^6\ve^4+19200\,\sqrt{3}\tau^6\ve^3-74304\,\tau^5\ve^4-4608\,\tau^6\ve^2+
18432\,\sqrt{3}\tau^5\ve^3+183600\,\tau^4\ve^4-320\,\sqrt{3}\tau^6\,\ve-65088\,\sqrt{3}\tau^4\ve^3+\\
179424\,\tau^3\ve^4+164\,\tau^6-960\,\sqrt{3}\tau^5\,\ve+22464\,\tau^4\ve^2-67584\,\sqrt{3}\tau^3\ve^3
+3456\,\tau^2\ve^4+156\,\tau^5+48\,\sqrt{3}\tau^4\,\ve+21456\,\tau^3\ve^2+1728\,\sqrt{3}\tau^2\ve^3+\\
7776\,\tau\,\ve^4-327\,\tau^4-112\,\sqrt{3}\tau^3\,\ve
-4320\,\tau^2\ve^2+10368\,\sqrt{3}\tau\,\ve^3+11664\,\ve^4-74\,\tau^3-384\,\sqrt{3}\tau^2\,\ve-3888\,\tau\,\ve^2+135\,\tau^2+54\,\tau)\tau\Bigr]\Big/\\
\Bigl[18\,(\tau-1)^2(2\,\tau+1)^2(\tau+2)\left( -6\,\ve+\sqrt{3} \right)^2\left( 6\,\ve+\sqrt{3} \right)^2\Bigr],
\end{multline*}
}
{\small
\begin{multline*}
\varpi_{1,4}(\tau,\ve)=
\Bigl[(-63072\,\tau^4\ve^4+19584\,\sqrt{3}\tau^4\ve^3-267552\,\tau^3\ve^4-5184\,\tau^4\ve^2+
81408\,\sqrt{3}\tau^3\ve^3-389304\,\tau^2\ve^4-192\,\sqrt{3}\tau^4\,\ve-20160\,\tau^3\ve^2+\\
118176\,\sqrt{3}\tau^2\ve^3-233712\,\tau\,\ve^4+132\,\tau^4-896\,\sqrt{3}\tau^3\,\ve-29520\,\tau^2\ve^2
+71712\,\sqrt{3}\tau\,\ve^3-48600\,\ve^4+488\,\tau^3-1072\,\sqrt{3}\tau^2\,\ve-18576\,\tau\,\ve^2+\\
15552\,\sqrt{3}\ve^3+601\,\tau^2-384\,\sqrt{3}\tau\,\ve-3888\,\ve^2+297\,\tau+54)\tau^2\Bigr]\Big/
\Bigl[18\,(2\,\tau+1)^2(\tau+2)(\tau+1)\left( -6\,\ve+\sqrt{3} \right)^2\left( 6\,\ve+\sqrt{3} \right)^2\Bigr],
\end{multline*}
}
{\small
\begin{multline*}
\varpi_{1,5}(\tau,\ve)=
\Bigl[-49968\,\tau^5\ve^4+15936\,\sqrt{3}\tau^5\ve^3-219384\,\tau^4\ve^4-4896\,\tau^5\ve^2+
64992\,\sqrt{3}\tau^4\ve^3-349920\,\tau^3\ve^4+32\,\sqrt{3}\tau^5\,\ve-18000\,\tau^4\ve^2+\\
90720\,\sqrt{3}\tau^3\ve^3-270216\,\tau^2\ve^4+58\,\tau^5+16\,\sqrt{3}\tau^4\,\ve-22032\,\tau^3\ve^2
+51840\,\sqrt{3}\tau^2\ve^3-112752\,\tau\,\ve^4+191\,\tau^4-9072\,\tau^2\ve^2+10368\,\sqrt{3}\tau\,\ve^3-\\
23328\,\ve^4+189\,\tau^3+54\,\tau^2\Bigr]\Big/
\Bigl[18\,(\tau+2)\left( 6\,\ve+\sqrt{3} \right)^2
\left( -6\,\ve+\sqrt{3} \right)^2(2\,\tau+1)(\tau+1)\Bigr],
\end{multline*}
}
{\small
\begin{multline*}
\varpi_{1,6}(\tau,\ve)=
\Bigl[-50040\,\tau^6\ve^4+16032\,\sqrt{3}\tau^6\ve^3-221220\,\tau^5\ve^4-5040\,\tau^6\ve^2+
66864\,\sqrt{3}\tau^5\ve^3-368100\,\tau^4\ve^4+64\,\sqrt{3}\tau^6\,\ve-19944\,\tau^5\ve^2+\\
103728\,\sqrt{3}\tau^4\ve^3-356616\,\tau^3\ve^4+50\,\tau^6+256\,\sqrt{3}\tau^5\,\ve-29880\,\tau^4\ve^2
+88992\,\sqrt{3}\tau^3\ve^3-307152\,\tau^2\ve^4+179\,\tau^5+304\,\sqrt{3}\tau^4\,\ve-17712\,\tau^3\ve^2+\\
46656\,\sqrt{3}\tau^2\ve^3-194400\,\tau\,\ve^4+185\,\tau^4+96\,\sqrt{3}\tau^3\,\ve-2592\,\tau^2\ve^2+
10368\,\sqrt{3}\tau\,\ve^3-46656\,\ve^4+54\,\tau^3\Bigr]\Big/
\Bigl[18\,(\tau+1)\left( 6\,\ve+\sqrt{3} \right)^2\\
\left( -6\,\ve+\sqrt{3} \right)^2(2\,\tau+1)(\tau+2)\tau\Bigr],
\end{multline*}
}
{\small
\begin{multline*}
\varpi_{1,7}(\tau,\ve)=
\Bigl[3\,(-1512\,\tau^4\ve^4+480\,\sqrt{3}\tau^4\ve^3-3780\,\tau^3\ve^4-144\,\tau^4\ve^2+
1200\,\sqrt{3}\tau^3\ve^3+216\,\tau^2\ve^4-360\,\tau^3\ve^2+480\,\sqrt{3}\tau^2\ve^3+4752\,\tau\,\ve^4+\\
2\,\tau^4-288\,\tau^2\ve^2+1728\,\ve^4+5\,\tau^3+2\,\tau^2)\Bigr]\Big/\Bigl[2\,\left( 6\,\ve+\sqrt{3} \right)^2
\left( -6\,\ve+\sqrt{3} \right)^2(2\,\tau+1)(\tau+2)\Bigr],
\end{multline*}
}
with the corresponding intervals
$\mI_1=\Bigl[0,\frac{3\,\sqrt{3}\,\ve}{2\,\left(1-\sqrt{3}\,\ve \right)}\Bigr)$,
$\mI_2=\Bigl[\frac{3\,\sqrt{3}\,\ve}{2\,\left(1-\sqrt{3}\,\ve \right)},\frac{2\,\sqrt{3}\,\ve}{1-2\,\sqrt{3}\,\ve}\Bigr)$,
$\mI_3=\Bigl[\frac{2\,\sqrt{3}\,\ve}{1-2\,\sqrt{3}\,\ve},\frac{3\,\sqrt{3}\,\ve}{1-\sqrt{3}\,\ve}\Bigr)$,
$\mI_4=\Bigl[\frac{3\,\sqrt{3}\,\ve}{1-\sqrt{3}\,\ve},\frac{3\,\sqrt{3}\,\ve}{1-4\,\sqrt{3}\,\ve}\Bigr)$,
$\mI_5=\Bigl[\frac{3\,\sqrt{3}\,\ve}{1-4\,\sqrt{3}\,\ve},\frac{6\,\sqrt{3}\,\ve}{1-\sqrt{3}\,\ve}\Bigr)$,
$\mI_6=\Bigl[\frac{6\,\sqrt{3}\,\ve}{1-\sqrt{3}\,\ve},\frac{6\,\sqrt{3}\,\ve}{1-4\,\sqrt{3}\,\ve}\Bigr)$
and
$\mI_7=\bigl[\frac{6\,\sqrt{3}\,\ve}{1-4\,\sqrt{3}\,\ve},1\bigr)$.

For $\ve \in \bigl[ \sqrt{3}/21,\sqrt{3}/12 \bigr)$,
$$\mu^A_{_{CS}}(\tau,\ve)=\sum_{j=1}^4 \varpi_{2,j}(\tau,\ve)\,\I(\tau \in \mI_j)$$
where
$\varpi_{2,j}(\tau,\ve)=\varpi_{1,j}(\tau,\ve)$ for
$j=1,2,3,4$ with the corresponding intervals
$\mI_1=\Bigl[0,\frac{3\,\sqrt{3}\,\ve}{2\,\left(1-\sqrt{3}\,\ve \right)}\Bigr)$,
$\mI_2=\Bigl[\frac{3\,\sqrt{3}\,\ve}{2\,\left(1-\sqrt{3}\,\ve \right)},\frac{2\,\sqrt{3}\,\ve}{1-2\,\sqrt{3}\,\ve}\Bigr)$,
$\mI_3=\Bigl[\frac{2\,\sqrt{3}\,\ve}{1-2\,\sqrt{3}\,\ve},\frac{3\,\sqrt{3}\,\ve}{1-\sqrt{3}\,\ve}\Bigr)$,
and
$\mI_4=\Bigl[\frac{3\,\sqrt{3}\,\ve}{1-\sqrt{3}\,\ve},1\Bigr)$.

For $\ve \in \bigl[ \sqrt{3}/12,2\,\sqrt{3}/15 \bigr)$,
$$\mu^A_{_{CS}}(\tau,\ve)=\sum_{j=1}^2 \varpi_{3,j}(\tau,\ve)\,\I(\tau \in \mI_j)$$
where
$$\varpi_{3,1}(\tau,\ve)=\frac{(\tau^2+5\,\tau+9)\tau^2}{18\,(\tau+2)^2},$$
{\small
\begin{multline*}
\varpi_{3,2}(\tau,\ve)=
-\Big[(-216\,\tau^6\ve^4+288\,\sqrt{3}\tau^6\ve^3-1836\,\tau^5\ve^4-432\,\tau^6\ve^2+
2160\,\sqrt{3}\tau^5\ve^3-6786\,\tau^4\ve^4+96\,\sqrt{3}\tau^6\,\ve-2808\,\tau^5\ve^2+\\
6600\,\sqrt{3}\tau^4\ve^3-13401\,\tau^3\ve^4-24\,\tau^6+528\,\sqrt{3}\tau^5\,\ve-
6876\,\tau^4\ve^2+9588\,\sqrt{3}\tau^3\ve^3-16074\,\tau^2\ve^4-108\,\tau^5+
1000\,\sqrt{3}\tau^4\,\ve-6498\,\tau^3\ve^2+\\
7500\,\sqrt{3}\tau^2\ve^3-10611\,\tau\,\ve^4-154\,\tau^4+484\,\sqrt{3}\tau^3\,\ve-
2178\,\tau^2\ve^2+2484\,\sqrt{3}\tau\,\ve^3-2916\,\ve^4-25\,\tau^3-92\,\sqrt{3}\tau^2\,\ve+162\,\tau\,\ve^2+23\,\tau^2\\
-36\,\sqrt{3}\tau\,\ve+9\,\tau)\tau\Big]\Big/\Big[2\,(2\,\tau-1)(2\,\tau+1)^2(\tau+2)^2\left( -3\,\ve+\sqrt{3} \right)^4\Big],
\end{multline*}
}
with the corresponding intervals
$\mI_1=\Bigl[0,\frac{3\,\sqrt{3}\,\ve}{2\,\left(1-\sqrt{3}\,\ve \right)}\Bigr)$,
and
$\mI_2=\Bigl[\frac{3\,\sqrt{3}\,\ve}{2\,\left(1-\sqrt{3}\,\ve \right)},1\Bigr)$.

For $\ve \in \bigl[ 2\,\sqrt{3}/15,\sqrt{3}/3 \bigr)$,
we have
$\mu^A_{_{CS}}(\tau,\ve)=\varpi_{3,1}(\tau,\ve)\,\I(\tau \in [0,1])$.

\end{document}